# Introduction into Calculus over Division Ring

Aleks Kleyn

*E-mail address*: Aleks_Kleyn@MailAPS.org
*URL*: http://sites.google.com/site/alekskleyn/
*URL*: http://arxiv.org/a/kleyn_a_1
*URL*: http://AleksKleyn.blogspot.com/

Abstract. Based on twin representations of division ring in an Abelian group I consider $D$-vector spaces over division ring. Morphism of $D$-vector spaces is linear map of $D$-vector spaces.

I consider derivative of function $f$ of normed division ring as linear map the most close to function $f$. I explore expansion of map into Taylor series and method to find solution of differential equation.

The norm in $D$-vector space allows considering of continuous mapping of $D$-vector spaces. Differential of mapping $f$ of $D$-vector spaces is defined as linear mapping the most close to map $f$.

# Contents

CHAPTER 1

# Preface

## 1.1. Preface to Version 1

> The way of charms and disappointments.
> However I would put it another way. When you understand that the reason for your disappointments is in your expectations, you start to examine surrounding landscape. You found yourself on the path where nobody went before. New impressions are stronger then impressions from travel along well known and trodden road.

Author is unknown. The Adventures of an Explorer.

In mathematics like in everyday life there are statements that after studying become obvious. It may appear surprising when you first time meet that these statements are false in new conditions.[1.1]

When I started to study linear algebra over division ring I was ready to different surprises related with noncommutativity of product. However large number of statements similar to statements of linear algebra over field relaxed my attention. Nevertheless I met first deep problem in linear algebra: definition of polylinear form and tensor product assumes that product is commutative. The firm belief that solving of particular problems gives key to solving of general problem, was main motive, prompting me to consider calculus over division ring. This book introduces technique of differentiation of mapping over division ring. Since the possibility of linear approximation of mapping is at the heart of calculus, this book is natural consequence of book [3]. In calculus like in linear algebra we see statements that we know well. However we see a lot of new statements.

The first considered definition of derivative was based on classic definition of derivative. However considering of derivative of function $x^2$, I discovered that I cannot express this derivative analytically. To answer the question of existence of derivative I decided to conduct calculating experiment. The best candidate for calculating experiment is division ring of quaternions. Simple calculations make it

[1.1] I experienced such sensation when I considered Rolle's theorem in the calculus over division ring. After numerous attempts to change statement of Rolle's theorem, I came to a conclusion that it is false even for complex numbers. It is enough to consider function $w = z^2$. $z = 0$ is the only point where derivative equal 0. However we can connect points $z = 0.5$ and $z = -0.5$ by curve that does not pass through point $z = 0$.

clear that sometimes derivative does not exist. It causes that the set of differentiable functions is very small, and such theory of differentiation is not interesting.[1.2] Analysis of calculating experiment shows that an expression that depends on increment of argument appears even derivative does not exist. We know such tool in functional analysis. This is the Fréchet derivative and the Gâteaux derivative in the space of linear operators.

When I defined the Gâteaux $D\star$-derivative[1.3], I put attention that differential does not depend on using of the Gâteaux $D\star$-derivative or the Gâteaux $\star D$-derivative. This caused me to change text and I started from definition of differential.

The concept of $D_*{}^*$-vector space is powerful tool for study of linear spaces over division ring. However the Gâteaux differential does not fit in frame of $D_*{}^*$-vector space. It follows from example 5.2.2 that $D$-vector space is the most appropriate model for studying of the Gâteaux differential.

The considering of $D_*{}^*$-linear maps leads to the concept of twin representations of division ring in Abelian group. Abelian group, in which we define twin representation of division ring $D$, is called $D$-vector space. Since we define multiplication over elements of division ring $D$ from left, as well as right, then homomorphisms of $D$-vector spaces cannot keep this operation. This leads to concept of additive map which is morphism of $D$-vector space.[1.4]

Consideration of the Gâteaux differential as additive map leads to new structure for definition of derivative. In contrast to the Gâteaux $D_*{}^*$-derivative components of the Gâteaux differential do not depend on direction. However components of the Gâteaux differential are defined non-uniquely and their number is uncertain. It is evident that when we consider calculus over field, the Gâteaux $D\star$-derivative and components of the Gâteaux differential coincide with derivative. These structures appear as result of different constructions. When product is commutative these structures coincide. However noncomutativity breaks expected relationship. In spite of difference in behavior, these structures are supplemental.

December, 2008

## 1.2. Preface to Version 3

Discussing the notion of derivative which I introduced in this book, I realized that the subject of the book is very complex. This is the reason why in this part of preface I want to consider in more detail the concept of differentiation.

When we study functions of one variable, the derivative in selected point is a number. When we study function of multiple variables, we realize that it is not enough to use number. The derivative becomes vector or gradient. When we study

[1.2]Search in papers and books shows that this is a deep problem. Mathematicians try to solve this problem suggesting different definitions of derivative. These definitions are based on different properties of derivative of function of real or complex variable. These definitions are very interesting. However my task is the study of common rules for differentiation of maps of complete division ring of characteristic 0.

[1.3]In this book I use notations introduced in [3].

[1.4]Design of additive mapping which I made proving theorem 4.1.4, is instructive methodically. From the time when I started to learn differential geometry, I accepted the fact that we can everything express using tensor, that we enumerate components of tensor by indexes and this relationship is permanent. Once I saw that it was not the case, I found myself unprepared to recognize it It got month before I was able to write down expression that I saw by my eyes. The problem was because I tried to write down this expression in tensor form.

maps of vector spaces, this is a first time that we tell about derivative as operator. However since this operator is linear, then we can represent derivative as product over matrix. Again we express the derivative as a set of numbers.

Surely, such behavior of derivative weaken our attention. When we consider objects which are more complex than fields or vector spaces, we still try to see an object which can be written as a factor before an increment and which does not depend on the increment. Quaternion algebra is the simplest model of continuous division ring. However, an attempt to define derivative as linear map is doomed to failure. Only certain linear functions are differentiable in this sense.

The question arises as to whether there exists an alternative method if the set of linear maps restricts our ability to study infinitesimal behavior of map? The answer on this question is affirmative. Let us consider the quaternion algebra. The quaternion algebra is normed space. We know two types of derivative in normed space. Strong derivative or the Fréchet derivative is analogue of the derivative which we know very well. When I wrote about attempt to find derivative as linear map as increment of argument, I meant the Fréchet derivative. Besides strong derivative there exits weak derivative or the Gâteaux derivative. The main idea is that derivative may depend on direction.

Vadim Komkov defines the Gâteaux derivative for quaternion algebra ([9], p. 322) at the point $q$ in the direction of $q_1$ to be the quantity $q_2$ such that, for $\epsilon > 0$

$$f(q + \epsilon q_1) - f(q) = \epsilon(q_1 q_2)$$

He made this definition according to definition [1]-3.1.2, page 177. This definition does not have difference from the Fréchet derivative. This phenomenon remains in general case. Every time when we try to extract the Gâteaux derivative as factor, the identification of the Gâteaux derivative and the Fréchet derivative occurs.

A linear map of division ring is homomorphism of its additive group and homomorphism of its multiplicative group in one map. Homomorphism of multiplicative group is restrictive factor. This is why it is mind to get rid of homomorphism of multiplicatie group and to consider only homomorphism of additive group. It is evident that the set of homomorphisms of additive group is more wide than the set of linear maps. However, because we defined two operations in division ring, additive map has a special structure. This is why I selected these map into separate set and called them additive maps.

Even in this case some elements of linearity remain. The center of division ring is field. As a consequence, additive map is linear over center. When division ring has characteristic 0 and does not have discrete topology, we can embed real field into center of division ring. As a consequence, the definition of derivative as additive map appears equivalent to the definition of the Gâteaux derivative.

Qualitative change happens exactly here. In additive map we cannot extract an argument as factor. We cannot extract increment from derivative. Is this idea fundamentally new? I do not think so. Derivative is a map. It is the map of differential of an argument to differential of function. Even additive map of increment of argument allows approximate the increment of function by polynomial of first power. The structure of polynomial over division ring is different from structure of polynomial over field. I consider some properties of polynomial in section 8.2. Using obtained theorems I explore Taylor series expansion of map and method to find solution of differential equation.

In the beginning I also tried to define derivative of map over division ring as Fréchet derivative. However, I saw that this derivative does not satisfy to basic rules of differentiation, and I realized that this is wrong way. To make sure that the problem is fundamental, and is not inherent in the logic of my research, I made computation in the division ring of quaternions. That is where I discovered that the derivative is dependent on the direction.

At this time I studied maps of $D$-vector spaces, and idea of additive map became leading one in this research. In preface to version 1 (section 1.1) I wrote about impressions of explorer who found oneself for the first time in unfamiliar land. Expecting to see familiar landscapes, explorer does not notice initially the peculiarity of new places. I experience such change of impressions during this year when I wrote the book.

I did not realize right away that derivative which I discovered is the Gâteaux derivative. I could see difference between my derivative and the Gâteaux derivative. I had to come a long way before realizing that two definitions are equivalent. Like explorer, who tried to see familiar landscapes, I considered possibility to extract increment from derivative. As a consequence, the Gâteaux $D\star$-derivative emerged.

If we consider new concepts, what will survive in the future? Surely, the Gâteaux derivative and the Gâteaux partial derivative will be basic concepts. Components of the Gâteaux derivative are important for considering of structure of the Gâteaux derivative. I do not think that I will consider the Gâteaux $D\star$-derivative in other cases except the Gâteaux derivative of first order. Primary advantage of the Gâteaux $D\star$-derivative is that it allows you to explore conditions when map is continuous.

May, 2009

## 1.3. Preface to Version 4

Quaternion algebra is the simplest example of division ring. This is why I verify theorems in quaternion algebra to see how they work. Quaternion algebra is similar to complex field and so it is natural to search some parallel.

In real field any additive map automatically turns out to be linear over real field. This statement is caused by statement that real field is completion of rational field and is corollary of theorem 3.1.3.

However the statement changes for complex field. Not every additive map of complex field is linear over complex field. Conjugation is the simplest example of such map. As soon as I discovered this extremely interesting statement, I returned to question about analytical representation of additive map.

Exploring additive maps, I realized that I too abruptly expanded the set of linear maps while switching from field to division ring. The reason for this was not a very clear understanding of how to overcome the noncommutativity of the product. However during the process of research it became more evident that any additive map is linear over some field. For the first time I used this concept when constructing the tensor product and the concept took shape in a subsequent study.

During constructing the tensor product of division rings $D_1$, ..., $D_n$ I assumed existence of field $F$ such that additive map of division ring $D_i$ is linear over this field for any $i$. If all division rings have characteristic 0, then according to theorem 3.1.3 such field always exists. However dependence of tensor product on selected field $F$ arises here. To get rid of this dependence, I assume that field $F$ is maximal field that possesses stated property.

If $D_1 = ... = D_n = D$, then such field is center $Z(D)$ of division ring $D$. If the product in division ring $D$ is commutative[1.5], then $Z(D) = D$. So, starting with additive map, I arrive at concept of linear map, which is generalization of linear map over field.

My decision to explore the Gâteaux derivative as linear map was influenced by following.

- Linear map of field is homogeneous polynomial of degree one.
- For continuous field derivative of function of field is homogeneous polynomial of degree one that is good approximation to change of function. The same time we explore the derivative as linear map of field. At the transition to division ring we can save the same link between derivative and additive map.

Research in area of complex numbers and quaternions revealed one more interesting phenomenon. In spite of the fact that complex field is extension of real field, the structure of linear map over complex field is different from the structure of linear map over real field. This difference leads to the statement that conjugation of complex numbers is additive map but not linear map over complex field.

Similarly, the structure of linear map over division ring of quaternions is different from the structure of linear map over complex field. The source of difference is statement that the center of quaternion algebra has more simple structure than complex field. This difference leads to the statement that conjugation of quaternion satisfies to equation

$$\overline{p} = -\frac{1}{2}(p + ipi + jpj + kpk)$$

Consequently, the problem to find a mapping satisfying to a theorem similar to the Riemann theorem (theorem 7.1.1), is a nontrivial task for quaternions.

August, 2009

## 1.4. Preface to Version 5

Recently, in internet, I found the book [11], where Hamilton found similar description of calculus in quaternion algebra. Although this book was written more than 100 years ago, I recommend to read this book everybody who are interested in problems of calculus over quaternion algebra.

December, 2009

## 1.5. Conventions

Convention 1.5.1. *Function and map are synonyms. However according to tradition, correspondence between either rings or vector spaces is called map and a map of either real field or quaternion algebra is called function.* □

Convention 1.5.2. *In any expression where we use index I assume that this index may have internal structure. For instance, considering the algebra $A$ we enumerate coordinates of $a \in A$ relative to basis $\overline{\overline{e}}$ by an index $i$. This means that $a$ is a vector. However, if $a$ is matrix, then we need two indexes, one enumerates rows, another enumerates columns. In the case, when index has structure, we begin*

[1.5] In other words, division ring is field.

*the index from symbol $\cdot$ in the corresponding position. For instance, if I consider the matrix $a^i_j$ as an element of a vector space, then I can write the element of matrix as $a^{\cdot i}_{\cdot j}$.* □

CONVENTION 1.5.3. *I assume sum over index $s$ in expression like*

$$a_{s\cdot 0}xa_{s\cdot 1}$$

□

CONVENTION 1.5.4. *We can consider division ring $D$ as $D$-vector space of dimension 1. According to this statement, we can explore not only homomorphisms of division ring $D_1$ into division ring $D_2$, but also linear maps of division rings.* □

CONVENTION 1.5.5. *In spite of noncommutativity of product a lot of statements remain to be true if we substitute, for instance, right representation by left representation or right vector space by left vector space. To keep this symmetry in statements of theorems I use symmetric notation. For instance, I consider $D\star$-vector space and $\star D$-vector space. We can read notation $D\star$-vector space as either $D$-star-vector space or left vector space. We can read notation $D\star$-linear dependent vectors as either $D$-star-linear dependent vectors or vectors that are linearly dependent from left.* □

CONVENTION 1.5.6. *Let $A$ be free algebra with finite or countable basis. Considering expansion of element of algebra $A$ relative basis $\overline{\overline{e}}$ we use the same root letter to denote this element and its coordinates. In expression $a^2$, it is not clear whether this is component of expansion of element $a$ relative basis, or this is operation $a^2 = aa$. To make text clearer we use separate color for index of element of algebra. For instance,*

$$a = a^i e_i$$

□

CONVENTION 1.5.7. *It is very difficult to draw the line between the module and the algebra. Especially since sometimes in the process of constructing, we must first prove that the set $A$ is a module, and then we prove that this set is an algebra. Therefore, to write the element of the module, we will also use the convention 1.5.6.* □

CONVENTION 1.5.8. *The identification of the vector and matrix of its coordinates can lead to ambiguity in the equation*

$$a = a^*{}_* e \tag{1.5.1}$$

*where $\overline{\overline{e}}$ is a basis of vector space. Therefore, we write the equation (1.5.1) in the following form*

$$\overline{a} = a^*{}_* e$$

*in order to see where we wrote vector.* □

CONVENTION 1.5.9. *If free finite dimensional algebra has unit, then we identify the vector of basis $\overline{e}_0$ with unit of algebra.* □

Without a doubt, the reader may have questions, comments, objections. I will appreciate any response.

# Part 1

# Linear Algebra over Division Ring

CHAPTER 2

# Linear Map of ${}_*{}^*D$-Vector Space

## 2.1. Linear Map of ${}_*{}^*D$-Vector Space

In this subsection we assume $V$, $W$ are ${}_*{}^*D$-vector spaces.

Definition 2.1.1. Let us denote by $\mathcal{L}({}_*{}^*D;V;W)$ set of linear maps

$$A:V\to W$$

of ${}_*{}^*D$-vector space $V$ into ${}_*{}^*D$-vector space $W$. Let us denote by $\mathcal{L}(D_*{}^*;V;W)$ set of linear maps

$$A:V\to W$$

of $D_*{}^*$-vector space $V$ into $D_*{}^*$-vector space $W$. □

We can consider division ring $D$ as ${}_*{}^*D$-vector space of dimension 1. Correspondingly we can consider set $\mathcal{L}(D_*{}^*;D;W)$ and $\mathcal{L}({}_*{}^*D;V;D)$.

Definition 2.1.2. Let us denote by $\mathcal{L}(*T;S;R)$ set of $\star T$-representations of division ring $S$ in additive group of division ring $R$. Let us denote by $\mathcal{L}(T*;S;R)$ set of $T\star$-representations of division ring $S$ in additive group of division ring $R$. □

Theorem 2.1.3. *Suppose $V$, $W$ are ${}_*{}^*D$-vector spaces. Then set $\mathcal{L}({}_*{}^*D;V;W)$ is an Abelian group relative composition law*

$$(A+B)_*{}^*x=A_*{}^*x+B_*{}^*x \tag{2.1.1}$$

*Linear map $A+B$ is called* **sum of linear maps** *$A$ and $B$.*

Proof. We need to show that mapping

$$A+B:V\to W$$

defined by equation (2.1.1) is linear map of ${}_*{}^*D$-vector spaces. According to the definition [3]-4.4.2

$$A_*{}^*(x_*{}^*a)=(A_*{}^*x)_*{}^*a$$
$$B_*{}^*(x_*{}^*a)=(B_*{}^*x)_*{}^*a$$

We see that

$$\begin{aligned}(A+B)_*{}^*(x_*{}^*a)&=A_*{}^*(x_*{}^*a)+B_*{}^*(x_*{}^*a)\\&=(A_*{}^*x)_*{}^*a+(B_*{}^*x)_*{}^*a\\&=(A_*{}^*x+B_*{}^*x)_*{}^*a\\&=((A+B)_*{}^*x)_*{}^*a\end{aligned}$$

We need to show also that this operation is commutative.

$$\begin{aligned}(A+B)_*{}^*x &= A_*{}^*x + B_*{}^*x\\ &= B_*{}^*x + A_*{}^*x\\ &= (B+A)_*{}^*x\end{aligned}$$

□

THEOREM 2.1.4. *Let* $\overline{\overline{e}}_V = (e_{V\cdot \boldsymbol{i}}, \boldsymbol{i} \in \boldsymbol{I})$ *be a basis in* ${}_*{}^*D$*-vector space* $V$ *and* $\overline{\overline{e}}_W = (e_{W\cdot \boldsymbol{j}}, \boldsymbol{j} \in \boldsymbol{J})$ *be a basis in* ${}_*{}^*D$*-vector space* $W$. *Let* $A = ({}^{\boldsymbol{j}}A_{\boldsymbol{i}})$, $\boldsymbol{i} \in \boldsymbol{I}$, $\boldsymbol{j} \in \boldsymbol{J}$, *be arbitrary matrix. Then map*

$$A : V \to W$$

*defined by equation*

$$b = A_*{}^*a$$

*relative to selected bases is a linear map of* ${}_*{}^*D$*-vector spaces.*

PROOF. Theorem 2.1.4 is inverse statement to theorem [3]-4.4.3. Suppose $A_*{}^*v = e_*{}^*A_*{}^*v$. Then

$$\begin{aligned}A_*{}^*(v_*{}^*a) &= e_*{}^*A_*{}^*v_*{}^*a\\ &= (A_*{}^*v)_*{}^*a\end{aligned}$$

□

THEOREM 2.1.5. *Let* $\overline{\overline{e}}_V = (e_{V\cdot \boldsymbol{i}}, \boldsymbol{i} \in \boldsymbol{I})$ *be a basis in* ${}_*{}^*D$*-vector space* $V$ *and* $\overline{\overline{e}}_W = (e_{W\cdot \boldsymbol{j}}, \boldsymbol{j} \in \boldsymbol{J})$ *be a basis in* ${}_*{}^*D$*-vector space* $W$. *Suppose linear map* $A$ *has matrix* $A = ({}^{\boldsymbol{j}}A_{\boldsymbol{i}})$, $\boldsymbol{i} \in \boldsymbol{I}$, $\boldsymbol{j} \in \boldsymbol{J}$, *relative to selected bases. Let* $m \in D$. *Then matrix*

$${}^{\boldsymbol{j}}(mA)_{\boldsymbol{i}} = m\ {}^{\boldsymbol{j}}A_{\boldsymbol{i}}$$

*defines linear map*

$$mA : V \to W$$

*which we call* **product of linear map** $A$ **over scalar**.

PROOF. The statement of the theorem is corollary of the theorem 2.1.4. □

THEOREM 2.1.6. *Set* $\mathcal{L}({}_*{}^*D; V; W)$ *is* $D_*{}^*$*-vector space.*

PROOF. Theorem 2.1.3 states that $\mathcal{L}({}_*{}^*D; V; W)$ is an Abelian group. It follows from theorem 2.1.5 that element of division ring $D$ defines $T\star$-transformation on the Abelian group $\mathcal{L}({}_*{}^*D; V; W)$. From theorems 2.1.4, [3]-4.1.1, and [3]-4.1.3 it follows that set $\mathcal{L}({}_*{}^*D; V; W)$ is $D\star$-vector space.

Writing elements of basis $D\star$-vector space $\mathcal{L}({}_*{}^*D; V; W)$ as ${}_*$-rows or ${}^*$-rows, we represent $D\star$-vector space $\mathcal{L}({}_*{}^*D; V; W)$ as $D^*{}_*$- or ${}_*{}^*D$-vector space. I want to stress that choice between $D_*{}^*$- and $D^*{}_*$-linear combination in $D\star$-vector space $\mathcal{L}({}_*{}^*D; V; W)$ does not depend on type of vector spaces $V$ and $W$.

To select the type of vector space $\mathcal{L}({}_*{}^*D; V; W)$ I draw attention to the following observation. let $V$ and $W$ be ${}_*{}^*D$-vector spaces. Suppose $\mathcal{L}({}_*{}^*D; V; W)$ is $D_*{}^*$- vector space. Then we can represent the operation of ${}^*$-row of ${}_*{}^*D$-linear maps ${}^jA$ on ${}_*$-row of vectors $f_i$ as matrix

$$\begin{pmatrix}{}^1A\\ ...\\ {}^mA\end{pmatrix}{}_*{}^*\begin{pmatrix}f_1 & ... & f_n\end{pmatrix}=\begin{pmatrix}{}^1A_*{}^*f_1 & ... & {}^1A_*{}^*f_n\\ ... & ... & ...\\ {}^mA_*{}^*f_1 & ... & {}^mA_*{}^*f_n\end{pmatrix}$$

This notation is coordinated with matrix notation of action of $D_*{}^*$-linear combination $a_*{}^*A$ of ${}_*{}^*D$-linear maps $A$. □

We can also define $\star D$-product of ${}_*{}^*D$-linear map $A$ over scalar. However in general we cannot carry this $\star T$-representation of division ring $D$ in $D\star$-vector space $\mathcal{L}({}_*{}^*D;V;W)$ into ${}_*{}^*D$-vector space $W$. Indeed, in case of $\star D$-product we get

$$(Am)_*{}^*v=(Am)_*{}^*f_*{}^*v=e_*{}^*(Am)_*{}^*v$$

Since product in division ring is noncommutative, we cannot express this expression as product of $A_*{}^*v$ over $m$.

Ambiguity of notation

$$m_*{}^*A_*{}^*v$$

is corollary of theorem 2.1.6. We may assume that meaning of this notation is clear from the text. To make notation more clear we will use brackets. Expression

$$w=[m_*{}^*A]_*{}^*v$$

means that $D_*{}^*$-linear composition of ${}_*{}^*D$-linear maps ${}^iA$ maps the vector $v$ to the vector $w$. Expression

$$w=B_*{}^*A_*{}^*v$$

means that ${}_*{}^*$-product of ${}_*{}^*D$-linear mappings $A$ and $B$ maps the vector $v$ to the vector $w$.

## 2.2. 1-Form on ${}_*{}^*D$-Vector Space

Definition 2.2.1. **1-form** on ${}_*{}^*D$-vector space $V$ is linear map

$$b:V\to D \tag{2.2.1}$$

□

We can write value of 1-form $b$, defined for vector $a$, as

$$b(a)=<b,a>$$

Theorem 2.2.2. *Set $\mathcal{L}({}_*{}^*D;V;D)$ is $D\star$-vector space.*

Proof. ${}_*{}^*D$-vector space of dimension 1 is equivalent to division ring $D$ □

Theorem 2.2.3. *Let $\overline{\overline{e}}$ be a basis in ${}_*{}^*D$-vector space $V$. 1-form $b$ has presentation*

$$<b,a>=b_*{}^*a \tag{2.2.2}$$

*relative to selected basis, where vector $a$ has expansion*

$$a=e_*{}^*a \tag{2.2.3}$$

*and*

$$b_i=<b,e_i> \tag{2.2.4}$$

Proof. Because $b$ is 1-form, it follows from (2.2.3) that

$$< b, a >=< b, e_*{}^*a >=< b, e_i > {}^ia \tag{2.2.5}$$

(2.2.2) follows from (2.2.5) and (2.2.4). □

Theorem 2.2.4. *Let $\overline{\overline{e}}$ be a basis in ${}_*{}^*D$-vector space $V$. 1-form (2.2.1) is uniquely defined by values (2.2.4) into which 1-form $b$ maps vectors of basis.*

Proof. Statement follows from theorems 2.2.3 and [3]-4.3.3. □

Theorem 2.2.5. *Let $\overline{\overline{e}}$ be a basis in ${}_*{}^*D$-vector space $V$. The set of 1-forms ${}^id$ such that*

$$< {}^jd, e_i >= {}^j\delta_i \tag{2.2.6}$$

*is basis $\overline{\overline{d}}$ of $D_*{}^*$-vector space $\mathcal{L}({}_*{}^*D; V; D)$.*

Proof. Since we assume $b_i = {}^j\delta_i$, then according to theorem 2.2.4 there exists the 1-form ${}^jd$ for given $j$. If we assume that there exists 1-form

$$b = b_*{}^*d = 0$$

then

$$b_j < {}^jd, e_i >= 0$$

According to equation (2.2.6),

$$b_i = b_j\ {}^j\delta_i = 0$$

Therefore, 1-forms ${}^jd$ are linearly independent. □

Definition 2.2.6. Let $V$ be ${}_*{}^*D$-vector space. $D_*{}^*$-vector space

$$V^* = \mathcal{L}({}_*{}^*D; V; D)$$

is called **dual space of ${}_*{}^*D$-vector space** $V$. Let $\overline{\overline{e}}$ be a basis in ${}_*{}^*D$-vector space $V$. basis $\overline{\overline{d}}$ of $D_*{}^*$-vector space $V^*$, satisfying to equation (2.2.6), is called **basis dual to basis** $\overline{\overline{e}}$. □

Theorem 2.2.7. *Let $A$ be passive transformation of basis manifold $\mathcal{B}(V, GL_n{}_*{}^*{}^D)$. Let basis*

$$\overline{\overline{e}}' = \overline{\overline{e}}_*{}^*A \tag{2.2.7}$$

*be image of basis $\overline{\overline{e}}$. Let $B$ be passive transformation of basis manifold $\mathcal{B}(V^*, GL_n{}_*{}^*{}^D)$ such, that basis*

$$\overline{\overline{d}}' = B_*{}^*\overline{\overline{d}} \tag{2.2.8}$$

*is dual to basis. Then*

$$B = A^{-1_*{}^*} \tag{2.2.9}$$

Proof. From equations (2.2.6), (2.2.7), (2.2.8) it follows

$$\begin{aligned} {}^j\delta_i &= < {}^jd', e'_i > \\ &= {}^jB_l\ < {}^ld, e_k >\ {}^kA_i \\ &= {}^jB_l\ {}^l\delta_k\ {}^kA_i \\ &= {}^jB_k\ {}^kA_i \end{aligned} \tag{2.2.10}$$

Equation (2.2.9) follows from equation (2.2.10). □

## 2.3. Twin Representations of Division Ring

Theorem 2.3.1. *We can introduce structure of $D\star$-vector space in any ${}_*{}^*D$-vector space defining $D\star$-product of vector over scalar using equation*

$$mv = m\delta_*{}^*v$$

Proof. We verify directly that the map

$$f : D \to V^\star$$

defined by equation

$$f(m) = m\delta$$

defines $T\star$-representation of the ring $D$. □

We can formulate the theorem 2.3.1 by other way.

Theorem 2.3.2. *Suppose we defined an effective $\star T$-representation $f$ of ring $D$ on the Abelian group $V$. Then we can uniquely define an effective $T\star$-representation $h$ of ring $D$ on the Abelian group $V$ such that diagram*

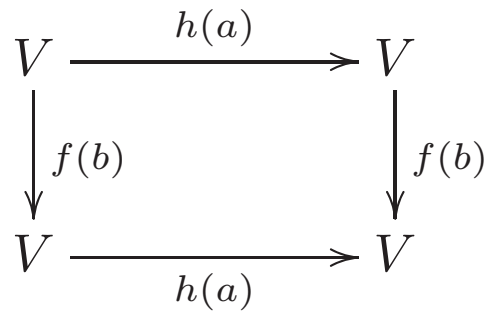

*is commutative for any $a$, $b \in D$.* □

We call representations $f$ and $h$ **twin representations of the division ring** $D$.

Theorem 2.3.3. *In vector space $V$ over division ring $D$ we can define $D\star$-product and $\star D$-product of vector over scalar. According to theorem 2.3.2 these operations satisfy equation*

$$(am)b = a(mb) \tag{2.3.1}$$

Equation (2.3.1) represents **associative law for twin representations**. This allows us writing of such expressions without using of brackets.

Proof. In section 2.1 there is definition of $D\star$-product of ${}_*{}^*D$-linear map $A$ over scalar. According to theorem 2.3.1 $T\star$-representation of division ring $D$ in $D\star$-vector space $\mathcal{L}({}_*{}^*D; V; W)$ can be carried into ${}_*{}^*D$-vector space $W$ according to rule

$$[mA]_*{}^*v = [m\delta]_*{}^*(A_*{}^*v) = m(A_*{}^*v)$$

□

The analogy of the vector space over field goes so far that we can assume an existence of the concept of basis which serves for $D\star$-product and $\star D$-product of vector over skalar.

Theorem 2.3.4. *Basis manifold of ${}_*{}^*D$-vector space $V$ and basis manifold of ${}^*{}_*D$-vector space $V$ are different*

$$\mathcal{B}(V, D_*{}^*) \neq \mathcal{B}(V, {}^*{}_*D)$$

Proof. To prove this theorem we use the standard representation of a matrix. Without loss of generality we prove theorem in coordinate vector space $D^n$.

Let $\overline{\overline{e}} = (e_i = (\delta_i^j), \ i, \ j \in I, \ |I| = n)$ be the set of vectors of vector space $D^n$. $\overline{\overline{e}}$ is evidently basis of ${}_*{}^*D$-vector space and basis of ${}^*{}_*D$-vector space. For arbitrary set of vectors $(f_i, i \in I, |I| = n)$ ${}_*{}^*D$-coordinate matrix

$$f = \begin{pmatrix} f_1^1 & ... & f_1^n \\ ... & ... & ... \\ f_n^1 & ... & f_n^n \end{pmatrix} \tag{2.3.2}$$

relative to basis $\overline{\overline{e}}$ coincide with ${}^*{}_*D$-coordinate matrix relative to basis $\overline{\overline{e}}$.

Let the set of vectors $(f_i, i \in I, |I| = n)$ be basis of ${}_*{}^*D$-vector space. According to theorem [3]-4.9.3 matrix (2.3.2) is ${}_*{}^*$-nonsingular matrix.

Let the set of vectors $(f_i, i \in I, |I| = n)$ be basis of ${}^*{}_*D$-vector space. According to theorem [3]-4.9.3 matrix (2.3.2) is ${}^*{}_*$-nonsingular matrix.

Therefore, if the set of vectors $(f_i, i \in I, |I| = n)$ is basis ${}_*{}^*D$-of vector space and basis of ${}^*{}_*D$-vector space, their coordinate matrix (2.3.2) is ${}_*{}^*$-nonsingular and ${}^*{}_*$-nonsingular matrix. The statement follows from theorem [3]-4.8.9. □

From theorem 2.3.4, it follows that there exists basis $\overline{\overline{e}}$ of $D_*{}^*$-vector space $V$ which is not basis of ${}^*{}_*D$-vector space $V$.

## 2.4. $D$-Vector Space

For many problems we may confine ourselves to considering of $D\star$-vector space or $\star D$-vector space. However there are problems where we forced to reject simple model and consider both structures of vector space at the same time.

Twin representations of division ring $D$ in Abelian group $V$ determine the structure of **$D$-vector space**. From this definition, it follows that if $v, \ w \in V$, then for any $a, \ b, \ c, \ d \in D$

$$avb + cwd \in V \tag{2.4.1}$$

The expression (2.4.1) is called **linear combination of vectors** of $D$-vector space $V$.

Definition 2.4.1. Vectors $a_i, \ i \in I,$ of $D$-vector space $V$ are **linearly independent** if $b^ic^i = 0, \ i = 1, \ ..., \ n,$ follows from the equation

$$b^i a_i c^i = 0$$

Otherwise vectors $a_i$ are **linearly dependent**. □

Definition 2.4.2. We call set of vectors

$$\overline{\overline{e}} = (e_i, i \in I)$$

**basis for $D$-vector space** if vectors $e_i$ are linearly independent and adding to this system any other vector we get a new system which is linearly dependent. □

Theorem 2.4.3. *A basis of $D$-vector space $V$ is basis of $D^*{}_*$-vector space $V$.*

Proof. Let vectors $a_i, \ i \in I$ are linear dependent in $D^*{}_*$-vector space. Then there exists ${}_*$-row

$$\begin{pmatrix} c^i & i \in I \end{pmatrix} \neq 0 \tag{2.4.2}$$

such that

$$c^*{}_*a = 0 \tag{2.4.3}$$

From the equation (2.4.3), it folows that

$$c^i a_i 1 = 0 \tag{2.4.4}$$

According to the definition 2.4.1, vectors $a_i$ are linear dependent in $D$-vector space.

Let $D^*{}_*$-dimension of vector space $V$ is equal to $\boldsymbol{n}$. Let $\overline{\overline{e}}$ be basis of $D^*{}_*$-vector space. Let

$$\begin{pmatrix} 1 & ... & 0 \\ ... & ... & ... \\ 0 & ... & 1 \end{pmatrix}$$

be coordinate matrix of basis $\overline{\overline{e}}$. It is evident that we cannot present vector $e_{\boldsymbol{n}}$ as $D$-linear combination of vectors $e_{\boldsymbol{1}}$, ..., $e_{\boldsymbol{n-1}}$. □

Whether the set of vectors $(e_{\boldsymbol{i}}, \boldsymbol{i} \in \boldsymbol{I})$ is a basis of ${}_*{}^*D$-vector space or a basis of $D^*{}_*$-vector space is of no concern to us. So we will use index in standard format to represent coordinates of vector. Based on the theorem 2.4.3, considering a basis of $D$-vector space $V$ we will consider it as a basis of $D^*{}_*$-vector space $V$.

CHAPTER 3

# Linear Map of Division Ring

## 3.1. Linear Map of Division Ring

According to the remark [3]-8.2.2, we consider division ring $D$ as algebra over field $F \subset Z(D)$.

DEFINITION 3.1.1. Let $D_1$, $D_2$ be division rings. Let $F$ be field such that $F \subset Z(D_1)$, $F \subset Z(D_2)$. Linear map

$$f : D_1 \to D_2$$

of $F$-vector space $D_1$ into $F$-vector space $D_2$ is called **linear map of division ring $D_1$ into division ring $D_2$.** □

According to definition 3.1.1, linear map $f$ of division ring $D_1$ into division ring $D_2$ holds

$$f(a+b) = f(a) + f(b) \quad a, b \in D$$

$$f(pa) = pf(a) \quad p \in F$$

THEOREM 3.1.2. *Let map*

$$f : D_1 \to D_2$$

*is linear map of division ring $D_1$ of characteristic* 0 *into division ring $D_2$ of characteristic* 0. *Then*[3.1]

$$f(nx) = nf(x)$$

*for any integer $n$.*

PROOF. We prove the theorem by induction on $n$. Statement is obvious for $n = 1$ because

$$f(1x) = f(x) = 1f(x)$$

Let statement is true for $n = k$. Then

$$f((k+1)x) = f(kx + x) = f(kx) + f(x) = kf(x) + f(x) = (k+1)f(x)$$

□

[3.1] Let

$$f : D_1 \to D_2$$

be a linear map of division ring $D_1$ of characteristic 2 into division ring $D_2$ of characteristic 3. Then for any $a \in D_1$, $2a = 0$, although $2f(a) \neq 0$. Therefore, if we assume that characteristic of division ring $D_1$ is greater than 0, then must demand that characteristic of division ring $D_1$ equal to characteristic of division ring $D_2$.

Theorem 3.1.3. *Let map*

$$f : D_1 \to D_2$$

*be linear map of division ring $D_1$ of characteristic 0 into division ring $D_2$ of characteristic 0. Then*

$$f(ax) = af(x)$$

*for any rational $a$.*

Proof. Let $a = \frac{p}{q}$. Assume $y = \frac{1}{q}x$. Then according to the theorem 3.1.2

$$f(x) = f(qy) = qf(y) = qf\left(\frac{1}{q}x\right) \tag{3.1.1}$$

From equation (3.1.1) it follows

$$\frac{1}{q}f(x) = f\left(\frac{1}{q}x\right) \tag{3.1.2}$$

From equation (3.1.2) it follows

$$f\left(\frac{p}{q}x\right) = pf\left(\frac{1}{q}x\right) = \frac{p}{q}f(x)$$

□

We cannot extend the statement of theorem 3.1.3 for arbitrary subfield of center $Z(D)$ of division ring $D$.

Theorem 3.1.4. *Let division ring $D$ is algebra over field $F \subset Z(D)$. If $F \neq Z(D)$, then there exists linear map*

$$f : D \to D$$

*which is not linear over field $Z(D)$.*

Proof. To prove the theorem it is enough to consider the complex field $C$ because $C = Z(C)$. Because the complex field is algebra over real field, then the function

$$z \to \overline{z}$$

is linear. However the equation

$$\overline{az} = a\overline{z}$$

is not true. □

Based on theorem 3.1.4, the question arises. Why do we consider linear maps over field $F \neq Z(D)$, if this leads us to sharp expansion of the set of linear maps? The answer to this question is a rich experience of the theory of complex function.

Theorem 3.1.5. *Let maps*

$$f : D_1 \to D_2$$

$$g : D_1 \to D_2$$

*be linear maps of division ring $D_1$ into division ring $D_2$. Then map $f + g$ is linear.*

PROOF. Statement of theorem follows from chain of equations

$$\begin{aligned}(f+g)(x+y)=&f(x+y)+g(x+y)=f(x)+f(y)+g(x)+g(y)\\=&(f+g)(x)+(f+g)(y)\\(f+g)(px)=&f(px)+g(px)=pf(x)+pg(x)=p(f(x)+g(x))\\=&p(f+g)(x)\end{aligned}$$

□

THEOREM 3.1.6. *Let map*

$$f:D_1\rightarrow D_2$$

*be linear map of division ring $D_1$ into division ring $D_2$. Then maps $af$, $fb$, $a$, $b\in R_2$ are linear.*

PROOF. Statement of theorem follows from chain of equations

$$\begin{aligned}(af)(x+y)=&a(f(x+y))=a(f(x)+f(y))=af(x)+af(y)\\=&(af)(x)+(af)(y)\\(af)(px)=&a(f(px))=a(pf(x))=p(af(x))\\=&p(af)(x)\\(fb)(x+y)=&(f(x+y))b=(f(x)+f(y))b=f(x)b+f(y)b\\=&(fb)(x)+(fb)(y)\\(fb)(px)=&(f(px))b=(pf(x))b=p(f(x)b)\\=&p(fb)(x)\end{aligned}$$

□

DEFINITION 3.1.7. Denote $\mathcal{L}(D_1;D_2)$ set of linear maps

$$f:D_1\rightarrow D_2$$

of division ring $D_1$ into division ring $D_2$. □

THEOREM 3.1.8. *We may represent linear map*

$$f:D_1\rightarrow D_2$$

*of division ring $D_1$ into division ring $D_2$ as*

$$f(x)=f_{k\cdot s_k\cdot 0}\ G_k(x)\ f_{k\cdot s_k\cdot 1}\tag{3.1.3}$$

*where $(G_k,k\in K)$ is set of additive maps of division ring $D_1$ into division ring $D_2$.*[3.2] *Expression $f_{k\cdot s_k\cdot p}$, $p=0,\ 1,$ in equation (3.1.3) is called* **component of linear map** $f$.

PROOF. The statement of theorem follows from theorems 3.1.5 and 3.1.6. □

If in the theorem 3.1.8 $|K|=1$, then the equation (3.1.3) has form

$$f(x)=f_{s\cdot 0}\ G(x)\ f_{s\cdot 1}\tag{3.1.4}$$

and map $f$ is called **linear map generated by map** $G$. Map $G$ is called **generator of linear map**.

[3.2] Here and in the following text we assume sum over index that is used in product few times. Equation (3.1.3) is recursive definition and there is hope that it is possible to simplify it.

Theorem 3.1.9. *Let $D_1$, $D_2$ be division rings of characteristic* 0. *Let $F$, $F \subset Z(D_1)$, $F \subset Z(D_2)$, be field. Let $G$ be linear map. Let $\overline{\overline{e}}$ be basis of division ring $D_2$ over field $F$.* **Standard representation of linear map** (3.1.4) *has form*[3.3]

$$f(x) = f_G^{ij}\ e_i\ G(x)\ e_j \tag{3.1.5}$$

*Expression $f_G^{ij}$ in equation* (3.1.5) *is called* **standard component of linear map** *$f$.*

Proof. Components of linear map $f$ have expansion

$$f_{s\cdot p} = f^i_{s\cdot p}\ e_i \tag{3.1.6}$$

relative to basis $\overline{\overline{e}}$. If we substitute (3.1.6) into (3.1.4), we get

$$f(x) = f^i_{s\cdot 0}\ e_i\ G(x)\ f^j_{s\cdot 1}\ e_j \tag{3.1.7}$$

If we substitute expression

$$f_G^{ij} = f^i_{s\cdot 0}\ f^j_{s\cdot 1}$$

into equation (3.1.7) we get equation (3.1.5). □

Theorem 3.1.10. *Let $D_1$, $D_2$ be division rings of characteristic* 0. *Let $F$, $F \subset Z(D_1)$, $F \subset Z(D_2)$, be field Let $G$ be linear map. Let $\overline{\overline{e}}_1$ be basis of division ring $D_1$ over field $F$. Let $\overline{\overline{e}}_2$ be basis of division ring $D_2$ over field $F$. Let $C_{2\cdot}{}^{p}_{kl}$ be structural constants of division ring $D_2$. Then it is possible to represent linear map* (3.1.4) *generated by linear map $G$ as*

$$\begin{aligned} &f(a) = e_{2\cdot j}\ f^j_i\,a^i && f^j_k \in F \\ &a = e_{1\cdot i}\,a^i && a^i \in F \quad a \in D_1 \end{aligned} \tag{3.1.8}$$

$$f^j_i = G^l_i f^{kr}_G C_{2\cdot}{}^{p}_{kl} C_{2\cdot}{}^{j}_{pr} \tag{3.1.9}$$

Proof. Consider map

$$\begin{aligned} G : D_1 \to D_2 \quad & a = e_{1\cdot i}a^i \to G(a) = e_{2\cdot j}G^j_i a^i \\ & a^i \in F \quad G^j_i \in F \end{aligned} \tag{3.1.10}$$

According to the theorem [3]-4.4.3, linear map $f(a)$ relative to bases $\overline{\overline{e}}_1$ and $\overline{\overline{e}}_2$ has form (3.1.8). From equations (3.1.5) and (3.1.10), it follows

$$f(a) = G^l_i a^i f^{kj}_G e_{2\cdot k} e_{2\cdot l} e_{2\cdot j} \tag{3.1.11}$$

From equations (3.1.8) and (3.1.11), it follows

$$e_{2\cdot j}\ f^j_i\ a^i = G^l_i a^i f^{kr}_G e_{2\cdot k} e_{2\cdot l} e_{2\cdot r} = G^l_i a^i f^{kr}_G C_{2\cdot}{}^{p}_{kl} C_{2\cdot}{}^{j}_{pr}\ e_{2\cdot j} \tag{3.1.12}$$

Since vectors $e_{2\cdot r}$ are linear independent over field $F$ and values $a^k$ are arbitrary, then equation (3.1.9) follows from equation (3.1.12). □

Considering map

$$f : D \to D \tag{3.1.13}$$

we assume $G(x) = x$.

[3.3] Representation of linear map of of division ring using components of linear map is ambiguous. We can increase or decrease number of summands using algebraic operations. Since dimension of division ring $D_2$ over field $F$ is finite, standard representation of linear map guarantees finiteness of set of items in the representation of map.

Theorem 3.1.11. *Let $D$ be division ring of characteristic 0. Linear map (3.1.13) has form*

$$f(x) = f_{s\cdot 0}\ x\ f_{s\cdot 1} \tag{3.1.14}$$

Theorem 3.1.12. *Let $D$ be division ring of characteristic 0. Let $\overline{\overline{e}}$ be the basis of division ring $D$ over field $F \subset Z(D)$. Standard representation of linear map (3.1.14) of division ring has form*

$$f(x) = f^{ij}\ e_i x e_j \tag{3.1.15}$$

Theorem 3.1.13. *Let $D$ be division ring of characteristic 0. Let $\overline{\overline{e}}$ be basis of division ring $D$ over field $F \subset Z(D)$. Then it is possible to represent linear map (3.1.13) as*

$$f(a) = e_j f^j_i a^i \qquad f^j_k \in F \tag{3.1.16}$$

$$a = e_i a^i \qquad a^i \in F \quad a \in D$$

$$f^j_i = f^{kr} C^p_{ki} C^j_{pr} \tag{3.1.17}$$

Theorem 3.1.14. *Consider matrix*

$$\mathcal{C} = \left(\mathcal{C}^{\cdot j}_{i \cdot kr}\right) = \left(C^p_{ki} C^j_{pr}\right) \tag{3.1.18}$$

*whose rows and columns are indexed by ${}^{\cdot j}_{i}$ and ${}_{\cdot kr}$, respectively. If $\det \mathcal{C} \neq 0$, then, for given coordinates of linear transformation $f^j_i$, the system of linear equations (3.1.17) with standard components of this transformation $f^{kr}$ has the unique solution. If $\det \mathcal{C} = 0$, then the equation*

$$\operatorname{rank}\left(\mathcal{C}^{\cdot j}_{i \cdot kr} \quad f^j_i\right) = \operatorname{rank} \mathcal{C} \tag{3.1.19}$$

*is the condition for the existence of solutions of the system of linear equations (3.1.17). In such case the system of linear equations (3.1.17) has infinitely many solutions and there exists linear dependence between values $f^j_i$.*

Proof. Equation (3.1.14) is special case of equation (3.1.4) when $G(x) = x$. Theorem 3.1.12 is special case of theorem 3.1.9 when $G(x) = x$. Theorem 3.1.13 is special case of theorem 3.1.10 when $G(x) = x$. The statement of the theorem 3.1.14 is corollary of the theory of linear equations over field. □

Theorem 3.1.15. *Standard components of the identity map have the form*

$$f^{kr} = \delta^k_0 \delta^r_0 \tag{3.1.20}$$

Proof. The equation (3.1.20) is corollary of the equation

$$x = e_0\ x\ e_0$$

Let us show that the standard components (3.1.20) of a linear transformation satisfy to the equation

$$\delta^j_i = f^{kr} C^p_{ki} C^j_{pr} \tag{3.1.21}$$

which follows from the equation (3.1.17) if $f = \delta$. From equations (3.1.20), (3.1.21), it follows that

$$\delta^j_i = C^p_{0i} C^j_{p0} \tag{3.1.22}$$

The equation (3.1.22) is true, because, from equations

$$e_j e_0 = e_0\ e_j = e_j$$

it follows that

$$C^{j}_{0r} = \delta^{j}_{r} \quad C^{j}_{r0} = \delta^{j}_{r}$$

If $\det \mathcal{C} \neq 0$, then the solution (3.1.20) is unique. If $\det \mathcal{C} = 0$, then the system of linear equations (3.1.21) has infinitely many solutions. However, we are looking for at least one solution. □

THEOREM 3.1.16. *If $\det \mathcal{C} \neq 0$, then standard components of the zero map*

$$z : A \to A \quad z(x) = 0$$

*are defined uniquely and have form $z^{ij} = 0$. If $\det \mathcal{C} = 0$, then the set of standard components of the zero map forms a vector space.*

PROOF. The theorem is true because standard components $z^{ij}$ are solution of homogeneous system of linear equations

$$0 = z^{kr} C^{p}_{ki} C^{j}_{pr}$$

□

REMARK 3.1.17. Consider equation

$$a^{kr}\ e_k\ x\ e_r = b^{kr}\ e_k\ x\ e_r \tag{3.1.23}$$

From the theorem 3.1.16, it follows that only when condition $\det \mathcal{C} \neq 0$ is true, from the equation (3.1.23), it follows that

$$a^{kr} = b^{kr} \tag{3.1.24}$$

Otherwise, we must assume equality

$$a^{kr} = b^{kr} + z^{kr} \tag{3.1.25}$$

Despite this, in case $\det \mathcal{C} = 0$, we also use standard representation because in general it is very hard to show the set of linear independent vectors. If we want to define operation over linear maps in standard representation, then as well as in the case of the theorem 3.1.15 we choose one element from the set of possible representations. □

THEOREM 3.1.18. *Expression*

$$f^{r}_{k} = f^{ij} C^{p}_{ik} C^{r}_{pj}$$

*is tensor over field $F$*

$${}_{i}f'^{j} = {}_{i}A^{k}\ {}_{k}f^{l}\ {}_{l}A^{-1j} \tag{3.1.26}$$

PROOF. $D$-linear map has form (3.1.16) relative to basis $\overline{\overline{e}}$. Let $\overline{\overline{e}}'$ be another basis. Let

$$e'_i = e_j A^{j}_{i} \tag{3.1.27}$$

be transformation map basis $\overline{\overline{e}}$ to basis $\overline{\overline{e}}'$. Since linear map $f$ is the same, then

$$f(x) = e'_l f'^{l}_{k} x'^{k} \tag{3.1.28}$$

Let us substitute [3]-(8.2.8), (3.1.27) into equation (3.1.28)

$$f(x) = e_j A^{j}_{l} f'^{l}_{k} A^{-1\cdot k}_{\ \ i} x^{i} \tag{3.1.29}$$

Because vectors $e_j$ are linear independent and components of vector $x^i$ are arbitrary, the equation (3.1.26) follows from equation (3.1.29). Therefore, expression $f^r_k$ is tensor over field $F$. □

Definition 3.1.19. The set

$$\ker f = \{x \in D_1 : f(x) = 0\}$$

is called **kernel of linear map**

$$f : D_1 \to D_2$$

of division ring $D_1$ into division ring $D_2$. □

Theorem 3.1.20. *Kernel of linear map*

$$f : D_1 \to D_2$$

*is subgroup of additive group of division ring* $D_1$*.*

Proof. Let $a, b \in \ker f$. Then

$$f(a) = 0$$
$$f(b) = 0$$
$$f(a+b) = f(a) + f(b) = 0$$

Therefore, $a + b \in \ker f$. □

Definition 3.1.21. The linear map

$$f : D_1 \to D_2$$

of division ring $D_1$ into division ring $D_2$ is called **singular**, when

$$\ker f \neq \{0\}$$

□

Theorem 3.1.22. *Let $D$ be division ring of characteristic 0. Let $\overline{\overline{e}}$ be basis of division ring $D$ over center $Z(D)$ of division ring $D$. Let*

$$f : D \to D \quad f(x) = f_{s\cdot 0}\ x\ f_{s\cdot 1} \tag{3.1.30}$$
$$= f^{ij}\ e_i x e_j \tag{3.1.31}$$
$$g : D \to D \quad g(x) = g_{t\cdot 0}\ x\ g_{t\cdot 1} \tag{3.1.32}$$
$$= g^{ij}\ e_i x e_j \tag{3.1.33}$$

*be linear maps of division ring $D$. Map*

$$h(x) = gf(x) = g(f(x)) \tag{3.1.34}$$

*is linear map*

$$h(x) = h_{ts\cdot 0}\ x\ h_{ts\cdot 1} \tag{3.1.35}$$
$$= h^{pr}\ e_p\ x\ e_r \tag{3.1.36}$$

*where*

$$h_{ts\cdot 0} = g_{t\cdot 0}\ f_{s\cdot 0} \tag{3.1.37}$$
$$h_{ts\cdot 1} = f_{s\cdot 1}\ g_{t\cdot 1} \tag{3.1.38}$$
$$h^{pr} = g^{ij} f^{kl} C^p_{ik} C^r_{lj} \tag{3.1.39}$$

Proof. Map (3.1.34) is linear because

$$h(x+y)=g(f(x+y))=g(f(x)+f(y))=g(f(x))+g(f(y))=h(x)+h(y)$$
$$h(ax)=g(f(ax))=g(af(x))=ag(f(x))=ah(x)$$

If we substitute (3.1.30) and (3.1.32) into (3.1.34), we get

$$h(x)=g_{t\cdot 0}\ f(x)\ g_{t\cdot 1}=g_{t\cdot 0}\ f_{s\cdot 0}\ x\ f_{s\cdot 1}\ g_{t\cdot 1} \tag{3.1.40}$$

Comparing (3.1.40) and (3.1.35), we get (3.1.37), (3.1.38).

If we substitute (3.1.31) and (3.1.33) into (3.1.34), we get

$$\begin{aligned} h(x)=&g^{ij}\ e_i\ f(x)\ e_j\\ =&g^{ij}\ e_i\ f^{kl}\ e_k\ x\ e_l e_j\\ =&g^{ij} f^{kl} C^p_{ik} C^r_{lj}\ e_p\ x\ e_r \end{aligned} \tag{3.1.41}$$

Comparing (3.1.41) and (3.1.36), we get (3.1.39). □

Definition 3.1.23. Let commutative ring $P$ be subring of center $Z(R)$ of ring $R$. Map

$$f:R\rightarrow R$$

of ring $R$ is called **projective over commutative ring** $P$, if

$$f(px)=f(x)$$

for any $p\in P$. Set

$$Px=\{px:p\in P,x\in R\}$$

is called **direction** $x$ **over commutative ring**[3.4] $P$. □

Example 3.1.24. If map $f$ of ring $R$ is linear over commutative ring $P$, then map

$$g(x)=x^{-1}f(x)$$

is projective over commutative ring $P$. □

## 3.2. Polylinear Map of Division Ring

Definition 3.2.1. Let $R_1$, ..., $R_n$, $P$ be rings of characteristic 0. Let $S$ be module over ring $P$. Let $F$ be commutative ring which is for any $i$ is subring of center of ring $R_i$. Map

$$f:R_1\times...\times R_n\rightarrow S$$

is called **polylinear over commutative ring** $F$, if

$$f(p_1,...,p_i+q_i,...,p_n)=f(p_1,...,p_i,...,p_n)+f(p_1,...,q_i,...,p_n)$$
$$f(a_1,...,ba_i,...,a_n)=bf(a_1,...,a_i,...,a_n)$$

for any $i$, $1\le i\le n$, and any $p_i$, $q_i\in R_i$, $b\in F$. Let us denote $\mathcal{L}(R_1,...,R_n;S)$ set of polylinear maps of rings $R_1$, ..., $R_n$ into module $S$. □

[3.4] Direction over commutative ring $P$ is subset of ring $R$. However we denote direction $Px$ by element $x\in R$ when this does not lead to ambiguity. We tell about direction over commutative ring $Z(R)$ when we do not show commutative ring $P$ explicitly.

Theorem 3.2.2. *Let $D$ be division ring of characteristic* 0. *Polylinear map*

$$f : D^n \to D, d = f(d_1, ..., d_n) \tag{3.2.1}$$

*has form*

$$d = f^n_{s\cdot 0}\ \sigma_s(d_1)\ f^n_{s\cdot 1}\ ...\ \sigma_s(d_n)\ f^n_{s\cdot n} \tag{3.2.2}$$

*$\sigma_s$ is a transposition of set of variables $\{d_1, ..., d_n\}$*

$$\sigma_s = \begin{pmatrix} d_1 & ... & d_n \\ \sigma_s(d_1) & ... & \sigma_s(d_n) \end{pmatrix}$$

Proof. We prove statement by induction on $n$.

When $n = 1$ the statement of theorem is corollary of theorem 3.1.11. In such case we may identify[3.5]

$$f^1_{s\cdot p} = f_{s\cdot p} \quad p = 0, 1$$

Let statement of theorem be true for $n = k-1$. Then it is possible to represent map (3.2.1) as

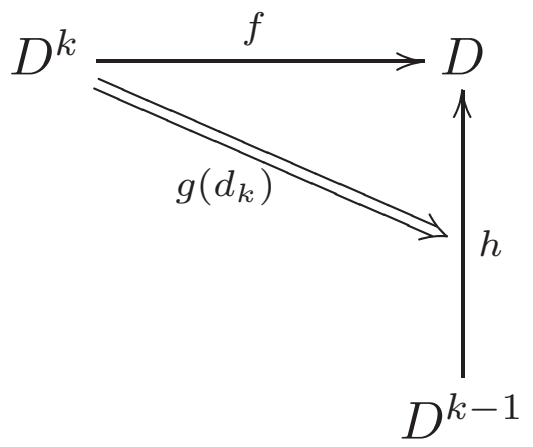

$$d = f(d_1, ..., d_k) = g(d_k)(d_1, ..., d_{k-1})$$

According to statement of induction pollinear map $h$ has form

$$d = h^{k-1}_{t\cdot 0}\ \sigma_t(d_1)\ h^{k-1}_{t\cdot 1}\ ...\ \sigma_t(d_{k-1})\ h^{k-1}_{t\cdot k-1}$$

According to construction $h = g(d_k)$. Therefore, expressions $h_{t\cdot p}$ are functions of $d_k$. Since $g(d_k)$ is linear map of $d_k$, then only one expression $h_{t\cdot p}$ is linear map of $d_k$, and rest expressions ${}_{t\cdot q}h$ do not depend on $d_k$.

Without loss of generality, assume $p = 0$. According to equation (3.1.14) for given $t$

$$h^{k-1}_{t\cdot 0} = g_{tr\cdot 0}\ d_k\ g_{tr\cdot 1}$$

Assume $s = tr$. Let us define transposition $\sigma_s$ according to rule

$$\sigma_s = \sigma_{tr} = \begin{pmatrix} d_k & d_1 & ... & d_{k-1} \\ d_k & \sigma_t(d_1) & ... & \sigma_t(d_{k-1}) \end{pmatrix}$$

Suppose

$$f^k_{tr\cdot q+1} = h^{k-1}_{t\cdot q} \quad q = 1, ..., k-1$$

$$f^k_{tr\cdot q} = g_{tr\cdot q} \qquad q = 0, 1$$

[3.5] In representation (3.2.2) we will use following rules.

- If range of any index is set consisting of one element, then we will omit corresponding index.
- If $n = 1$, then $\sigma_s$ is identical transformation. We will not show such transformation in the expression.

We proved step of induction. □

DEFINITION 3.2.3. Expression $f^n_{s\cdot p}$ in equation (3.2.2) is called **component of polylinear map** $f$. □

THEOREM 3.2.4. *Let $D$ be division ring of characteristic 0. Let $\overline{\overline{e}}$ be basis in division ring $D$ over field $F \subset Z(D)$.* **Standard representation of polylinear map** *of division ring has form*

$$f(d_1, ..., d_n) = f_t^{i_0 ... i_n}\ e_{i_0}\ \sigma_t(d_1)\ e_{i_1}\ ...\ \sigma_t(d_n)\ e_{i_n} \tag{3.2.3}$$

*Index $t$ enumerates every possible transpositions $\sigma_t$ of the set of variables $\{d_1, ..., d_n\}$. Expression $f_t^{i_0 ... i_n}$ in equation (3.2.3) is called* **standard component of polylinear map** $f$.

PROOF. Components of polylinear map $f$ have expansion

$$f^n_{s\cdot p} = e_i f^{ni}_{s\cdot p} \tag{3.2.4}$$

relative to basis $\overline{\overline{e}}$. If we substitute (3.2.4) into (3.2.2), we get

$$d = f^{nj_1}_{s\cdot 0}\ e_{j_1}\ \sigma_s(d_1)\ f^{nj_2}_{s\cdot 1}\ e_{j_2}\ ...\ \sigma_s(d_n)\ f^{nj_n}_{s\cdot n}\ e_{j_n} \tag{3.2.5}$$

Consider expression

$$f_t^{j_0 ... j_n} = f^{nj_1}_{s\cdot 0}\ ... f^{nj_n}_{s\cdot n} \tag{3.2.6}$$

The right-hand side is supposed to be the sum of the terms with the index $s$, for which the transposition $\sigma_s$ is the same. Each such sum has a unique index $t$. If we substitute expression (3.2.6) into equation (3.2.5) we get equation (3.2.3). □

THEOREM 3.2.5. *Let $\overline{\overline{e}}$ be basis of division ring $D$ over field $F \subset Z(D)$. Polylinear map (3.2.1) can be represented as $D$-valued form of degree $n$ over field $F \subset Z(D)$*[3.6]

$$f(a_1, ..., a_n) = a_1^{i_1} ... a_n^{i_n} f_{i_1 ... i_n} \tag{3.2.7}$$

*where*

$$a_j = e_i a_j^i$$

$$f_{i_1 ... i_n} = f(e_{i_1}, ..., e_{i_n}) \tag{3.2.8}$$

*and values ${}_{i_1 ... i_n}f$ are coordinates of $D$-valued covariant tensor over field $F$.*

PROOF. According to the definition 3.2.1, the equation (3.2.7) follows from the chain of equations

$$f(a_1, ..., a_n) = f(e_{i_1} a_1^{i_1}, ..., e_{i_n} a_n^{i_n}) = a_1^{i_1} ... a_n^{i_n} f(e_{i_1}, ..., e_{i_n})$$

Let $\overline{\overline{e}}'$ be another basis. Let

$$e'_i = e_j A_i^j \tag{3.2.9}$$

[3.6] We proved the theorem by analogy with theorem in [2], p. 107, 108

be transformation, mapping basis $\overline{\overline{e}}$ into basis $\overline{\overline{e}}'$. From equations (3.2.9) and (3.2.8) it follows

$$f'_{i_1...i_n} = f(e'_{i_1}, ..., e'_{i_n}) \tag{3.2.10}$$
$$= f(e_{j_1}A^{j_1}_{i_1}, ..., e'_{j_n}A^{j_n}_{i_n})$$
$$= A^{j_1}_{i_1}...A^{j_n}_{i_n} f(e_{j_1}, ..., e_{j_n})$$
$$= A^{j_1}_{i_1}...A^{j_n}_{i_n} f_{j_1...j_n}$$

From equation (3.2.10) the tensor law of transformation of coordinates of polylinear map follows. From equation (3.2.10) and theorem [3]-8.2.3 it follows that value of the map $f(a_1, ..., a_n)$ does not depend from choice of basis. □

Polylinear map (3.2.1) is **symmetric**, if

$$f(d_1, ..., d_n) = f(\sigma(d_1), ..., \sigma(d_n))$$

for any transposition $\sigma$ of set $\{d_1, ..., d_n\}$.

THEOREM 3.2.6. *If polyadditive map $f$ is symmetric, then*

$$f_{i_1,...,i_n} = f_{\sigma(i_1),...,\sigma(i_n)} \tag{3.2.11}$$

PROOF. Equation (3.2.11) follows from equation

$$a_1^{i_1}...a_n^{i_n} f_{i_1...i_n} = f(a_1, ..., a_n)$$
$$= f(\sigma(a_1), ..., \sigma(a_n))$$
$$= a_1^{i_1}...a_n^{i_n} f_{\sigma(i_1)...\sigma(i_n)}$$

□

Polylinear map (3.2.1) is **skew symmetric**, if

$$f(d_1, ..., d_n) = |\sigma| f(\sigma(d_1), ..., \sigma(d_n))$$

for any transposition $\sigma$ of set $\{d_1, ..., d_n\}$. Here

$$|\sigma| = \begin{cases} 1 & \text{transposition } \sigma \text{ even} \\ -1 & \text{transposition } \sigma \text{ odd} \end{cases}$$

THEOREM 3.2.7. *If polylinear map $f$ is skew symmetric, then*

$$f_{i_1,...,i_n} = |\sigma| f_{\sigma(i_1),...,\sigma(i_n)} \tag{3.2.12}$$

PROOF. Equation (3.2.12) follows from equation

$$a_1^{i_1}...a_n^{i_n} f_{i_1...i_n} = f(a_1, ..., a_n)$$
$$= |\sigma| f(\sigma(a_1), ..., \sigma(a_n))$$
$$= a_1^{i_1}...a_n^{i_n} |\sigma| f_{\sigma(i_1)...\sigma(i_n)}$$

□

THEOREM 3.2.8. *Coordinates of the polylinear over field $F$ map (3.2.1) and its components relative basis $\overline{\overline{e}}$ satisfy to the equation*

$$f_{j_1...j_n} = f_t^{i_0...i_n} C^{k_1}_{i_0\sigma_t(j_1)} C^{l_1}_{k_1 i_1} ... B^{k_n}_{l_{n-1}\sigma_t(j_n)} C^{l_n}_{k_n i_n} e_{l_n} \tag{3.2.13}$$

$$f^p_{j_1...j_n} = f_t^{i_0...i_n} C^{k_1}_{i_0\sigma_t(j_1)} C^{l_1}_{k_1 i_1} ... C^{k_n}_{l_{n-1}\sigma_t(j_n)} C^{p}_{k_n i_n} \tag{3.2.14}$$

Proof. In equation (3.2.3), we assume

$$d_i = e_{j_i} d_i^{j_i}$$

Then equation (3.2.3) gets form

$$\begin{aligned} f(d_1, ..., d_n) &= f_t^{i_0...i_n}\ e_{i_0}\sigma_t(d_1^{j_1}e_{j_1})e_{i_1}...\sigma_t(d_n^{j_n}e_{j_n})e_{i_n} \\ &= d_1^{j_1}...d_n^{j_n} f_t^{i_0...i_n} e_{i_0}\sigma_t(e_{j_1})e_{i_1}...\sigma_t(e_{j_n})e_{i_n} \\ &= d_1^{j_1}...d_n^{j_n} f_t^{i_0...i_n} C^{k_1}_{i_0\sigma_t(j_1)} C^{l_1}_{k_1 i_1} \\ &\quad ...C^{k_n}_{l_{n-1}\sigma_t(j_n)} C^{l_n}_{k_n i_n} e_{l_n} \end{aligned} \tag{3.2.15}$$

From equation (3.2.7) it follows that

$$f(a_1, ..., a_n) = e_p f^p_{i_1...i_n} a_1^{i_1}...a_n^{i_n} \tag{3.2.16}$$

Equation (3.2.13) follows from comparison of equations (3.2.15) and (3.2.7). Equation (3.2.14) follows from comparison of equations (3.2.15) and (3.2.16). □

CHAPTER 4

# Linear Map of $D$-Vector Spaces

## 4.1. Linear Map of $D$-Vector Spaces

Considering linear map of $D$-vector spaces we assume that division ring $D$ is finite dimensional algebra over field $F$.

Definition 4.1.1. Let field $F$ be subring of center $Z(D)$ of division ring $D$. Suppose $V$ and $W$ are $D$-vector spaces. We call map

$$A:V\rightarrow W$$

of $D$-vector space $V$ into $D$-vector space $W$ **linear map** if

$$A(x+y)=A(x)+A(y)\quad x,y\in V$$
$$A(px)=pA(x)\qquad\qquad p\in F$$

Let us denote $\mathcal{L}(D;V;W)$ set of linear maps

$$A:V\rightarrow W$$

of $D$-vector space $V$ into $D$-vector space $W$. □

It is evident that linear map of $D_*{}^*$-vector space as well linear map of ${}^*{}_*D$-vector space are linear maps. Set of morphisms of $D$-vector space is wider then set of morphisms of $D_*{}^*$-vector space. To consider linear map of vector space, we will follow method used in section [3]-4.4.

Theorem 4.1.2. *Let $D$ be division ring of characteristic 0. Let $V$, $W$ be $D$-vector spaces. Linear map*

$$A:V\rightarrow W$$

*relative to basis $\overline{\overline{e}}_V$ of ${}_*{}^*D$-vector space $V$ and basis $\overline{\overline{e}}_W$ of ${}_*{}^*D$-vector space $W$ has form*

$$A(v)=e_{W\cdot j}\ A^j_i(v^i)\qquad v=e_{V*}{}^*v \tag{4.1.1}$$

*where $A^j_i(v^i)$ linearly depends on one variable $v^i$ and does not depends on the rest of coordinates of vector $v$.*

Proof. According to definition 4.1.1

$$A(v)=A(e_{V*}{}^*v)=A\left(\sum_i e_{V\cdot i}\ v^i\right)=\sum_i A(e_{V\cdot i}\ v^i) \tag{4.1.2}$$

For any given $i$ vector $A(e_{V\cdot i}\ v^i)\in W$ has only expansion

$$A(e_{V\cdot i}\ v^i)=e_{W\cdot j}\ A^j_i(v^i)\qquad A(e_{V\cdot i}\ v^i)=e_{W*}{}^*A_i(v^i) \tag{4.1.3}$$

relative to basis $\overline{\overline{e}}_W$ of ${}_*{}^*D$-vector space. Let us substitute (4.1.3) into (4.1.2). We get (4.1.1). □

DEFINITION 4.1.3. The linear map

$$A_i^j : D \to D$$

is called **partial linear map** of variable $v^i$. □

We can write linear map as product of matrices

$$A(v) = \begin{pmatrix} e_{W\cdot 1} & ... & e_{W\cdot m} \end{pmatrix} {}_*{}^* \begin{pmatrix} A_i^1(v^i) \\ ... \\ A_i^m(v^i) \end{pmatrix} \tag{4.1.4}$$

Let us define product of matrices

$$\begin{pmatrix} A_1^1 & ... & A_n^1 \\ ... & ... & ... \\ A_1^m & ... & A_n^m \end{pmatrix} {}_\circ{}^\circ \begin{pmatrix} v^1 \\ ... \\ v^n \end{pmatrix} = \begin{pmatrix} A_i^1(v^i) \\ ... \\ A_i^m(v^i) \end{pmatrix} \tag{4.1.5}$$

where $A = \begin{pmatrix} A_i^j \end{pmatrix}$ is matrix of partial linear maps. Using the equation (4.1.5), we can write the equation (4.1.4) in the form

$$A(v) = \begin{pmatrix} e_{W\cdot 1} & ... & e_{W\cdot m} \end{pmatrix} {}_*{}^* \begin{pmatrix} A_1^1 & ... & A_n^1 \\ ... & ... & ... \\ A_1^m & ... & A_n^m \end{pmatrix} {}_\circ{}^\circ \begin{pmatrix} v^1 \\ ... \\ v^n \end{pmatrix} \tag{4.1.6}$$

THEOREM 4.1.4. *Let $D$ be division ring of characteristic 0. Linear map*

$$A : v \in V \to w \in W \quad w = A(v) \tag{4.1.7}$$

*relative to basis $\overline{\overline{e}}_V$ of $D$-vector space $V$ and basis $\overline{\overline{e}}_W$ of $D$-vector space $W$*[4.1] *has form*

$$v = e_{V*}{}^*v$$

$$w = e_{W*}{}^*w \tag{4.1.8}$$

$$w^j = A_i^j(v^i) = A_{s\cdot 0\cdot}{}_i^j\ v^i\ A_{s\cdot 1\cdot}{}_i^j \tag{4.1.9}$$

PROOF. According to theorem 4.1.2 we can write linear map $A(v)$ as (4.1.1). Because for given indexes $i$, $j$ partial linear map $A_i^j(v^i)$ is linear with respect to variable $v^i$, then according to (3.1.14) it is possible to represent expression $A_i^j(v^i)$ as

$$A_i^j(v^i) = A_{s\cdot 0\cdot}{}_i^j\ v^i\ A_{s\cdot 1\cdot}{}_i^j \tag{4.1.10}$$

[4.1] Coordinate representation of map (4.1.7) depends on choice of basis. Equations change form if, for instance, we choose ${}^*{}_*D$-basis $\overline{\overline{r}}_*$ in $D$-vector space $W$.

where index $s$ is numbering items. Range of index $s$ depends on indexes $i$ and $j$. Combining equations (4.1.2) and (4.1.10), we get

$$
\begin{aligned}
A(v) &= e_{W\cdot j}A^{j}_{i}(v^{i}) &&= e_{W\cdot j}A_{s\cdot 0\cdot}{}^{j}_{i}\ v^{i}\ A_{s\cdot 1\cdot}{}^{j}_{i} \\
A(v) &= e_{W*}{}^{*}A_{i}(v^{i}) &&= e_{W*}{}^{*}(A_{s\cdot 0\cdot i}\ v^{i}\ A_{s\cdot 1\cdot i})
\end{aligned}
\tag{4.1.11}
$$

In equation (4.1.11), we also summarize on the index $i$. Equation (4.1.9) follows from comparison of equations (4.1.8) and (4.1.11). □

Definition 4.1.5. Expression $A_{s\cdot p\cdot}{}^{j}_{i}$ in equation (4.1.3) is called **component of linear map** $A$. □

Theorem 4.1.6. *Let $D$ be division ring of characteristic 0. Let $\overline{\overline{e}}_{V}$ be a basis of $D$-vector space $V$, $\overline{\overline{e}}_{U}$ be a basis of $D$-vector space $U$, and $\overline{\overline{e}}_{W}$ be a basis of $D$-vector space $W$. Suppose diagram of maps*

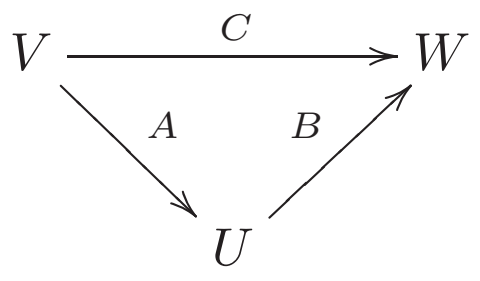

*is commutative diagram where linear map $A$ has presentation*

$$
u = A(v) = e_{U\cdot j}A^{j}_{i}(v^{i}) = e_{U\cdot j}A_{s\cdot 0\cdot}{}^{j}_{i}\ v^{i}\ A_{s\cdot 1\cdot}{}^{j}_{i}
\tag{4.1.12}
$$

*relative to selected bases and linear map $B$ has presentation*

$$
w = B(u) = e_{W\cdot k}B^{k}_{j}(u^{j}) = e_{W\cdot k}B_{t\cdot 0\cdot}{}^{k}_{j}\ u^{j}\ B_{t\cdot 1\cdot}{}^{k}_{j}
\tag{4.1.13}
$$

*relative to selected bases. Then map $C$ is linear map and has presentation*

$$
w = C(v) = e_{W\cdot k}C^{k}_{i}(v^{i}) = e_{W\cdot k}C_{u\cdot 0\cdot}{}^{k}_{i}\ v^{i}\ C_{u\cdot 1\cdot}{}^{k}_{i}
\tag{4.1.14}
$$

*relative to selected bases, where*[4.2]

$$
\begin{aligned}
&C^{k}_{i}(v^{i}) && = B^{k}_{j}(A^{j}_{i}(v^{i})) \\
C_{u\cdot 0\cdot}{}^{k}_{i} &= C_{ts\cdot 0\cdot}{}^{k}_{i} && = B_{t\cdot 0\cdot}{}^{k}_{j}\ A_{s\cdot 0\cdot}{}^{j}_{i} \\
C_{u\cdot 1\cdot}{}^{k}_{i} &= C_{ts\cdot 1\cdot}{}^{k}_{i} && = A_{s\cdot 1\cdot}{}^{j}_{i}\ B_{t\cdot 1\cdot}{}^{k}_{j}
\end{aligned}
\tag{4.1.15}
$$

Proof. The map $C$ is linear map because

$$
\begin{aligned}
C(a+b) &= B(A(a+b)) \\
&= B(A(a)+A(b)) \\
&= B(A(a))+B(A(b)) \\
&= C(a)+C(b) \\
C(ab) &= B(A(ab)) = B(aA(b)) \\
&= aB(A(b)) = aC(b) \\
& a,b \in V \\
& a \in F
\end{aligned}
$$

Equation (4.1.14) follows from substituting (4.1.12) into (4.1.13). □

---

[4.2] Index $u$ appeared composite index, $u = st$. However it is possible that some items in (4.1.15) may be summed together.

Theorem 4.1.7. *For linear map $A$ there exists linear map $B$ such, that*

$$A(axb) = B(x)$$

$$B_{s\cdot 0\cdot}{}^{j}_{i} = A_{s\cdot 0\cdot}{}^{j}_{i}\ a$$

$$B_{s\cdot 1\cdot}{}^{j}_{i} = b\ A_{s\cdot 1\cdot}{}^{j}_{i}$$

Proof. The map $B$ is linear map because

$$B(x+y) = A(a(x+y)b) = A(axb+ayb) = A(axb)+A(ayb) = B(x)+B(y)$$

$$B(cx) = A(acxb) = A(caxb) = cA(axb) = cB(x)$$

$$x, y \in V, c \in F$$

According to equation (4.1.9)

$$B_{s\cdot 0\cdot}{}^{j}_{i}\ v^i\ B_{s\cdot 1\cdot}{}^{j}_{i} = A_{s\cdot 0\cdot}{}^{j}_{i}\ (a\ v^i\ b)A_{s\cdot 1\cdot}{}^{j}_{i} = (A_{s\cdot 0\cdot}{}^{j}_{i}\ a)v^i\ (b\ A_{s\cdot 1\cdot}{}^{j}_{i})$$

□

Theorem 4.1.8. *Let $D$ be division ring of characteristic* 0. *Let*

$$A : V \to W$$

*linear map of $D$-vector space $V$ into $D$-vector space $W$. Then* $A(0) = 0$.

Proof. Corollary of equation

$$A(a+0) = A(a) + A(0)$$

□

Definition 4.1.9. The set

$$\ker f = \{x \in V : f(x) = 0\}$$

is called **kernel of linear map**

$$A : V \to W$$

of $D$-vector space $V$ into $D$-vector space $W$. □

Definition 4.1.10. The linear map

$$A : V \to W$$

of $D$-vector space $V$ into $D$-vector space $W$ is called **singular**, if

$$\ker A = V$$

□

## 4.2. Polylinear Map of $D$-Vector Spaces

Definition 4.2.1. Let field $F$ be subring of center $Z(D)$ of division ring $D$ of characteristic 0. Suppose $V_1$, ..., $V_n$, $W_1$, ..., $W_m$ are $D$-vector spaces. We call map

$$A : V_1 \times ... \times V_n \to W_1 \times ... \times W_m$$
$$w_1 \times ... \times w_m = A(v_1, ..., v_n) \tag{4.2.1}$$

**polylinear map** of $\times$-$D$-vector space $V_1\times...\times V_n$ into $\times$-$D$-vector space $W_1\times ...\times W_m$, if

$$A(p_1,...,p_i+q_i,...,p_n)=A(p_1,...,p_i,...,p_n)+A(p_1,...,q_i,...,p_n)$$
$$A(p_1,...,ap_i,...,p_n)=aA(p_1,...,p_i,...,p_n)$$
$$1\le i\le n\quad p_i,q_i\in V_i\quad a\in F$$

$\square$

Definition 4.2.2. Let us denote $\mathcal{L}(D;V_1,...,V_n;W_1,...,W_m)$ set of polylinear maps of $\times$-$D$-vector space $V_1\times...\times V_n$ into $\times$-$D$-vector space $W_1\times...\times W_m$. $\square$

Theorem 4.2.3. *Let $D$ be division ring of characteristic $0$. For each $k\in K=[1,n]$ let $\overline{\overline{e}}_{V_k}$ be basis in $D$-vector space $V_k$ and*

$$v_k=v_k{}^*{}_*e_{V_k}\quad v_k\in V_k$$

*For each $l$, $1\le l\le m$, let $\overline{\overline{e}}_{W_l}$ be basis in $D$-vector space $W_l$ and*

$$w_l=w_l{}^*{}_*e_{W_l}\quad w_l\in W_l$$

*Polylinear map (4.2.1) relative to bases $\overline{\overline{e}}_{V_1}$, ..., $\overline{\overline{e}}_{V_n}$, $\overline{\overline{e}}_{W_1}$, ..., $\overline{\overline{e}}_{W_m}$ has form*

$$\begin{aligned}w_l^{\boldsymbol{j}}&=A_{l\cdot}{}^{\boldsymbol{j}}_{\boldsymbol{i}_1...\boldsymbol{i}_n}(v_1^{\boldsymbol{i}_1},...,v_n^{\boldsymbol{i}_n})\\&=A^n_{s\cdot 0\cdot l\cdot}{}^{\boldsymbol{j}}_{\boldsymbol{i}_1...\boldsymbol{i}_n}\ \sigma_s(v_1^{\boldsymbol{i}_1})\ A^n_{s\cdot 1.l\cdot}{}^{\boldsymbol{j}}_{\boldsymbol{i}_1...\boldsymbol{i}_n}\ ...\ \sigma_s(v_n^{\boldsymbol{i}_n})\ A^n_{s\cdot n\cdot l\cdot}{}^{\boldsymbol{j}}_{\boldsymbol{i}_1...\boldsymbol{i}_n}\end{aligned}\tag{4.2.2}$$

*Range $S$ of index $s$ depends on values of indexes $\boldsymbol{i}_1$, ..., $\boldsymbol{i}_n$. $\sigma_s$ is a transposition of set of variables $\{v_1^{\boldsymbol{i}_1},...,v_n^{\boldsymbol{i}_n}\}$.*

Proof. Since we may consider map $A$ into $\times$-$D$-vector space $W_1\times...\times W_m$ componentwise, then we may confine to considering of map

$$A_l:V_1\times...\times V_n\to W_l\quad w_l=A_l(v_1,...,v_n)\tag{4.2.3}$$

We prove statement by induction on $n$.

When $n=1$ the statement of theorem is statement of theorem 4.1.4. In such case we may identify[4.3]

$$A^1_{s\cdot p\cdot}{}^{\boldsymbol{j}}_{\boldsymbol{i}}=A_{s\cdot p\cdot}{}^{\boldsymbol{j}}_{\boldsymbol{i}}\quad p=0,1$$

[4.3] In representation (4.2.2) we will use following rules.

- If range of any index is set consisting of one element, then we will omit corresponding index.
- If $n=1$, then $\sigma_s$ is identical transformation. We will not show such transformation in the expression.

Let statement of theorem be true for $n = k-1$. Then it is possible to represent map (4.2.3) as

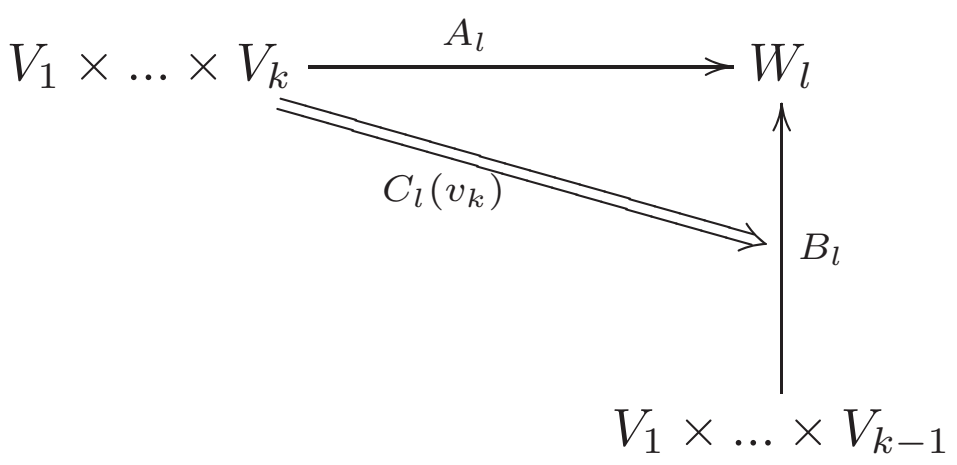

$$w_l = A_l(v_1, ..., v_k) = C_l(v_k)(v_1, ..., v_{k-1})$$

According to statement of induction polylinear map $B_l$ has form

$$w_l^{j} = B^{k-1}{}_{t\cdot 0\cdot l\cdot}{}^{j}{}_{i_1...i_{k-1}}\ \sigma_t(v_1^{i_1})\ B^{k-1}{}_{t\cdot 1.l\cdot}{}^{j}{}_{i_1...i_{k-1}}\ ...\ \sigma_t(v_{k-1}^{i_{k-1}})\ B^{k-1}{}_{t\cdot k-1.l\cdot}{}^{j}{}_{i_1...i_{k-1}}$$

According to construction $B_l = C_l(v_k)$. Therefore, expressions $B^{k-1}{}_{t\cdot p\cdot l\cdot}{}^{j}{}_{i_1...i_{k-1}}$ are functions of $v_k$. Since $C_l(v_k)$ is linear map of $v_k$, then only one expression $B^{k-1}{}_{t\cdot p\cdot l\cdot}{}^{j}{}_{i_1...i_{k-1}}$ is linear map $v_k$, and rest expressions $B^{k-1}{}_{t\cdot q\cdot l\cdot}{}^{j}{}_{i_1...i_{k-1}}$ do not depend on $v_k$.

Without loss of generality, assume $p = 0$. According to theorem 4.1.4

$$B^{k-1}{}_{t\cdot 0\cdot l\cdot}{}^{j}{}_{i_1...i_{k-1}} = C^{k}{}_{tr\cdot 0\cdot l\cdot}{}^{j}{}_{i_k i_1...i_{k-1}}\ v_k^{i_k}\ C^{k}{}_{tr\cdot 1.l\cdot}{}^{j}{}_{i_k i_1...i_{k-1}}$$

Assume $s = tr$. Let us define transposition $\sigma_s$ according to rule

$$\sigma_s = \sigma(tr) = \begin{pmatrix} v_k^{i_k} & v_1^{i_1} & ... & v_{k-1}^{i_{k-1}} \\ v_k^{i_k} & \sigma_t(v_1^{i_1}) & ... & \sigma_t(v_{k-1}^{i_{k-1}}) \end{pmatrix}$$

Suppose

$$A^{k}{}_{tr\cdot q+1.l\cdot}{}^{j}{}_{i_k i_1...i_{k-1}} = B^{k-1}{}_{t\cdot q\cdot l\cdot}{}^{j}{}_{i_1...i_{k-1}} \qquad q = 1, ..., k-1$$
$$A^{k}{}_{tr\cdot q\cdot l\cdot}{}^{j}{}_{i_k i_1...i_{k-1}} = C^{k}{}_{tr\cdot q\cdot l\cdot}{}^{j}{}_{i_k i_1...i_{k-1}} \qquad q = 0, 1$$

We proved step of induction. □

Definition 4.2.4. Expression $A^{n}{}_{s\cdot p\cdot l\cdot}{}^{j}{}_{i_1...i_n}$ in equation (4.2.2) is called **component of polylinear map** $A$. □

# Part 2

# Calculus over Division Ring

CHAPTER 5

# Differentiable Maps

## 5.1. Topological Division Ring

Definition 5.1.1. Division ring $D$ is called **topological division ring**[5.1] if $D$ is topological space and the algebraic operations defined in $D$ are continuous in the topological space $D$. □

According to definition, for arbitrary elements $a, b \in D$ and for arbitrary neighborhoods $W_{a-b}$ of the element $a-b$, $W_{ab}$ of the element $ab$ there exists neighborhoods $W_a$ of the element $a$ and $W_b$ of the element $b$ such that $W_a - W_b \subset W_{a-b}$, $W_a W_b \subset W_{ab}$. For any $a \neq 0$ and for arbitrary neighborhood $W_{a^{-1}}$ there exists neighborhood $W_a$ of the element $a$, satisfying the condition $W_a^{-1} \subset W_{a^{-1}}$.

Definition 5.1.2. **Absolute value on division ring** $D$[5.2] is a map

$$d \in D \to |d| \in R$$

which satisfies the following axioms

- $|a| \geq 0$
- $|a| = 0$ if, and only if, $a = 0$
- $|ab| = |a|\ |b|$
- $|a+b| \leq |a| + |b|$

Division ring $D$, endowed with the structure defined by a given absolute value on $D$, is called **valued division ring**. □

Invariant distance on additive group of division ring $D$

$$d(a, b) = |a - b|$$

defines topology of metric space, compatible with division ring structure of $D$.

Definition 5.1.3. Let $D$ be valued division ring. Element $a \in D$ is called **limit of a sequence** $\{a_n\}$

$$a = \lim_{n \to \infty} a_n$$

if for every $\epsilon \in R$, $\epsilon > 0$, there exists positive integer $n_0$ depending on $\epsilon$ and such, that $|a_n - a| < \epsilon$ for every $n > n_0$. □

Theorem 5.1.4. *Let $D$ be valued division ring of characteristic* 0 *and let $d \in D$. Let $a \in D$ be limit of a sequence $\{a_n\}$. Then*

$$\lim_{n \to \infty} (a_n d) = ad$$

$$\lim_{n \to \infty} (d a_n) = da$$

[5.1] I made definition according to definition from [6], chapter 4

[5.2] I made definition according to definition from [4], IX, §3.2

Proof. Statement of the theorem is trivial, however I give this proof for completeness sake. Since $a \in D$ is limit of the sequence $\{a_n\}$, then according to definition 5.1.3 for given $\epsilon \in R$, $\epsilon > 0$, there exists positive integer $n_0$ such, that $|a_n - a| < \epsilon/|d|$ for every $n > n_0$. According to definition 5.1.2 the statement of theorem follows from inequalities

$$|a_n d - ad| = |(a_n - a)d| = |a_n - a||d| < \epsilon/|d||d| = \epsilon$$
$$|da_n - da| = |d(a_n - a)| = |d||a_n - a| < |d|\epsilon/|d| = \epsilon$$

for any $n > n_0$. □

Definition 5.1.5. Let $D$ be valued division ring. The sequence $\{a_n\}$, $a_n \in D$ is called **fundamental** or **Cauchy sequence**, if for every $\epsilon \in R$, $\epsilon > 0$, there exists positive integer $n_0$ depending on $\epsilon$ and such, that $|a_p - a_q| < \epsilon$ for every $p$, $q > n_0$. □

Definition 5.1.6. Valued division ring $D$ is called **complete** if any fundamental sequence of elements of division ring $D$ converges, i.e. has limit in division ring $D$. □

So far everything was good. Ring of characteristic 0 contains ring of integers. Division ring of characteristic 0 contains field of rational numbers. However this does not mean that complete division ring $D$ of characteristic 0 contains field of real numbers. For instance, absolute value $|x| = 1$ determines discrete topology on valued division ring and is not interesting for us, because any fundamental sequence is constant map.

Later on, speaking about valued division ring of characteristic 0, we will assume that homeomorphism of field of rational numbers $Q$ into division ring $D$ is defined.

Theorem 5.1.7. *Complete division ring $D$ of characteristic 0 contains as subfield an isomorphic image of the field $R$ of real numbers. It is customary to identify it with $R$.*

Proof. Consider fundamental sequence of rational numbers $\{p_n\}$. Let $p'$ be limit of this sequence in division ring $D$. Let $p$ be limit of this sequence in field $R$. Since immersion of field $Q$ into division ring $D$ is homeomorphism, then we may identify $p' \in D$ and $p \in R$. □

Theorem 5.1.8. *Let $D$ be complete division ring of characteristic 0 and let $d \in D$. Then any real number $p \in R$ commute with $d$.*

Proof. Let us represent real number $p \in R$ as fundamental sequence of rational numbers $\{p_n\}$. Statement of theorem follows from chain of equations

$$pd = \lim_{n\to\infty}(p_n d) = \lim_{n\to\infty}(dp_n) = dp$$

based on statement of theorem 5.1.4. □

Theorem 5.1.9. *Let $D$ be complete division ring of characteristic 0. Then field of real numbers $R$ is subfield of center $Z(D)$ of division ring $D$.*

Proof. Corollary of theorem 5.1.8. □

Definition 5.1.10. Let $D$ be complete division ring of characteristic 0. Set of elements $d \in D$, $|d| = 1$, is called **unit sphere in division ring $D$**. □

Definition 5.1.11. Let $D_1$ be complete division ring of characteristic 0 with absolute value $|x|_1$. Let $D_2$ be complete division ring of characteristic 0 with absolute value $|x|_2$. A map

$$f : D_1 \to D_2$$

is called **continuous**, if for every as small as we please $\epsilon > 0$ there exist such $\delta > 0$, that

$$|x' - x|_1 < \delta$$

implies

$$|f(x') - f(x)|_2 < \epsilon$$

□

Theorem 5.1.12. *Let $D$ be complete division ring of characteristic 0. Since index $s$ in expansion (3.1.14) of linear map*

$$f : D \to D$$

*belongs to finite range, then linear map $f$ is continuous.*

Proof. Suppose $x' = x + a$. Then

$$f(x') - f(x) = f(x+a) - f(x) = f(a) = {}_{s\cdot 0}f\ a\ {}_{s\cdot 1}f$$

$$|f(x') - f(x)| = |{}_{s\cdot 0}f\ a\ {}_{s\cdot 1}f| < (|{}_{s\cdot 0}f|\ |{}_{s\cdot 1}f|)|a|$$

Let us denote $F = |{}_{s\cdot 0}f|\ |{}_{s\cdot 1}f|$. Then

$$|f(x') - f(x)| < F|a|$$

Suppose $\epsilon > 0$ and let us assume $a = \frac{\epsilon}{F}e$. Then $\delta = |a| = \frac{\epsilon}{F}$. According to definition 5.1.11, linear map $f$ is continuous. □

Similarly, since index $s$ belongs to countable range, then for continuity of additive map $f$ it is necessary that series $|{}_{s\cdot 0}f|\ |{}_{s\cdot 1}f|$ converges. Since index $s$ belongs to continues range, then for continuity of additive map $f$ it is necessary that integral $\int |{}_{s\cdot 0}f|\ |{}_{s\cdot 1}f| ds$ exist.

Definition 5.1.13. Let

$$f : D_1 \to D_2$$

be map of complete division ring $D_1$ of characteristic 0 with absolute value $|x|_1$ into complete division ring $D_2$ of characteristic 0 with absolute value $|y|_2$. Value

$$\|f\| = \sup \frac{|f(x)|_2}{|x|_1} \tag{5.1.1}$$

is called **norm of map** $f$. □

Theorem 5.1.14. *Let*

$$f : D_1 \to D_2$$

*be linear map of complete division ring $D_1$ into complete division ring $D_2$. Since $\|f\| < \infty$, then map $f$ is continuous.*

Proof. Because map $f$ is additive, then according to definition 5.1.13

$$|f(x)-f(y)|_2=|f(x-y)|_2\le\|f\|\ |x-y|_1$$

Let us assume arbitrary $\epsilon>0$. Assume $\delta=\frac{\epsilon}{\|f\|}$. Then

$$|f(x)-f(y)|_2\le\|f\|\ \delta=\epsilon$$

follows from inequality

$$|x-y|_1<\delta$$

According to definition 5.1.11, map $f$ is continuous. □

Theorem 5.1.15. *Let map*

$$f:D_1\to D_2$$

*be continuous linear map of division ring $D_1$ of characteristic 0 into division ring $D_2$ of characteristic 0. Then*

$$f(ax)=af(x)$$

*for any real number $a$.*

Proof. According to the theorem 5.1.9, real field is subfield of the center $Z(D)$ of division ring $D$. Let us represent real number $p\in R$ as fundamental sequence of rational numbers $\{p_n\}$. The statement of theorem follows from theorems 3.1.3, 5.1.4 and chain of equations

$$f(ph)=\lim_{n\to\infty}f(p_nh)=\lim_{n\to\infty}(p_nf(h))=pf(h) \tag{5.1.2}$$

□

Theorem 5.1.16. *Let $D$ be complete division ring of characteristic 0. Either continuous map $f$ of division ring which is projective over field $P$, does not depend on direction over field $P$, or value $f(0)$ is not defined.*

Proof. According to definition 3.1.23, map $f$ is constant on direction $Pa$. Since $0\in Pa$, then we may assume

$$f(0)=f(a)$$

based on continuity. However this leads to uncertainty of value of map $f$ in direction 0, when map $f$ has different values for different directions $a$. □

If projective over field $R$ map $f$ is continuous, then we say that **map $f$ is continuous in direction over field** $R$. Since for any $a\in D$, $a\ne 0$ we may choose $a_1=|a|^{-1}a$, $f(a_1)=f(a)$, then it is possible to make definition more accurate.

Definition 5.1.17. Let $D$ be complete division ring of characteristic 0. Projective over field $R$ function $f$ is continuous in direction over field $R$, if for every as small as we please $\epsilon>0$ there exists such $\delta>0$, that

$$|x'-x|_1<\delta\quad |x'|_1=|x|_1=1$$

implies

$$|f(x')-f(x)|_2<\epsilon$$

□

THEOREM 5.1.18. *Let $D$ be complete division ring of characteristic $0$. Projective over field $R$ function $f$ is continuous in direction over field $R$ iif this function is continuous on unit sphere of division ring $D$.*

PROOF. Corollary of definitions 5.1.11, 3.1.23, 5.1.17. □

## 5.2. Differentiable Map of Division Ring

DEFINITION 5.2.1. Let $D$ be valued division ring.[5.3] Function

$$f : D \to D$$

is called **$D\star$-differentiable in the Fréchet sense** on the set $U \subset D$[5.4], if at every point $x \in U$ the increment of the function $f$ can be represented as

$$f(x+h) - f(x) = h\frac{df(x)}{d_* x} + o(h) \tag{5.2.1}$$

where $o$ is such continuous map

$$o : D \to D$$

that

$$\lim_{h \to 0} \frac{|o(h)|}{|h|} = 0 \tag{5.2.2}$$

□

According to definition 5.2.1, **the Fréchet $D\star$-derivative of map $f$ at point** $x$

$$\frac{df(x)}{d_* x} \in {}_*\mathcal{L}(D, D)$$

generates a homomorphism

$$\Delta f = \Delta x \frac{df(x)}{d_* x}$$

which maps increment of argument into increment of function.

If we multiply both sides of equation (5.2.1) by $h^{-1}$, then we get equation

$$h^{-1}(f(x+h) - f(x)) = \frac{df(x)}{d_* x} + h^{-1}o(h) \tag{5.2.3}$$

Alternative definition of the Fréchet $D\star$-derivative

$$\frac{df(x)}{d_* x} = \lim_{h \to 0} (h^{-1}(f(x+h) - f(x))) \tag{5.2.4}$$

follows from equations (5.2.3) and (5.2.2).

---

[5.3] I made definition according to definition [1]-3.1.1, page 177.

[5.4] Like in remark [3]-4.4.5 we may define $D\star$-differentiability of map

$$f : S \to D$$

of valued division ring $S$ into valued division ring $D$ according to rule

$$f(x+h) - f(x) = F(h)\frac{df(x)}{d_* x} + o(h)$$

where

$$F : S \to D$$

is homomorphism of division rings. However, based on isomorphism theorem, we may confine ourselves by case of maps of division ring $D$ into $D$.

In contrast to derivative of map over field $F$ the Fréchet $D\star$-derivative of map of division ring $D$ is not linear map. Indeed, the Fréchet $D\star$-derivative is $\star D$-linear map; however, in general, it is not $D\star$-linear map. In fact,

$$f(x+h)a - f(x)a = h\frac{df}{d_*x}a + o(h)$$

$$af(x+h) - af(x) = ah\frac{df}{d_*x} + o(h) \tag{5.2.5}$$

However, in general,

$$ah\frac{df}{d_*x} \neq ha\frac{df}{d_*x}$$

Similar problem arises when we differentiate product of maps. According to definition 5.2.1

$$f(x+h)g(x+h) - f(x)g(x) = h\frac{df(x)g(x)}{d_*x} + o(h) \tag{5.2.6}$$

It is possible to represent expression in left side as

$$\begin{aligned} &f(x+h)g(x+h) - f(x)g(x)\\ =&f(x+h)g(x+h) - f(x)g(x+h) + f(x)g(x+h) - f(x)g(x)\\ =&(f(x+h) - f(x))g(x+h) + f(x)(g(x+h) - g(x)) \end{aligned} \tag{5.2.7}$$

According to definition 5.2.1, it is possible to represent (5.2.7) as

$$\begin{aligned} &f(x+h)g(x+h) - f(x)g(x)\\ =&\left(h\frac{df(x)}{d_*x} + o(h)\right)\left(g(x) + h\frac{dg(x)}{d_*x} + o(h)\right)\\ &+f(x)\left(h\frac{dg(x)}{d_*x} + o(h)\right) \end{aligned} \tag{5.2.8}$$

Since, in general,

$$f(x)h\frac{dg}{d_*x} \neq hf(x)\frac{dg}{d_*x}$$

then it is natural to anticipate that map $f(x)g(x)$ is not $D\star$-differentiable in the Fréchet sense.

Thus definition the Fréchet $D\star$-derivative is extremely restricted and does not satisfy standard definition of differentiation. To find solution of problem of differentiation, Consider this problem on the other hand.

Example 5.2.2. Consider increment of map $f(x) = x^2$.

$$\begin{aligned} f(x+h) - f(x) &= (x+h)^2 - x^2\\ &= xh + hx + h^2\\ &= xh + hx + o(h) \end{aligned}$$

As can be easily seen, the component of the increment of the function $f(x) = x^2$ that is linearly dependent on the increment of the argument, is of the form

$$xh + hx$$

Since product is non commutative, we cannot represent increment of map $f(x+h) - f(x)$ as $Ah$ or $hA$ where $A$ does not depend on $h$. It results in the unpredictable behavior of the increment of the function $f(x) = x^2$ when the increment of the argument converges to 0. However, since infinitesimal $h$ is infinitesimal like $h = ta$, $a \in D$, $t \in R$, $t \to 0$, the answer becomes more definite

$$(xa + ax)t$$

□

DEFINITION 5.2.3. Let $D$ be complete division ring of characteristic 0.[5.5] The map

$$f : D \to D$$

is called **differentiable in the Gâteaux sense** on the set $U \subset D$, if at every point $x \in U$ the increment of the function $f$ can be represented as

$$f(x + a) - f(x) = \partial f(x)(a) + o(a) = \frac{\partial f(x)}{\partial x}(a) + o(a) \tag{5.2.9}$$

where **the Gâteaux derivative** $\partial f(x)$ **of map** $f$ is linear map of increment $a$ and $o : D \to D$ is such continuous map that

$$\lim_{a \to 0} \frac{|o(a)|}{|a|} = 0$$

□

REMARK 5.2.4. According to definition 5.2.3 for given $x$, the Gâteaux derivative $\partial f(x) \in \mathcal{L}(D; D)$. Therefore, the Gâteaux derivative of map $f$ is map

$$\partial f : D \to \mathcal{L}(D; D)$$

Expressions $\partial f(x)$ and $\dfrac{\partial f(x)}{\partial x}$ are different notations for the same function. We will use notation $\dfrac{\partial f(x)}{\partial x}$ to underline that this is the Gâteaux derivative with respect to variable $x$. □

THEOREM 5.2.5. *It is possible to represent the Gâteaux derivative of map $f$ as*

$$\partial f(x)(a) = \frac{\partial_{s \cdot 0} f(x)}{\partial x} \; a \; \frac{\partial_{s \cdot 1} f(x)}{\partial x} \tag{5.2.10}$$

PROOF. Corollary of definition 5.2.3 and theorem 3.1.11. □

DEFINITION 5.2.6. Expression $\dfrac{\partial_{s \cdot p} f(x)}{\partial x}$, $p = 0, 1,$ is called **component of the Gâteaux derivative of map** $f(x)$. □

THEOREM 5.2.7. *Let $D$ be division ring of characteristic 0. Let function*

$$f : D \to D$$

*be differentiable in the Gâteaux sense Then*

$$\partial f(x)(ah) = a\partial f(x)(h) \tag{5.2.11}$$

*for any real number $a$.*

PROOF. Corollary of theorem 5.1.15 and definition 5.2.3. □

[5.5] I made definition according to definition [1]-3.1.2, page 177.

Combining equation (5.2.11) and definition 5.2.3, we get known definition of the Gâteaux differential

$$\partial f(x)(a)=\lim_{t\to 0,\ t\in R}(t^{-1}(f(x+ta)-f(x))) \tag{5.2.12}$$

Definitions of the Gâteaux derivative (5.2.9) and (5.2.12) are equivalent. Using this equivalence we tell that map $f$ is called differentiable in the Gâteaux sense on the set $U\subset D$, if at every point $x\in U$ the increment of the function $f$ can be represented as

$$f(x+ta)-f(x)=t\partial f(x)(a)+o(t) \tag{5.2.13}$$

where $o:R\to D$ is such continuous map that

$$\lim_{t\to 0}\frac{|o(t)|}{|t|}=0$$

Since infinitesimal $a$ in the equation (5.2.10) is differential $dx$, then equation (5.2.10) becomes definition of **the Gâteaux differential**

$$\partial f(x)(dx)=\frac{\partial_{s\cdot 0}f(x)}{\partial x}dx\frac{\partial_{s\cdot 1}f(x)}{\partial x} \tag{5.2.14}$$

Theorem 5.2.8. *Let $D$ be division ring of characteristic $0$. Let $\overline{\overline{e}}$ be basis of division ring $D$ over center $Z(D)$ of division ring $D$.* **Standard representation of the Gâteaux derivative** (5.2.10) *of map*

$$f:D\to D$$

*has form*

$$\partial f(x)(a)=\frac{\partial^{ij}f(x)}{\partial x}\ e_i\ a\ e_j \tag{5.2.15}$$

*Expression* $\frac{\partial^{ij}f(x)}{\partial x}$ *in equation* (5.2.15) *is called* **standard component of the Gâteaux derivative** *of map $f$.*

Proof. Statement of theorem is corollary of theorem 3.1.12. □

Theorem 5.2.9. *Let $D$ be division ring of characteristic $0$. Let $\overline{\overline{e}}$ be basis of division ring $D$ over field $F\subset Z(D)$. Then it is possible to represent the Gâteaux differential of map*

$$f:D\to D$$

*as*

$$\partial f(x)(a)=a^i\ \frac{\partial f^j}{\partial x^i}e_j \tag{5.2.16}$$

*where $a\in D$ has expansion*

$$a=a^i\ e_i\qquad a^i\in F$$

*relative to basis $\overline{\overline{e}}$ and Jacobian matrix of map $f$ has form*

$$\frac{\partial f^j}{\partial x^i}=\frac{\partial^{kr}f(x)}{\partial x}\ C^p_{ki}\ C^j_{pr} \tag{5.2.17}$$

Proof. Statement of theorem is corollary of theorem 3.1.13. □

To find structure similar to derivative in commutative calculus we need extract derivative from differential of map. For this purpose we need take the increment of argument outside brackets.[5.6] It is possible to take $a$ outside brackets using equations

$$\frac{\partial_{s\cdot 0}f(x)}{\partial x}a\frac{\partial_{s\cdot 1}f(x)}{\partial x}=\frac{\partial_{s\cdot 0}f(x)}{\partial x}a\frac{\partial_{s\cdot 1}f(x)}{\partial x}a^{-1}a$$
$$\frac{\partial_{s\cdot 0}f(x)}{\partial x}a\frac{\partial_{s\cdot 1}f(x)}{\partial x}=aa^{-1}\frac{\partial_{s\cdot 0}f(x)}{\partial x}a\frac{\partial_{s\cdot 1}f(x)}{\partial x}$$

Definition 5.2.10. Let $D$ be complete division ring of characteristic 0 and $a\in D$. We define **the Gâteaux $D\star$-derivative** $\frac{\partial f(x)(a)}{\partial_* x}$ **of map** $f:D\rightarrow D$ using equation

$$\partial f(x)(a)=a\frac{\partial f(x)(a)}{\partial_* x}\tag{5.2.18}$$

We define **the Gâteaux $\star D$-derivative** $\frac{\partial f(x)(a)}{{}_*\partial x}$ **of map** $f:D\rightarrow D$ using equation

$$\partial f(x)(a)=\frac{\partial f(x)(a)}{{}_*\partial x}a\tag{5.2.19}$$

□

For representation of the Gâteaux $D\star$-derivative we also use notation

$$\frac{\partial f(x)(a)}{\partial_* x}=\frac{\partial}{\partial_* x}f(x)(a)=\partial_* f(x)(a)$$

For representation of the Gâteaux $\star D$-derivative we also use notation

$$\frac{\partial f(x)(a)}{{}_*\partial x}=\frac{\partial}{{}_*\partial x}f(x)(a)={}_*\partial f(x)(a)$$

Based on duality principle [3]-4.3.8, we will study $D\star$-derivative, keeping in mind that dual statement is true for $\star D$-derivative.

Because product is not commutative, concept of fraction is not defined in division ring. We need explicitly show from which side denominator has an effect on numerator. Symbol $*$ in denominator of fraction implements this function. According to definition

$$\frac{\partial f}{\partial_* x}=(\partial_* x)^{-1}\partial f$$
$$\frac{\partial f}{{}_*\partial x}=\partial f({}_*\partial x)^{-1}$$

Thus, we may represent the Gâteaux $D\star$-derivative of function $f$ as ratio of change of function to change of argument. We do not extend this remark to components of the Gâteaux differential. In this case we consider notation $\frac{\partial}{\partial_* x}$ and $\frac{\partial}{{}_*\partial x}$ not as fraction, but as symbol of operator.

[5.6] It is possible when all $\frac{{}_{s\cdot 0}\partial f(x)}{\partial x}=e$ or all $\frac{{}_{s\cdot 1}\partial f(x)}{\partial x}=e$. Therefore definition [1]-3.1.2, page 177, leads us either to definition of the Fréchet $D\star$-derivative, or to definition of the Fréchet $\star D$-derivative.

Then we can write equation (5.2.18) as

$$\partial f(x)(dx) = dx(\partial_* x)^{-1}\partial f(x)(dx) \tag{5.2.20}$$

It is easy to see that we wrote $D\star$-differential of variable $x$ in denominator of fraction in equation (5.2.20), however we multiply fraction by differential of variable $x$.

THEOREM 5.2.11. *Let $D$ be complete division ring of characteristic $0$. The Gâteaux $D\star$-derivative is projective over field of real numbers $R$.*

PROOF. Corollary of theorems 5.2.7 and example 3.1.24. □

From theorem 5.2.11 it follows

$$\frac{\partial f(x)(ra)}{\partial_* x} = \frac{\partial f(x)(a)}{\partial_* x} \tag{5.2.21}$$

for every $r \in R$, $r \neq 0$ and $a \in D$, $a \neq 0$. Therefore the Gâteaux $D\star$-derivative is well defined in direction $a$ over field $R$, $a \in D$, $a \neq 0$, and does not depend on the choice of value in this direction.

THEOREM 5.2.12. *Let $D$ be complete division ring of characteristic $0$ and $a \neq 0$. The Gâteaux $D\star$-derivative and the Gâteaux $\star D$-derivative of map $f$ of division ring $D$ are bounded by relationship*

$$\frac{\partial f(x)(a)}{{}_*\partial x} = a\frac{\partial f(x)(a)}{\partial_* x}a^{-1} \tag{5.2.22}$$

PROOF. From equations (5.2.18) and (5.2.19) it follows

$$\frac{\partial f(x)(a)}{{}_*\partial x} = \partial f(x)(a)a^{-1} = a\frac{\partial f(x)(a)}{\partial_* x}a^{-1}$$

□

THEOREM 5.2.13. *Let $D$ be complete division ring of characteristic $0$. The Gâteaux derivative satisfies to relationship*

$$\partial(f(x)g(x))(a) = \partial f(x)(a)\ g(x) + f(x)\ \partial g(x)(a) \tag{5.2.23}$$

PROOF. Equation (5.2.23) follows from chain of equations

$$\begin{aligned}
\partial(f(x)g(x))(a) &= \lim_{t\to 0}(t^{-1}(f(x+ta)g(x+ta) - f(x)g(x)))\\
&= \lim_{t\to 0}(t^{-1}(f(x+ta)g(x+ta) - f(x)g(x+ta)))\\
&+ \lim_{t\to 0}(t^{-1}(f(x)g(x+ta) - f(x)g(x)))\\
&= \lim_{t\to 0}(t^{-1}(f(x+ta) - f(x)))g(x)\\
&+ f(x)\lim_{t\to 0}(t^{-1}(g(x+ta) - g(x)))
\end{aligned}$$

based on definition (5.2.12). □

THEOREM 5.2.14. *Let $D$ be complete division ring of characteristic $0$. Suppose the Gâteaux derivative of map $f : D \to D$ has expansion*

$$\partial f(x)(a) = \frac{\partial_{s\cdot 0} f(x)}{\partial x}a\frac{\partial_{s\cdot 1} f(x)}{\partial x} \tag{5.2.24}$$

*Suppose the Gâteaux derivative of map $g : D \to D$ has expansion*

$$\partial g(x)(a) = \frac{\partial_{t\cdot 0} g(x)}{\partial x}a\frac{\partial_{t\cdot 1} g(x)}{\partial x} \tag{5.2.25}$$

*Components of the Gâteaux derivative of map $f(x)g(x)$ have form*

$$\frac{\partial_{s\cdot 0}f(x)g(x)}{\partial x}=\frac{\partial_{s\cdot 0}f(x)}{\partial x}\qquad\frac{\partial_{t\cdot 0}f(x)g(x)}{\partial x}=f(x)\frac{\partial_{t\cdot 0}g(x)}{\partial x}\tag{5.2.26}$$

$$\frac{\partial_{s\cdot 1}f(x)g(x)}{\partial x}=\frac{\partial_{s\cdot 1}f(x)}{\partial x}g(x)\qquad\frac{\partial_{t\cdot 1}f(x)g(x)}{\partial x}=\frac{\partial_{t\cdot 1}g(x)}{\partial x}\tag{5.2.27}$$

Proof. Let us substitute (5.2.24) and (5.2.25) into equation (5.2.23)

$$\begin{aligned}\partial(f(x)g(x))(a)&=\partial f(x)(a)\ g(x)+f(x)\ \partial g(x)(a)\\&=\frac{\partial_{s\cdot 0}f(x)}{\partial x}a\frac{\partial_{s\cdot 1}f(x)}{\partial x}g(x)+f(x)\frac{\partial_{t\cdot 0}g(x)}{\partial x}a\frac{\partial_{t\cdot 1}g(x)}{\partial x}\end{aligned}\tag{5.2.28}$$

Based (5.2.28), we define equations (5.2.26), (5.2.27). □

Theorem 5.2.15. *Let $D$ be complete division ring of characteristic 0. The Gâteaux $D\star$-derivative satisfy to relationship*

$$\frac{\partial f(x)g(x)}{\partial_* x}(a)=\frac{\partial f(x)(a)}{\partial_* x}g(x)+a^{-1}f(x)a\frac{\partial g(x)(a)}{\partial_* x}\tag{5.2.29}$$

Proof. Equation (5.2.29) follows from chain of equations

$$\begin{aligned}\frac{\partial f(x)g(x)}{\partial_* x}(a)&=a^{-1}\partial f(x)g(x)(a)\\&=a^{-1}(\partial f(x)(a)g(x)+f(x)\partial g(x)(a))\\&=a^{-1}\partial f(x)(a)g(x)+a^{-1}f(x)aa^{-1}\partial g(x)(a)\\&=\frac{\partial f(x)(a)}{\partial_* x}g(x)+a^{-1}f(x)a\frac{\partial g(x)(a)}{\partial_* x}\end{aligned}$$

□

Theorem 5.2.16. *Let $D$ be complete division ring of characteristic 0. Either the Gâteaux $D\star$-derivative does not depend on direction, or the Gâteaux $D\star$-derivative in direction 0 is not defined.*

Proof. Statement of theorem is corollary of theorem 5.2.11 and theorem 5.1.16. □

Theorem 5.2.17. *Let $D$ be complete division ring of characteristic 0. Let unit sphere of division ring $D$ be compact. If the Gâteaux $D\star$-derivative $\dfrac{\partial f(x)(a)}{\partial_* x}$ exists in point $x$ and is continuous in direction over field $R$, then there exist norm $\|\partial f(x)\|$ of the Gâteaux $D\star$-differential.*

Proof. From definition 5.2.10 it follows

$$|\partial f(x)(a)|=|a|\left|\frac{\partial f(x)(a)}{\partial_* x}\right|\tag{5.2.30}$$

From theorems 5.1.18, 5.2.11 it follows, that the Gâteaux $D\star$-derivative is continuous on unit sphere. Since unit sphere is compact, then range the Gâteaux $D\star$-derivative of function $f$ at point $x$ is bounded

$$\left|\frac{\partial f(x)(a)}{\partial_* x}\right|<F=\sup\left|\frac{\partial f(x)(a)}{\partial_* x}\right|$$

According to definition 5.1.13

$$\|\partial f(x)\|=F$$

□

Theorem 5.2.18. *Let $D$ be complete division ring of characteristic 0. Let unit sphere of division ring $D$ be compact. If the Gâteaux $D\star$-derivative $\frac{\partial f(x)(a)}{\partial_* x}$ exists in point $x$ and is continuous in direction over field $R$, then function $f$ is continuous at point $x$.*

Proof. From theorem 5.2.17 it follows

$$|\partial f(x)(a)| \le \|\partial f(x)\||a| \tag{5.2.31}$$

From (5.2.9), (5.2.31) it follows

$$|f(x+a)-f(x)| < |a|\ \|\partial f(x)\| \tag{5.2.32}$$

Let us assume arbitrary $\epsilon > 0$. Assume

$$\delta = \frac{\epsilon}{\|\partial f(x)\|}$$

Then from inequality

$$|a| < \delta$$

it follows

$$|f(x+a)-f(x)| \le \|\partial f(x)\|\ \delta = \epsilon$$

According to definition 5.1.11, map $f$ is continuous at point $x$. □

Theorem 5.2.18 has interesting generalization. If unit sphere of division ring $D$ is not compact, then we may consider compact set of directions instead of unit sphere. In this case we may speak about continuity of function $f$ along selected set of directions.

## 5.3. Table of Derivatives of Map of Division Ring

Theorem 5.3.1. *Let $D$ be complete division ring of characteristic 0. Then for any $b \in D$*

$$\partial(b)(a) = 0 \tag{5.3.1}$$

Proof. Immediate corollary of definition 5.2.3. □

Theorem 5.3.2. *Let $D$ be complete division ring of characteristic 0. Then for any $b,\ c \in D$*

$$\partial(bf(x)c)(a) = b\partial f(x)(a)c \tag{5.3.2}$$

$$\frac{\partial_{s\cdot 0} bf(x)c}{\partial x} = b\frac{\partial_{s\cdot 0} f(x)}{\partial x} \tag{5.3.3}$$

$$\frac{\partial_{s\cdot 1} bf(x)c}{\partial x} = \frac{\partial_{s\cdot 1} f(x)}{\partial x}c \tag{5.3.4}$$

$$\frac{\partial bf(x)c}{\partial_* x}(a) = a^{-1}ba\frac{\partial f(x)(a)}{\partial_* x}c \tag{5.3.5}$$

Proof. Immediate corollary of equations (5.2.23), (5.2.26), (5.2.27), (5.2.29) because $\partial b = \partial c = 0$. □

THEOREM 5.3.3. *Let $D$ be complete division ring of characteristic $0$. Then for any $b$, $c \in D$*

$$\partial(bxc)(h) = bhc \tag{5.3.6}$$

$$\frac{\partial_{1\cdot 0} bxc}{\partial x} = b \tag{5.3.7}$$

$$\frac{\partial_{1\cdot 1} bxc}{\partial x} = c \tag{5.3.8}$$

$$\frac{\partial bxc}{\partial_* x}(h) = h^{-1}bhc \tag{5.3.9}$$

PROOF. Corollary of theorem 5.3.2, when $f(x) = x$. □

THEOREM 5.3.4. *Let $D$ be complete division ring of characteristic $0$. Then for any $b \in D$*

$$\begin{aligned} &\partial(xb - bx)(h) = hb - bh \\ &\frac{\partial_{1\cdot 0}(xb - bx)}{\partial x} = 1 \qquad \frac{\partial_{1\cdot 1}(xb - bx)}{\partial x} = b \\ &\frac{\partial_{2\cdot 0}(xb - bx)}{\partial x} = -b \qquad \frac{\partial_{2\cdot 1}(xb - bx)}{\partial x} = 1 \\ &\frac{\partial(xb - bx)}{\partial_* x}(h) = h^{-1}bhc \end{aligned} \tag{5.3.10}$$

PROOF. Corollary of theorem 5.3.2, when $f(x) = x$. □

THEOREM 5.3.5. *Let $D$ be complete division ring of characteristic $0$. Then*[5.7]

$$\begin{aligned} \partial(x^2)(a) &= xa + ax \\ \frac{\partial x^2}{\partial_* x}(a) &= a^{-1}xa + x \end{aligned} \tag{5.3.11}$$

$$\begin{aligned} &\frac{\partial_{1\cdot 0} x^2}{\partial x} = x \qquad \frac{\partial_{1\cdot 1} x^2}{\partial x} = e \\ &\frac{\partial_{2\cdot 0} x^2}{\partial x} = e \qquad \frac{\partial_{2\cdot 1} x^2}{\partial x} = x \end{aligned} \tag{5.3.12}$$

PROOF. (5.3.11) follows from example 5.2.2 and definition 5.2.10. (5.3.12) follows from example 5.2.2 and equation (5.2.14). □

[5.7]The statement of the theorem is similar to example VIII, [11], p. 451. If product is commutative, then the equation (5.3.11) gets form

$$\partial(x^2)(h) = 2hx$$

$$\frac{dx^2}{dx} = 2x$$

Theorem 5.3.6. *Let $D$ be complete division ring of characteristic $0$. Then*[5.8]

$$\partial(x^{-1})(h) = -x^{-1}hx^{-1} \tag{5.3.13}$$
$$\frac{\partial x^{-1}}{\partial_* x}(h) = -h^{-1}x^{-1}hx^{-1}$$
$$\frac{\partial_{1\cdot 0}x^{-1}}{\partial x} = -x^{-1} \qquad \frac{\partial_{1\cdot 1}x^{-1}}{\partial x} = x^{-1}$$

Proof. Let us substitute $f(x) = x^{-1}$ in definition (5.2.12).

$$\begin{aligned}
\partial f(x)(h) &= \lim_{t\to 0,\ t\in R}(t^{-1}((x+th)^{-1}-x^{-1}))\\
&= \lim_{t\to 0,\ t\in R}(t^{-1}((x+th)^{-1}-x^{-1}(x+th)(x+th)^{-1}))\\
&= \lim_{t\to 0,\ t\in R}(t^{-1}(1-x^{-1}(x+th))(x+th)^{-1})\\
&= \lim_{t\to 0,\ t\in R}(t^{-1}(1-1-x^{-1}th)(x+th)^{-1})\\
&= \lim_{t\to 0,\ t\in R}(-x^{-1}h(x+th)^{-1})
\end{aligned} \tag{5.3.14}$$

Equation (5.3.13) follows from chain of equations (5.3.14). □

Theorem 5.3.7. *Let $D$ be complete division ring of characteristic $0$. Then*[5.9]

$$\partial(xax^{-1})(h) = hax^{-1} - xax^{-1}hx^{-1} \tag{5.3.15}$$
$$\frac{\partial xax^{-1}}{\partial_* x}(h) = ax^{-1} - h^{-1}xax^{-1}hx^{-1}$$
$$\frac{\partial_{1\cdot 0}x^{-1}}{\partial x} = 1 \qquad \frac{\partial_{1\cdot 1}x^{-1}}{\partial x} = ax^{-1}$$
$$\frac{\partial_{2\cdot 0}x^{-1}}{\partial x} = -xax^{-1} \qquad \frac{\partial_{2\cdot 1}x^{-1}}{\partial x} = x^{-1}$$

Proof. Equation (5.3.15) is corollary of equations (5.2.23), (5.3.6), (5.3.15). □

[5.8]The statement of the theorem is similar to example IX, [11], p. 451. If product is commutative, then the equation (5.3.13) gets form

$$\partial(x^{-1})(h) = -hx^{-2}$$
$$\frac{dx^{-1}}{dx} = -x^{-2}$$

[5.9]If product is commutative, then

$$y = xax^{-1} = a$$

Accordingly, the derivative is 0.

CHAPTER 6

# Differentiable Maps of $D$-Vector Space

## 6.1. Topological $D$-Vector Space

Definition 6.1.1. Given a topological division ring $D$ and $D$-vector space $V$ such that $V$ has a topology compatible with the structure of the additive group of $V$ and maps

$$(a, v) \in D \times V \to av \in V$$
$$(v, a) \in V \times D \to va \in V$$

are continuous, then $V$ is called a **topological $D$-vector space**[6.1]. □

Definition 6.1.2. The map

$$f : F \to A$$

of set $F$ to an arbitrary algebra $A$ is called **$A$-valued function**. The map

$$f : F \to V$$

of set $F$ to $D$-vector space $V$ is called **$D$-vector function**. □

We consider following definitions in this section for topological $D$-vector space. However definitions do not change, if we consider topological $D^*{}_*$-vector space.

Definition 6.1.3. Let $V$ be $D$-vector space of dimension $n$. Let $\overline{\overline{e}}$ be basis of $D^*{}_*$-vector space $V$. We represent an arbitrary map

$$f : V \to A$$

$D$-vector space $V$ to set $A$ as **map**

$$f' : D^n \to A$$

of $n$ **$D$-valued variables**; function $f'$ is defined by equation

$$f'(a^1, ..., a^n) = f(a^*{}_* e)$$

□

Definition 6.1.4. A map

$$f : V \to A$$

of topological $D$-vector space $V$ of dimension $n$ to topological space $A$ is called **continuous with respect to set of the arguments** if for arbitrary neighbourhood $U$ of image of element $\overline{a} = a^*{}_* e \in V$ for every $a_i \in D$ there exist such neighbourhood $V_i$, that

$$f'(V_1, ..., V_n) \subset U$$

□

[6.1] I made definition according to definition from [5], p. TVS I.1

THEOREM 6.1.5. *Continuous map*

$$f : V \to A$$

*of topological $D$-vector space $V$ of dimension $n$ to topological space $A$ is continuous with respect to set of arguments.*

PROOF. Let $b = f(a^*{}_*e) \in A$ and $U$ be neigborhood of point $b$. Since map $f$ is continuous, there exists such neighborhood $V$ of vector $\overline{a} = a^*{}_*e$ that $f(V) \subseteq U$. Since additive operation is continues, then there exists such neighborhood $E_i$ of vector $e_i$ and neighborhood $W_i$ of $a_i \in D$ that $W_*{}^*E \subset V$. □

DEFINITION 6.1.6. **Norm on $D$-vector space $\overline{V}$** over non-discrete valued division ring $D$[6.2] is a map

$$v \in V \to \|v\| \in R$$

which satisfies the following axioms

- $\|v\| \geq 0$
- $\|v\| = 0$ if, and only if, $v = 0$
- $\|v + w\| \leq \|v\| + \|w\|$
- $\|av\| = |a|\|v\|$ for all $a \in D$ and all $v \in V$

If $D$ is a non-discrete valued division ring, $D$-vector space $V$, endowed with the structure defined by a given norm on $V$, is called **normed $D$-vector space**. □

Invariant distance on additive group of $D$-vector space $V$

$$d(a, b) = \|a - b\|$$

defines topology of metric space compatible with structure of $D$-vector space $V$.

DEFINITION 6.1.7. Let $D$ be complete division ring of characteristic 0. Let

$$f : V \to W$$

be map of normed $D$-vector space $V$ with norm $\|x\|_1$ to normed $D$-vector space $\overline{W}$ with norm $\|\overline{y}\|_2$. Value

$$\|f\| = \sup \frac{\|f(x)\|_2}{\|x\|_1}$$

is called **norm of map $f$**. □

## 6.2. Differentiable Maps of $D$-Vector Space

DEFINITION 6.2.1. The map

$$f : V \to W$$

of normed $D$-vector space $V$ with norm $\|x\|_1$ to normed $D$-vector space $W$ with norm $\|x\|_2$[6.3] is called **differentiable in the Gâteaux sense on the set $U \subset V$**, if at every point $x \in U$ the increment of the map $f$ can be represented as

$$f(x + a) - f(x) = \partial f(x)(a) + o(a) \tag{6.2.1}$$

[6.2] I made definition according to definition from [4], IX, §3.3

[6.3] When dimension of $D$-vector space $V$ equal 1, we can identify it with division ring $D$. Then we speak about map of division ring $D$ into $D$-vector space $W$. When dimension of $D$-vector space $W$ equal 1, we can identify it with division ring $D$. Then we speak about map of $D$-vector space $V$ into division ring $D$. In both cases we will omit corresponding index, because this index has single value.

where **the Gâteaux derivative $\partial f(x)$ of map** $f$ is linear map of increment $a$ and $o: V \rightarrow W$ is such continuous map that

$$\lim_{a\rightarrow 0}\frac{\|o(a)\|_2}{\|a\|_1}=0$$

□

REMARK 6.2.2. According to definition 6.2.1 for given $x$, the Gâteaux derivative $\partial f(x)\in\mathcal{L}(D;V;W)$. Therefore, the Gâteaux derivative of map $f$ is map

$$\partial f: V\rightarrow\mathcal{L}(D;V;W)$$

□

THEOREM 6.2.3. *Let $D$ be division ring of characteristic 0. Let $V$, $W$ be $D$-vector spaces. Let $\overline{\overline{e}}_V$ be basis of $D^*{}_*$-vector space $V$ and $h\in V$*

$$\overline{h}=h^*{}_*e_V$$

*Let $\overline{\overline{e}}_W$ be basis of $D^*{}_*$-vector space $W$. It is possible to represent the Gâteaux derivative of map $f$ as*

$$\partial f(x)(h)=\partial_{s\cdot 0\cdot i}f^j(x)\ h^i\ \partial_{s\cdot 1\cdot i}f^j(x)\ e_{W\cdot j}=\frac{\partial_{s\cdot 0}f^j(x)}{\partial x^i}\ h^i\ \frac{\partial_{s\cdot 1}f^j(x)}{\partial x^i}\ e_{W\cdot j} \tag{6.2.2}$$

PROOF. Corollary of definition 6.2.1 and theorem 4.1.4. □

DEFINITION 6.2.4. Expression $\partial_{s\cdot p\cdot i}f^j(x)=\dfrac{\partial_{s\cdot p}f^j(x)}{\partial x^i}$, $p=0,\ 1,$ is called **component of the Gâteaux derivative** of map $f(x)$. □

THEOREM 6.2.5. *Let complete division ring $D$ of characteristic 0 is algebra over field $F\subset Z(D)$. The Gâteaux derivative of the map*

$$f: V\rightarrow W$$

*of normed $D$-vector space $V$ to normed $D$-vector space $W$ satisfies to equation*

$$\partial f(x)(ph)=p\partial f(x)(h)\quad p\in F$$

PROOF. Corollary of definitions 4.1.1, 6.2.1. □

THEOREM 6.2.6. *Let $D$ be complete division ring of characteristic 0. The Gâteaux derivative of the function*

$$f: V\rightarrow W$$

*of normed $D$-vector space $V$ to normed $D$-vector space $W$ satisfies to equation*

$$\partial f(x)(ra)=r\partial f(x)(a) \tag{6.2.3}$$

*for every $r\in R$, $r\neq 0$ and $a\in V$, $a\neq 0$.*

PROOF. Corollary of theorems 5.1.9, 6.2.5. □

Combining equation (6.2.3) and definition 6.2.1, we get known definition of the Gâteaux derivative

$$\partial f(x)(a)=\lim_{t\rightarrow 0}(t^{-1}(f(x+ta)-f(x))) \tag{6.2.4}$$

Definitions of the Gâteaux derivative (6.2.1) and (6.2.4) are equivalent. Using this equivalence we tell that map $f$ is called differentiable in the Gâteaux sense on the set $U \subset D$, if at every point $x \in U$ the increment of the function $f$ can be represented as

$$f(x+ta)-f(x)=t\partial f(x)(a)+o(t) \tag{6.2.5}$$

where $o: R \to W$ is such continuous map that

$$\lim_{t\to 0}\frac{\|o(t)\|_2}{|t|}=0$$

Since infinitesimal $a$ is differential $dx$, then equation (6.2.2) becomes definition of **the Gâteaux differential**

$$\partial f(x)(dx)=\frac{\partial_{s\cdot 0}f^j(x)}{\partial x^i}\ dx^i\ \frac{\partial_{s\cdot 1}f^j(x)}{\partial x^i}\ e_{W\cdot j} \tag{6.2.6}$$

Definition 6.2.7. Let $D$ be division ring of characteristic 0. Let $V$, $W$ be $D$-vector spaces. Let $\overline{\overline{e}}_V$ be basis of $D^*{}_*$-vector space $V$ and $h \in V$

$$\overline{h}=h^*{}_*e_V$$

Let $\overline{\overline{e}}_W$ be basis of $D^*{}_*$-vector space $W$. If we consider the Gâteaux derivative of map $f$ with respect to variable $v^i$ assuming that other coordinates of vector $v$ are fixed, then corresponding additive map

$$\frac{\partial f^j(v)}{\partial v^i}(h^i)=\lim_{t\to 0}(t^{-1}(f^j(v^1,...,v^i+th^i,...,v^n)-f^j(v))) \tag{6.2.7}$$

is called **the Gâteaux partial derivative** of map $f^j$ with respect to variable $v^i$. □

Theorem 6.2.8. *Let $D$ be division ring of characteristic 0. Let $\overline{\overline{e}}_V$ be basis of $D$-vector space $V$ and $h \in V$*

$$\overline{h}=h^*{}_*e_V$$

*Let $\overline{\overline{e}}_W$ be basis of $D$-vector space $W$. Let the Gâteaux partial derivatives of map*

$$f: V \to W$$

*are continuous on the set $U \in V$. Then, on the set $U$, the Gâteaux derivative of map $f$ and the Gâteaux partial derivatives hold*

$$\partial f(x)(h)=\left(\frac{\partial f^j(x)}{\partial x^i}(h^i)\right)e_{W\cdot j} \tag{6.2.8}$$

Proof. From theorem 4.1.2, it follows that we can expand the Gâteaux derivative of map $f$ into sum (6.2.8) of partial linear maps. Our task is to clarify when these partial linear maps become the Gâteaux partial derivatives.

According to definitions 6.2.4,

$$\begin{aligned}\partial f(x)(h)&=\lim_{t\to 0}(t^{-1}(f(x+th)-f(x)))\\&=\lim_{t\to 0}(t^{-1}(f^j(x^1+th^1,x^2+th^2,...,x^n+th^n)\\&-f^j(x^1,x^2+th^2,...,x^n+th^n)\\&+f^j(x^1,x^2+th^2,...,x^n+th^n)-...\\&-f^j(x^1,...,x^n)))e_{W\cdot j}\end{aligned} \tag{6.2.9}$$

From equations (6.2.7) and (6.2.9), up to infinitesimal with respect to $t$, it follows that

$$\begin{aligned}\partial f(x)(h) = \lim_{t\to 0}\Bigg(&\frac{\partial f^{j}(x^{1},x^{2}+th^{2},...,x^{n}+th^{n})}{\partial x^{1}}(h^{1})\\ &+...+\frac{\partial f^{j}(x^{1},...,x^{n})}{\partial x^{n}}(h^{n})\Bigg)e_{W\cdot j}\end{aligned} \tag{6.2.10}$$

Therefore, for the equation (6.2.8) to follow from the equation (6.2.10) it is necessary that the Gâteaux partial derivatives be continuous. □

Definition 6.2.9. Let map

$$f:V\to W$$

of normed $D$-vector space $V$ to normed $D$-vector space $W$ be differentiable in the Gâteaux sense on the set $U\times W$. Assume we selected basis $\overline{\overline{e}}_V$ in ${D^*}_*$-vector space $V$. Assume we selected basis $\overline{\overline{e}}_W$ in ${D^*}_*$-vector space $W$. We define **the Gâteaux ${D^*}_*$-derivative $\frac{\partial f(x)(a)}{(\partial^*{}_*)x}$ of map $f$.** using equation

$$\partial f(x)(a) = a^*{}_*\frac{\partial f(x)(a)}{(\partial^*{}_*)x} = a^{i}\,\frac{\partial f^{j}(x)(a)}{(\partial^*{}_*)x^{i}}\,e_{W\cdot j} \tag{6.2.11}$$

Expression $\frac{\partial f^{j}(x)(a)}{(\partial^*{}_*)x^{i}}$ in equation (6.2.11) is called **the Gâteaux partial ${D^*}_*$-derivative of map $f^{j}$ with respect to variable $x^{i}$.** □

For representation of the Gâteaux ${D^*}_*$-derivative we also use notation

$$\begin{aligned}\frac{\partial f(x)(a)}{(\partial^*{}_*)x} &= \frac{\partial}{(\partial^*{}_*)x}f(x)(a) &&= (\partial^*{}_*)f(x)(a)\\ \frac{\partial f^{j}(x)(a^{i})}{(\partial^*{}_*)x^{i}} &= \frac{\partial}{(\partial^*{}_*)x^{i}}f^{j}(x)(a) &&= (\partial_i{}^*{}_*)f^{j}(x)(a)\end{aligned}$$

Definition 6.2.10. Let map

$$f:V\to W$$

of normed $D$-vector space $V$ to normed $D$-vector space $W$ be differentiable in the Gâteaux sense on the set $U\times W$. Assume we selected basis $\overline{\overline{e}}_V$ in ${}_*{}^*D$-vector space $V$. Assume we selected basis $\overline{\overline{e}}_W$ in ${}_*{}^*D$-vector space $W$. We define **the Gâteaux ${}_*{}^*D$-derivative $\frac{\partial f(x)(a)}{({}_*{}^*\partial)x}$ of map $f$.** using equation

$$\partial f(x)(a) = \frac{\partial f(x)(a)}{({}_*{}^*\partial)x}{}_*{}^*\overline{a} = e_{W\cdot j}\,\frac{\partial f^{j}(x)(a)}{({}_*{}^*\partial)x^{i}}\,a^{i} \tag{6.2.12}$$

Expression $\frac{\partial f^{j}(x)(a)}{({}_*{}^*\partial)x^{i}}$ in equation (6.2.12) is called **the Gâteaux partial ${}_*{}^*D$-derivative of map $f^{j}$ with respect to variable $x^{i}$.** □

For representation of the Gâteaux ${}_*{}^*D$-derivative we also use notation

$$\begin{aligned}\frac{\partial f(x)(a)}{({}_*{}^*\partial)x} &= \frac{\partial}{({}_*{}^*\partial)x}f(x)(a) &&= ({}_*{}^*\partial)f(x)(a)\\ \frac{\partial f^{j}(x)(a)}{({}_*{}^*\partial)x^{i}} &= \frac{\partial}{({}_*{}^*\partial)x^{i}}f^{j}(x)(a) &&= ({}_*{}^*\partial_i)\ f^{j}(x)(a)\end{aligned}$$

THEOREM 6.2.11. *Let $D$ be division ring of characteristic* 0. *Let function*

$$f:U\to V$$

*of normed $D$-vector space $U$ to normed $D$-vector space $V$ be differentiable in the Gâteaux sense at point $x$. Then*

$$\partial f(x)(0)=0$$

PROOF. Corollary of definitions 6.2.1 and theorem 4.1.8. □

THEOREM 6.2.12. *Let $D$ be division ring of characteristic* 0. *Let $U$ be normed $D$-vector space with norm $\|x\|_1$. Let $V$ be normed $D$-vector space with norm $\|y\|_2$. Let $W$ be normed $D$-vector space with norm $\|z\|_3$. Let map*

$$f:U\to V$$

*be differentiable in the Gâteaux sense at point $x$ and norm of the Gâteaux derivative of map $f$ be finite*

$$\|\partial f(x)\|=F\le\infty \tag{6.2.13}$$

*Let map*

$$g:V\to W$$

*be differentiable in the Gâteaux sense at point*

$$y=f(x) \tag{6.2.14}$$

*and norm of the Gâteaux derivative of map $g$ be finite*

$$\|\partial g(y)\|=G\le\infty \tag{6.2.15}$$

*Map*

$$gf(x)=g(f(x)) \tag{6.2.16}$$

*is differentiable in the Gâteaux sense at point $x$*

$$\partial gf(x)(a)=\partial g(f(x))(\partial f(x)(a)) \tag{6.2.17}$$

$$\frac{\partial gf(x)(a)}{(\partial^*{}_*)x}=\frac{\partial f(x)(a)}{(\partial^*{}_*)x}{}_*{}^*\frac{\partial g(f(x))(\partial f(x)(a))}{(\partial^*{}_*)f(x)} \tag{6.2.18}$$

$$\frac{\partial (gf)^{k}(x)(a)}{(\partial^*{}_*)x^{i}}=\frac{\partial f^{j}(x)(a)}{(\partial^*{}_*)x^{i}}\frac{\partial g^{k}(f(x))(\partial f(x)(a))}{(\partial^*{}_*)f^{j}(x)} \tag{6.2.19}$$

$$\frac{\partial_{st\cdot 0}(gf)^{k}(x)}{\partial x^{i}}=\frac{\partial_{s\cdot 0}g^{k}(f(x))}{\partial f^{j}(x)}\frac{\partial_{t\cdot 0}f^{j}(x)}{\partial x^{i}} \tag{6.2.20}$$

$$\frac{\partial_{st\cdot 1}(gf)^{k}(x)}{\partial x^{i}}=\frac{\partial_{t\cdot 1}f^{j}(x)}{\partial x^{i}}\frac{\partial_{s\cdot 1}g^{k}(f(x))}{\partial f^{j}(x)} \tag{6.2.21}$$

PROOF. According to definition 6.2.1

$$g(y+b)-g(y)=\partial g(y)(b)+o_1(b) \tag{6.2.22}$$

where $o_1:V\to W$ is such continuous map that

$$\lim_{b\to 0}\frac{\|o_1(b)\|_3}{\|b\|_2}=0$$

According to definition 6.2.1

$$f(x+a)-f(x)=\partial f(x)(a)+o_2(a) \tag{6.2.23}$$

where $o_2:U\to V$ is such continuous map that

$$\lim_{a\to 0}\frac{\|o_2(a)\|_2}{\|a\|_1}=0$$

According to (6.2.23) increment of vector $x\in U$ in direction $a$ leads to increment of vector $y$ in direction

$$b=\partial f(x)(a)+o_2(a) \tag{6.2.24}$$

Using (6.2.14), (6.2.24) in equation (6.2.22), we get

$$\begin{aligned}&g(f(x+a))-g(f(x))\\ =&g(f(x)+\partial f(x)(a)+o_2(a))-g(f(x))\\ =&\partial g(f(x))(\partial f(x)(a)+o_2(a))-o_1(\partial f(x)(a)+o_2(a))\end{aligned} \tag{6.2.25}$$

According to definitions 6.2.1, 4.1.1 from equation (6.2.25) it follows

$$\begin{aligned}&g(f(x+a))-g(f(x))\\ =&\partial g(f(x))(\partial f(x)(a))+\partial g(f(x))(o_2(a))-o_1(\partial f(x)(a)+o_2(a))\end{aligned} \tag{6.2.26}$$

Let us find limit

$$\lim_{a\to 0}\frac{\|\partial g(f(x))(o_2(a))-o_1(\partial f(x)(a)+o_2(a))\|_3}{\|a\|_1}$$

According to definition 6.1.6

$$\begin{aligned}&\lim_{a\to 0}\frac{\|\partial g(f(x))(o_2(a))-o_1(\partial f(x)(a)+o_2(a))\|_3}{\|a\|_1}\\ \le&\lim_{a\to 0}\frac{\|\partial g(f(x))(o_2(a))\|_3}{\|a\|_1}+\lim_{a\to 0}\frac{\|o_1(\partial f(x)(a)+o_2(a))\|_3}{\|a\|_1}\end{aligned} \tag{6.2.27}$$

From (6.2.15) it follows that

$$\lim_{a\to 0}\frac{\|\partial g(f(x))(o_2(a))\|_3}{\|a\|_1}\le G\lim_{a\to 0}\frac{\|o_2(a)\|_2}{\|a\|_1}=0 \tag{6.2.28}$$

From (6.2.13) it follows that

$$\begin{aligned}&\lim_{a\to 0}\frac{\|o_1(\partial f(x)(a)+o_2(a))\|_3}{\|a\|_1}\\ =&\lim_{a\to 0}\frac{\|o_1(\partial f(x)(a)+o_2(a))\|_3}{\|\partial f(x)(a)+o_2(a)\|_2}\lim_{a\to 0}\frac{\|\partial f(x)(a)+o_2(a)\|_2}{\|a\|_1}\\ \le&\lim_{a\to 0}\frac{\|o_1(\partial f(x)(a)+o_2(a))\|_3}{\|\partial f(x)(a)+o_2(a)\|_2}\lim_{a\to 0}\frac{\|\partial f(x)\|\|a\|_1+\|o_2(a)\|_2}{\|a\|_1}\\ =&\lim_{a\to 0}\frac{\|o_1(\partial f(x)(a)+o_2(a))\|_3}{\|\partial f(x)(a)+o_2(a)\|_2}\|\partial f(x)\|\end{aligned}$$

According to the theorem 6.2.11

$$\lim_{a\to 0}(\partial f(x)(a)+o_2(a))=0$$

Therefore,

$$\lim_{a\to 0}\frac{\|o_1(\partial f(x)(a)+o_2(a))\|_3}{\|a\|_1}=0 \tag{6.2.29}$$

From equations (6.2.27), (6.2.28), (6.2.29) it follows that

$$\lim_{a\to 0}\frac{\|\partial g(f(x))(o_2(a))-o_1(\partial f(x)(a)+o_2(a))\|_3}{\|a\|_1}=0 \tag{6.2.30}$$

According to definition 6.2.1

$$gf(x+a)-gf(x)=\partial gf(x)(a)+o(a) \tag{6.2.31}$$

where $o:U\to W$ is such continuous map that

$$\lim_{a\to 0}\frac{\|o(a)\|_3}{\|a\|_1}=0$$

Equation (6.2.17) follows from (6.2.26), (6.2.30), (6.2.31).

From equation (6.2.17) and the definition 6.2.9 it follows that

$$\begin{aligned}a^*{}_*\frac{\partial gf(x)(a)}{(\partial^*{}_*)x}&=\partial f(x)(a)^*{}_*\frac{\partial g(f(x))(\partial f(x)(a)}{(\partial^*{}_*)f(x)}\\&=a^*{}_*\frac{\partial f(x)(a)}{(\partial^*{}_*)x}{}^*{}_*\frac{\partial g(f(x))(\partial f(x)(a)}{(\partial^*{}_*)f(x)}\end{aligned} \tag{6.2.32}$$

Since increment $a$ is arbitrary, then equation (6.2.18) follows from equation (6.2.32).

From equation (6.2.17) and the theorem 6.2.3 it follows that

$$\begin{aligned}&\frac{\partial_{st\cdot 0}gf^{\boldsymbol{k}}(x)}{\partial x^{\boldsymbol{i}}}\ a^{\boldsymbol{i}}\ \frac{\partial_{st\cdot 1}gf^{\boldsymbol{k}}(x)}{\partial x^{\boldsymbol{i}}}\ e_{W\cdot\boldsymbol{k}}\\=&\frac{\partial_{s\cdot 0}g^{\boldsymbol{k}}(f(x))}{\partial f^{\boldsymbol{j}}(x)}\ \partial f^{\boldsymbol{j}}(x)(a)\ \frac{\partial_{s\cdot 1}g^{\boldsymbol{k}}(f(x))}{\partial f^{\boldsymbol{j}}(x)}\ e_{W\cdot\boldsymbol{k}}\\=&\frac{\partial_{s\cdot 0}g^{\boldsymbol{k}}(f(x))}{\partial f^{\boldsymbol{j}}(x)}\frac{\partial_{t\cdot 0}f^{\boldsymbol{j}}(x)}{\partial x^{\boldsymbol{i}}}\ a^{\boldsymbol{i}}\ \frac{\partial_{t\cdot 1}f^{\boldsymbol{j}}(x)}{\partial x^{\boldsymbol{i}}}\ \frac{\partial_{s\cdot 1}g^{\boldsymbol{k}}(f(x))}{\partial f^{\boldsymbol{j}}(x)}\ e_{W\cdot\boldsymbol{k}}\end{aligned} \tag{6.2.33}$$

(6.2.20), (6.2.21) follow from equation (6.2.33). □

CHAPTER 7

# Quaternion Algebra

## 7.1. Linear Function of Complex Field

THEOREM 7.1.1 (the Cauchy-Riemann equations). *Let us consider complex field $C$ as two-dimensional algebra over real field. Let*

$$e_{C\cdot 0}=1 \quad e_{C\cdot 1}=i \tag{7.1.1}$$

*be the basis of algebra $C$. Then in this basis product has form*

$$e_{C\cdot 1}^2=-e_{C\cdot 0} \tag{7.1.2}$$

*and structural constants have form*

$$\begin{array}{ll} C_{C\cdot}{}^{0}_{00}=1 & C_{C\cdot}{}^{1}_{01}=\ \ 1\\ C_{C\cdot}{}^{1}_{10}=1 & C_{C\cdot}{}^{0}_{11}=-1 \end{array} \tag{7.1.3}$$

*Matrix of linear function*

$$y^i=x^jf^i_j$$

*of complex field over real field satisfies relationship*

$$f^0_0=f^1_1 \tag{7.1.4}$$

$$f^1_0=-f^0_1 \tag{7.1.5}$$

PROOF. Equations (7.1.2) and (7.1.3) follow from equation $i^2=-1$. Using equation (3.1.17) we get relationships

$$f^0_0=f^{kr}C_{C\cdot}{}^{p}_{k0}C_{C\cdot}{}^{0}_{pr}=f^{0r}C_{C\cdot}{}^{0}_{00}C_{C\cdot}{}^{0}_{0r}+f^{1r}C_{C\cdot}{}^{1}_{10}C_{C\cdot}{}^{0}_{1r}=f^{00}-f^{11} \tag{7.1.6}$$

$$f^1_0=f^{kr}C_{C\cdot}{}^{p}_{k0}C_{C\cdot}{}^{1}_{pr}=f^{0r}C_{C\cdot}{}^{0}_{00}C_{C\cdot}{}^{1}_{0r}+f^{1r}C_{C\cdot}{}^{1}_{10}C_{C\cdot}{}^{1}_{1r}=f^{01}+f^{10} \tag{7.1.7}$$

$$f^0_1=f^{kr}C_{C\cdot}{}^{p}_{k1}C_{C\cdot}{}^{0}_{pr}=f^{0r}C_{C\cdot}{}^{1}_{01}C_{C\cdot}{}^{0}_{1r}+f^{1r}C_{C\cdot}{}^{0}_{11}C_{C\cdot}{}^{0}_{0r}=-f^{01}-f^{10} \tag{7.1.8}$$

$$f^1_1=f^{kr}C_{C\cdot}{}^{p}_{k1}C_{C\cdot}{}^{1}_{pr}=f^{0r}C_{C\cdot}{}^{1}_{01}C_{C\cdot}{}^{1}_{1r}+f^{1r}C_{C\cdot}{}^{0}_{11}C_{C\cdot}{}^{1}_{0r}=f^{00}-f^{11} \tag{7.1.9}$$

(7.1.4) follows from equations (7.1.6) and (7.1.9). (7.1.5) follows from equations (7.1.7) and (7.1.8). □

THEOREM 7.1.2 (the Cauchy-Riemann equations). *Since matrix*

$$\begin{pmatrix} \dfrac{\partial y^0}{\partial x^0} & \dfrac{\partial y^0}{\partial x^1}\\ \dfrac{\partial y^1}{\partial x^0} & \dfrac{\partial y^1}{\partial x^1} \end{pmatrix}$$

*is Jacobian matrix of map of complex variable*

$$x=x^0+x^1i\rightarrow y=y^0(x^0,x^1)+y^1(x^0,x^1)i$$

*over real field, then*

$$\begin{aligned}\frac{\partial y^1}{\partial x^0} &= -\frac{\partial y^0}{\partial x^1}\\ \frac{\partial y^0}{\partial x^0} &= \frac{\partial y^1}{\partial x^1}\end{aligned} \tag{7.1.10}$$

Proof. The statement of theorem is corollary of theorem 7.1.1. □

Theorem 7.1.3. *Derivative of function of complex variable satisfyes to equation*

$$\frac{\partial y}{\partial x^0} + i\frac{\partial y}{\partial x^1} = 0 \tag{7.1.11}$$

Proof. Equation

$$\frac{\partial y^0}{\partial x^0} + i\frac{\partial y^1}{\partial x^0} + i\frac{\partial y^0}{\partial x^1} - \frac{\partial y^1}{\partial x^1} = 0$$

follows from equations (7.1.10). □

Equation (7.1.11) is equivalent to equation

$$\begin{pmatrix}1 & i\end{pmatrix}\begin{pmatrix}\frac{\partial y^0}{\partial x^0} & \frac{\partial y^0}{\partial x^1}\\ \frac{\partial y^1}{\partial x^0} & \frac{\partial y^1}{\partial x^1}\end{pmatrix}\begin{pmatrix}1\\ i\end{pmatrix} = 0 \tag{7.1.12}$$

## 7.2. Quaternion Algebra

In this paper I explore the set of quaternion algebras defined in [8].

Definition 7.2.1. Let $F$ be field. Extension field $F(i,j,k)$ is called **the quaternion algebra** $E(F,a,b)$ **over the field** $F$[7.1] if multiplication in algebra $E$ is defined according to rule

(7.2.1)

| | $i$ | $j$ | $k$ |
|---|---|---|---|
| $i$ | $a$ | $k$ | $aj$ |
| $j$ | $-k$ | $b$ | $-bi$ |
| $k$ | $-aj$ | $bi$ | $-ab$ |

where $a$, $b \in F$, $ab \neq 0$. □

Elements of the algebra $E(F,a,b)$ have form

$$x = x^0 + x^1 i + x^2 j + x^3 k$$

where $x^i \in F$, $i = 0, 1, 2, 3$. Quaternion

$$\overline{x} = x^0 - x^1 i - x^2 j - x^3 k$$

is called conjugate to the quaternion $x$. We define **the norm of the quaternion** $x$ using equation

$$|x|^2 = x\overline{x} = (x^0)^2 - a(x^1)^2 - b(x^2)^2 + ab(x^3)^2 \tag{7.2.2}$$

From equation (7.2.2), it follows that $E(F,a,b)$ is algebra with division only when $a < 0$, $b < 0$. In this case we can renorm basis such that $a = -1$, $b = -1$.

[7.1] I follow definition from [8].

We use symbol $E(F)$ to denote the quaternion division algebra $E(F,-1,-1)$ over the field $F$. We will use notation $H = E(R,-1,-1)$. Multiplication in quaternion algebra $E(F)$ is defined according to rule

(7.2.3)

| | $i$ | $j$ | $k$ |
|---|---|---|---|
| $i$ | $-1$ | $k$ | $-j$ |
| $j$ | $-k$ | $-1$ | $i$ |
| $k$ | $j$ | $-i$ | $-1$ |

In algebra $E(F)$, the norm of the quaternion has form

$$|x|^2 = x\overline{x} = (x^0)^2 + (x^1)^2 + (x^2)^2 + (x^3)^2 \tag{7.2.4}$$

In this case inverse element has form

$$x^{-1} = |x|^{-2}\overline{x} \tag{7.2.5}$$

The inner automorphism of quaternion algebra $H$[7.2]

$$\begin{aligned} &p \rightarrow qpq^{-1} \\ &q(ix + jy + kz)q^{-1} = ix' + jy' + kz' \end{aligned} \tag{7.2.6}$$

describes the rotation of the vector with coordinates $x$, $y$, $z$. The norm of quaternion $q$ is irrelevant, although usually we assume $|q| = 1$. If $q$ is written as sum of scalar and vector

$$q = \cos\alpha + (ia + jb + kc)\sin\alpha \quad a^2 + b^2 + c^2 = 1$$

then (7.2.6) is a rotation of the vector $(x, y, z)$ about the vector $(a, b, c)$ through an angle $2\alpha$.

## 7.3. Linear Function of Quaternion Algebra

THEOREM 7.3.1. *Let*

$$e_0 = 1 \quad e_1 = i \quad e_2 = j \quad e_3 = k \tag{7.3.1}$$

*be basis of quaternion algebra $H$. Then in the basis (7.3.1), structural constants have form*

$$\begin{array}{llll} C^0_{00} = 1 & C^1_{01} = \ \ 1 & C^2_{02} = \ \ 1 & C^3_{03} = \ \ 1 \\ C^1_{10} = 1 & C^0_{11} = -1 & C^3_{12} = \ \ 1 & C^2_{13} = -1 \\ C^2_{20} = 1 & C^3_{21} = -1 & C^0_{22} = -1 & C^1_{23} = \ \ 1 \\ C^3_{30} = 1 & C^2_{31} = \ \ 1 & C^1_{32} = -1 & C^0_{33} = -1 \end{array}$$

PROOF. Value of structural constants follows from multiplication table (7.2.3). □

Since calculations in this section get a lot of space, I put in one place references to theorems in this section.

**Theorem 7.3.2, page 66:** the definition of coordinates of linear map of quaternion algebra $H$ using standard components of this map.

[7.2] See [10], p. 643.

**Equation** (7.3.22)**, page 69:** matrix form of dependence of coordinates of linear map of quaternion algebra $H$ from standard components of this map.

**Equation** (7.3.23)**, page 69:** matrix form of dependence of standard components of linear map of quaternion algebra $H$ from coordinates of this map.

**Theorem 7.3.4, page 72:** dependence of standard components of a linear map of quaternion algebra $H$ from coordinates of this map.

THEOREM 7.3.2. *Standard components of linear function of quaternion algebra $H$ relative to basis* (7.3.1) *and coordinates of corresponding linear map satisfy relationship*

$$\left\{\begin{aligned} f_0^0 &= f^{00} - f^{11} - f^{22} - f^{33} \\ f_1^1 &= f^{00} - f^{11} + f^{22} + f^{33} \\ f_2^2 &= f^{00} + f^{11} - f^{22} + f^{33} \\ f_3^3 &= f^{00} + f^{11} + f^{22} - f^{33} \end{aligned}\right. \tag{7.3.2}$$

$$\left\{\begin{aligned} f_0^1 &= \ \ f^{01} + f^{10} + f^{23} - f^{32} \\ f_1^0 &= -f^{01} - f^{10} + f^{23} - f^{32} \\ f_2^3 &= -f^{01} + f^{10} - f^{23} - f^{32} \\ f_3^2 &= \ \ f^{01} - f^{10} - f^{23} - f^{32} \end{aligned}\right. \tag{7.3.3}$$

$$\left\{\begin{aligned} f_0^2 &= \ \ f^{02} - f^{13} + f^{20} + f^{31} \\ f_1^3 &= \ \ f^{02} - f^{13} - f^{20} - f^{31} \\ f_2^0 &= -f^{02} - f^{13} - f^{20} + f^{31} \\ f_3^1 &= -f^{02} - f^{13} + f^{20} - f^{31} \end{aligned}\right. \tag{7.3.4}$$

$$\left\{\begin{aligned} f_0^3 &= \ \ f^{03} + f^{12} - f^{21} + f^{30} \\ f_1^2 &= -f^{03} - f^{12} - f^{21} + f^{30} \\ f_2^1 &= \ \ f^{03} - f^{12} - f^{21} - f^{30} \\ f_3^0 &= -f^{03} + f^{12} - f^{21} - f^{30} \end{aligned}\right. \tag{7.3.5}$$

PROOF. Using equation (3.1.17) we get relationships

$$\begin{aligned} f_0^0 &= f^{kr} C_{k0}^p C_{pr}^0 \\ &= f^{00} C_{00}^0 C_{00}^0 + f^{11} C_{10}^1 C_{11}^0 + f^{22} C_{20}^2 C_{22}^0 + f^{33} C_{30}^3 C_{33}^0 \\ &= f^{00} - f^{11} - f^{22} - f^{33} \end{aligned} \tag{7.3.6}$$

$$\begin{aligned}
f_0^1 &= f^{kr} C_{k0}^p C_{pr}^1 \\
&= f^{01} C_{00}^0 C_{01}^1 + f^{10} C_{10}^1 C_{10}^1 + f^{23} C_{20}^2 C_{23}^1 + f^{32} C_{30}^3 C_{32}^1 \\
&= f^{01} + f^{10} + f^{23} - f^{32}
\end{aligned} \tag{7.3.7}$$

$$\begin{aligned}
f_0^2 &= f^{kr} C_{k0}^p C_{pr}^2 \\
&= f^{02} C_{00}^0 C_{02}^2 + f^{13} C_{10}^1 C_{13}^2 + f^{20} C_{20}^2 C_{20}^2 + f^{31} C_{30}^3 C_{31}^2 \\
&= f^{02} - f^{13} + f^{20} + f^{31}
\end{aligned} \tag{7.3.8}$$

$$\begin{aligned}
f_0^3 &= f^{kr} C_{k0}^p C_{pr}^3 \\
&= f^{03} C_{00}^0 C_{03}^3 + f^{12} C_{10}^1 C_{12}^3 + f^{21} C_{20}^2 C_{21}^3 + f^{30} C_{30}^3 C_{30}^3 \\
&= f^{03} + f^{12} - f^{21} + f^{30}
\end{aligned} \tag{7.3.9}$$

$$\begin{aligned}
f_1^0 &= f^{kr} C_{k1}^p C_{pr}^0 \\
&= f^{01} C_{01}^1 C_{11}^0 + f^{10} C_{11}^0 C_{00}^0 + f^{23} C_{21}^3 C_{33}^0 + f^{32} C_{31}^2 C_{22}^0 \\
&= -f^{01} - f^{10} + f^{23} - f^{32}
\end{aligned} \tag{7.3.10}$$

$$\begin{aligned}
f_1^1 &= f^{kr} C_{k1}^p C_{pr}^1 \\
&= f^{00} C_{01}^1 C_{10}^1 + f^{11} C_{11}^0 C_{01}^1 + f^{22} C_{21}^3 C_{32}^1 + f^{33} C_{31}^2 C_{23}^1 \\
&= f^{00} - f^{11} + f^{22} + f^{33}
\end{aligned} \tag{7.3.11}$$

$$\begin{aligned}
f_1^2 &= f^{kr} C_{k1}^p C_{pr}^2 \\
&= f^{03} C_{01}^1 C_{13}^2 + f^{12} C_{11}^0 C_{02}^2 + f^{21} C_{21}^3 C_{31}^2 + f^{30} C_{31}^2 C_{20}^2 \\
&= -f^{03} - f^{12} - f^{21} + f^{30}
\end{aligned} \tag{7.3.12}$$

$$\begin{aligned}
f_1^3 &= f^{kr} C_{k1}^p C_{pr}^3 \\
&= f^{02} C_{01}^1 C_{12}^3 + f^{13} C_{11}^0 C_{03}^3 + f^{20} C_{21}^3 C_{30}^3 + f^{31} C_{31}^2 C_{21}^3 \\
&= f^{02} - f^{13} - f^{20} - f^{31}
\end{aligned} \tag{7.3.13}$$

$$\begin{aligned}
f_2^0 &= f^{kr} C_{k2}^p C_{pr}^0 \\
&= f^{02} C_{02}^2 C_{22}^0 + f^{13} C_{12}^3 C_{33}^0 + f^{20} C_{22}^0 C_{00}^0 + f^{31} C_{32}^1 C_{11}^0 \\
&= -f^{02} - f^{13} - f^{20} + f^{31}
\end{aligned} \tag{7.3.14}$$

(7.3.15)
$$\begin{aligned}f_2^1 &= f^{kr}C_{k2}^pC_{pr}^1\\ &= f^{03}C_{02}^2C_{23}^1+f^{12}C_{12}^3C_{32}^1+f^{21}C_{22}^0C_{01}^1+f^{30}C_{32}^1C_{10}^1\\ &= f^{03}-f^{12}-f^{21}-f^{30}\end{aligned}$$

(7.3.16)
$$\begin{aligned}f_2^2 &= f^{kr}C_{k2}^pC_{pr}^2\\ &= f^{00}C_{02}^2C_{20}^2+f^{11}C_{12}^3C_{31}^2+f^{22}C_{22}^0C_{02}^2+f^{33}C_{32}^1C_{13}^2\\ &= f^{00}+f^{11}-f^{22}+f^{33}\end{aligned}$$

(7.3.17)
$$\begin{aligned}f_2^3 &= f^{kr}C_{k2}^pC_{pr}^3\\ &= f^{01}C_{02}^2C_{21}^3+f^{10}C_{12}^3C_{30}^3+f^{23}C_{22}^0C_{03}^3+f^{32}C_{32}^1C_{12}^3\\ &= -f^{01}+f^{10}-f^{23}-f^{32}\end{aligned}$$

(7.3.18)
$$\begin{aligned}f_3^0 &= f^{kr}C_{k3}^pC_{pr}^0\\ &= f^{03}C_{03}^3C_{33}^0+f^{12}C_{13}^2C_{22}^0+f^{21}C_{23}^1C_{11}^0+f^{30}C_{33}^0C_{00}^0\\ &= -f^{03}+f^{12}-f^{21}-f^{30}\end{aligned}$$

(7.3.19)
$$\begin{aligned}f_3^1 &= f^{kr}C_{k3}^pC_{pr}^1\\ &= f^{02}C_{03}^3C_{32}^1+f^{13}C_{13}^2C_{23}^1+f^{20}C_{23}^1C_{10}^1+f^{31}C_{33}^0C_{01}^1\\ &= -f^{02}-f^{13}+f^{20}-f^{31}\end{aligned}$$

(7.3.20)
$$\begin{aligned}f_3^2 &= f^{kr}C_{k3}^pC_{pr}^2\\ &= f^{01}C_{03}^3C_{31}^2+f^{10}C_{13}^2C_{20}^2+f^{23}C_{23}^1C_{13}^2+f^{32}C_{33}^0C_{02}^2\\ &= f^{01}-f^{10}-f^{23}-f^{32}\end{aligned}$$

(7.3.21)
$$\begin{aligned}f_3^3 &= f^{kr}C_{k3}^pC_{pr}^3\\ &= f^{00}C_{03}^3C_{30}^3+f^{11}C_{13}^2C_{21}^3+f^{22}C_{23}^1C_{12}^3+f^{33}C_{33}^0C_{03}^3\\ &= f^{00}+f^{11}+f^{22}-f^{33}\end{aligned}$$

Equations (7.3.6), (7.3.11), (7.3.16), (7.3.21) form the system of linear equations (7.3.2).

Equations (7.3.7), (7.3.10), (7.3.17), (7.3.20) form the system of linear equations (7.3.3).

Equations (7.3.8), (7.3.13), (7.3.14), (7.3.19) form the system of linear equations (7.3.4).

Equations (7.3.9), (7.3.12), (7.3.15), (7.3.18) form the system of linear equations (7.3.5). □

THEOREM 7.3.3. *Consider quaternion algebra $H$ with the basis (7.3.1). Standard components of linear function over field $F$ and coordinates of this function over field $F$ satisfy relationship*

$$
\begin{pmatrix} f^{0}_{0} & f^{0}_{1} & f^{0}_{2} & f^{0}_{3} \\ f^{1}_{1} & -f^{0}_{1} & f^{1}_{3} & -f^{1}_{2} \\ f^{2}_{2} & -f^{2}_{3} & -f^{2}_{0} & f^{2}_{1} \\ f^{3}_{3} & f^{3}_{2} & -f^{3}_{1} & -f^{3}_{0} \end{pmatrix}
= \begin{pmatrix} 1 & -1 & -1 & -1 \\ 1 & -1 & 1 & 1 \\ 1 & 1 & -1 & 1 \\ 1 & 1 & 1 & -1 \end{pmatrix}
\begin{pmatrix} f^{00} & -f^{01} & -f^{02} & -f^{03} \\ f^{11} & f^{10} & f^{13} & -f^{12} \\ f^{22} & -f^{23} & f^{20} & f^{21} \\ f^{33} & f^{32} & -f^{31} & f^{30} \end{pmatrix} \tag{7.3.22}
$$

$$
\begin{pmatrix} f^{00} & -f^{01} & -f^{02} & -f^{03} \\ f^{11} & f^{10} & f^{13} & -f^{12} \\ f^{22} & -f^{23} & f^{20} & f^{21} \\ f^{33} & f^{32} & -f^{31} & f^{30} \end{pmatrix}
= \frac{1}{4} \begin{pmatrix} 1 & 1 & 1 & 1 \\ -1 & -1 & 1 & 1 \\ -1 & 1 & -1 & 1 \\ -1 & 1 & 1 & -1 \end{pmatrix}
\begin{pmatrix} f^{0}_{0} & f^{0}_{1} & f^{0}_{2} & f^{0}_{3} \\ f^{1}_{1} & -f^{0}_{1} & f^{1}_{3} & -f^{1}_{2} \\ f^{2}_{2} & -f^{2}_{3} & -f^{2}_{0} & f^{2}_{1} \\ f^{3}_{3} & f^{3}_{2} & -f^{3}_{1} & -f^{3}_{0} \end{pmatrix} \tag{7.3.23}
$$

*where*

$$
\begin{pmatrix} 1 & -1 & -1 & -1 \\ 1 & -1 & 1 & 1 \\ 1 & 1 & -1 & 1 \\ 1 & 1 & 1 & -1 \end{pmatrix}^{-1}
= \frac{1}{4} \begin{pmatrix} 1 & 1 & 1 & 1 \\ -1 & -1 & 1 & 1 \\ -1 & 1 & -1 & 1 \\ -1 & 1 & 1 & -1 \end{pmatrix}
$$

PROOF. Let us write the system of linear equations (7.3.2) as product of matrices

$$
\begin{pmatrix} f^{0}_{0} \\ f^{1}_{1} \\ f^{2}_{2} \\ f^{3}_{3} \end{pmatrix}
= \begin{pmatrix} 1 & -1 & -1 & -1 \\ 1 & -1 & 1 & 1 \\ 1 & 1 & -1 & 1 \\ 1 & 1 & 1 & -1 \end{pmatrix}
\begin{pmatrix} f^{00} \\ f^{11} \\ f^{22} \\ f^{33} \end{pmatrix} \tag{7.3.24}
$$

Let us write the system of linear equations (7.3.3) as product of matrices

$$\begin{pmatrix} f_{0}^{1} \\ f_{1}^{0} \\ f_{2}^{3} \\ f_{3}^{2} \end{pmatrix} = \begin{pmatrix} 1 & 1 & 1 & -1 \\ -1 & -1 & 1 & -1 \\ -1 & 1 & -1 & -1 \\ 1 & -1 & -1 & -1 \end{pmatrix} \begin{pmatrix} f^{01} \\ f^{10} \\ f^{23} \\ f^{32} \end{pmatrix} \tag{7.3.25}$$

From the equation (7.3.25), it follows that

$$\begin{aligned} \begin{pmatrix} -f_{1}^{0} \\ f_{1}^{0} \\ f_{2}^{3} \\ -f_{3}^{2} \end{pmatrix} &= \begin{pmatrix} -1 & -1 & -1 & 1 \\ -1 & -1 & 1 & -1 \\ -1 & 1 & -1 & -1 \\ -1 & 1 & 1 & 1 \end{pmatrix} \begin{pmatrix} f^{01} \\ f^{10} \\ f^{23} \\ f^{32} \end{pmatrix} \\ &= \begin{pmatrix} 1 & -1 & 1 & 1 \\ 1 & -1 & -1 & -1 \\ 1 & 1 & 1 & -1 \\ 1 & 1 & -1 & 1 \end{pmatrix} \begin{pmatrix} -f^{01} \\ f^{10} \\ -f^{23} \\ f^{32} \end{pmatrix} \end{aligned}$$

$$\begin{pmatrix} f_{1}^{0} \\ -f_{1}^{0} \\ -f_{3}^{2} \\ f_{2}^{3} \end{pmatrix} = \begin{pmatrix} 1 & -1 & -1 & -1 \\ 1 & -1 & 1 & 1 \\ 1 & 1 & -1 & 1 \\ 1 & 1 & 1 & -1 \end{pmatrix} \begin{pmatrix} -f^{01} \\ f^{10} \\ -f^{23} \\ f^{32} \end{pmatrix} \tag{7.3.26}$$

Let us write the system of linear equations (7.3.4) as product of matrices

$$\begin{pmatrix} f_{0}^{2} \\ f_{1}^{3} \\ f_{2}^{0} \\ f_{3}^{1} \end{pmatrix} = \begin{pmatrix} 1 & -1 & 1 & 1 \\ 1 & -1 & -1 & -1 \\ -1 & -1 & -1 & 1 \\ -1 & -1 & 1 & -1 \end{pmatrix} \begin{pmatrix} f^{02} \\ f^{13} \\ f^{20} \\ f^{31} \end{pmatrix} \tag{7.3.27}$$

From the equation (7.3.27), it follows that

$$\begin{pmatrix}-f_0^2\\-f_1^3\\f_2^0\\f_3^1\end{pmatrix}=\begin{pmatrix}-1&1&-1&-1\\-1&1&1&1\\-1&-1&-1&1\\-1&-1&1&-1\end{pmatrix}\begin{pmatrix}f^{02}\\f^{13}\\f^{20}\\f^{31}\end{pmatrix}$$
$$=\begin{pmatrix}1&1&-1&1\\1&1&1&-1\\1&-1&-1&-1\\1&-1&1&1\end{pmatrix}\begin{pmatrix}-f^{02}\\f^{13}\\f^{20}\\-f^{31}\end{pmatrix}$$

$$\begin{pmatrix}f_2^0\\f_3^1\\-f_0^2\\-f_1^3\end{pmatrix}=\begin{pmatrix}1&-1&-1&-1\\1&-1&1&1\\1&1&-1&1\\1&1&1&-1\end{pmatrix}\begin{pmatrix}-f^{02}\\f^{13}\\f^{20}\\-f^{31}\end{pmatrix}\tag{7.3.28}$$

Let us write the system of linear equations (7.3.5) as product of matrices

$$\begin{pmatrix}f_0^3\\f_1^2\\f_2^1\\f_3^0\end{pmatrix}=\begin{pmatrix}1&1&-1&1\\-1&-1&-1&1\\1&-1&-1&-1\\-1&1&-1&-1\end{pmatrix}\begin{pmatrix}f^{03}\\f^{12}\\f^{21}\\f^{30}\end{pmatrix}\tag{7.3.29}$$

From the equation (7.3.29), it follows that

$$\begin{pmatrix}-f_0^3\\f_1^2\\-f_2^1\\f_3^0\end{pmatrix}=\begin{pmatrix}-1&-1&1&-1\\-1&-1&-1&1\\-1&1&1&1\\-1&1&-1&-1\end{pmatrix}\begin{pmatrix}f^{03}\\f^{12}\\f^{21}\\f^{30}\end{pmatrix}$$
$$=\begin{pmatrix}1&1&1&-1\\1&1&-1&1\\1&-1&1&1\\1&-1&-1&-1\end{pmatrix}\begin{pmatrix}-f^{03}\\-f^{12}\\f^{21}\\f^{30}\end{pmatrix}$$

$$\begin{pmatrix} f_3^0 \\ -f_2^1 \\ f_1^2 \\ -f_0^3 \end{pmatrix} = \begin{pmatrix} 1 & -1 & -1 & -1 \\ 1 & -1 & 1 & 1 \\ 1 & 1 & -1 & 1 \\ 1 & 1 & 1 & -1 \end{pmatrix} \begin{pmatrix} -f^{03} \\ -f^{12} \\ f^{21} \\ f^{30} \end{pmatrix} \tag{7.3.30}$$

We join equations (7.3.24), (7.3.26), (7.3.28), (7.3.30) into equation (7.3.22). □

THEOREM 7.3.4. *Standard components of linear function of quaternion algebra $H$ relative to basis (7.3.1) and coordinates of corresponding linear map satisfy relationship*

$$\left\{\begin{array}{l} 4f^{00} = f_0^0 + f_1^1 + f_2^2 + f_3^3 \\ 4f^{11} = -f_0^0 - f_1^1 + f_2^2 + f_3^3 \\ 4f^{22} = -f_0^0 + f_1^1 - f_2^2 + f_3^3 \\ 4f^{33} = -f_0^0 + f_1^1 + f_2^2 - f_3^3 \end{array}\right. \tag{7.3.31}$$

$$\left\{\begin{array}{l} 4f^{10} = -f_1^0 + f_1^0 - f_3^2 + f_2^3 \\ 4f^{01} = -f_1^0 + f_1^0 + f_3^2 - f_2^3 \\ 4f^{32} = -f_1^0 - f_1^0 - f_3^2 - f_2^3 \\ 4f^{23} = f_1^0 + f_1^0 - f_3^2 - f_2^3 \end{array}\right. \tag{7.3.32}$$

$$\left\{\begin{array}{l} 4f^{20} = -f_2^0 + f_3^1 + f_0^2 - f_1^3 \\ 4f^{31} = f_2^0 - f_3^1 + f_0^2 - f_1^3 \\ 4f^{02} = -f_2^0 - f_3^1 + f_0^2 + f_1^3 \\ 4f^{13} = -f_2^0 - f_3^1 - f_0^2 - f_1^3 \end{array}\right. \tag{7.3.33}$$

$$\left\{\begin{array}{l} 4f^{30} = -f_3^0 - f_2^1 + f_1^2 + f_0^3 \\ 4f^{21} = -f_3^0 - f_2^1 - f_1^2 - f_0^3 \\ 4f^{12} = f_3^0 - f_2^1 - f_1^2 + f_0^3 \\ 4f^{03} = -f_3^0 + f_2^1 - f_1^2 + f_0^3 \end{array}\right. \tag{7.3.34}$$

PROOF. We get systems of linear equations (7.3.31), (7.3.32), (7.3.33), (7.3.34) as the product of matrices in equation (7.3.23). □

## 7.4. Differentiable Map of Division Ring of Quaternions

THEOREM 7.4.1. *Since matrix* $\left(\frac{\partial y^{i}}{\partial x^{j}}\right)$ *is Jacobian matrix of the map* $x \to y$ *of division ring of quaternions over real field, then*

$$\left\{\begin{aligned}\frac{\partial y^{0}}{\partial x^{0}}&=\frac{\partial^{00}y}{\partial x}-\frac{\partial^{11}y}{\partial x}-\frac{\partial^{22}y}{\partial x}-\frac{\partial^{33}y}{\partial x}\\ \frac{\partial y^{1}}{\partial x^{1}}&=\frac{\partial^{00}y}{\partial x}-\frac{\partial^{11}y}{\partial x}+\frac{\partial^{22}y}{\partial x}+\frac{\partial^{33}y}{\partial x}\\ \frac{\partial y^{2}}{\partial x^{2}}&=\frac{\partial^{00}y}{\partial x}+\frac{\partial^{11}y}{\partial x}-\frac{\partial^{22}y}{\partial x}+\frac{\partial^{33}y}{\partial x}\\ \frac{\partial y^{3}}{\partial x^{3}}&=\frac{\partial^{00}y}{\partial x}+\frac{\partial^{11}y}{\partial x}+\frac{\partial^{22}y}{\partial x}-\frac{\partial^{33}y}{\partial x}\end{aligned}\right. \tag{7.4.1}$$

$$\left\{\begin{aligned}\frac{\partial y^{1}}{\partial x^{0}}&=\frac{\partial^{01}y}{\partial x}+\frac{\partial^{10}y}{\partial x}+\frac{\partial^{23}y}{\partial x}-\frac{\partial^{32}y}{\partial x}\\ \frac{\partial y^{0}}{\partial x^{1}}&=-\frac{\partial^{01}y}{\partial x}-\frac{\partial^{10}y}{\partial x}+\frac{\partial^{23}y}{\partial x}-\frac{\partial^{32}y}{\partial x}\\ \frac{\partial y^{3}}{\partial x^{2}}&=-\frac{\partial^{01}y}{\partial x}+\frac{\partial^{10}y}{\partial x}-\frac{\partial^{23}y}{\partial x}-\frac{\partial^{32}y}{\partial x}\\ \frac{\partial y^{2}}{\partial x^{3}}&=\frac{\partial^{01}y}{\partial x}-\frac{\partial^{10}y}{\partial x}-\frac{\partial^{23}y}{\partial x}-\frac{\partial^{32}y}{\partial x}\end{aligned}\right. \tag{7.4.2}$$

$$\left\{\begin{aligned}\frac{\partial y^{2}}{\partial x^{0}}&=\frac{\partial^{02}y}{\partial x}-\frac{\partial^{13}y}{\partial x}+\frac{\partial^{20}y}{\partial x}+\frac{\partial^{31}y}{\partial x}\\ \frac{\partial y^{3}}{\partial x^{1}}&=\frac{\partial^{02}y}{\partial x}-\frac{\partial^{13}y}{\partial x}-\frac{\partial^{20}y}{\partial x}-\frac{\partial^{31}y}{\partial x}\\ \frac{\partial y^{0}}{\partial x^{2}}&=-\frac{\partial^{02}y}{\partial x}-\frac{\partial^{13}y}{\partial x}-\frac{\partial^{20}y}{\partial x}+\frac{\partial^{31}y}{\partial x}\\ \frac{\partial y^{1}}{\partial x^{3}}&=-\frac{\partial^{02}y}{\partial x}-\frac{\partial^{13}y}{\partial x}+\frac{\partial^{20}y}{\partial x}-\frac{\partial^{31}y}{\partial x}\end{aligned}\right. \tag{7.4.3}$$

$$\left\{\begin{aligned}\frac{\partial y^{3}}{\partial x^{0}}&=\frac{\partial^{03}y}{\partial x}+\frac{\partial^{12}y}{\partial x}-\frac{\partial^{21}y}{\partial x}+\frac{\partial^{30}y}{\partial x}\\ \frac{\partial y^{2}}{\partial x^{1}}&=-\frac{\partial^{03}y}{\partial x}-\frac{\partial^{12}y}{\partial x}-\frac{\partial^{21}y}{\partial x}+\frac{\partial^{30}y}{\partial x}\\ \frac{\partial y^{1}}{\partial x^{2}}&=\frac{\partial^{03}y}{\partial x}-\frac{\partial^{12}y}{\partial x}-\frac{\partial^{21}y}{\partial x}-\frac{\partial^{30}y}{\partial x}\\ \frac{\partial y^{0}}{\partial x^{3}}&=-\frac{\partial^{03}y}{\partial x}+\frac{\partial^{12}y}{\partial x}-\frac{\partial^{21}y}{\partial x}-\frac{\partial^{30}y}{\partial x}\end{aligned}\right. \tag{7.4.4}$$

$$
\begin{cases}
4\dfrac{\partial^{00}y}{\partial x} = \dfrac{\partial y^0}{\partial x^0} + \dfrac{\partial y^1}{\partial x^1} + \dfrac{\partial y^2}{\partial x^2} + \dfrac{\partial y^3}{\partial x^3}\\
4\dfrac{\partial^{11}y}{\partial x} = -\dfrac{\partial y^0}{\partial x^0} - \dfrac{\partial y^1}{\partial x^1} + \dfrac{\partial y^2}{\partial x^2} + \dfrac{\partial y^3}{\partial x^3}\\
4\dfrac{\partial^{22}y}{\partial x} = -\dfrac{\partial y^0}{\partial x^0} + \dfrac{\partial y^1}{\partial x^1} - \dfrac{\partial y^2}{\partial x^2} + \dfrac{\partial y^3}{\partial x^3}\\
4\dfrac{\partial^{33}y}{\partial x} = -\dfrac{\partial y^0}{\partial x^0} + \dfrac{\partial y^1}{\partial x^1} + \dfrac{\partial y^2}{\partial x^2} - \dfrac{\partial y^3}{\partial x^3}
\end{cases}
\tag{7.4.5}
$$

$$
\begin{cases}
4\dfrac{\partial^{10}y}{\partial x} = -\dfrac{\partial y^0}{\partial x^1} + \dfrac{\partial y^1}{\partial x^0} - \dfrac{\partial y^2}{\partial x^3} + \dfrac{\partial y^3}{\partial x^2}\\
4\dfrac{\partial^{01}y}{\partial x} = -\dfrac{\partial y^0}{\partial x^1} + \dfrac{\partial y^1}{\partial x^0} + \dfrac{\partial y^2}{\partial x^3} - \dfrac{\partial y^3}{\partial x^2}\\
4\dfrac{\partial^{32}y}{\partial x} = -\dfrac{\partial y^0}{\partial x^1} - \dfrac{\partial y^1}{\partial x^0} - \dfrac{\partial y^2}{\partial x^3} - \dfrac{\partial y^3}{\partial x^2}\\
4\dfrac{\partial^{23}y}{\partial x} = \dfrac{\partial y^0}{\partial x^1} + \dfrac{\partial y^1}{\partial x^0} - \dfrac{\partial y^2}{\partial x^3} - \dfrac{\partial y^3}{\partial x^2}
\end{cases}
\tag{7.4.6}
$$

$$
\begin{cases}
4\dfrac{\partial^{20}y}{\partial x} = -\dfrac{\partial y^0}{\partial x^2} + \dfrac{\partial y^1}{\partial x^3} + \dfrac{\partial y^2}{\partial x^0} - \dfrac{\partial y^3}{\partial x^1}\\
4\dfrac{\partial^{31}y}{\partial x} = +\dfrac{\partial y^0}{\partial x^2} - \dfrac{\partial y^1}{\partial x^3} + \dfrac{\partial y^2}{\partial x^0} - \dfrac{\partial y^3}{\partial x^1}\\
4\dfrac{\partial^{02}y}{\partial x} = -\dfrac{\partial y^0}{\partial x^2} - \dfrac{\partial y^1}{\partial x^3} + \dfrac{\partial y^2}{\partial x^0} + \dfrac{\partial y^3}{\partial x^1}\\
4\dfrac{\partial^{13}y}{\partial x} = -\dfrac{\partial y^0}{\partial x^2} - \dfrac{\partial y^1}{\partial x^3} - \dfrac{\partial y^2}{\partial x^0} - \dfrac{\partial y^3}{\partial x^1}
\end{cases}
\tag{7.4.7}
$$

$$
\begin{cases}
4\dfrac{\partial^{30}y}{\partial x} = -\dfrac{\partial y^0}{\partial x^3} - \dfrac{\partial y^1}{\partial x^2} + \dfrac{\partial y^2}{\partial x^1} + \dfrac{\partial y^3}{\partial x^0}\\
4\dfrac{\partial^{21}y}{\partial x} = -\dfrac{\partial y^0}{\partial x^3} - \dfrac{\partial y^1}{\partial x^2} - \dfrac{\partial y^2}{\partial x^1} - \dfrac{\partial y^3}{\partial x^0}\\
4\dfrac{\partial^{12}y}{\partial x} = \dfrac{\partial y^0}{\partial x^3} - \dfrac{\partial y^1}{\partial x^2} - \dfrac{\partial y^2}{\partial x^1} + \dfrac{\partial y^3}{\partial x^0}\\
4\dfrac{\partial^{03}y}{\partial x} = -\dfrac{\partial y^0}{\partial x^3} + \dfrac{\partial y^1}{\partial x^2} - \dfrac{\partial y^2}{\partial x^1} + \dfrac{\partial y^3}{\partial x^0}
\end{cases}
\tag{7.4.8}
$$

PROOF. The statement of theorem is corollary of theorem 7.3.2. □

THEOREM 7.4.2. *Quaternionic map*

$$f(x) = \overline{x}$$

*has the Gâteaux derivative*

$$\partial(\overline{x})(h) = -\frac{1}{2}(h + ihi + jhj + khk) \tag{7.4.9}$$

PROOF. Jacobian matrix of the map $f$ has form

$$\begin{pmatrix}1&0&0&0\\0&-1&0&0\\0&0&-1&0\\0&0&0&-1\end{pmatrix}\tag{7.4.10}$$

From equations (7.4.5) it follows that

$$\frac{\partial^{00}y}{\partial x}=\frac{\partial^{11}y}{\partial x}=\frac{\partial^{22}y}{\partial x}=\frac{\partial^{33}y}{\partial x}=-\frac{1}{2}\tag{7.4.11}$$

From equations (7.4.6), (7.4.7), (7.4.8) it follows that

$$\frac{\partial^{ij}y}{\partial x}=0\quad i\ne j\tag{7.4.12}$$

Equation (7.4.9) follows from equations (5.2.15), (7.4.11), (7.4.12). □

THEOREM 7.4.3. *The norm of quaternion satisfies equations*[7.3]

$$\partial(|x|^2)\circ dx=x\overline{dx}+\overline{x}dx\tag{7.4.13}$$

$$\partial(|x|)(h)=\frac{1}{2\sqrt{x\overline{dx}+\overline{x}dx}}\tag{7.4.14}$$

PROOF. theorem 7.4.2 and equations (7.2.4), (5.2.23). Since $x\overline{dx}+\overline{x}dx\in R$, than the equation (7.4.14) follows from the equation (7.4.13). □

THEOREM 7.4.4. *In division ring of quaternions over real field $H$ with basis (7.3.1) standard $D\star$-representation of the Gâteaux differential of map $x^2$ has form*

$$\partial x^2(a)=(x+x_1)a+x_2ai+x_3aj+x_4ak$$

PROOF. The theorem follows from the definition 5.2.10. □

THEOREM 7.4.5. *We can identify quaternion*

$$a=a^0+a^1i+a^2j+a^3k\tag{7.4.15}$$

*and matrix*

$$J_a=\begin{pmatrix}a^0&-a^1&-a^2&-a^3\\a^1&a^0&-a^3&a^2\\a^2&a^3&a^0&-a^1\\a^3&-a^2&a^1&a^0\end{pmatrix}\tag{7.4.16}$$

PROOF. The product of quaternions (7.4.15) and

$$x=x^0+x^1i+x^2j+x^3k$$

has form

$$\begin{aligned}ax=&a^0x^0-a^1x^1-a^2x^2-a^3x^3+(a^0x^1+a^1x^0+a^2x^3-a^3x^2)i\\&+(a^0x^2+a^2x^0+a^3x^1-a^1x^3)j+(a^0x^3+a^3x^0+a^1x^2-a^2x^1)k\end{aligned}$$

[7.3] The statement of the theorem is similar to example X, [11], p. 455.

Therefore, function $f_a(x) = ax$ has Jacobian matrix (7.4.16). It is evident that $f_a \circ f_b = f_{ab}$. Similar equation is true for matrices

$$\begin{pmatrix} a^0 & -a^1 & -a^2 & -a^3 \\ a^1 & a^0 & -a^3 & a^2 \\ a^2 & a^3 & a^0 & -a^1 \\ a^3 & -a^2 & a^1 & a^0 \end{pmatrix} \begin{pmatrix} b^0 & -b^1 & -b^2 & -b^3 \\ b^1 & b^0 & -b^3 & b^2 \\ b^2 & b^3 & b^0 & -b^1 \\ b^3 & -b^2 & b^1 & b^0 \end{pmatrix}$$

$$= \begin{pmatrix} \begin{matrix} a^0b^0 - a^1b^1 \\ -a^2b^2 - a^3b^3 \end{matrix} & \begin{matrix} -a^0b^1 - a^1b^0 \\ -a^2b^3 + a^3b^2 \end{matrix} & \begin{matrix} -a^0b^2 + a^1b^3 \\ -a^2b^0 - a^3b^1 \end{matrix} & \begin{matrix} -a^0b^3 - a^1b^2 \\ +a^2b^1 - a^3b^0 \end{matrix} \\ \\ \begin{matrix} a^0b^1 + a^1b^0 \\ +a^2b^3 - a^3b^2 \end{matrix} & \begin{matrix} a^0b^0 - a^1b^1 \\ -a^2b^2 - a^3b^3 \end{matrix} & \begin{matrix} -a^0b^3 - a^1b^2 \\ +a^2b^1 - a^3b^0 \end{matrix} & \begin{matrix} a^0b^2 - a^1b^3 \\ +a^2b^0 + a^3b^1 \end{matrix} \\ \\ \begin{matrix} a^0b^2 - a^1b^3 \\ +a^2b^0 + a^3b^1 \end{matrix} & \begin{matrix} a^0b^3 + a^1b^2 \\ -a^2b^1 + a^3b^0 \end{matrix} & \begin{matrix} a^0b^0 - a^1b^1 \\ -a^2b^2 - a^3b^3 \end{matrix} & \begin{matrix} -a^0b^1 - a^1b^0 \\ -a^2b^3 - a^3b^2 \end{matrix} \\ \\ \begin{matrix} a^0b^3 + a^1b^2 \\ -a^2b^1 + a^3b^0 \end{matrix} & \begin{matrix} -a^0b^2 + a^1b^3 \\ -a^2b^0 - a^3b^1 \end{matrix} & \begin{matrix} a^0b^1 + a^1b^0 \\ +a^2b^3 - a^3b^2 \end{matrix} & \begin{matrix} a^0b^0 - a^1b^1 \\ -a^2b^2 - a^3b^3 \end{matrix} \end{pmatrix}$$

Therefore, we can identify the quaternion $a$ and the matrix $J_a$. $\square$

CHAPTER 8

# Derivative of Second Order of Map of Division Ring

## 8.1. Derivative of Second Order of Map of Division Ring

Let $D$ be valued division ring of characteristic 0. Let

$$f : D \to D$$

function differentiable in the Gâteaux sense. According to remark 5.2.4 the Gâteaux derivative is map

$$\partial f : D \to \mathcal{L}(D; D)$$

According to theorems 3.1.5, 3.1.6 and definition 5.1.13, set $\mathcal{L}(D; D)$ is normed $D$-vector space. Therefore, we may consider the question, if map $\partial f$ is differentiable in the Gâteaux sense.

According to definition 6.2.1,

$$\partial f(x + a_2)(a_1) - \partial f(x)(a_1) = \partial(\partial f(x)(a_1))(a_2) + o_2(a_2) \tag{8.1.1}$$

where $o_2 : D \to \mathcal{L}(D; D)$ is such continuous map, that

$$\lim_{a_2 \to 0} \frac{\|o_2(a_2)\|}{|a_2|} = 0$$

According to definition 6.2.1, map $\partial(\partial f(x)(a_1))(a_2)$ is linear map of variable $a_2$. From equation (8.1.1), it follows that map $\partial(\partial f(x)(a_1))(a_2)$ is linear map of variable $a_1$.

Definition 8.1.1. Polylinear map

$$\partial^2 f(x)(a_1; a_2) = \frac{\partial^2 f(x)}{\partial x^2}(a_1; a_2) = \partial(\partial f(x)(a_1))(a_2) \tag{8.1.2}$$

is called **the Gâteaux derivative of second order** of map $f$. □

Remark 8.1.2. According to definition 8.1.1, for given $x$ the Gâteaux derivative of second order $\partial^2 f(x) \in \mathcal{L}(D, D; D)$. Therefore, the Gâteaux derivative of second order of map $f$ is map

$$\partial^2 f : D \to \mathcal{L}(D, D; D)$$

Theorem 8.1.3. *It is possible to represent the Gâteaux derivative of second order of map $f$ as*

$$\partial^2 f(x)(a_1; a_2) = \frac{\partial^2_{s \cdot 0} f(x)}{\partial x^2} \sigma_s(a_1) \frac{\partial^2_{s \cdot 1} f(x)}{\partial x^2} \sigma_s(a_2) \frac{\partial^2_{s \cdot 2} f(x)}{\partial x^2} \tag{8.1.3}$$

Proof. Corollary of definition 8.1.1 and theorem 4.2.3. □

Definition 8.1.4. Expression[8.1] $\dfrac{\partial^2_{s\cdot p} f(x)}{\partial x^2}$, $p=0,\ 1,\ 2,$ is called **component of the Gâteaux derivative of second order** of map $f(x)$. □

By induction, assuming that we defined the Gâteaux derivative $\partial^{n-1}f(x)$ of order $n-1$, we define

$$\partial^n f(x)(a_1;...;a_n)=\frac{\partial^n f(x)}{\partial x^n}(a_1;...;a_n)=\partial(\partial^{n-1}f(x)(a_1;...;a_{n-1}))(a_n) \tag{8.1.4}$$

**the Gâteaux derivative of order** $n$ of map $f$. We also assume $\partial^0 f(x)=f(x)$.

## 8.2. Taylor Series

Consider polynomial in one variable over division ring $D$ of power $n$, $n>0$. We want to explore the structure of monomial $p_k(x)$ of polynomial of power $k$.

It is evident that monomial of power 0 has form $a_0$, $a_0\in D$. Let $k>0$. Let us prove that

$$p_k(x)=p_{k-1}(x)xa_k$$

where $a_k\in D$. Actually, last factor of monomial $p_k(x)$ is either $a_k\in D$, or has form $x^l$, $l\ge 1$. In the later case we assume $a_k=1$. Factor preceding $a_k$ has form $x^l$, $l\ge 1$. We can represent this factor as $x^{l-1}x$. Therefore, we proved the statement.

In particular, monomial of power 1 has form $p_1(x)=a_0xa_1$.

Without loss of generality, we assume $k=n$.

Theorem 8.2.1. *For any $m>0$ the following equation is true*

$$\begin{aligned}&\partial^m(f(x)x)(h_1;...;h_m)\\&=\partial^m f(x)(h_1;...;h_m)x+\partial^{m-1}f(x)(h_1;...;h_{m-1})h_m\\&+\partial^{m-1}f(x)(\widehat{h_1};...;h_{m-1};h_m)h_1+...\\&+\partial^{m-1}f(x)(h_1;...;\widehat{h_{m-1}};h_m)h_{m-1}\end{aligned} \tag{8.2.1}$$

*where symbol $\widehat{h^i}$ means absense of variable $h^i$ in the list.*

Proof. For $m=1$, this is corollary of equation (5.2.23)

$$\partial(f(x)x)(h_1)=\partial f(x)(h_1)x+f(x)h_1$$

Assume, (8.2.1) is true for $m-1$. Then

$$\begin{aligned}&\partial^{m-1}(f(x)x)(h_1;...;h_{m-1})\\&=\partial^{m-1}f(x)(h_1;...;h_{m-1})x+\partial^{m-2}f(x)(h_1;...;h_{m-2})h_{m-1}\\&+\partial^{m-2}f(x)(\widehat{h_1},...,h_{m-2},h_{m-1})h_1+...\\&+\partial^{m-2}f(x)(h_1,...,\widehat{h_{m-2}},h_{m-1})h_{m-2}\end{aligned} \tag{8.2.2}$$

[8.1] We suppose

$$\frac{\partial^2_{s\cdot p} f(x)}{\partial x^2}=\frac{\partial^2_{s\cdot p} f(x)}{\partial x\partial x}$$

Using equations (5.2.23) and (5.3.2) we get

$$
\begin{aligned}
&\partial^m(f(x)x)(h_1;...;h_{m-1};h_m)\\
=&\partial^m f(x)(h_1;...;h_{m-1};h_m)x+\partial^{m-1}f(x)(h_1;...;h_{m-2};h_{m-1})h_m\\
&+\partial^{m-1}f(x)(h_1;...;h_{m-2};\widehat{h_{m-1}};h_m)h_{m-1}\\
&+\partial^{m-2}f(x)(\widehat{h_1},...,h_{m-2},h_{m-1};h_m)h_1+...\\
&+\partial^{m-2}f(x)(h_1,...,\widehat{h_{m-2}},h_{m-1};h_m)h_{m-2}
\end{aligned}
\tag{8.2.3}
$$

The difference between equations (8.2.1) and (8.2.3) is only in form of presentation. We proved the theorem. □

Theorem 8.2.2. *The Gâteaux derivative $\partial^m p_n(x)(h_1,...,h_m)$ is symmetric polynomial with respect to variables $h_1$, ..., $h_m$.*

Proof. To prove the theorem we consider algebraic properties of the Gâteaux derivative and give equivalent definition. We start from construction of monomial. For any monomial $p_n(x)$ we build symmetric polynomial $r_n(x)$ according to following rules

- If $p_1(x)=a_0xa_1$, then $r_1(x_1)=a_0x_1a_1$
- If $p_n(x)=p_{n-1}(x)a_n$, then

  $$r_n(x_1,...,x_n)=r_{n-1}(x_{[1},...,x_{n-1})x_{n]}a_n$$

  where square brackets express symmetrization of expression with respect to variables $x_1$, ..., $x_n$.

It is evident that

$$p_n(x)=r_n(x_1,...,x_n)\quad x_1=...=x_n=x$$

We define the Gâteaux derivative of power $k$ according to rule

$$\partial^k p_n(x)(h_1,...,h_k)=r_n(h_1,...,h_k,x_{k+1},...,x_n)\quad x_{k+1}=x_n=x \tag{8.2.4}$$

According to construction, polynomial $r_n(h_1,...,h_k,x_{k+1},...,x_n)$ is symmetric with respect to variables $h_1$, ..., $h_k$, $x_{k+1}$, ..., $x_n$. Therefore, polynomial (8.2.4) is symmetric with respect to variables $h_1$, ..., $h_k$.

For $k=1$, we will prove that definition (8.2.4) of the Gâteaux derivative coincides with definition (5.2.18).

For $n=1$, $r_1(h_1)=a_0h_1a_1$. This expression coincides with expression of the Gâteaux derivative in theorem 5.3.3.

Let the statement be true for $n-1$. The following equation is true

$$r_n(h_1,x_2,...,x_n)=r_{n-1}(h_1,x_{[2},...,x_{n-1})x_{n]}a_n+r_{n-1}(x_2,...,x_n)h_1a_n \tag{8.2.5}$$

Assume $x_2=...=x_n=x$. According to suggestion of induction, from equations (8.2.4), (8.2.5) it follows that

$$r_n(h_1,x_2,...,x_n)=\partial p_{n-1}(x)(h_1)xa_n+p_{n-1}(x)h_1a_n$$

According to theorem 8.2.1

$$r_n(h_1,x_2,...,x_n)=\partial p_n(x)(h_1)$$

This proves the equation (8.2.4) for $k=1$.

Let us prove now that definition (8.2.4) of the Gâteaux derivative coincides with definition (8.1.4) for $k>1$.

Let equation (8.2.4) be true for $k-1$. Consider arbitrary monomial of polynomial $r_n(h_1,...,h_{k-1},x_k,...,x_n)$. Identifying variables $h_1$, ..., $h_{k-1}$ with elements of division ring $D$, we consider polynomial

$$R_{n-k}(x_k,...,x_n)=r_n(h_1,...,h_{k-1},x_k,...,x_n) \tag{8.2.6}$$

Assume $P_{n-k}(x)=R_{n-k}(x_k,...,x_n)$, $x_k=...=x_n=x$. Therefore

$$P_{n-k}(x)=\partial^{k-1}p_n(x)(h_1;...;h_{k-1})$$

According to definition of the Gâteaux derivative (8.1.4)

$$\begin{aligned}\partial P_{n-k}(x)(h_k)=&\partial(\partial^{k-1}p_n(x)(h_1;...;h_{k-1}))(h_k)\\=&\partial^k p_n(x)(h_1;...;h_{k-1};h_k)\end{aligned} \tag{8.2.7}$$

According to definition (8.2.4) of the Gâteaux derivative

$$\partial P_{n-k}(x)(h_k)=R_{n-k}(h_k,x_{k+1},...,x_n)\quad x_{k+1}=x_n=x \tag{8.2.8}$$

According to definition (8.2.6), from equation (8.2.8) it follows that

$$\partial P_{n-k}(x)(h_k)=r_n(h_1,...,h_k,x_{k+1},...,x_n)\quad x_{k+1}=x_n=x \tag{8.2.9}$$

From comparison of equations (8.2.7) and (8.2.9) it follows that

$$\partial^k p_n(x)(h_1;...;h_k)=r_n(h_1,...,h_k,x_{k+1},...,x_n)\quad x_{k+1}=x_n=x$$

Therefore equation (8.2.4) is true for any $k$ and $n$.

We proved the statement of theorem. □

THEOREM 8.2.3. *For any $n\ge 0$ following equation is true*

$$\partial^{n+1}p_n(x)(h_1;...;h_{n+1})=0 \tag{8.2.10}$$

PROOF. Since $p_0(x)=a_0$, $a_0\in D$, then for $n=0$ theorem is corollary of theorem 5.3.1. Let statement of theorem is true for $n-1$. According to theorem 8.2.1 when $f(x)=p_{n-1}(x)$ we get

$$\begin{aligned}\partial^{n+1}p_n(x)(h_1;...;h_{n+1})=&\partial^{n+1}(p_{n-1}(x)xa_n)(h_1;...;h_{n+1})\\=&\partial^{n+1}p_{n-1}(x)(h_1;...;h_m)xa_n\\&+\partial^n p_{n-1}(x)(h_1;...;h_{m-1})h_ma_n\\&+\partial^n p_{n-1}(x)(\widehat{h_1};...;h_{m-1};h_m)h_1a_n+...\\&+\partial^{m-1}p_{n-1}(x)(h_1;...;\widehat{h_{m-1}};h_m)h_{m-1}a_n\end{aligned}$$

According to suggestion of induction all monomials are equal 0. □

THEOREM 8.2.4. *If $m<n$, then following equation is true*

$$\partial^m p_n(0)(h)=0 \tag{8.2.11}$$

PROOF. For $n=1$ following equation is true

$$\partial^0 p_1(0)=a_0xa_1=0$$

Assume that statement is true for $n-1$. Then according to theorem 8.2.1

$$\begin{aligned}&\partial^m(p_{n-1}(x)xa_n)(h_1;...;h_m)\\=&\partial^m p_{n-1}(x)(h_1;...;h_m)xa_n+\partial^{m-1}p_{n-1}(x)(h_1;...;h_{m-1})h_ma_n\\&+\partial^{m-1}p_{n-1}(x)(\widehat{h_1};...;h_{m-1};h_m)h_1a_n+...\\&+\partial^{m-1}p_{n-1}(x)(h_1;...;\widehat{h_{m-1}};h_m)h_{m-1}a_n\end{aligned}$$

First term equal 0 because $x=0$. Because $m-1<n-1$, then rest terms equal 0 according to suggestion of induction. We proved the statement of theorem. □

When $h_1=...=h_n=h$, we assume

$$\partial^n f(x)(h)=\partial^n f(x)(h_1;...;h_n)$$

This notation does not create ambiguity, because we can determine function according to number of arguments.

THEOREM 8.2.5. *For any $n>0$ following equation is true*

$$\partial^n p_n(x)(h)=n!p_n(h) \tag{8.2.12}$$

PROOF. For $n=1$ following equation is true

$$\partial p_1(x)(h)=\partial(a_0xa_1)(h)=a_0ha_1=1!p_1(h)$$

Assume the statement is true for $n-1$. Then according to theorem 8.2.1

$$\begin{aligned}\partial^n p_n(x)(h)=&\partial^n p_{n-1}(x)(h)xa_n+\partial^{n-1}p_{n-1}(x)(h)ha_n\\&+...+\partial^{n-1}p_{n-1}(x)(h)ha_n\end{aligned} \tag{8.2.13}$$

First term equal 0 according to theorem 8.2.3. The rest $n$ terms equal, and according to suggestion of induction from equation (8.2.13) it follows

$$\partial^n p_n(x)(h)=n\partial^{n-1}p_{n-1}(x)(h)ha_n=n(n-1)!p_{n-1}(h)ha_n=n!p_n(h)$$

Therefore, statement of theorem is true for any $n$. □

Let $p(x)$ be polynomial of power $n$.[8.2]

$$p(x)=p_0+p_{1i_1}(x)+...+p_{ni_n}(x)$$

We assume sum by index $i_k$ which enumerates terms of power $k$. According to theorem 8.2.3, 8.2.4, 8.2.5

$$\partial^k p(0)(x)=k!p_{ki_k}(x)$$

Therefore, we can write

$$p(x)=p_0+(1!)^{-1}\partial p(0)(x)+(2!)^{-1}\partial^2 p(0)(x)+...+(n!)^{-1}\partial^n p(0)(x)$$

This representation of polynomial is called **Taylor polynomial**. If we consider substitution of variable $x=y-y_0$, then considered above construction remain true for polynomial

$$p(y)=p_0+p_{1i_1}(y-y_0)+...+p_{ni_n}(y-y_0)$$

Therefore

$$p(y)=p_0+(1!)^{-1}\partial p(y_0)(y-y_0)+(2!)^{-1}\partial^2 p(y_0)(y-y_0)+...+(n!)^{-1}\partial^n p(y_0)(y-y_0)$$

[8.2] I consider Taylor polynomial for polynomials by analogy with construction of Taylor polynomial in [7], p. 246.

Assume that function $f(x))$ is differentiable in the Gâteaux sense at point $x_0$ up to any order.[8.3]

Theorem 8.2.6. *If function $f(x)$ holds*

(8.2.14) $$f(x_0) = \partial f(x_0)(h) = ... = \partial^n f(x_0)(h) = 0$$

*then for $t \to 0$ expression $f(x+th)$ is infinitesimal of order higher then $n$ with respect to $t$*

$$f(x_0 + th) = o(t^n)$$

Proof. When $n = 1$ this statement follows from equation (5.2.13).

Let statement be true for $n-1$. Map

$$f_1(x) = \partial f(x)(h)$$

satisfies to condition

$$f_1(x_0) = \partial f_1(x_0)(h) = ... = \partial^{n-1} f_1(x_0)(h) = 0$$

According to suggestion of induction

$$f_1(x_0 + th) = o(t^{n-1})$$

Then equation (5.2.12) gets form

$$o(t^{n-1}) = \lim_{t\to 0,\ t\in R} (t^{-1} f(x+th))$$

Therefore,

$$f(x+th) = o(t^n)$$

□

Let us form polynomial

(8.2.15) $$p(x) = f(x_0) + (1!)^{-1}\partial f(x_0)(x-x_0) + ... + (n!)^{-1}\partial^n f(x_0)(x-x_0)$$

According to theorem 8.2.6

$$f(x_0 + t(x-x_0)) - p(x_0 + t(x-x_0)) = o(t^n)$$

Therefore, polynomial $p(x)$ is good approximation of map $f(x)$.

If map $f(x)$ has the Gâteaux derivative of any order, the passing to the limit $n \to \infty$, we get expansion into series

$$f(x) = \sum_{n=0}^{\infty} (n!)^{-1} \partial^n f(x_0)(x-x_0)$$

which is called **Taylor series**.

---

[8.3] I explore construction of Taylor series by analogy with construction of Taylor series in [7], p. 248, 249.

## 8.3. Integral

Concept of integral has different aspect. In this section we consider integration as operation inverse to differentiation. As a matter of fact, we consider procedure of solution of ordinary differential equation

$$\partial f(x)(h) = F(x;h)$$

Example 8.3.1. I start from example of differential equation over real field.

$$y' = 3x^2 \tag{8.3.1}$$

$$x_0 = 0 \quad y_0 = 0 \tag{8.3.2}$$

Differentiating one after another equation (8.3.1), we get the chain of equations

$$y'' = 6x \tag{8.3.3}$$

$$y''' = 6 \tag{8.3.4}$$

$$y^{(n)} = 0 \qquad n > 3 \tag{8.3.5}$$

From equations (8.3.1), (8.3.2), (8.3.3), (8.3.4), (8.3.5) it follows expansion into Taylor series

$$y = x^3$$

□

Example 8.3.2. Consider similar differential equation over division ring

$$\partial(y)(h) = hx^2 + xhx + x^2h \tag{8.3.6}$$

$$x_0 = 0 \quad y_0 = 0 \tag{8.3.7}$$

Differentiating one after another equation (8.3.6), we get the chain of equations

$$\begin{aligned}&\partial^2(y)(h_1;h_2)\\ =&h_1h_2x + h_1xh_2 + h_2h_1x + xh_1h_2 + h_2xh_1 + xh_2h_1\end{aligned} \tag{8.3.8}$$

$$\begin{aligned}&\partial^3(y)(h_1;h_2;h_3)\\ =&h_1h_2h_3 + h_1h_3h_2 + h_2h_1h_3 + h_3h_1h_2 + h_2h_3h_1 + h_3h_2h_1\end{aligned} \tag{8.3.9}$$

$$\partial^n(y)(h_1;...;h_n) = 0 \quad n > 3 \tag{8.3.10}$$

From equations (8.3.6), (8.3.7), (8.3.8), (8.3.9), (8.3.10) expansion into Taylor series follows

$$y = x^3$$

□

Remark 8.3.3. Differential equation

$$\partial(y)(h) = 3hx^2 \tag{8.3.11}$$

$$x_0 = 0 \quad y_0 = 0 \tag{8.3.12}$$

also leads to answer $y = x^3$. it is evident that this map does not satisfies differential equation. However, contrary to theorem 8.2.2 second derivative is not symmetric polynomial. This means that equation (8.3.11) does not possess a solution. □

Example 8.3.4. It is evident that, if function satisfies to differential equation

$$\partial(y)(h) = f_{s\cdot 0}\ h\ f_{s\cdot 1} \tag{8.3.13}$$

then The Gâteaux derivative of second order

$$\partial^2 f(x)(h_1; h_2) = 0$$

In that case, if initial condition is $y(0) = 0$, then differential equation (8.3.13) has solution

$$y = f_{s\cdot 0} x f_{s\cdot 1}$$

□

## 8.4. Exponent

In this section we consider one of possible models of exponent.

In a field we can define exponent as solution of differential equation

$$y' = y \tag{8.4.1}$$

It is evident that we cannot write such equation for division ring. However we can use equation

$$\partial(y)(h) = y'h \tag{8.4.2}$$

From equations (8.4.1), (8.4.2) it follows

$$\partial(y)(h) = yh \tag{8.4.3}$$

This equation is closer to our goal, however there is the question: in which order we should multiply $y$ and $h$? To answer this question we change equation

$$\partial(y)(h) = \frac{1}{2}(yh + hy) \tag{8.4.4}$$

Hence, our goal is to solve differential equation (8.4.4) with initial condition $y(0) = 1$.

For the statement and proof of the theorem 8.4.1 I introduce following notation. Let

$$\sigma = \begin{pmatrix} y & h_1 & ... & h_n \\ \sigma(y) & \sigma(h_1) & ... & \sigma(h_n) \end{pmatrix}$$

be transposition of the tuple of variables

$$\begin{pmatrix} y & h_1 & ... & h_n \end{pmatrix}$$

Let $p_\sigma(h_i)$ be position that variable $h_i$ gets in the tuple

$$\begin{pmatrix} \sigma(y) & \sigma(h_1) & ... & \sigma(h_n) \end{pmatrix}$$

For instance, if transposition $\sigma$ has form

$$\begin{pmatrix} y & h_1 & h_2 & h_3 \\ h_2 & y & h_3 & h_1 \end{pmatrix}$$

then following tuples equal

$$\begin{aligned}\begin{pmatrix}\sigma(y) & \sigma(h_1) & \sigma(h_2) & \sigma(h_3)\end{pmatrix} &= \begin{pmatrix}h_2 & y & h_3 & h_1\end{pmatrix}\\ &= \begin{pmatrix}p_\sigma(h_2) & p_\sigma(y) & p_\sigma(h_3) & p_\sigma(h_1)\end{pmatrix}\end{aligned}$$

THEOREM 8.4.1. *If function $y$ is solution of differential equation (8.4.4) then the Gâteaux derivative of order $n$ of function $y$ has form*

$$\partial^n(y)(h_1,...,h_n) = \frac{1}{2^n}\sum_\sigma \sigma(y)\sigma(h_1)...\sigma(h_n) \tag{8.4.5}$$

*where sum is over transpositions*

$$\sigma = \begin{pmatrix} y & h_1 & ... & h_n\\ \sigma(y) & \sigma(h_1) & ... & \sigma(h_n)\end{pmatrix}$$

*of the set of variables $y$, $h_1$, ..., $h_n$. Transposition $\sigma$ has following properties*

(1) *If there exist $i$, $j$, $i \ne j$, such that $p_\sigma(h_i)$ is situated in product (8.4.5) on the left side of $p_\sigma(h_j)$ and $p_\sigma(h_j)$ is situated on the left side of $p_\sigma(y)$, then $i<j$.*
(2) *If there exist $i$, $j$, $i \ne j$, such that $p_\sigma(h_i)$ is situated in product (8.4.5) on the right side of $p_\sigma(h_j)$ and $p_\sigma(h_j)$ is situated on the right side of $p_\sigma(y)$, then $i>j$.*

PROOF. We prove this statement by induction. For $n=1$ the statement is true because this is differential equation (8.4.4). Let the statement be true for $n=k-1$. Hence

$$\partial^{k-1}(y)(h_1,...,h_{k-1}) = \frac{1}{2^{k-1}}\sum_\sigma \sigma(y)\sigma(h_1)...\sigma(h_{k-1}) \tag{8.4.6}$$

where the sum is over transposition

$$\sigma = \begin{pmatrix} y & h_1 & ... & h_{k-1}\\ \sigma(y) & \sigma(h_1) & ... & \sigma(h_{k-1})\end{pmatrix}$$

of the set of variables $y$, $h_1$, ..., $h_{k-1}$. Transposition $\sigma$ satisfies to conditions (1), (2) in theorem. According to definition (8.1.4) the Gâteaux derivative of order $k$ has form

$$\begin{aligned}\partial^k(y)(h_1,...,h_k) =&\partial(\partial^{k-1}(y)(h_1,...,h_{k-1}))(h_k)\\ =&\frac{1}{2^{k-1}}\partial\left(\sum_\sigma \sigma(y)\sigma(h_1)...\sigma(h_{k-1})\right)(h_k)\end{aligned} \tag{8.4.7}$$

From equations (8.4.4), (8.4.7) it follows that

$$\begin{aligned}&\partial^k(y)(h_1,...,h_k)\\ =&\frac{1}{2^{k-1}}\frac{1}{2}\left(\sum_\sigma \sigma(yh_k)\sigma(h_1)...\sigma(h_{k-1}) + \sum_\sigma \sigma(h_ky)\sigma(h_1)...\sigma(h_{k-1})\right)\end{aligned} \tag{8.4.8}$$

It is easy to see that arbitrary transposition $\sigma$ from sum (8.4.8) forms two transpositions

(8.4.9)
$$\begin{aligned}\tau_1 &= \begin{pmatrix} y & h_1 & ... & h_{k-1} & h_k \\ \tau_1(y) & \tau_1(h_1) & ... & \tau_1(h_{k-1}) & \tau_1(h_k)\end{pmatrix}\\ &= \begin{pmatrix} h_k y & h_1 & ... & h_{k-1} \\ \sigma(h_k y) & \sigma(h_1) & ... & \sigma(h_{k-1})\end{pmatrix}\\ \tau_2 &= \begin{pmatrix} y & h_1 & ... & h_{k-1} & h_k \\ \tau_2(y) & \tau_2(h_1) & ... & \tau_2(h_{k-1}) & \tau_2(h_k)\end{pmatrix}\\ &= \begin{pmatrix} y h_k & h_1 & ... & h_{k-1} \\ \sigma(y h_k) & \sigma(h_1) & ... & \sigma(h_{k-1})\end{pmatrix}\end{aligned}$$

From (8.4.8) and (8.4.9) it follows that

(8.4.10)
$$\begin{aligned}&\partial^k(y)(h_1,...,h_k)\\ =&\frac{1}{2^k}\left(\sum_{\tau_1}\tau_1(y)\tau_1(h_1)...\tau_1(h_{k-1})\tau_1(h_k)\right.\\ &\left.+\sum_{\tau_2}\tau_2(y)\tau_2(h_1)...\tau_2(h_{k-1})\tau_2(h_k)\right)\end{aligned}$$

In expression (8.4.10) $p_{\tau_1}(h_k)$ is written immediately before $p_{\tau_1}(y)$. Since $k$ is smallest value of index then transposition $\tau_1$ satisfies to conditions (1), (2) in the theorem. In expression (8.4.10) $p_{\tau_2}(h_k)$ is written immediately after $p_{\tau_2}(y)$. Since $k$ is largest value of index then transposition $\tau_2$ satisfies to conditions (1), (2) in the theorem.

It remains to show that in the expression (8.4.10) we get all transpositions $\tau$ that satisfy to conditions (1), (2) in the theorem. Since $k$ is largest index then according to conditions (1), (2) in the theorem $\tau(h_k)$ is written either immediately before or immediately after $\tau(y)$. Therefore, any transposition $\tau$ has either form $\tau_1$ or form $\tau_2$. Using equation (8.4.9), we can find corresponding transposition $\sigma$ for given transposition $\tau$. Therefore, the statement of theorem is true for $n = k$. We proved the theorem. □

THEOREM 8.4.2. *The solution of differential equation (8.4.4) with initial condition $y(0) = 1$ is exponent $y = e^x$ that has following Taylor series expansion*

(8.4.11)
$$e^x = \sum_{n=0}^{\infty}\frac{1}{n!}x^n$$

PROOF. The Gâteaux derivative of order $n$ has $2^n$ items. In fact, the Gâteaux derivative of order 1 has 2 items, and each differentiation increase number of items twice. From initial condition $y(0) = 1$ and theorem 8.4.1 it follows that the Gâteaux derivative of order $n$ of required solution has form

(8.4.12)
$$\partial^n(y(0))(h,...,h) = 1$$

Taylor series expansion (8.4.11) follows from equation (8.4.12). □

Theorem 8.4.3. *The equation*

$$e^{a+b} = e^a e^b \tag{8.4.13}$$

*is true iff*

$$ab = ba \tag{8.4.14}$$

Proof. To prove the theorem it is enough to consider Taylor series

$$e^a = \sum_{n=0}^{\infty} \frac{1}{n!} a^n \tag{8.4.15}$$

$$e^b = \sum_{n=0}^{\infty} \frac{1}{n!} b^n \tag{8.4.16}$$

$$e^{a+b} = \sum_{n=0}^{\infty} \frac{1}{n!} (a+b)^n \tag{8.4.17}$$

Let us multiply expressions (8.4.15) and (8.4.16). The sum of monomials of order 3 has form

$$\frac{1}{6}a^3 + \frac{1}{2}a^2 b + \frac{1}{2}ab^2 + \frac{1}{6}b^3 \tag{8.4.18}$$

and in general does not equal expression

$$\frac{1}{6}(a+b)^3 = \frac{1}{6}a^3 + \frac{1}{6}a^2 b + \frac{1}{6}aba + \frac{1}{6}ba^2 + \frac{1}{6}ab^2 + \frac{1}{6}bab + \frac{1}{6}b^2 a + \frac{1}{6}b^3 \tag{8.4.19}$$

The proof of statement that (8.4.13) follows from (8.4.14) is trivial. □

The meaning of the theorem 8.4.3 becomes more clear if we recall that there exist two models of design of exponent. First model is the solution of differential equation (8.4.4). Second model is exploring of one parameter group of transformations. For field both models lead to the same function. I cannot state this now for general case. This is the subject of separate research. However if we recall that quaternion is analogue of transformation of three dimensional space then the statement of the theorem becomes evident.

CHAPTER 9

# Derivative of Second Order of Map of $D$-Vector Space

## 9.1. Derivative of Second Order of Map of $D$-Vector Space

Let $D$ be valued division ring of characteristic 0. Let $V$ be normed $D$-vector space with norm $\|y\|_V$. Let $W$ be normed $D$-vector space with norm $\|z\|_W$. Let

$$f : V \to W$$

function differentiable in the Gâteaux sense. According to remark 6.2.2 the Gâteaux derivative is map

$$\partial f : V \to \mathcal{L}(D; V; W)$$

According to theorems 3.1.5, and definition 6.1.7 set $\mathcal{L}(D; V; W)$ is normed $D$-vector space. Therefore, we may consider the question, if map $\partial f$ is differentiable in the Gâteaux sense.

According to the definition 6.2.1,

$$\partial f(x + a_2)(a_1) - \partial f(x)(a_1) = \partial(\partial f(x)(a_1))(a_2) + o_2(a_2) \tag{9.1.1}$$

where $o_2 : V \to \mathcal{L}(D; V; W)$ is such continuous map, that

$$\lim_{a_2 \to 0} \frac{\|o_2(a_2)\|_W}{\|a_2\|_V} = 0$$

According to the definition 6.2.1, the map $\partial(\partial f(x)(a_1))(a_2)$ is linear map of variable $a_2$. From equation (9.1.1), it follows that map $\partial(\partial f(x)(a_1))(a_2)$ is linear map of variable $a_1$.

Definition 9.1.1. Polylinear map

$$\partial^2 f(x)(a_1; a_2) = \partial(\partial f(x)(a_1))(a_2) \tag{9.1.2}$$

is called **the Gâteaux derivative of second order** of map $f$. □

Remark 9.1.2. According to definition 9.1.1 for given $x$ the Gâteaux derivative of second order $\partial^2 f(x) \in \mathcal{L}(D; V, V; W)$. Therefore, the Gâteaux derivative of second order of map $f$ is map

$$\partial^2 f : V \to \mathcal{L}(D; V, V; W)$$

Theorem 9.1.3. *Let $\overline{\overline{e}}_V$ be basis of $D$-vector space $V$. Let $\overline{\overline{e}}_W$ be basis of $D$-vector space $W$.* **The Gâteaux differential of second order** *of map $f$ relative to basis $\overline{\overline{e}}_V$ and basis $\overline{\overline{e}}_W$ has form*

$$\partial^2 f(x)(a_1; a_2) = \frac{\partial^2_{s\cdot 0} f^j(x)}{\partial x^{i_1} \partial x^{i_2}} \sigma_s(a_1^{i_1}) \frac{\partial^2_{s\cdot 1} f^j(x)}{\partial x^{i_1} \partial x^{i_2}} \sigma_s(a_2^{i_2}) \frac{\partial^2_{s\cdot 2} f^j(x)}{\partial x^{i_1} \partial x^{i_2}} e_{W\cdot j} \tag{9.1.3}$$

Proof. Corollary of definition 9.1.1 and theorem 4.2.3. □

Definition 9.1.4. Expression $\frac{\partial^2_{s\cdot p} f^j(x)}{\partial x^{i_1}\partial x^{i_2}}$, $p = 0, 1, 2,$ is called **component of the Gâteaux derivative of second order** of map $f(x)$. □

Definition 9.1.5. Let $D$ be division ring of characteristic 0. Let $\overline{\overline{e}}_V$ be basis of $D$-vector space $V$. Let $\overline{\overline{e}}_W$ be basis of $D$-vector space $W$. The Gâteaux partial derivative

$$\frac{\partial^2 f^k(v)}{\partial v^j \partial v^i}(h^i, h^j) = \frac{\partial}{\partial v^j}\left(\frac{\partial f^k(v)}{\partial v^i}(h^i)\right)(h^j) \tag{9.1.4}$$

is called **the Gâteaux mixed partial derivative** of map $f^k$ with respect to variables $v^i$, $v^j$. □

Theorem 9.1.6. *Let $D$ be division ring of characteristic 0.*[9.1] *Let $\overline{\overline{e}}_V$ be basis of $D$-vector space $V$. Let $\overline{\overline{e}}_W$ be basis of $D$-vector space $W$. Let the Gâteaux partial derivatives of map*

$$f : V \rightarrow W$$

*are continuous and differentiable in the Gâteaux sense on the set $U \subset V$. Let the Gâteaux mixed partial derivatives of second order is continuous on the set $U \subset V$. Then on the set $U$ the Gâteaux mixed partial derivatives satisfy equation*

$$\frac{\partial^2 f^k(v)}{\partial v^j \partial v^i}(h^i, h^j) = \frac{\partial^2 f^k(v)}{\partial v^i \partial v^j}(h^j, h^i) \tag{9.1.5}$$

Proof. To prove theorem we will use notation $f^k(v) = f^k(v^i, v^j)$, assuming, that rest variables have given values. The statement of theorem follows from following chain of equations

$$\begin{aligned}
&\frac{\partial^2 f(v^i, v^j)}{\partial v^i \partial v^j}(h^j, h^i)\\
=&\frac{\partial}{\partial v^i}\left(\frac{\partial f(v^i, v^j)}{\partial v^j}(h^j)\right)(h^i)\\
=&\lim_{t\to 0} t^{-1}\left(\frac{\partial f(v^i + h^i t, v^j)}{\partial v^j}(h^j) - \frac{\partial f(v^i, v^j)}{\partial v^j}(h^j)\right)\\
=&\lim_{t\to 0} t^{-2}\left(f(v^i + h^i t, v^j + h^j t) - f(v^i + h^i t, v^j) - f(v^i, v^j + h^j t) + f(v^i, v^j)\right)\\
=&\lim_{t\to 0} t^{-2}\left(f(v^i + h^i t, v^j + h^j t) - f(v^i, v^j + h^j t) - f(v^i + h^i t, v^j) + f(v^i, v^j)\right)\\
=&\lim_{t\to 0} t^{-1}\left(\frac{\partial f(v^i, v^j + h^j t)}{\partial v^i}(h^i) - \frac{\partial f(v^i, v^j)}{\partial v^i}(h^i)\right)\\
=&\frac{\partial}{\partial v^j}\left(\frac{\partial f(v^i, v^j)}{\partial v^i}(h^i)\right)(h^j)\\
=&\frac{\partial^2 f(v^i, v^j)}{\partial v^j \partial v^i}(h^i, h^j)
\end{aligned}$$

□

By induction, assuming that the Gâteaux derivative $\partial^{n-1} f(x)$ of order $n-1$ is defined, we define

$$\partial^n f(x)(a_1; ...; a_n) = \partial(\partial^{n-1} f(x)(a_1; ...; a_{n-1}))a_n$$

[9.1] I proved the theorem the same way like it is done for theorem in [7], p. 405, 406.

**the Gâteaux derivative of order** $n$ of map $f$. We also assume $\partial^0 f(x) = f(x)$.

## 9.2. Homogeneous Map

DEFINITION 9.2.1. The map

$$f : V \to W$$

of $D$-vector space $V$ into $D$-vector space $W$ is called **homogeneous map of degree** $k$ **over field** $F$, if for any $a \in F$

$$f(av) = a^k f(v) \tag{9.2.1}$$

□

THEOREM 9.2.2 (Euler's theorem). *Let map*

$$f : V \to W$$

*of $D$-vector space $V$ into $D$-vector space $W$ be homogeneous map of degree $k$ over field $F \subset Z(D)$.*

$$\partial f(v)(v) = k f(v) \tag{9.2.2}$$

PROOF. Let us differentiate equation (9.2.1) with respect to $a$[9.2]

$$\frac{d}{da}(f(av))\Delta a = \frac{d}{da}(a^k f(v))\Delta a \tag{9.2.3}$$

From theorem 6.2.12 it follows

$$\begin{aligned} \frac{d}{da}(f(av))\Delta a &= \partial f(av)(\partial(av)(\Delta a)) \\ &= \partial f(av)\left(\frac{d(av)}{da}\Delta a\right) \\ &= \partial f(av)(v\Delta a) \end{aligned} \tag{9.2.4}$$

From theorem 6.2.5 and equation (9.2.4) it follows

$$\frac{d}{da}(f(av))\Delta a = \partial f(av)(v)\Delta a \tag{9.2.5}$$

It is evident that

$$\frac{d}{da}(a^k f(v)) = k a^{k-1} f(v) \tag{9.2.6}$$

Let us substitute (9.2.5) and (9.2.6) in (9.2.3)

$$\partial f(av)(v)\Delta a = k a^{k-1} f(v)\Delta a \tag{9.2.7}$$

Because $\Delta a \neq 0$, then (9.2.2) follows from (9.2.7) under the condition $a = 1$. □

Let map

$$f : V \to W$$

from $D$-vector space $V$ into $D$-vector space $W$ be homogeneous map of degree $k$ over field $F \subset Z(D)$. From theorem 6.2.5 it follows

$$\partial f(av)(av) = a\partial f(av)(v) \tag{9.2.8}$$

---

[9.2] For map of field it is true that

$$\partial f(x)(h) = \frac{df(x)}{dx}h$$

From equations (9.2.2), (9.2.1) it follows

$$\partial f(av)(av) = kf(av) = ka^k f(v) = a^k \partial f(v)(v) \tag{9.2.9}$$

From equations (9.2.8), (9.2.9) it follows

$$\partial f(av)(v) = a^{k-1} \partial f(v)(v) \tag{9.2.10}$$

However from equation (9.2.10) it does not follow that map

$$\partial f : V \to \mathcal{A}(D; V; W)$$

is homogeneous map of degree $k-1$ over field $F$.

CHAPTER 10

CHAPTER 11

# Index

CHAPTER 12

# Special Symbols and Notations

# Введение в математический анализ над телом

Александр Клейн

*E-mail address*: Aleks_Kleyn@MailAPS.org
*URL*: http://sites.google.com/site/alekskleyn/
*URL*: http://arxiv.org/a/kleyn_a_1
*URL*: http://AleksKleyn.blogspot.com/

Аннотация. В книге рассматриваются парные представления тела в абелевой группе и $D$-векторные пространства над телом. Морфизмы $D$-векторных пространств являются линейными отображениями $D$-векторных пространств.

Я изучаю производную функции $f$ нормированых тел как линейное отображение, наиболее близкое к функции $f$. Я изучаю разложение отображения в ряд Тейлора и метод решения дифференциального уравнения.

Норма в $D$-векторном пространстве позволяет рассматривать непрерывные отображения $D$-векторных пространств. Дифференциал отображения $f$ $D$-векторных пространств определён как линейное отображение, наиболее близкое к отображению $f$.

# Оглавление

Глава 1

# Предисловие

## 1.1. Предисловие к изданию 1

> Дорога очарований и разочарований.
> Но я бы сказал наоборот. Когда ты понимаешь, что причина твоих разочарований - в твоих ожиданиях, ты начинаешь присматриваться к окружающему ландшафту. Ты оказался на тропе, по которой до тебя никто не ходил. И новые впечатления куда сильнее, чем если бы ты шёл по дороге многократно протоптанной и хорошо изученной.
>
> Автор неизвестен. Записки путешественника.

В математике, как и в обычной жизни, есть утверждения, которые после изучения становятся очевидными. И когда эти утверждения в новых условиях не работают, это в первый момент вызывает удивление.[1.1]

Когда я начал изучать линейную алгебру над телом, я был готов к разным сюрпризам, связанными с некоммутативностью умножения. Но большое количество утверждений, похожих на утверждения линейной алгебры над полем, повидимому ослабили моё внимание. И тем не менее, уже в линейной алгебре я встретил первую серьёзную проблему: определение полилинейных форм и тензорного произведения подразумевает коммутативность произведения. Твёрдая уверенность, что решение конкретных задач даст ключ к решению общей задачи, была основным мотивом, побуждающим меня изучать математический анализ над телом. Предлагаемая книга вводит в методы дифференцирования отображений над телом и является естественным следствием книги [3], поскольку в основе математического анализа лежит возможность линейного приближения к отображению. В анализе, как и в линейной алгебре, мы встречаем утверждения, которые нам хорошо знакомы. Однако есть немало новых и неожиданных утверждений.

Первое определение производной, которое я рассмотрел, было основано на классическом определении производной. Однако при попытке рассмотреть

[1.1]Подобного рода ощущение я испытал, когда пытался сформулировать теорему Ролля при изучении математического анализа над телом. После многократных попыток изменить формулировку теоремы Ролля я понял, что это утверждение неверно даже в случае комплексных чисел. Достаточно рассмотреть функцию $w = z^2$. $z = 0$ - единственная точка, где производная равна 0. Однако точки $z = 0.5$ и $z = -0.5$ могут быть соединены кривой, не проходящей через точку $z = 0$.

производную функции $x^2$ я обнаружил, что не могу эту производную выразить аналитически. Чтобы ответить на вопрос, существует ли производная, я решил провести вычислительный эксперимент. Самый лучший кандидат для вычислительного эксперимента - это тело кватернионов. Несложные расчёты показали, что в некоторых случаях производная действительно не существует. Это приводит к тому, что множество дифференцируемых функций крайне бедно, а сама теория дифференцирования не представляет серьёзного интереса.[1.2] Анализ вычислительного эксперимента обнаруживает, что даже если производная не существует, возникает выражение, зависящее от приращения аргумента. Подобная конструкция известна в функциональном анализе. Речь идёт о производных Фреше и Гато в пространстве линейных операторов.

Когда я определил $D\star$-производную[1.3] Гато, я обратил внимание, что дифференциал не зависит от того, использую ли я $D\star$-производную Гато или $\star D$-производную Гато. Поэтому я изменил характер изложения и начал с определения дифференциала.

Концепция $D_*{}^*$-векторного пространства является мощным инструментом изучения линейных пространств над телом. Однако дифференциал Гато не укладывается в рамках $D_*{}^*$-векторного пространства. Пример 5.2.2 приводит к выводу, что $D$-векторное пространство - наиболее адекватная модель для изучения дифференциала Гато.

Изучение $D_*{}^*$-линейных отображений приводит к концепции парных представлений тела в абелевой группе. Абелева группа, в которой определено парное представление тела $D$, называется $D$-векторным пространством. Поскольку умножение на элементы тела $D$ определено и слева, и справа, гомоморфизмы $D$-векторных пространств не могут сохранять эту операцию. Это приводит к концепции аддитивного отображения, которое является морфизмом $D$-векторных пространств.[1.4]

Рассмотрение дифференциала Гато как аддитивное отображение приводит к новой конструкции для определения производной. В отличии от $D_*{}^*$-производной Гато компоненты дифференциала Гато не зависят от направления. Однако компоненты дифференциала Гато определены неоднозначно и неопределенно их количество. Очевидно, что, когда мы рассматриваем математический анализ над полем, $D\star$-производная Гато и компоненты дифференциала Гато совпадают с производной. Эти конструкции возникают в результате различных построений. Когда произведение коммутативно эти конструкции совпадают. Но как только произведение становится не коммутативным, ожидаемая

[1.2]Поиск в статьях и книгах показал, что это действительно серьёзная проблема. Математики пытаются решать эту задачу, предлагая разные определения производной, основанные на тех или иных свойствах производной функции действительного или комплексного переменного. Эти конструкции очень интересны. Но моя задача рассмотреть общие закономерности для дифференцирования отображений полного тела характеристики 0.

[1.3]В этой книге я использую обозначения, введенные в [3].

[1.4]Построение аддитивного отображения, выполненное при доказательстве теоремы 4.1.4, весьма поучительно методически. С того самого момента, когда я начал изучать дифференциальную геометрию, я признал факт, что всё сводится к тензорам, что компоненты тензора нумеруются с помощью индексов и эта связь неразрывна. В тот момент, когда это оказалось не так, я оказался не готов признать это. Понадобился месяц прежде, чем я смог записать выражение, которое я видел глазами, но которое я не мог записать до этого, так как пытался выразить его в тензорной или операторной форме.

связь разрушается. Несмотря на различие в поведении, эти конструкции взаимно дополняют друг друга.

Декабрь, 2008

## 1.2. Предисловие к изданию 3

Обсуждая понятие производной, приведенное в этой книге, я понял, что это очень сложная тема. Поэтому в этой части введения я хочу рассмотреть более подробно концепцию дифференцирования.

Когда мы изучаем функции одной переменной, то производная в заданной точке является числом. Когда мы изучаем функцию нескольких переменных, выясняется, что числа недостаточно. Производная становится вектором или градиентом. При изучении отображений векторных пространств мы впервые говорим о производной как об операторе. Но так как этот оператор линеен, то мы можем представить производную как умножение на матрицу, т. е. опять мы выражаем производную как набор чисел.

Без сомнения, подобное поведение производной ослабляет наше внимание. Когда мы переходим к объектам, более сложным чем поля или векторные пространства, мы по прежнему пытаемся увидеть объект, который можно записать как множитель перед приращением и который от приращения не зависит. Алгебра кватернионов - это самая простая модель непрерывного тела. Однако попытка определить производную как линейное отображение обречена на неудачу. Только линейная функция, да и то не любая, дифференцируема в этом смысле.

Если множество линейных отображений ограничивает нашу возможность изучения поведения отображений в малом, существует ли альтернатива. Ответ на этот вопрос положительный. Рассмотрим алгебру кватернионов. Алгебра кватернионов - это нормированное пространство. Мы знаем два типа производных в нормированном пространстве. Сильная производная или производная Фреше является аналогом той производной, к которой мы привыкли. Когда я писал о попытке найти производную как линейное отображение приращения аргумента, я имел в виду производную Фреше. Помимо сильной производной существует слабая производная или производная Гато. Основная идея состоит в том, что производная может зависеть от направления.

Вадим Комков даёт определение производной Гато для тела кватернионов ([9], с. 322) следующим образом. Положим $q_2$ является производной Гато в направлении $q_1$, если для $\epsilon > 0$ справедливо

$$f(q + \epsilon q_1) - f(q) = \epsilon(q_1 q_2)$$

Это определение сделано согласно определению [1]-3.1.2, стр. 256 и не сильно отличается от производной Фреше. Это явление сохраняется и в общем случае. Всякий раз, когда мы пытаемся выделить производную Гато как множитель, происходит отождествление производной Гато и производной Фреше.

Линейное отображение тела - это гомоморфизм его аддитивной группы и гомоморфизм его мультипликативной группы в одном отображении. Ограничивающим фактором является гомоморфизм мультипликативной группы. Поэтому есть смысл отложить в сторону гомоморфизм мультипликативной группы и рассматривать только гомоморфизм аддитивной группы. Очевидно, что множество гомоморфизмов аддитивной группы гораздо шире, чем множество линейных отображений. Однако, так как в теле определены две операции, аддитивное отображение имеет специальную структуру. Поэтому я выделил эти отображения в отдельный класс и назвал их аддитивными отображениями.

Даже в этом случае сохраняются некоторые элементы линейности. Центр тела является полем. В результате аддитивное отображение оказывается линейным над центром. Если тело имеет характеристику 0 и не имеет дискретную топологию, то мы можем вложить поле действительных чисел в центр тела. В результате определение производной как аддитивного отображения оказывается эквивалентным определению производной Гато.

Именно здесь происходит качественный скачок. В аддитивном отображении мы не можем выделить аргумент как множитель. Мы не можем выделить приращение из производной. Является ли эта концепция принципиально новой? Я не думаю. Производная - это отображение. Это отображение дифференциала аргумента в дифференциал функции. Даже аддитивная функция приращения аргумента позволяет аппроксимировать приращение функции полиномом первой степени. Структура полинома над телом отличается от структуры полинома над полем. Я рассматриваю некоторые свойства полинома в разделе 8.2. Опираясь на полученные результаты, я изучаю разложение отображения в ряд Тейлора и метод решения дифференциального уравнения.

Я вначале тоже пытался определить производную отображения над телом как производную Фреше. Но увидев, что эта производная не удовлетворяет основным правилам дифференцирования, я понял, что это неверный путь. Чтобы убедиться, что проблема фундаментальна, а не в логике моих построений, я выполнил расчёт в теле кватернионов. Именно здесь я обнаружил, что производная зависит от направления.

В это время я изучал отображения $D$-векторных пространств, и концепция аддитивного отображения стала ведущей идеей в этом исследовании. В предисловии к изданию 1 (раздел 1.1) я писал о впечатлениях путешественника, оказавшегося впервые в незнакомых краях. Желая увидеть знакомые пейзажи, путешественник не сразу замечает необычность новых мест. Такую смену впечатлений я испытал на протяжении этого года пока я писал книгу.

Я не сразу осознал, что найденная мною производная является производной Гато. Я даже видел различия между моей производной и производной Гато. Потребовался длинный путь, пока я понял, что два определения эквивалентны. Подобно путешественнику, который пытался увидеть знакомые пейзажи, я изучал варианты как можно выделить приращение из производной. В результате появилась $D\star$-производная Гато.

Что уцелеет из введенных понятий в будущем? Без сомнения основными понятиями будет собственно производная Гато и частные производные Гато. Компоненты производной Гато важны для исследования структуры производной Гато. Я не думаю, что $D\star$-производная Гато будет меня интересовать за пределами производных Гато первого порядка. Основное преимущество $D\star$-производной Гато состоит в том, что она позволяет изучить условия, когда отображение непрерывно.

Май, 2009

## 1.3. Предисловие к изданию 4

Алгебра кватернионов - самый простой пример тела. Поэтому я проверяю полученные мной теоремы в алгебре кватернионов, чтобы увидеть как они работают. Алгебра кватернионов похожа на поле комплексных чисел и поэтому естественно искать некоторые параллели.

В поле действительных чисел любое аддитивное отображение автоматически оказывается линейным над полем действительных чисел. Это связано с тем, что поле действительных чисел является пополнением поля рациональных чисел и является следствием теоремы 3.1.3.

Однако в случае комплексных чисел ситуация меняется. Не всякое аддитивное отображение поля комплексных чисел оказывается линейным над полем комплексных чисел. Операция сопряжения - простейший пример такого отображения. Когда я обнаружил этот крайне интересный факт, я вернулся к вопросу об аналитическом представлении аддитивного отображения.

Изучая аддитивные отображения, я понял, что я слишком резко расширил множество линейных отображений при переходе от поля к телу. Причиной этому было не вполне ясное понимание как преодолеть некоммутативность произведения. Однако в процессе построений становилось всё более очевидным, что любое аддитивное отображение линейно над некоторым полем. Впервые эта концепция появилась при построении тензорных произведений и отчётливо проявилась в последующем исследовании.

При построении тензорного произведения тел $D_1$, ..., $D_n$ я предполагаю существование поля $F$, над которым аддитивное преобразование тела $D_i$ линейно для любого $i$. Если все тела имеют характеристику 0, то согласно теореме 3.1.3 такое поле всегда существует. Однако здесь возникает зависимость тензорного произведения от выбранного поля $F$. Чтобы избавиться от этой зависимости, я предполагаю, что поле $F$ - максимальное поле, обладающее указанным свойством.

Если $D_1 = ... = D_n = D$, то такое поле является центром $Z(D)$ тела $D$. Если произведение в теле $D$ коммутативно[1.5], то $Z(D) = D$. Следовательно, начав с аддитивного отображения, я пришёл к концепции линейного отображения, которая является обобщением линейного отображения над полем.

На моё решение изучать производную Гато как линейное отображение повлияли следующие обстоятельства.

- Линейное отображение поля - это в тоже время однородный многочлен первой степени.
- Для непрерывных полей производная функции поля - это однородный многочлен первого порядка, аппроксимирующий изменение функции. Одновременно мы рассматриваем производную как линейное отображение поля. При переходе к телу мы можем сохранить ту же связь производной и аддитивного отображения.

Исследование в области комплексных чисел и кватернионов проявили ещё одно интересное явление. Несмотря на то, что поле комплексных чисел является расширением поля действительных чисел, структура линейного отображения над полем комплексных чисел отличается от структуры линейного отображения над полем действительных чисел. Это различие приводит к тому, что

[1.5] Иными словами, тело является полем.

операция сопряжения комплексных чисел является аддитивным отображением, но не является линейным отображением над полем комплексных чисел.

Аналогично, структура линейного отображения над телом кватернионов отличается от структуры линейного отображения над полем комплексных чисел. Причина различия в том, что центр алгебры кватернионов имеет более простую структуру, чем поле комплексных чисел. Это отличие приводит к тому, что операция сопряжения кватерниона удовлетворяет равенству

$$\overline{p} = -\frac{1}{2}(p + ipi + jpj + kpk)$$

Вследствие этого задача найти отображения, удовлетворяющие теореме подобной теореме Римана (теорема 7.1.1), является нетривиальной задачей для кватернионов.

Август, 2009

## 1.4. Предисловие к изданию 5

Недавно в интернете я нашёл книгу [11], в которой Гамильтон нашёл аналогичное описание дифференциального исчисления в алгебре кватернионов. И хотя эта книга написана более чем 100 лет назад, я рекомендую прочесть эту книгу всем, кому интересны проблемы анализа в алгебре кватернионов.

Декабрь, 2009

## 1.5. Соглашения

Соглашение 1.5.1. *Функция и отображение - синонимы. Однако существует традиция соответствие между кольцами или векторными пространствами называть отображением, а отображение поля действительных чисел или алгебры кватернионов называть функцией. Я тоже следую этой традиции, хотя встречается текст, в котором неясно, какому термину надо отдать предпочтение.* □

Соглашение 1.5.2. *В любом выражении, где появляется индекс, я предполагаю, что этот индекс может иметь внутреннюю структуру. Например, при рассмотрении алгебры $A$ координаты $a \in A$ относительно базиса $\overline{\overline{e}}$ пронумерованы индексом $i$. Это означает, что $a$ является вектором. Однако, если $a$ является матрицей, нам необходимо два индекса, один нумерует строки, другой - столбцы. В том случае, когда мы уточняем структуру индекса, мы будем начинать индекс с символа $\cdot$ в соответствующей позиции. Например, если я рассматриваю матрицу $a^i_j$ как элемент векторного пространства, то я могу записать элемент матрицы в виде $a^{\cdot i}_{\cdot j}$.* □

Соглашение 1.5.3. *В выражении вида*

$$a_{s\cdot 0} x a_{s\cdot 1}$$

*предполагается сумма по индексу $s$.* □

Соглашение 1.5.4. *Тело $D$ можно рассматривать как $D$-векторное пространство размерности 1. Соответственно этому, мы можем изучать не только гомоморфизм тела $D_1$ в тело $D_2$, но и линейное отображение тел.* □

Соглашение 1.5.5. *Несмотря на некоммутативность произведения многие утверждения сохраняются, если заменить например правое представление на левое представление или правое векторное пространство на левое векторное пространство. Чтобы сохранить эту симметрию в формулировках теорем я пользуюсь симметричными обозначениями. Например, я рассматриваю $D\star$-векторное пространство и $\star D$-векторное пространство. Запись $D\star$-векторное пространство можно прочесть как $D$-star-векторное пространство либо как левое векторное пространство. Запись $D\star$-линейно зависимые векторы можно прочесть как $D$-star-линейно зависимые векторы либо как векторы, линейно зависимые слева.* □

Соглашение 1.5.6. *Пусть $A$ - свободная алгебра с конечным или счётным базисом. При разложении элемента алгебры $A$ относительно базиса $\overline{\overline{e}}$ мы пользуемся одной и той же корневой буквой для обозначения этого элемента и его координат. В выражении $a^2$ не ясно - это компонента разложения элемента $a$ относительно базиса или это операция возведения в степень. Для облегчения чтения текста мы будем индекс элемента алгебры выделять цветом. Например,*

$$a = a^{i} e_{i}$$

□

Соглашение 1.5.7. *Очень трудно провести границу между модулем и алгеброй. Тем более, иногда в процессе построения мы должны сперва доказать, что множество $A$ является модулем, а потом мы доказываем, что это множество является алгеброй. Поэтому для записи координат элемента модуля мы также будем пользоваться соглашением 1.5.6.* □

Соглашение 1.5.8. *Отождествление вектора и матрицы его координат может привести к неоднозначности в равенстве*

$$a = a^*{}_* e \tag{1.5.1}$$

*где $\overline{\overline{e}}$ - базис векторного пространства. Поэтому мы будем записывать равенство (1.5.1) в виде*

$$\overline{a} = a^*{}_* e$$

*чтобы видеть, где записан вектор.* □

Соглашение 1.5.9. *Если свободная конечномерная алгебра имеет единицу, то мы будем отождествлять вектор базиса $\overline{e}_{0}$ с единицей алгебры.* □

Без сомнения, у читателя могут быть вопросы, замечания, возражения. Я буду признателен любому отзыву.

**Часть 1**

# Линейная алгебра над телом

Глава 2

# Линейное отображение ${}_*{}^*D$-векторного пространства

## 2.1. Линейное отображение ${}_*{}^*D$-векторного пространства

В этом разделе мы положим, что $V$, $W$ - это ${}_*{}^*D$-векторные пространства.

Определение 2.1.1. Обозначим $\mathcal{L}({}_*{}^*D;V;W)$ множество линейных отображений

$$A:V\to W$$

${}_*{}^*D$-векторного пространства $V$ в ${}_*{}^*D$-векторное пространство $W$. Обозначим $\mathcal{L}(D_*{}^*;V;W)$ множество линейных отображений

$$A:V\to W$$

${}_*{}^*D$-векторного пространства $V$ в ${}_*{}^*D$-векторное пространство $W$. □

Мы можем рассматривать тело $D$ как одномерное ${}_*{}^*D$-векторное пространство. Соответственно мы можем рассматривать множества $\mathcal{L}({}_*{}^*D;D;W)$ и $\mathcal{L}({}_*{}^*D;V;D)$.

Определение 2.1.2. Обозначим $\mathcal{L}(*T;S;R)$ множество $\star T$-представлений тела $S$ в аддитивной группе тела $R$. Обозначим $\mathcal{L}(T*;S;R)$ множество $T\star$-представлений тела $S$ в аддитивной группе тела $R$. □

Теорема 2.1.3. *Предположим, что $V$, $W$ - ${}_*{}^*D$-векторные пространства. Тогда множество $\mathcal{L}({}_*{}^*D;V;W)$ является абелевой группой относительно закона композиции*

$$(A+B)_*{}^*x=A_*{}^*x+B_*{}^*x \tag{2.1.1}$$

*Линейное отображение $A+B$ называется* **суммой линейных отображений** *$A$ и $B$.*

Доказательство. Нам надо показать, что отображение

$$A+B:V\to W$$

определённое равенством (2.1.1), - это линейное отображение ${}_*{}^*D$-векторных пространств. Согласно определению [3]-4.4.2

$$A_*{}^*(x_*{}^*a)=(A_*{}^*x)_*{}^*a$$
$$B_*{}^*(x_*{}^*a)=(B_*{}^*x)_*{}^*a$$

Мы видим, что

$$\begin{aligned}(A+B)_*{}^*(x_*{}^*a) &= A_*{}^*(x_*{}^*a)+B_*{}^*(x_*{}^*a)\\ &=(A_*{}^*x)_*{}^*a+(B_*{}^*x)_*{}^*a\\ &=(A_*{}^*x+B_*{}^*x)_*{}^*a\\ &=((A+B)_*{}^*x)_*{}^*a\end{aligned}$$

Нам надо показать так же, что эта операция коммутативна.

$$\begin{aligned}(A+B)_*{}^*x &= A_*{}^*x+B_*{}^*x\\ &=B_*{}^*x+A_*{}^*x\\ &=(B+A)_*{}^*x\end{aligned}$$

□

Теорема 2.1.4. *Пусть* $\overline{\overline{e}}_V=(e_{V\cdot\boldsymbol{i}},\boldsymbol{i}\in\boldsymbol{I})$ *- базис в* ${}_*{}^*D$*-векторном пространстве* $V$ *и* $\overline{\overline{e}}_W=(e_{W\cdot\boldsymbol{j}},\boldsymbol{j}\in\boldsymbol{J})$ *- базис в* ${}_*{}^*D$*-векторном пространстве* $W$*. Пусть* $A=({}^{\boldsymbol{j}}A_{\boldsymbol{i}})$, $\boldsymbol{i}\in\boldsymbol{I}$, $\boldsymbol{j}\in\boldsymbol{J}$, *- произвольная матрица. Тогда отображение*

$$A:V\to W$$

*определённое равенством*

$$b=A_*{}^*a$$

*относительно выбранных базисов, является линейным отображением* ${}_*{}^*D$*-векторных пространств.*

Доказательство. Теорема 2.1.4 является обратным утверждением теореме [3]-4.4.3. Предположим $A_*{}^*v=e_*{}^*A_*{}^*v$. Тогда

$$\begin{aligned}A_*{}^*(v_*{}^*a) &= e_*{}^*A_*{}^*v_*{}^*a\\ &=(A_*{}^*v)_*{}^*a\end{aligned}$$

□

Теорема 2.1.5. *Пусть* $\overline{\overline{e}}_V=(e_{V\cdot\boldsymbol{i}},\boldsymbol{i}\in\boldsymbol{I})$ *- базис в* ${}_*{}^*D$*-векторном пространстве* $V$ *и* $\overline{\overline{e}}_W=(e_{W\cdot\boldsymbol{j}},\boldsymbol{j}\in\boldsymbol{J})$ *- базис в* ${}_*{}^*D$*-векторном пространстве* $W$*. Предположим линейное отображение* $A$ *имеет матрицу* $A=({}^{\boldsymbol{j}}A_{\boldsymbol{i}})$, $\boldsymbol{i}\in\boldsymbol{I}$, $\boldsymbol{j}\in\boldsymbol{J}$, *относительно выбранных базисов. Пусть* $m\in D$*. Тогда матрица*

$${}^{\boldsymbol{j}}(mA)_{\boldsymbol{i}}=m\,{}^{\boldsymbol{j}}A_{\boldsymbol{i}}$$

*определяет линейное отображение*

$$mA:V\to W$$

*которое мы будем называть* **произведением линейного отображения** $A$ **на скаляр**.

Доказательство. Утверждение теоремы является следствием теоремы 2.1.4. □

Теорема 2.1.6. *Множество* $\mathcal{L}({}_*{}^*D;V;W)$ *является* $D_*{}^*$*-векторным пространством.*

Доказательство. Теорема 2.1.3 утверждает, что $\mathcal{L}({}_*{}^*D;V;W)$ - абелевая группа. Из теоремы 2.1.5 следует, что элемент тела $D$ порождает $T\star$-преобразование на абелевой группе $\mathcal{L}({}_*{}^*D;V;W)$. Из теорем 2.1.4, [3]-4.1.1 и [3]-4.1.3 следует, что множество $\mathcal{L}({}_*{}^*D;V;W)$ является $D\star$-векторным пространством.

Выписывая элементы базиса $D\star$-векторного пространства $\mathcal{L}({}_*{}^*D;V;W)$ в виде ${}_*$-строки или ${}^*$-строки, мы представим $D\star$-векторное пространство $\mathcal{L}({}_*{}^*D;V;W)$ как $D^*{}_*$- или $D_*{}^*$-векторное пространство. При этом надо иметь в виду, что выбор между $D_*{}^*$- и $D^*{}_*$-линейной зависимостью в $D\star$-векторном пространстве $\mathcal{L}({}_*{}^*D;V;W)$ не зависит от типа векторных пространств $V$ и $W$.

Для того, чтобы выбрать тип векторного пространства $\mathcal{L}({}_*{}^*D;V;W)$ мы обратим внимание на следующее обстоятельство. Допустим $V$ и $W$ - ${}_*{}^*D$-векторные пространства. Допустим $\mathcal{L}({}_*{}^*D;V;W)$ - $D_*{}^*$-векторное пространство. Тогда действие ${}^*$-строки ${}_*{}^*D$-линейных отображений ${}^jA$ на ${}_*$-строку векторов $f_i$ можно представить в виде матрицы

$$\begin{pmatrix}{}^1A\\ ...\\ {}^mA\end{pmatrix}{}_*{}^*\begin{pmatrix}f_1 & ... & f_n\end{pmatrix}=\begin{pmatrix}{}^1A_*{}^*f_1 & ... & {}^1A_*{}^*f_n\\ ... & ... & ...\\ {}^mA_*{}^*f_1 & ... & {}^mA_*{}^*f_n\end{pmatrix}$$

Эта запись согласуется с матричной записью действия $D_*{}^*$-линейной комбинации $a_*{}^*A$ ${}_*{}^*D$-линейных отображений $A$. □

Мы можем также определить $\star D$-произведение ${}_*{}^*D$-линейного отображения $A$ на скаляр. Однако, вообще говоря, это $\star T$-представление тела $D$ в $D\star$-векторном пространстве $\mathcal{L}({}_*{}^*D;V;W)$ не может быть перенесено в ${}_*{}^*D$-векторное пространство $W$. Действительно, в случае $\star D$-произведения мы имеем

$$(Am)_*{}^*v=(Am)_*{}^*f_*{}^*v=e_*{}^*(Am)_*{}^*v$$

Поскольку произведение в теле не коммутативно, мы не можем выразить полученное выражение как произведение $A_*{}^*v$ на скаляр $m$.

Следствием теоремы 2.1.6 является неоднозначность записи

$$m_*{}^*A_*{}^*v$$

Мы можем предположить, что смысл этой записи ясен из контекста. Однако желательно неоднозначность избежать. Мы будем пользоваться для этой цели скобками. Выражение

$$w=[m_*{}^*A]_*{}^*v$$

означает, что $D_*{}^*$-линейная комбинация ${}_*{}^*D$-линейных отображений ${}^iA$ отображает вектор $v$ в вектор $w$. Выражение

$$w=B_*{}^*A_*{}^*v$$

означает, что ${}_*{}^*$-произведение ${}_*{}^*D$-линейных отображений $A$ и $B$ отображает вектор $v$ в вектор $w$.

## 2.2. 1-форма на ${}_*{}^*D$-векторном пространстве

Определение 2.2.1. **1-форма** на ${}_*{}^*D$-векторном пространстве $V$ - это линейное отображение

$$b:V\rightarrow D\tag{2.2.1}$$

□

Мы можем значение 1-формы $b$, определённое для вектора $a$, записывать в виде

$$b(a) =< b, a >$$

Теорема 2.2.2. *Множество $\mathcal{L}({}_*{}^*D; V; D)$ является $D\star$-векторным пространством.*

Доказательство. ${}_*{}^*D$-векторное пространство размерности 1 эквивалентно телу $D$ □

Теорема 2.2.3. *Пусть $\overline{\overline{e}}$ - базис в ${}_*{}^*D$-векторном пространстве $V$. 1-форма $b$ имеет представление*

$$< b, a >= b_*{}^*a \tag{2.2.2}$$

*относительно выбранного базиса, где вектор $a$ имеет разложение*

$$a = e_*{}^*a \tag{2.2.3}$$

*и*

$$b_i =< b, e_i > \tag{2.2.4}$$

Доказательство. Так как $b$ - 1-форма, из (2.2.3) следует, что

$$< b, a >=< b, e_*{}^*a >=< b, e_i > \ {}^ia \tag{2.2.5}$$

(2.2.2) следует из (2.2.5) и (2.2.4). □

Теорема 2.2.4. *Пусть $\overline{\overline{e}}$ - базис в ${}_*{}^*D$-векторном пространстве $V$. 1-форма (2.2.1) однозначно определена значениями (2.2.4), в которые 1-форма $b$ отображает вектора базиса.*

Доказательство. Утверждение является следствием теорем 2.2.3 и [3]-4.3.3. □

Теорема 2.2.5. *Пусть $\overline{\overline{e}}$ - базис в ${}_*{}^*D$-векторном пространстве $V$. Множество 1-форм ${}^id$ таких, что*

$$< {}^jd, e_i >= {}^j\delta_i \tag{2.2.6}$$

*является базисом $\overline{\overline{d}}$ $D_*{}^*$-векторного пространства $\mathcal{L}({}_*{}^*D; V; D)$.*

Доказательство. 1-форма ${}^jd$ существует для данного $j$ согласно теореме 2.2.4, если положить $b_i = {}^j\delta_i$ . Если предположить, что существует 1-форма

$$b = b_*{}^*d = 0$$

то

$$b_j < {}^jd, e_i >= 0$$

Согласно равенству (2.2.6)

$$b_i = b_j\ {}^j\delta_i = 0$$

Следовательно, 1-формы ${}^jd$ линейно независимы. □

Определение 2.2.6. Пусть $V$ - ${}_*{}^*D$-векторное пространство. $D_*{}^*$-векторное пространство

$$V^* = \mathcal{L}({}_*{}^*D; V; D)$$

называется **дуальным пространством к ${}_*{}^*D$-векторному пространству** $V$. Пусть $\overline{\overline{e}}$ - базис в ${}_*{}^*D$-векторном пространстве $V$. базис $\overline{\overline{d}}$ $D_*{}^*$-векторного пространства $V^*$, удовлетворяющий равенству (2.2.6), называется **базисом, дуальным базису** $\overline{\overline{e}}$. □

Теорема 2.2.7. *Пусть $A$ - пассивное преобразование многообразия базисов $\mathcal{B}(V, GL_n^{{}_*{}^*D})$. Допустим базис*

$$\overline{\overline{e}}' = \overline{\overline{e}}_*{}^*A \tag{2.2.7}$$

*является образом базиса $\overline{\overline{e}}$. Пусть $B$ - пассивное преобразование многообразия базисов $\mathcal{B}(V^*, GL_n^{{}_*{}^*D})$ такое, что базис*

$$\overline{\overline{d}}' = B_*{}^*\overline{\overline{d}} \tag{2.2.8}$$

*дуален базису. Тогда*

$$B = A^{-1_*{}^*} \tag{2.2.9}$$

Доказательство. Из равенств (2.2.6), (2.2.7), (2.2.8) следует

$$\begin{aligned} {}^j\delta_i &= <{}^jd', e'_i> \\ &={}^jB_l\ <{}^ld, e_k>\ {}^kA_i \\ &={}^jB_l\ {}^l\delta_k\ {}^kA_i \\ &={}^jB_k\ {}^kA_i \end{aligned} \tag{2.2.10}$$

Равенство (2.2.9) следует из равенства (2.2.10). □

## 2.3. Парные представления тела

Теорема 2.3.1. *В любом ${}_*{}^*D$-векторном пространстве можно определить структуру $D\star$-векторного пространства, определив $D\star$-произведение вектора на скаляр равенством*

$$mv = m\delta_*{}^*v$$

Доказательство. Непосредственная проверка доказывает, что отображение

$$f : D \to V^\star$$

определённое равенством

$$f(m) = m\delta$$

определяет $T\star$-представление кольца $D$. □

Мы можем иначе сформулировать теорему 2.3.1.

Теорема 2.3.2. *Если мы определили эффективное $\star T$-представление $f$ тела $D$ на абелевой группе $V$, то мы можем однозначно определить эффективное $T\star$-представление $h$ тела $D$ на абелевой группе $V$ такое, что диаграмма*

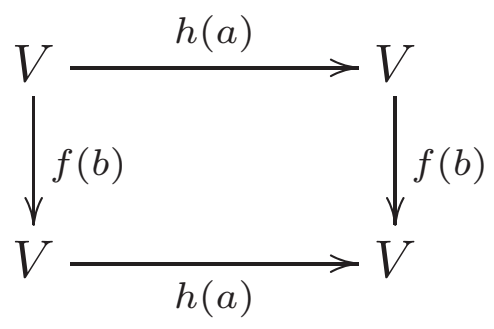

*коммутативна для любых $a$, $b \in D$.* □

Мы будем называть представления $f$ и $h$ **парными представлениями тела** $D$.

Теорема 2.3.3. *В векторном пространстве $V$ над телом $D$ мы можем определить $D\star$-произведение и $\star D$-произведение вектора на скаляр. Согласно теореме 2.3.2 эти операции удовлетворяют равенству*

$$(am)b = a(mb) \tag{2.3.1}$$

Равенство (2.3.1) представляет **закон ассоциативности для парных представлений**. Это позволяет нам писать подобные выражения не пользуясь скобками.

Доказательство. В разделе 2.1 дано определение $D\star$-произведения ${}_*{}^*D$-линейного отображения $A$ на скаляр. Согласно теореме 2.3.1 $T\star$-представление тела $D$ в $D\star$-векторном пространстве $\mathcal{L}({}_*{}^*D; V; W)$ может быть перенесено в ${}_*{}^*D$-векторном пространство $W$ согласно правилу

$$[mA]_*{}^*v = [m\delta]_*{}^*(A_*{}^*v) = m(A_*{}^*v)$$

□

Аналогия с векторными пространствами над полем заходит столь далеко, что мы можем предположить существование концепции базиса, который годится для $D\star$-произведения и $\star D$-произведения вектора на скаляр.

Теорема 2.3.4. *Многообразие базисов $D_*{}^*$-векторного пространства $V$ и многообразие базисов ${}^*{}_*D$-векторного пространства $V$ отличны*

$$\mathcal{B}(V, D_*{}^*) \neq \mathcal{B}(V, {}^*{}_*D)$$

Доказательство. При доказательстве этой теоремы мы будем пользоваться стандартным представлением матрицы. Не нарушая общности, мы проведём доказательство в координатном векторном пространстве $D^n$.

Пусть $\overline{\overline{e}} = (e_i = (\delta_i^j), i, j \in I, |I| = n)$ множество векторов векторного пространства $D^n$. Очевидно, $\overline{\overline{e}}$ является одновременно базисом $D_*{}^*$-векторного пространства и базисом ${}^*{}_*D$-векторного пространства. Для произвольного множества векторов $(f_i, i \in I, |I| = n)$ $D_*{}^*$-координатная матрица

$$f = \begin{pmatrix} f_1^1 & ... & f_1^n \\ ... & ... & ... \\ f_n^1 & ... & f_n^n \end{pmatrix} \tag{2.3.2}$$

относительно базиса $\overline{\overline{e}}$ совпадает с ${}^*{}_*D$-координатной матрицей относительно базиса $\overline{\overline{e}}$.

Если множество векторов $(f_{\boldsymbol{i}}, \boldsymbol{i} \in \boldsymbol{I}, |\boldsymbol{I}| = \boldsymbol{n})$ - базис ${}_*{}^*D$-векторного пространства, то согласно [3]-4.9.3 матрица (2.3.2) - ${}_*{}^*$-невырожденная матрица.

Если множество векторов $(f_{\boldsymbol{i}}, \boldsymbol{i} \in \boldsymbol{I}, |\boldsymbol{I}| = \boldsymbol{n})$ - базис ${}^*{}_*D$-векторного пространства, то согласно теореме [3]-4.9.3 матрица (2.3.2) - ${}^*{}_*$-невырожденная матрица.

Следовательно, если множество векторов $(f_{\boldsymbol{i}}, \boldsymbol{i} \in \boldsymbol{I}, |\boldsymbol{I}| = \boldsymbol{n})$ порождают базис $D_*{}^*$-векторного пространства и базис ${}^*{}_*D$-векторного пространства, их координатная матрица (2.3.2) является ${}_*{}^*$-невырожденной и ${}^*{}_*$-невырожденной матрицей. Утверждение следует из теоремы [3]-4.8.9. □

Из теоремы 2.3.4 следует, что существует базис $\overline{\overline{e}}$ $D_*{}^*$-векторного пространства $V$, который не является базисом ${}^*{}_*D$-векторного пространства $V$.

## 2.4. $D$-векторное пространство

При изучении многих задач мы вполне можем ограничиться рассмотрением $D\star$-векторного пространства либо $\star D$-векторного пространства. Однако есть задачи, в которых мы вынуждены отказаться от простой модели и одновременно рассматривать обе структуры векторного пространства.

Парное представление тела $D$ в абелевой группе $V$ определяет структуру $D$ **-векторного пространства**. Из этого определения следует, что если $v$, $w \in V$, то для любых $a$, $b$, $c$, $d \in D$

$$avb + cwd \in V \tag{2.4.1}$$

Выражение (2.4.1) называется **линейной комбинацией векторов** $D$-векторного пространства $V$.

Определение 2.4.1. Векторы $a_i$, $i \in I$, $D$-векторного пространства $V$ **линейно независимы**, если из равенства

$$b^i a_i c^i = 0$$

следует $b^i c^i = 0$, $i = 1, ..., n$, В противном случае, векторы $a_i$ **линейно зависимы**. □

Определение 2.4.2. Мы называем множество векторов

$$\overline{\overline{e}} = (e_{\boldsymbol{i}}, \boldsymbol{i} \in \boldsymbol{I})$$

**базисом $D$-векторного пространства**, если векторы $e_{\boldsymbol{i}}$ линейно независимы и, добавив к этой системе любого другого вектора, мы получим новую систему, которая является линейно зависимой. □

Теорема 2.4.3. *Базис $D$-векторного пространства $V$ является базисом $D^*{}_*$-векторного пространства $V$.*

Доказательство. Пусть векторы $a_i$, $i \in I$ линейно зависимы в $D^*{}_*$-векторном пространстве. Тогда существует ${}_*$-строка

$$\left(c^i \quad i \in I\right) \neq 0 \tag{2.4.2}$$

такая, что

$$c^*{}_*a = 0 \tag{2.4.3}$$

Из равенства (2.4.3), следует, что

$$c^i a_i 1 = 0 \tag{2.4.4}$$

Согласно определению 2.4.1, векторы $a_i$ линейно зависимы в $D$-векторном пространстве.

Пусть $D^*{}_*$-размерность векторного пространства $V$ равна $n$. Пусть $\overline{\overline{e}}$ - базис $D^*{}_*$-векторного пространства. Пусть

$$\begin{pmatrix}1 & ... & 0\\ ... & ... & ...\\ 0 & ... & 1\end{pmatrix}$$

координатная матрица базиса $\overline{\overline{e}}$. Очевидно, что мы не можем представить вектор $e_n$ как $D$-линейную комбинацию векторов $e_1$, ..., $e_{n-1}$. □

Для нас несущественно, является ли множество векторов $(e_i, i \in I)$ базисом ${}_*{}^*D$-векторного пространства или базисом $D^*{}_*$-векторного пространства. Поэтому при записи координат вектора мы индексы будем записывать в стандартном формате. Рассматривая базис $D$-векторного пространства $V$ мы, опираясь на теорему 2.4.3, будем рассматривать его как базис $D^*{}_*$-векторного пространства $V$.

Глава 3

# Линейное отображение тела

## 3.1. Линейное отображение тела

Согласно замечанию [3]-8.2.2, мы рассматриваем тело $D$ как алгебру над полем $F \subset Z(D)$.

Определение 3.1.1. Пусть $D_1$, $D_2$ - тела. Пусть $F$ - поле такое, что $F \subset Z(D_1)$, $F \subset Z(D_2)$. Линейное отображение

$$f : D_1 \to D_2$$

$F$-векторного пространства $D_1$ в $F$-векторное пространство $D_2$ называется **линейным отображением тела $D_1$ в тело $D_2$**. □

Согласно определению 3.1.1, линейное отображение $f$ тела $D_1$ в тело $D_2$ удовлетворяет свойству

$$f(a+b) = f(a) + f(b) \quad a, b \in D$$
$$f(pa) = pf(a) \quad p \in F$$

Теорема 3.1.2. *Пусть отображение*

$$f : D_1 \to D_2$$

*является линейным отображением тела $D_1$ характеристики 0 в тело $D_2$ характеристики 0. Тогда*[3.1]

$$f(nx) = nf(x)$$

*для любого целого $n$.*

Доказательство. Мы докажем теорему индукцией по $n$. При $n = 1$ утверждение очевидно, так как

$$f(1x) = f(x) = 1f(x)$$

Допустим уравнение справедливо при $n = k$. Тогда

$$f((k+1)x) = f(kx + x) = f(kx) + f(x) = kf(x) + f(x) = (k+1)f(x)$$

□

[3.1] Допустим

$$f : D_1 \to D_2$$

является линейным отображением тела $D_1$ характеристики 2 в тело $D_2$ характеристики 3. Тогда для любого $a \in D_1$, $2a = 0$, хотя $2f(a) \neq 0$. Следовательно, если мы предположим, что характеристика тела $D_1$ больше 0, то мы должны потребовать, чтобы характеристика тела $D_1$ равнялась характеристике тела $D_2$.

Теорема 3.1.3. *Пусть отображение*

$$f : D_1 \to D_2$$

*является линейным отображением тела* $D_1$ *характеристики* 0 *в тело* $D_2$ *характеристики* 0*. Тогда*

$$f(ax) = af(x)$$

*для любого рационального* $a$.

Доказательство. Запишем $a$ в виде $a = \frac{p}{q}$. Положим $y = \frac{1}{q}x$. Тогда согласно теореме 3.1.2

$$f(x) = f(qy) = qf(y) = qf\left(\frac{1}{q}x\right) \tag{3.1.1}$$

Из равенства (3.1.1) следует

$$\frac{1}{q}f(x) = f\left(\frac{1}{q}x\right) \tag{3.1.2}$$

Из равенства (3.1.2) следует

$$f\left(\frac{p}{q}x\right) = pf\left(\frac{1}{q}x\right) = \frac{p}{q}f(x)$$

□

Мы не можем распространить утверждение теоремы 3.1.3 на произвольное подполе центра $Z(D)$ тела $D$.

Теорема 3.1.4. *Пусть тело* $D$ *является алгеброй над полем* $F \subset Z(D)$*. Если* $F \neq Z(D)$*, то существует линейное отображение*

$$f : D \to D$$

*которое нелинейно над полем* $Z(D)$*.*

Доказательство. Для доказательства теоремы достаточно рассмотреть поле комплексных чисел $C$ так как $C = Z(C)$. Так как поле комплексных чисел является алгеброй над полем вещественных чисел, то функция

$$z \to \overline{z}$$

линейна. Однако равенство

$$\overline{az} = a\overline{z}$$

неверно. □

Из теоремы 3.1.4 возникает вопрос. Для чего мы рассматриваем линейные отображения над полем $F \neq Z(D)$, если это приводит к резкому расширению множества линейных отображений? Ответом на этот вопрос служит богатый опыт теории функции комплексного переменного.

Теорема 3.1.5. *Пусть отображения*

$$f : D_1 \to D_2$$

$$g : D_1 \to D_2$$

*являются линейными отображениями тела* $D_1$ *в тело* $D_2$*. Тогда отображение* $f + g$ *также является линейным.*

Доказательство. Утверждение теоремы следует из цепочки равенств

$$\begin{aligned}(f+g)(x+y) =&f(x+y)+g(x+y)=f(x)+f(y)+g(x)+g(y)\\ =&(f+g)(x)+(f+g)(y)\\ (f+g)(px) =&f(px)+g(px)=pf(x)+pg(x)=p(f(x)+g(x))\\ =&p(f+g)(x)\end{aligned}$$

□

Теорема 3.1.6. *Пусть отображение*

$$f:D_1\to D_2$$

*является линейным отображением тела $D_1$ в тело $D_2$. Тогда отображения $af$, $fb$, $a,\ b\in R_2$, также являются линейными.*

Доказательство. Утверждение теоремы следует из цепочки равенств

$$\begin{aligned}(af)(x+y) =&a(f(x+y))=a(f(x)+f(y))=af(x)+af(y)\\ =&(af)(x)+(af)(y)\\ (af)(px) =&a(f(px))=a(pf(x))=p(af(x))\\ =&p(af)(x)\\ (fb)(x+y) =&(f(x+y))b=(f(x)+f(y))b=f(x)b+f(y)b\\ =&(fb)(x)+(fb)(y)\\ (fb)(px) =&(f(px))b=(pf(x))b=p(f(x)b)\\ =&p(fb)(x)\end{aligned}$$

□

Определение 3.1.7. Обозначим $\mathcal{L}(D_1;D_2)$ множество линейных отображений

$$f:D_1\to D_2$$

тела $D_1$ в тело $D_2$. □

Теорема 3.1.8. *Мы можем представить линейное отображение*

$$f:D_1\to D_2$$

*тела $D_1$ в тело $D_2$ в виде*

$$f(x)=f_{k\cdot s_k\cdot 0}\ G_k(x)\ f_{k\cdot s_k\cdot 1} \tag{3.1.3}$$

*где $(G_k,k\in K)$ - множество аддитивных отображений тела $D_1$ в тело $D_2$.*[3.2] *Выражение $f_{k\cdot s_k\cdot p}$, $p=0,\ 1,$ в равенстве (3.1.3) называется* **компонентой линейного отображения** *$f$.*

Доказательство. Утверждение теоремы следует из теорем 3.1.5 и 3.1.6. □

[3.2] Здесь и в дальнейшем мы будем предполагать сумму по индексу, который встречается в произведении несколько раз. Равенство (3.1.3) является рекурсивным определением и есть надежда, что мы можем его упростить.

Если в теореме 3.1.8 $|K| = 1$, то равенство (3.1.3) имеет вид

$$f(x) = f_{s\cdot 0}\ G(x)\ f_{s\cdot 1} \tag{3.1.4}$$

и отображение $f$ называется **линейным отображением, порождённым отображением** $G$. Отображение $G$ мы будем называть **образующей линейного отображения**.

Теорема 3.1.9. *Пусть $D_1$, $D_2$ - тела характеристики* 0. *Пусть $F$, $F \subset Z(D_1)$, $F \subset Z(D_2)$, - поле. Пусть $G$ - линейное отображение. Пусть $\overline{\overline{e}}$ - базис тела $D_2$ над полем $F$.* **Стандартное представление линейного отображения** (3.1.4) *имеет вид*[3.3]

$$f(x) = f_G^{ij}\ e_i\ G(x)\ e_j \tag{3.1.5}$$

*Выражение $f_G^{ij}$ в равенстве* (3.1.5) *называется* **стандартной компонентой линейного отображения** $f$.

Доказательство. Компоненты линейного отображения $f$ имеют разложение

$$f_{s\cdot p} = f_{s\cdot p}^{i}\ e_i \tag{3.1.6}$$

относительно базиса $\overline{\overline{e}}$. Если мы подставим (3.1.6) в (3.1.4), мы получим

$$f(x) = f_{s\cdot 0}^{i}\ e_i\ G(x)\ f_{s\cdot 1}^{j}\ e_j \tag{3.1.7}$$

Подставив в равенство (3.1.7) выражение

$$f_G^{ij} = f_{s\cdot 0}^{i}\ f_{s\cdot 1}^{j}$$

мы получим равенство (3.1.5). □

Теорема 3.1.10. *Пусть $D_1$, $D_2$ - тела характеристики* 0. *Пусть $F$, $F \subset Z(D_1)$, $F \subset Z(D_2)$, - поле. Пусть $G$ - линейное отображение. Пусть $\overline{\overline{e}}_1$ - базис тела $D_1$ над полем $F$. Пусть $\overline{\overline{e}}_2$ - базис тела $D_2$ над полем $F$. Пусть $C_{2\cdot}{}^{p}_{kl}$ - структурные константы тела $D_2$. Тогда линейное отображение* (3.1.4), *порождённое линейным отображением $G$, можно записать в виде*

$$\begin{aligned} f(a) &= e_{2\cdot j}\ f^j_i\ a^i & f^j_k &\in F \\ a &= e_{1\cdot i}\ a^i & a^i &\in F \quad a \in D_1 \end{aligned} \tag{3.1.8}$$

$$f^j_i = G^l_i f_G^{kr} C_{2\cdot}{}^{p}_{kl} C_{2\cdot}{}^{j}_{pr} \tag{3.1.9}$$

Доказательство. Выберем отображение

$$\begin{aligned} G : D_1 \to D_2 \quad a = e_{1\cdot i}a^i \to G(a) = e_{2\cdot j}G^j_i a^i \\ a^i \in F \quad G^j_i \in F \end{aligned} \tag{3.1.10}$$

Согласно теореме [3]-4.4.3 линейное отображение $f(a)$ относительно базисов $\overline{\overline{e}}_1$ и $\overline{\overline{e}}_2$ принимает вид (3.1.8). Из равенств (3.1.5) и (3.1.10) следует

$$f(a) = G^l_i a^i f_G^{kj} e_{2\cdot k} e_{2\cdot l} e_{2\cdot j} \tag{3.1.11}$$

[3.3] Представление линейного отображения тела с помощью компонент линейного отображения неоднозначно. Чисто алгебраическими методами мы можем увеличить либо уменьшить число слагаемых. Если размерность тела $D_2$ над полем $F$ конечна, то стандартное представление линейного отображения гарантирует конечность множества слагаемых в представлении отображения.

Из равенств (3.1.8) и (3.1.11) следует

$$e_{2\cdot j}\ f^{\ j}_{i}\ a^i = G^l_i a^i f^{kr}_G e_{2\cdot k} e_{2\cdot l} e_{2\cdot r} = G^l_i a^i f^{kr}_G C_{2\cdot}{}^{p}_{kl} C_{2\cdot}{}^{j}_{pr}\ e_{2\cdot j} \tag{3.1.12}$$

Так как векторы $e_{2\cdot r}$ линейно независимы над полем $F$ и величины $a^k$ произвольны, то из равенства (3.1.12) следует равенство (3.1.9). □

При рассмотрении отображений

$$f : D \to D \tag{3.1.13}$$

мы будем полагать $G(x) = x$.

Теорема 3.1.11. *Пусть $D$ является телом характеристики $0$. Линейное отображение (3.1.13) имеет вид*

$$f(x) = f_{s\cdot 0}\ x\ f_{s\cdot 1} \tag{3.1.14}$$

Теорема 3.1.12. *Пусть $D$ - тело характеристики $0$. Пусть $\overline{\overline{e}}$ - базис тела $D$ над полем $F \subset Z(D)$. Стандартное представление линейного отображения (3.1.14) тела имеет вид*

$$f(x) = f^{ij}\ e_i x e_j \tag{3.1.15}$$

Теорема 3.1.13. *Пусть $D$ является телом характеристики $0$. Пусть $\overline{\overline{e}}$ - базис тела $D$ над полем $F \subset Z(D)$. Тогда линейное отображение (3.1.13) можно записать в виде*

$$\begin{aligned} f(a) &= e_j f^{\ j}_{i} a^i & f^{\ j}_{k} &\in F \\ a &= e_i a^i & a^i &\in F \quad a \in D \end{aligned} \tag{3.1.16}$$

$$f^{\ j}_{i} = f^{kr} C^{p}_{ki} C^{j}_{pr} \tag{3.1.17}$$

Теорема 3.1.14. *Рассмотрим матрицу*

$$\mathcal{C} = \left(\mathcal{C}^{\cdot j}_{i\cdot kr}\right) = \left(C^{p}_{ki} C^{j}_{pr}\right) \tag{3.1.18}$$

*строки которой пронумерованы индексом $\cdot j \atop i$ и столбцы пронумерованы индексом $\cdot kr$ Если $\det \mathcal{C} \neq 0$, то для заданных координат линейного преобразования $f^{\ j}_{i}$ система линейных уравнений (3.1.17) относительно стандартных компонент этого преобразования $f^{kr}$ имеет единственное решение. Если $\det \mathcal{C} = 0$, то условием существования решения системы линейных уравнений (3.1.17) является равенство*

$$\operatorname{rank}\left(\mathcal{C}^{\cdot j}_{i\cdot kr}\quad f^{\ j}_{i}\right) = \operatorname{rank}\mathcal{C} \tag{3.1.19}$$

*В этом случае система линейных уравнений (3.1.17) имеет бесконечно много решений и существует линейная зависимость между величинами $f^{\ j}_{i}$.*

Доказательство. Равенство (3.1.14) является частным случаем равенства (3.1.4) при условии $G(x) = x$. Теорема 3.1.12 является частным случаем теоремы 3.1.9 при условии $G(x) = x$. Теорема 3.1.13 является частным случаем теоремы 3.1.10 при условии $G(x) = x$. Утверждение теоремы 3.1.14 является следствием теории линейных уравнений над полем. □

Теорема 3.1.15. *Стандартные компоненты тождественного отображения имеют вид*

$$f^{kr} = \delta^k_0 \delta^r_0 \tag{3.1.20}$$

Доказательство. Равенство (3.1.20) является следствием равенства

$$x = e_0\ x\ e_0$$

Убедимся, что стандартные компоненты (3.1.20) линейного преобразования удовлетворяют уравнению

$$\delta_i^j = f^{kr} C_{ki}^p C_{pr}^j \tag{3.1.21}$$

которое следует из уравнения (3.1.17) если $f = \delta$. Из равенств (3.1.20), (3.1.21) следует

$$\delta_i^j = C_{0i}^p C_{p0}^j \tag{3.1.22}$$

Равенство (3.1.22) верно, так как из равенств

$$e_j e_0 = e_0\ e_j = e_j$$

следует

$$C_{0r}^j = \delta_r^j \quad C_{r0}^j = \delta_r^j$$

Если $\det \mathcal{C} \neq 0$, то решение (3.1.20) единственно. Если $\det \mathcal{C} = 0$, то система линейных уравнений (3.1.21) имеет бесконечно много решений. Однако нас интересует по крайней мере одно. □

Теорема 3.1.16. *Если $\det \mathcal{C} \neq 0$, то стандартные компоненты нулевого отображения*

$$z : A \to A \quad z(x) = 0$$

*определены однозначно и имеют вид $z^{ij} = 0$. Если $\det \mathcal{C} = 0$, то множество стандартных компонент нулевого отображения порождает векторное пространство.*

Доказательство. Теорема верна, поскольку стандартные компоненты $z^{ij}$ являются решением однородной системы линейных уравнений

$$0 = z^{kr} C_{ki}^p C_{pr}^j$$

□

Замечание 3.1.17. Рассмотрим равенство

$$a^{kr}\ e_k\ x\ e_r = b^{kr}\ e_k\ x\ e_r \tag{3.1.23}$$

Из теоремы 3.1.16 следует, что только при условии $\det \mathcal{C} \neq 0$ из равенства (3.1.23) следует

$$a^{kr} = b^{kr} \tag{3.1.24}$$

В противном случае мы должны предполагать равенство

$$a^{kr} = b^{kr} + z^{kr} \tag{3.1.25}$$

Несмотря на это, мы и в случае $\det \mathcal{C} = 0$ будем пользоваться стандартным представлением так как в общем случае указать множество линейно независимых векторов - задача достаточно сложная. Если мы хотим определить операцию над линейными отображениями, записанными в стандартном представлении, то также как в случае теоремы 3.1.15 мы будем выбирать представитель из множества возможных представлений. □

Теорема 3.1.18. *Выражение*

$$f_k^r = f^{ij} C_{ik}^p C_{pj}^r$$

*является тензором над полем $F$*

$$f'^{j}_{i} = A_i^k f_k^l A^{-1\cdot j}_{l} \tag{3.1.26}$$

Доказательство. $D$-линейное отображение относительно базиса $\overline{\overline{e}}$ имеет вид (3.1.16). Пусть $\overline{\overline{e}}'$ - другой базис. Пусть

$$e'_i = e_j A_i^j \tag{3.1.27}$$

преобразование, отображающее базис $\overline{\overline{e}}$ в базис $\overline{\overline{e}}'$. Так как линейное отображение $f$ не меняется, то

$$f(x) = e'_l f'^{l}_{k} x'^k \tag{3.1.28}$$

Подставим [3]-(8.2.8), (3.1.27) в равенство (3.1.28)

$$f(x) = e_j A_l^j f'^{l}_{k} A^{-1\cdot k}_{i} x^i \tag{3.1.29}$$

Так как векторы $e_j$ линейно независимы и компоненты вектора $x^i$ произвольны, то равенство (3.1.26) следует из равенства (3.1.29). Следовательно, выражение $f_k^r$ является тензором над полем $F$. □

Определение 3.1.19. Множество

$$\ker f = \{x \in D_1 : f(x) = 0\}$$

называется **ядром линейного отображения**

$$f : D_1 \to D_2$$

тела $D_1$ в тело $D_2$. □

Теорема 3.1.20. *Ядро линейного отображения*

$$f : D_1 \to D_2$$

*является подгруппой аддитивной группы тела $D_1$.*

Доказательство. Пусть $a$, $b \in \ker f$. Тогда

$$f(a) = 0$$
$$f(b) = 0$$
$$f(a + b) = f(a) + f(b) = 0$$

Следовательно, $a + b \in \ker f$. □

Определение 3.1.21. Линейное отображение

$$f : D_1 \to D_2$$

тела $D_1$ в тело $D_2$ называется **вырожденным**, если

$$\ker f \neq \{0\}$$

□

Теорема 3.1.22. *Пусть $D$ - тело характеристики 0. Пусть $\overline{\overline{e}}$ - базис тела $D$ над центром $Z(D)$ тела $D$. Пусть*

$$f:D\to D \quad f(x)=f_{s\cdot 0}\ x\ f_{s\cdot 1} \tag{3.1.30}$$
$$= f^{ij}\ e_i x e_j \tag{3.1.31}$$
$$g:D\to D \quad g(x)=g_{t\cdot 0}\ x\ g_{t\cdot 1} \tag{3.1.32}$$
$$= g^{ij}\ e_i x e_j \tag{3.1.33}$$

*линейные отображения тела $D$. Отображение*

$$h(x)=gf(x)=g(f(x)) \tag{3.1.34}$$

*является линейным отображением*

$$h(x)=h_{ts\cdot 0}\ x\ h_{ts\cdot 1} \tag{3.1.35}$$
$$= h^{pr}\ e_p\ x\ e_r \tag{3.1.36}$$

*где*

$$h_{ts\cdot 0}=g_{t\cdot 0}\ f_{s\cdot 0} \tag{3.1.37}$$
$$h_{ts\cdot 1}=f_{s\cdot 1}\ g_{t\cdot 1} \tag{3.1.38}$$
$$h^{pr}=g^{ij}f^{kl}C^p_{ik}C^r_{lj} \tag{3.1.39}$$

Доказательство. Отображение (3.1.34) линейно так как

$$h(x+y)=g(f(x+y))=g(f(x)+f(y))=g(f(x))+g(f(y))=h(x)+h(y)$$
$$h(ax)=g(f(ax))=g(af(x))=ag(f(x))=ah(x)$$

Если мы подставим (3.1.30) и (3.1.32) в (3.1.34), то мы получим

$$h(x)=g_{t\cdot 0}\ f(x)\ g_{t\cdot 1}=g_{t\cdot 0}\ f_{s\cdot 0}\ x\ f_{s\cdot 1}\ g_{t\cdot 1} \tag{3.1.40}$$

Сравнивая (3.1.40) и (3.1.35), мы получим (3.1.37), (3.1.38).

Если мы подставим (3.1.31) и (3.1.33) в (3.1.34), то мы получим

$$\begin{aligned}h(x)=&g^{ij}\ e_i\ f(x)\ e_j\\ =&g^{ij}\ e_i\ f^{kl}\ e_k\ x\ e_l e_j\\ =&g^{ij}f^{kl}C^p_{ik}C^r_{lj}\ e_p\ x\ e_r\end{aligned} \tag{3.1.41}$$

Сравнивая (3.1.41) и (3.1.36), мы получим (3.1.39). □

Определение 3.1.23. Пусть коммутативное кольцо $P$ является подкольцом центра $Z(R)$ кольца $R$. Отображение

$$f:R\to R$$

кольца $R$ называется **проективным над коммутативным кольцом** $P$, если

$$f(px)=f(x)$$

для любого $p\in P$. Множество

$$Px=\{px:p\in P,x\in R\}$$

называется **направлением** $x$ **над коммутативным кольцом**[3.4] $P$. □

[3.4]Направление над коммутативным кольцом $P$ является подмножеством кольца $R$. Однако мы будем обозначать направление $Px$ элементом $x\in R$, когда это не приводит к неоднозначности. Мы будем говорить о направлении над коммутативным кольцом $Z(R)$, если мы явно не указываем коммутативное кольцо $P$.

Пример 3.1.24. Если отображение $f$ кольца $R$ линейно над коммутативным кольцом $P$, то отображение

$$g(x)=x^{-1}f(x)$$

проективно над коммутативным кольцом $P$. □

## 3.2. Полилинейное отображение тела

Определение 3.2.1. Пусть $R_1$, ..., $R_n$, $P$ - кольца характеристики 0. Пусть $S$ - модуль над кольцом $P$. Пусть $F$ - коммутативное кольцо, которое для любого $i$ является подкольцом центра кольца $R_i$. Отображение

$$f:R_1\times ...\times R_n\to S$$

называется **полилинейным над коммутативным кольцом** $F$, если

$$f(p_1,...,p_i+q_i,...,p_n)=f(p_1,...,p_i,...,p_n)+f(p_1,...,q_i,...,p_n)$$
$$f(a_1,...,ba_i,...,a_n)=bf(a_1,...,a_i,...,a_n)$$

для любого $i$, $1\le i\le n$, и любых $p_i$, $q_i\in R_i$, $b\in F$. Обозначим $\mathcal{L}(R_1,...,R_n;S)$ множество полилинейных отображений колец $R_1$, ..., $R_n$ в модуль $S$. □

Теорема 3.2.2. *Пусть $D$ - тело характеристики 0. Полилинейное отображение*

$$f:D^n\to D,d=f(d_1,...,d_n) \tag{3.2.1}$$

*имеет вид*

$$d=f^n_{s\cdot 0}\ \sigma_s(d_1)\ f^n_{s\cdot 1}\ ...\ \sigma_s(d_n)\ f^n_{s\cdot n} \tag{3.2.2}$$

*$\sigma_s$ - перестановка множества переменных $\{d_1,...,d_n\}$*

$$\sigma_s=\begin{pmatrix}d_1&...&d_n\\ \sigma_s(d_1)&...&\sigma_s(d_n)\end{pmatrix}$$

Доказательство. Мы докажем утверждение индукцией по $n$.

При $n=1$ доказываемое утверждение является следствием теоремы 3.1.11. При этом мы можем отождествить[3.5]

$$f^1_{s\cdot p}=f_{s\cdot p}\quad p=0,1$$

[3.5] В представлении (3.2.2) мы будем пользоваться следующими правилами.

- Если область значений какого-либо индекса - это множество, состоящее из одного элемента, мы будем опускать соответствующий индекс.
- Если $n=1$, то $\sigma_s$ - тождественное преобразование. Это преобразование можно не указывать в выражении.

Допустим, что утверждение теоремы справедливо при $n = k-1$. Тогда отображение (3.2.1) можно представить в виде

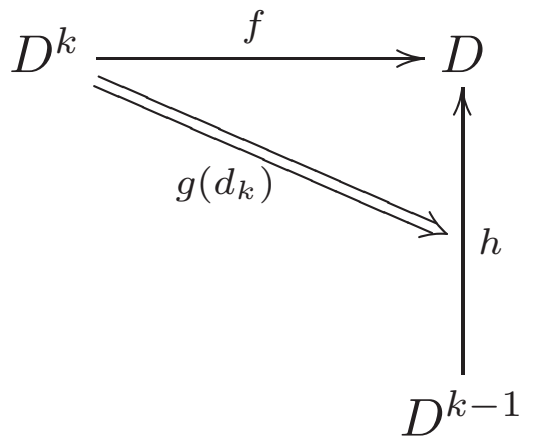

$$d = f(d_1, ..., d_k) = g(d_k)(d_1, ..., d_{k-1})$$

Согласно предположению индукции полилинейное отображение $h$ имеет вид

$$d = h^{k-1}_{t\cdot 0}\ \sigma_t(d_1)\ h^{k-1}_{t\cdot 1}\ ...\ \sigma_t(d_{k-1})\ h^{k-1}_{t\cdot k-1}$$

Согласно построению $h = g(d_k)$. Следовательно, выражения $h_{t\cdot p}$ являются функциями $d_k$. Поскольку $g(d_k)$ - линейная функция $d_k$, то только одно выражение $h_{t\cdot p}$ является линейной функцией переменной $d_k$, и остальные выражения ${}_{t\cdot q}h$ не зависят от $d_k$.

Не нарушая общности, положим $p = 0$. Согласно равенству (3.1.14) для заданного $t$

$$h^{k-1}_{t\cdot 0} = g_{tr\cdot 0}\ d_k\ g_{tr\cdot 1}$$

Положим $s = tr$ и определим перестановку $\sigma_s$ согласно правилу

$$\sigma_s = \sigma_{tr} = \begin{pmatrix} d_k & d_1 & ... & d_{k-1} \\ d_k & \sigma_t(d_1) & ... & \sigma_t(d_{k-1}) \end{pmatrix}$$

Положим

$$f^k_{tr\cdot q+1} = h^{k-1}_{t\cdot q} \qquad q = 1, ..., k-1$$
$$f^k_{tr\cdot q} = g_{tr\cdot q} \qquad q = 0, 1$$

Мы доказали шаг индукции. □

Определение 3.2.3. Выражение $f^n_{s\cdot p}$ в равенстве (3.2.2) называется **компонентой полилинейного отображения** $f$. □

Теорема 3.2.4. *Пусть $D$ - тело характеристики 0. Допустим $\overline{\overline{e}}$ - базис тела $D$ над полем $F \subset Z(D)$.* **Стандартное представление полилинейного отображения** *тела имеет вид*

(3.2.3) $$f(d_1, ..., d_n) = f^{i_0...i_n}_t\ e_{i_0}\ \sigma_t(d_1)\ e_{i_1}\ ...\ \sigma_t(d_n)\ e_{i_n}$$

*Индекс $t$ нумерует всевозможные перестановки $\sigma_t$ множества переменных $\{d_1, ..., d_n\}$. Выражение $f^{i_0...i_n}_t$ в равенстве (3.2.3) называется* **стандартной компонентой полилинейного отображения** *$f$.*

Доказательство. Компоненты полилинейного отображения $f$ имеют разложение

(3.2.4) $$f^n_{s\cdot p} = e_i f^{ni}_{s\cdot p}$$

относительно базиса $\overline{\overline{e}}$. Если мы подставим (3.2.4) в (3.2.2), мы получим

(3.2.5) $$d = f^{nj_1}_{s\cdot 0}\ e_{j_1}\ \sigma_s(d_1)\ f^{nj_2}_{s\cdot 1}\ e_{j_2}\ ...\ \sigma_s(d_n)\ f^{nj_n}_{s\cdot n}\ e_{j_n}$$

Рассмотрим выражение

$$f_t^{j_0...j_n} = f_{s\cdot 0}^{nj_1} \; ...f_{s\cdot n}^{nj_n} \tag{3.2.6}$$

В правой части подразумевается сумма тех слагаемых с индексом $s$, для которых перестановка $\sigma_s$ совпадает. Каждая такая сумма будет иметь уникальный индекс $t$. Подставив в равенство (3.2.5) выражение (3.2.6) мы получим равенство (3.2.3). □

Теорема 3.2.5. *Пусть $\overline{\overline{e}}$ - базис тела $D$ над полем $F \subset Z(D)$. Полилинейное отображение (3.2.1) можно представить в виде $D$-значной формы степени $n$ над полем $F \subset Z(D)$*[3.6]

$$f(a_1, ..., a_n) = a_1^{i_1}...a_n^{i_n} f_{i_1...i_n} \tag{3.2.7}$$

*где*

$$a_j = e_i a_j^i$$

$$f_{i_1...i_n} = f(e_{i_1}, ..., e_{i_n}) \tag{3.2.8}$$

*и величины $f_{i_1...i_n}$ являются координатами $D$-значного ковариантнго тензора над полем $F$.*

Доказательство. Согласно определению 3.2.1 равенство (3.2.7) следует из цепочки равенств

$$f(a_1, ..., a_n) = f(e_{i_1} a_1^{i_1}, ..., e_{i_n} a_n^{i_n}) = a_1^{i_1}...a_n^{i_n} f(e_{i_1}, ..., e_{i_n})$$

Пусть $\overline{\overline{e}}'$ - другой базис. Пусть

$$e'_i = e_j A_i^j \tag{3.2.9}$$

преобразование, отображающее базис $\overline{\overline{e}}$ в базис $\overline{\overline{e}}'$. Из равенств (3.2.9) и (3.2.8) следует

$$\begin{aligned} f'_{i_1...i_n} &= f(e'_{i_1}, ..., e'_{i_n}) \\ &= f(e_{j_1} A_{i_1}^{j_1}, ..., e'_{j_n} A_{i_n}^{j_n}) \\ &= A_{i_1}^{j_1}...A_{i_n}^{j_n} f(e_{j_1}, ..., e_{j_n}) \\ &= A_{i_1}^{j_1}...A_{i_n}^{j_n} f_{j_1...j_n} \end{aligned} \tag{3.2.10}$$

Из равенства (3.2.10) следует тензорный закон преобразования координат полилинейного отображения. Из равенства (3.2.10) и теоремы [3]-8.2.3 следует, что значение отображения $f(a_1, ..., a_n)$ не зависит от выбора базиса. □

Полилинейное отображение (3.2.1) **симметрично**, если

$$f(d_1, ..., d_n) = f(\sigma(d_1), ..., \sigma(d_n))$$

для любой перестановки $\sigma$ множества $\{d_1, ..., d_n\}$.

Теорема 3.2.6. *Если полилинейное отображение $f$ симметрично, то*

$$f_{i_1,...,i_n} = f_{\sigma(i_1),...,\sigma(i_n)} \tag{3.2.11}$$

[3.6] Теорема доказана по аналогии с теоремой в [2], с. 107, 108

Доказательство. Равенство (3.2.11) следует из равенства

$$\begin{aligned} a_1^{i_1}...a_n^{i_n}f_{i_1...i_n} =&f(a_1,...,a_n)\\ =&f(\sigma(a_1),...,\sigma(a_n))\\ =&a_1^{i_1}...a_n^{i_n}f_{\sigma(i_1)...\sigma(i_n)} \end{aligned}$$

□

Полилинейное отображение (3.2.1) **косо симметрично**, если

$$f(d_1,...,d_n)=|\sigma|f(\sigma(d_1),...,\sigma(d_n))$$

для любой перестановки $\sigma$ множества $\{d_1,...,d_n\}$. Здесь

$$|\sigma|=\begin{cases}1 & \text{перестановка } \sigma \text{ чётная}\\ -1 & \text{перестановка } \sigma \text{ нечётная}\end{cases}$$

Теорема 3.2.7. *Если полилинейное отображение $f$ косо симметрично, то*

(3.2.12) $$f_{i_1,...,i_n}=|\sigma|f_{\sigma(i_1),...,\sigma(i_n)}$$

Доказательство. Равенство (3.2.12) следует из равенства

$$\begin{aligned} a_1^{i_1}...a_n^{i_n}f_{i_1...i_n} =&f(a_1,...,a_n)\\ =&|\sigma|f(\sigma(a_1),...,\sigma(a_n))\\ =&a_1^{i_1}...a_n^{i_n}|\sigma|f_{\sigma(i_1)...\sigma(i_n)} \end{aligned}$$

□

Теорема 3.2.8. *Координаты полилинейного над полем $F$ отображения* (3.2.1) *и его компоненты относительно базиса $\overline{\overline{e}}$ удовлетворяют равенству*

(3.2.13) $$f_{j_1...j_n}=f_t^{i_0...i_n}C_{i_0\sigma_t(j_1)}^{k_1}C_{k_1i_1}^{l_1}...B_{l_{n-1}\sigma_t(j_n)}^{k_n}C_{k_ni_n}^{l_n}e_{l_n}$$

(3.2.14) $$f_{j_1...j_n}^{p}=f_t^{i_0...i_n}C_{i_0\sigma_t(j_1)}^{k_1}C_{k_1i_1}^{l_1}...C_{l_{n-1}\sigma_t(j_n)}^{k_n}C_{k_ni_n}^{p}$$

Доказательство. В равенстве (3.2.3) положим

$$d_i=e_{j_i}d_i^{j_i}$$

Тогда равенство (3.2.3) примет вид

(3.2.15) $$\begin{aligned} f(d_1,...,d_n)&=f_t^{i_0...i_n}\ e_{i_0}\sigma_t(d_1^{j_1}e_{j_1})e_{i_1}...\sigma_t(d_n^{j_n}e_{j_n})e_{i_n}\\ &=d_1^{j_1}...d_n^{j_n}f_t^{i_0...i_n}e_{i_0}\sigma_t(e_{j_1})e_{i_1}...\sigma_t(e_{j_n})e_{i_n}\\ &=d_1^{j_1}...d_n^{j_n}f_t^{i_0...i_n}C_{i_0\sigma_t(j_1)}^{k_1}C_{k_1i_1}^{l_1}\\ &\quad ...C_{l_{n-1}\sigma_t(j_n)}^{k_n}C_{k_ni_n}^{l_n}e_{l_n} \end{aligned}$$

Из равенства (3.2.7) следует

(3.2.16) $$f(a_1,...,a_n)=e_pf_{i_1...i_n}^{p}a_1^{i_1}...a_n^{i_n}$$

Равенство (3.2.13) следует из сравнения равенств (3.2.15) и (3.2.7). Равенство (3.2.14) следует из сравнения равенств (3.2.15) и (3.2.16). □

Глава 4

# Линейное отображение $D$-векторных пространств

## 4.1. Линейное отображение $D$-векторных пространств

При изучении линейного отображения $D$-векторных пространств мы будем рассматривать тело $D$ как конечно мерную алгебру над полем $F$.

Определение 4.1.1. Пусть поле $F$ является подкольцом центра $Z(D)$ тела $D$. Пусть $V$ и $W$ - $D$-векторные пространства. Мы будем называть отображение

$$A : V \to W$$

$D$-векторного пространства $V$ в $D$-векторное пространство $W$ **линейным отображением**, если

$$A(x+y) = A(x) + A(y) \quad x, y \in V$$
$$A(px) = pA(x) \quad p \in F$$

Обозначим $\mathcal{L}(D; V; W)$ множество линейных отображений

$$A : V \to W$$

$D$-векторного пространства $V$ в $D$-векторное пространство $W$. □

Очевидно, что линейное отображение $D_*{}^*$-векторного пространства, так же как линейное отображение ${}^*{}_*D$-векторного пространства являются линейными отображениями. Множество морфизмов $D$-векторного пространства шире, чем множество морфизмов $D_*{}^*$-векторного пространства. Чтобы рассмотреть линейное отображение векторных пространств, мы будем следовать методике, предложенной в разделе [3]-4.4.

Теорема 4.1.2. *Пусть $D$ - тело характеристики 0. Пусть $V$, $W$ - $D$-векторные пространства. Линейное отображение*

$$A : V \to W$$

*относительно базиса $\overline{\overline{e}}_V$ ${}_*{}^*D$-векторного пространства $V$ и базиса $\overline{\overline{e}}_W$ ${}_*{}^*D$-векторного пространства $W$ имеет вид*

$$A(v) = e_{W\cdot j}\ A^j_i(v^i) \qquad v = e_{V*}{}^*v \tag{4.1.1}$$

*где $A^j_i(v^i)$ линейно зависит от одной переменной $v^i$ и не зависит от остальных координат вектора $v$.*

Доказательство. Согласно определению 4.1.1

$$A(v) = A(e_{V*}{}^*v) = A\left(\sum_i e_{V\cdot i}\ v^i\right) = \sum_i A(e_{V\cdot i}\ v^i) \tag{4.1.2}$$

Для любого заданного $i$ вектор $A(e_{V\cdot i}\ v^i) \in W$ имеет единственное разложение

$$A(e_{V\cdot i}\ v^i) = e_{W\cdot j}\ A^j_i(v^i) \qquad A(e_{V\cdot i}\ v^i) = e_{W*}{}^*A_i(v^i) \tag{4.1.3}$$

относительно базиса $\overline{\overline{e}}_W{}_*{}^*D$-векторного пространства. Подставив (4.1.3) в (4.1.2), мы получим (4.1.1). □

Определение 4.1.3. Линейное отображение

$$A^j_i : D \to D$$

называется **частным линейным отображением** переменной $v^i$. □

Мы можем записать линейное отображение в виде произведения матриц

$$A(v) = \begin{pmatrix} e_{W\cdot 1} & ... & e_{W\cdot m}\end{pmatrix} {}_*{}^* \begin{pmatrix} A^1_i(v^i) \\ ... \\ A^m_i(v^i)\end{pmatrix} \tag{4.1.4}$$

Определим произведение матриц

$$\begin{pmatrix} A^1_1 & ... & A^1_n \\ ... & ... & ... \\ A^m_1 & ... & A^m_n \end{pmatrix} {}_\circ{}^\circ \begin{pmatrix} v^1 \\ ... \\ v^n \end{pmatrix} = \begin{pmatrix} A^1_i(v^i) \\ ... \\ A^m_i(v^i)\end{pmatrix} \tag{4.1.5}$$

где $A = \begin{pmatrix} A^j_i \end{pmatrix}$ - матрица частных линейных отображений. Используя равенство (4.1.5), мы можем записать равенство (4.1.4) в виде

$$A(v) = \begin{pmatrix} e_{W\cdot 1} & ... & e_{W\cdot m}\end{pmatrix} {}_*{}^* \begin{pmatrix} A^1_1 & ... & A^1_n \\ ... & ... & ... \\ A^m_1 & ... & A^m_n \end{pmatrix} {}_\circ{}^\circ \begin{pmatrix} v^1 \\ ... \\ v^n \end{pmatrix} \tag{4.1.6}$$

Теорема 4.1.4. *Пусть $D$ - тело характеристики $0$. Линейное отображение*

$$A : v \in V \to w \in W \quad w = A(v) \tag{4.1.7}$$

*относительно базиса $\overline{\overline{e}}_V$ в $D$-векторном пространстве $V$ и базиса $\overline{\overline{e}}_W$ в $D$-векторном пространстве $W$*[4.1] *имеет вид*

$$v = e_{V*}{}^*v$$

$$w = e_{W*}{}^*w \tag{4.1.8}$$

$$w^j = A^j_i(v^i) = A_{s\cdot 0\cdot}{}^j_i\ v^i\ A_{s\cdot 1\cdot}{}^j_i \tag{4.1.9}$$

[4.1] Координатная запись отображения (4.1.7) зависит от выбора базиса. Равенства изменят свой вид, если например мы выберем ${}^*{}_*D$-базис $\overline{\overline{r}}_*$ в $D$-векторном пространстве $W$.

Доказательство. Согласно теореме 4.1.2 линейное отображение $A(v)$ можно записать в виде (4.1.1). Так как для заданных индексов $i$, $j$ частное линейное отображение $A^j_i(v^i)$ линейно по переменной $v^i$, то согласно (3.1.14) выражение $A^j_i(v^i)$ можно представить в виде

$$A^j_i(v^i) = A_{s\cdot 0\cdot}{}^j_i\ v^i\ A_{s\cdot 1\cdot}{}^j_i \tag{4.1.10}$$

где индекс $s$ нумерует слагаемые. Множество значений индекса $s$ зависит от индексов $i$ и $j$. Комбинируя равенства (4.1.2) и (4.1.10), мы получим

$$\begin{aligned} A(v) &= e_{W\cdot j}A^j_i(v^i) &&= e_{W\cdot j}A_{s\cdot 0\cdot}{}^j_i\ v^i\ A_{s\cdot 1\cdot}{}^j_i \\ A(v) &= e_{W*}{}^*A_i(v^i) &&= e_{W*}{}^*(A_{s\cdot 0\cdot i}\ v^i\ A_{s\cdot 1\cdot i}) \end{aligned} \tag{4.1.11}$$

В равенстве (4.1.11) мы суммируем также по индексу $i$. Равенство (4.1.9) следует из сравнения равенств (4.1.8) и (4.1.11). □

Определение 4.1.5. Выражение $A_{s\cdot p\cdot}{}^j_i$ в равенстве (4.1.3) называется **компонентой линейного отображения** $A$. □

Теорема 4.1.6. *Пусть $D$ - тело характеристики 0. Пусть $\overline{\overline{e}}_V$ - базис в $D$-векторном пространстве $V$, $\overline{\overline{e}}_U$ - базис в $D$-векторном пространстве $U$, и $\overline{\overline{e}}_W$ - базис в $D$-векторном пространстве $W$. Предположим, что мы имеем коммутативную диаграмму отображений*

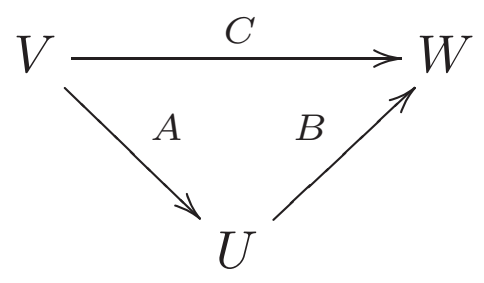

*где линейное отображение $A$ имеет представление*

$$u = A(v) = e_{U\cdot j}A^j_i(v^i) = e_{U\cdot j}A_{s\cdot 0\cdot}{}^j_i\ v^i\ A_{s\cdot 1\cdot}{}^j_i \tag{4.1.12}$$

*относительно заданных базисов и линейное отображение $B$ имеет представление*

$$w = B(u) = e_{W\cdot k}B^k_j(u^j) = e_{W\cdot k}B_{t\cdot 0\cdot}{}^k_j\ u^j\ B_{t\cdot 1\cdot}{}^k_j \tag{4.1.13}$$

*относительно заданных базисов. Тогда отображение $C$ линейно и имеет представление*

$$w = C(v) = e_{W\cdot k}C^k_i(v^i) = e_{W\cdot k}C_{u\cdot 0\cdot}{}^k_i\ v^i\ C_{u\cdot 1\cdot}{}^k_i \tag{4.1.14}$$

*относительно заданных базисов, где*[4.2]

$$\begin{aligned} & C^k_i(v^i) && = B^k_j(A^j_i(v^i)) \\ C_{u\cdot 0\cdot}{}^k_i &= C_{ts\cdot 0\cdot}{}^k_i && = B_{t\cdot 0\cdot}{}^k_j\ A_{s\cdot 0\cdot}{}^j_i \\ C_{u\cdot 1\cdot}{}^k_i &= C_{ts\cdot 1\cdot}{}^k_i && = A_{s\cdot 1\cdot}{}^j_i\ B_{t\cdot 1\cdot}{}^k_j \end{aligned} \tag{4.1.15}$$

---

[4.2] Индекс $u$ оказался составным индексом, $u = st$. Однако не исключено, что некоторые слагаемые в (4.1.15) могут быть объединены вместе.

Доказательство. Отображение $C$ линейно, так как

$$\begin{aligned}
C(a+b) &= B(A(a+b))\\
&= B(A(a)+A(b))\\
&= B(A(a))+B(A(b))\\
&= C(a)+C(b)\\
C(ab) &= B(A(ab)) = B(aA(b))\\
&= aB(A(b)) = aC(b)\\
& a,b \in V\\
& a \in F
\end{aligned}$$

Равенство (4.1.14) следует из подстановки (4.1.12) в (4.1.13). □

Теорема 4.1.7. *Для линейного отображения $A$ существует линейное отображение $B$ такое, что*

$$\begin{aligned}
A(axb) &= B(x)\\
B_{s\cdot 0\cdot}{}^j_i &= A_{s\cdot 0\cdot}{}^j_i\ a\\
B_{s\cdot 1\cdot}{}^j_i &= b\ A_{s\cdot 1\cdot}{}^j_i
\end{aligned}$$

Доказательство. Отображение $B$ линейно, так как

$$\begin{aligned}
B(x+y) &= A(a(x+y)b) = A(axb+ayb) = A(axb)+A(ayb) = B(x)+B(y)\\
B(cx) &= A(acxb) = A(caxb) = cA(axb) = cB(x)\\
& x,y \in V, c \in F
\end{aligned}$$

Согласно равенству (4.1.9)

$$B_{s\cdot 0\cdot}{}^j_i\ v^i\ B_{s\cdot 1\cdot}{}^j_i = A_{s\cdot 0\cdot}{}^j_i\ (a\ v^i\ b)A_{s\cdot 1\cdot}{}^j_i = (A_{s\cdot 0\cdot}{}^j_i\ a)v^i\ (b\ A_{s\cdot 1\cdot}{}^j_i)$$

□

Теорема 4.1.8. *Пусть $D$ - тело характеристики* 0. *Пусть*

$$A : V \to W$$

*линейное отображение $D$-векторного пространства $V$ в $D$-векторное пространство $W$. Тогда* $A(0) = 0$.

Доказательство. Следствие равенства

$$A(a+0) = A(a) + A(0)$$

□

Определение 4.1.9. Множество

$$\ker f = \{x \in V : f(x) = 0\}$$

называется **ядром линейного отображения**

$$A : V \to W$$

$D$-векторного пространства $V$ в $D$-векторное пространство $W$. □

Определение 4.1.10. Линейное отображение

$$A : V \to W$$

$D$-векторного пространства $V$ в $D$-векторное пространство $W$ называется **вырожденным**, если

$$\mathrm{ker} A = V$$

□

## 4.2. Полилинейное отображение $D$-векторного пространства

Определение 4.2.1. Пусть поле $F$ является подкольцом центра $Z(D)$ тела $D$ характеристики 0. Пусть $V_1$, ..., $V_n$, $W_1$, ..., $W_m$ - $D$-векторные пространства. Мы будем называть отображение

$$A : V_1 \times ... \times V_n \to W_1 \times ... \times W_m$$
$$w_1 \times ... \times w_m = A(v_1, ..., v_n) \tag{4.2.1}$$

**полилинейным отображением** $\times$-$D$-векторного пространства $V_1 \times ... \times V_n$ в $\times$-$D$-векторное пространство $W_1 \times ... \times W_m$, если

$$A(p_1, ..., p_i + q_i, ..., p_n) = A(p_1, ..., p_i, ..., p_n) + A(p_1, ..., q_i, ..., p_n)$$
$$A(p_1, ..., ap_i, ..., p_n) = aA(p_1, ..., p_i, ..., p_n)$$

$$1 \le i \le n \quad p_i, q_i \in V_i \quad a \in F$$

□

Определение 4.2.2. Обозначим $\mathcal{L}(D; V_1, ..., V_n; W_1, ..., W_m)$ множество полилинейных отображений $\times$-$D$-векторного пространства $V_1 \times ... \times V_n$ в $\times$-$D$-векторное пространство $W_1 \times ... \times W_m$. □

Теорема 4.2.3. *Пусть $D$ - тело характеристики 0. Для каждого $k \in K = [1, n]$ допустим $\overline{\overline{e}}_{V_k}$ - базис в $D$-векторном пространстве $V_k$ и*

$$v_k = v_k{}^*{}_* e_{V_k} \quad v_k \in V_k$$

*Для каждого $l$, $1 \le l \le m$, допустим $\overline{\overline{e}}_{W_l}$ - базис в $D$-векторном пространстве $W_l$ и*

$$w_l = w_l{}^*{}_* e_{W_l} \quad w_l \in W_l$$

*Полилинейное отображение (4.2.1) относительно базисов $\overline{\overline{e}}_{V_1}$, ..., $\overline{\overline{e}}_{V_n}$, $\overline{\overline{e}}_{W_1}$, ..., $\overline{\overline{e}}_{W_m}$ имеет вид*

$$\begin{aligned} w_l^{\boldsymbol{j}} &= A_{l\cdot}{}^{\boldsymbol{j}}_{\boldsymbol{i}_1...\boldsymbol{i}_n}(v_1^{\boldsymbol{i}_1}, ..., v_n^{\boldsymbol{i}_n}) \\ &= A^n_{s\cdot 0\cdot l\cdot}{}^{\boldsymbol{j}}_{\boldsymbol{i}_1...\boldsymbol{i}_n} \ \sigma_s(v_1^{\boldsymbol{i}_1}) \ A^n_{s\cdot 1.l\cdot}{}^{\boldsymbol{j}}_{\boldsymbol{i}_1...\boldsymbol{i}_n} \ ... \ \sigma_s(v_n^{\boldsymbol{i}_n}) \ A^n_{s\cdot n\cdot l\cdot}{}^{\boldsymbol{j}}_{\boldsymbol{i}_1...\boldsymbol{i}_n} \end{aligned} \tag{4.2.2}$$

*Область значений $S$ индекса $s$ зависит от значений индексов $\boldsymbol{i}_1$, ..., $\boldsymbol{i}_n$. $\sigma_s$ - перестановка множества переменных $\{v_1^{\boldsymbol{i}_1}, ..., v_n^{\boldsymbol{i}_n}\}$.*

Доказательство. Так как отображение $A$ в $\times$-$D$-векторное пространство $W_1 \times ... \times W_m$ можно рассматривать покомпонентно, то мы можем ограничиться рассмотрением отображения

$$A_l : V_1 \times ... \times V_n \to W_l \quad w_l = A_l(v_1, ..., v_n) \tag{4.2.3}$$

Мы докажем утверждение индукцией по $n$.

При $n = 1$ доказываемое утверждение является утверждением теоремы 4.1.4. При этом мы можем отождествить[4.3]

$$A^{1}_{s\cdot p\cdot}{}^{j}_{i} = A_{s\cdot p\cdot}{}^{j}_{i} \quad p = 0, 1$$

Допустим, что утверждение теоремы справедливо при $n = k - 1$. Тогда отображение (4.2.3) можно представить в виде

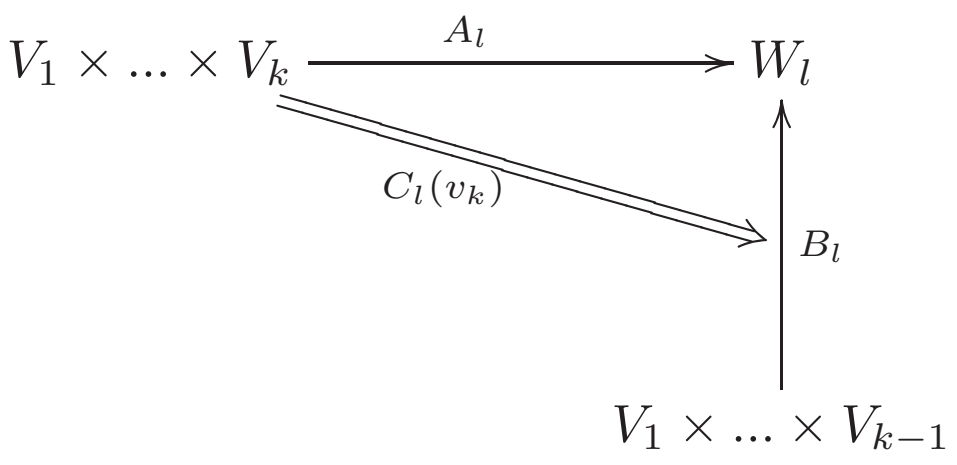

$$w_l = A_l(v_1, ..., v_k) = C_l(v_k)(v_1, ..., v_{k-1})$$

Согласно предположению индукции полилинейное отображение $B_l$ имеет вид

$$w_l^{j} = B^{k-1}_{t\cdot 0\cdot l\cdot}{}^{j}_{i_1...i_{k-1}}\ \sigma_t(v_1^{i_1})\ B^{k-1}_{t\cdot 1.l\cdot}{}^{j}_{i_1...i_{k-1}}\ ...\ \sigma_t(v_{k-1}^{i_{k-1}})\ B^{k-1}_{t\cdot k-1.l\cdot}{}^{j}_{i_1...i_{k-1}}$$

Согласно построению $B_l = C_l(v_k)$. Следовательно, выражения $B^{k-1}_{t\cdot p\cdot l\cdot}{}^{j}_{i_1...i_{k-1}}$ являются функциями $v_k$. Поскольку $C_l(v_k)$ - линейная функция $v_k$, то только одно выражение $B^{k-1}_{t\cdot p\cdot l\cdot}{}^{j}_{i_1...i_{k-1}}$ является линейной функцией $v_k$, и остальные выражения $B^{k-1}_{t\cdot q\cdot l\cdot}{}^{j}_{i_1...i_{k-1}}$ не зависят от $v_k$.

Не нарушая общности, положим $p = 0$. Согласно теореме 4.1.4

$$B^{k-1}_{t\cdot 0\cdot l\cdot}{}^{j}_{i_1...i_{k-1}} = C^{k}_{tr\cdot 0\cdot l\cdot}{}^{j}_{i_k i_1...i_{k-1}}\ v_k^{i_k}\ C^{k}_{tr\cdot 1.l\cdot}{}^{j}_{i_k i_1...i_{k-1}}$$

Положим $s = tr$ и определим перестановку $\sigma_s$ согласно правилу

$$\sigma_s = \sigma(tr) = \begin{pmatrix} v_k^{i_k} & v_1^{i_1} & ... & v_{k-1}^{i_{k-1}} \\ v_k^{i_k} & \sigma_t(v_1^{i_1}) & ... & \sigma_t(v_{k-1}^{i_{k-1}}) \end{pmatrix}$$

Положим

$$A^{k}_{tr\cdot q+1.l\cdot}{}^{j}_{i_k i_1...i_{k-1}} = B^{k-1}_{t\cdot q\cdot l\cdot}{}^{j}_{i_1...i_{k-1}} \qquad q = 1, ..., k-1$$

$$A^{k}_{tr\cdot q\cdot l\cdot}{}^{j}_{i_k i_1...i_{k-1}} = C^{k}_{tr\cdot q\cdot l\cdot}{}^{j}_{i_k i_1...i_{k-1}} \qquad q = 0, 1$$

Мы доказали шаг индукции. □

Определение 4.2.4. Выражение $A^{n}_{s\cdot p\cdot l\cdot}{}^{j}_{i_1...i_n}$ в равенстве (4.2.2) называется **компонентой полилинейного отображения** $A$. □

---

[4.3] В представлении (4.2.2) мы будем пользоваться следующими правилами.

- Если область значений какого-либо индекса - это множество, состоящее из одного элемента, мы будем опускать соответствующий индекс.
- Если $n = 1$, то $\sigma_s$ - тождественное преобразование. Это преобразование можно не указывать в выражении.

# Часть 2

# Математический анализ над телом

Глава 5

# Дифференцируемые отображения

## 5.1. Топологическое тело

Определение 5.1.1. Тело $D$ называется **топологическим телом**[5.1], если $D$ является топологическим пространством, и алгебраические операции, определённые в $D$, непрерывны в топологическом пространстве $D$. □

Согласно определению, для произвольных элементов $a, b \in D$ и для произвольных окрестностей $W_{a-b}$ элемента $a-b$, $W_{ab}$ элемента $ab$ существуют такие окрестности $W_a$ элемента $a$ и $W_b$ элемента $b$, что $W_a - W_b \subset W_{a-b}$, $W_a W_b \subset W_{ab}$. Если $a \neq 0$, то для произвольной окрестности $W_{a^{-1}}$ существует окрестность $W_a$ элемента $a$, удовлетворяющая условию $W_a^{-1} \subset W_{a^{-1}}$.

Определение 5.1.2. **Норма на теле** $D$[5.2] - это отображение

$$d \in D \to |d| \in R$$

такое, что

- $|a| \geq 0$
- $|a| = 0$ равносильно $a = 0$
- $|ab| = |a|\ |b|$
- $|a+b| \leq |a| + |b|$

Тело $D$, наделённое структурой, определяемой заданием на $D$ нормы, называется **нормированным телом**. □

Инвариантное расстояние на аддитивной группе тела $D$

$$d(a, b) = |a - b|$$

определяет топологию метрического пространства, согласующуюся со структурой тела в $D$.

Определение 5.1.3. Пусть $D$ - нормированное тело. Элемент $a \in D$ называется **предел последовательности** $\{a_n\}$

$$a = \lim_{n \to \infty} a_n$$

если для любого $\epsilon \in R$, $\epsilon > 0$, существует, зависящее от $\epsilon$, натуральное число $n_0$ такое, что $|a_n - a| < \epsilon$ для любого $n > n_0$. □

Теорема 5.1.4. *Пусть $D$ - нормированное тело характеристики $0$ и пусть $d \in D$. Пусть $a \in D$ - предел последовательности $\{a_n\}$. Тогда*

$$\lim_{n \to \infty} (a_n d) = ad$$

[5.1] Определение дано согласно определению из [6], глава 4

[5.2] Определение дано согласно определению из [4], гл. IX, §3, п°2

$$\lim_{n\to\infty}(da_n)=da$$

Доказательство. Утверждение теоремы тривиально, однако я привожу доказательство для полноты текста. Поскольку $a\in D$ - предел последовательности $\{a_n\}$, то согласно определению 5.1.3 для заданного $\epsilon\in R$, $\epsilon>0$, существует натуральное число $n_0$ такое, что $|a_n-a|<\epsilon/|d|$ для любого $n>n_0$. Согласно определению 5.1.2 утверждение теоремы следует из неравенств

$$|a_nd-ad|=|(a_n-a)d|=|a_n-a||d|<\epsilon/|d||d|=\epsilon$$
$$|da_n-da|=|d(a_n-a)|=|d||a_n-a|<|d|\epsilon/|d|=\epsilon$$

для любого $n>n_0$. □

Определение 5.1.5. Пусть $D$ - нормированное тело. Последовательность $\{a_n\}$, $a_n\in D$ называется **фундаментальной** или **последовательностью Коши**, если для любого $\epsilon\in R$, $\epsilon>0$, существует, зависящее от $\epsilon$, натуральное число $n_0$ такое, что $|a_p-a_q|<\epsilon$ для любых $p$, $q>n_0$. □

Определение 5.1.6. Нормированное тело $D$ называется **полным** если любая фундаментальная последовательность элементов данного тела сходится, т. е. имеет предел в этом теле. □

До сих пор всё было хорошо. Кольцо целых чисел содержится в кольце характеристики 0. Поле рациональных чисел содержится в теле характеристики 0. Однако это не означает, что поле действительных чисел содержится в полном теле $D$ характеристики 0. Например, мы можем определить норму на теле, тождественно равную 1. Эта норма определяет на нормированном теле дискретную топологию и не представляет для нас интереса, так как любая фундаментальная последовательность является постоянным отображением.

В дальнейшем, говоря о нормированном теле характеристики 0, мы будем предполагать, что определён гомеоморфизм поля рациональных чисел $Q$ в тело $D$.

Теорема 5.1.7. *Полное тело $D$ характеристики 0 содержит в качестве подполя изоморфный образ поля $R$ действительных чисел. Это поле обычно отождествляют с $R$.*

Доказательство. Рассмотрим фундаментальную последовательность рациональных чисел $\{p_n\}$. Пусть $p'$ - предел этой последовательности в теле $D$. Пусть $p$ - предел этой последовательности в поле $R$. Так как вложение поля $Q$ в тело $D$ гомеоморфно, то мы можем отождествить $p'\in D$ и $p\in R$. □

Теорема 5.1.8. *Пусть $D$ - полное тело характеристики 0 и пусть $d\in D$. Тогда любое действительное число $p\in R$ коммутирует с $d$.*

Доказательство. Мы можем представить действительное число $p\in R$ в виде фундаментальной последовательности рациональных чисел $\{p_n\}$. Утверждение теоремы следует из цепочки равенств

$$pd=\lim_{n\to\infty}(p_nd)=\lim_{n\to\infty}(dp_n)=dp$$

основанной на утверждении теоремы 5.1.4. □

Теорема 5.1.9. *Пусть $D$ - полное тело характеристики 0. Тогда поле действительных чисел $R$ является подполем центра $Z(D)$ тела $D$.*

Доказательство. Следствие теоремы 5.1.8. □

Определение 5.1.10. Пусть $D$ - полное тело характеристики 0. Множество элементов $d\in D$, $|d|=1$, называется **единичной сферой в теле** $D$. □

Определение 5.1.11. Пусть $D_1$ - полное тело характеристики 0 с нормой $|x|_1$. Пусть $D_2$ - полное тело характеристики 0 с нормой $|x|_2$. Отображение

$$f:D_1\rightarrow D_2$$

называется **непрерывной**, если для любого сколь угодно малого $\epsilon>0$ существует такое $\delta>0$, что

$$|x'-x|_1<\delta$$

влечёт

$$|f(x')-f(x)|_2<\epsilon$$

□

Теорема 5.1.12. *Пусть $D$ - полное тело характеристики 0. Если в разложении (3.1.14) линейного отображения*

$$f:D\rightarrow D$$

*индекс $s$ принимает конечное множество значений, то линейное отображение $f$ непрерывно.*

Доказательство. Положим $x'=x+a$. Тогда

$$f(x')-f(x)=f(x+a)-f(x)=f(a)={}_{s\cdot 0}f\ a\ {}_{s\cdot 1}f$$

$$|f(x')-f(x)|=|{}_{s\cdot 0}f\ a\ {}_{s\cdot 1}f|<(|{}_{s\cdot 0}f|\ |{}_{s\cdot 1}f|)|a|$$

Положим $F=|{}_{s\cdot 0}f|\ |{}_{s\cdot 1}f|$. Тогда

$$|f(x')-f(x)|<F|a|$$

Выберем $\epsilon>0$ и положим $a=\dfrac{\epsilon}{F}e$. Тогда $\delta=|a|=\dfrac{\epsilon}{F}$. Согласно определению 5.1.11 линейное отображение $f$ непрерывно. □

Аналогично, если индекс $s$ принимает счётное множество значений, то для непрерывности аддитивного отображения $f$ мы требуем, что бы ряд $|{}_{s\cdot 0}f|\ |{}_{s\cdot 1}f|$ сходился. Если индекс $s$ принимает непрерывное множество значений, то для непрерывности аддитивного отображения $f$ мы требуем, что бы интеграл $\int|{}_{s\cdot 0}f|\ |{}_{s\cdot 1}f|ds$ существовал.

Определение 5.1.13. Пусть

$$f:D_1\rightarrow D_2$$

отображение полного тела $D_1$ характеристики 0 с нормой $|x|_1$ в полное тело $D_2$ характеристики 0 с нормой $|y|_2$. Величина

$$\|f\|=\sup\frac{|f(x)|_2}{|x|_1}\tag{5.1.1}$$

называется **нормой отображения** $f$. □

Теорема 5.1.14. *Пусть*

$$f:D_1\rightarrow D_2$$

*линейное отображение полного тела $D_1$ в полное тело $D_2$. Отображение $f$ непрерывно, если $\|f\|<\infty$.*

Доказательство. Поскольку отображение $f$ аддитивно, то согласно определению 5.1.13

$$|f(x)-f(y)|_2=|f(x-y)|_2\le\|f\|\ |x-y|_1$$

Возьмём произвольное $\epsilon>0$. Положим $\delta=\frac{\epsilon}{\|f\|}$. Тогда из неравенства

$$|x-y|_1<\delta$$

следует

$$|f(x)-f(y)|_2\le\|f\|\ \delta=\epsilon$$

Согласно определению 5.1.11 отображение $f$ непрерывно. □

Теорема 5.1.15. *Пусть отображение*

$$f:D_1\to D_2$$

*является непрерывным линейным отображением тела $D_1$ характеристики $0$ в тело $D_2$ характеристики $0$. Тогда*

$$f(ax)=af(x)$$

*для любого действительного числа $a$.*

Доказательство. Согласно теореме 5.1.9, поле действительных чисел является подполем центра $Z(D)$ тела $D$. Мы можем представить действительное число $p\in R$ в виде фундаментальной последовательности рациональных чисел $\{p_n\}$. Утверждение теоремы следует из теорем 3.1.3, 5.1.4 и цепочки равенств

$$f(ph)=\lim_{n\to\infty}f(p_nh)=\lim_{n\to\infty}(p_nf(h))=pf(h) \tag{5.1.2}$$

□

Теорема 5.1.16. *Пусть $D$ - полное тело характеристики $0$. Либо непрерывное отображение $f$ тела, проективное над полем $P$, не зависит от направления над полем $P$, либо значение $f(0)$ не определено.*

Доказательство. Согласно определению 3.1.23, отображение $f$ постоянно на направлении $Pa$. Так как $0\in Pa$, то естественно положить по непрерывности

$$f(0)=f(a)$$

Однако это приводит к неопределённости значения отображения $f$ в направлении $0$, если отображение $f$ имеет разное значение для разных направлений $a$. □

Если проективное над полем $R$ отображение $f$ непрерывно, то мы будем говорить, что **отображение $f$ непрерывно по направлению над полем $R$**. Поскольку для любого $a\in D$, $a\ne0$ мы можем выбрать $a_1=|a|^{-1}a$, $f(a_1)=f(a)$, то мы можем сделать определение более точным.

Определение 5.1.17. Пусть $D$ - полное тело характеристики $0$. Проективная над полем $R$ функция $f$ непрерывна по направлению над полем $R$, если для любого сколь угодно малого $\epsilon>0$ существует такое $\delta>0$, что

$$|x'-x|_1<\delta\quad |x'|_1=|x|_1=1$$

влечёт

$$|f(x')-f(x)|_2<\epsilon$$

□

Теорема 5.1.18. *Пусть $D$ - полное тело характеристики 0. Проективная над полем $R$ функция $f$ непрерывна по направлению над полем $R$ тогда и только тогда, когда эта функция непрерывна на единичной сфере тела $D$.*

Доказательство. Следствие определений 5.1.11, 3.1.23, 5.1.17. □

## 5.2. Дифференцируемое отображение тела

Определение 5.2.1. Пусть $D$ - нормированное тело.[5.3] Функция

$$f : D \to D$$

называется **$D\star$-дифференцируемой по Фреше** на множестве $U \subset D$[5.4], если в каждой точке $x \in U$ изменение функции $f$ может быть представлено в виде

$$f(x+h) - f(x) = h\frac{df(x)}{d_* x} + o(h) \tag{5.2.1}$$

где $o$ - такое непрерывное отображение

$$o : D \to D$$

что

$$\lim_{h\to 0} \frac{|o(h)|}{|h|} = 0 \tag{5.2.2}$$

□

Согласно определению 5.2.1, **$D\star$-производная Фреше отображения $f$ в точке $x$**

$$\frac{df(x)}{d_* x} \in {}_*\mathcal{L}(D, D)$$

порождает гомоморфизм

$$\Delta f = \Delta x \frac{df(x)}{d_* x}$$

отображающий приращение аргумента в приращение функции.

Если мы умножим обе части равенства (5.2.1) на $h^{-1}$, то мы получим равенство

$$h^{-1}(f(x+h) - f(x)) = \frac{df(x)}{d_* x} + h^{-1}o(h) \tag{5.2.3}$$

---

[5.3] Определение дано согласно определению [1]-3.1.1, стр. 256.

[5.4] Так же как в замечании [3]-4.4.5 мы можем определить $D\star$-дифференцируемость отображения

$$f : S \to D$$

нормированного тела $S$ в нормированное тело $D$ согласно правилу

$$f(x+h) - f(x) = F(h)\frac{df(x)}{d_* x} + o(h)$$

где

$$F : S \to D$$

гомоморфизм тел. Однако, опираясь на теоремы об изоморфизмах, мы можем ограничиться случаем отображений тела $D$ в $D$.

Из равенств (5.2.3) и (5.2.2) следует альтернативное определение $D\star$-производной Фреше

$$\frac{df(x)}{d_* x}=\lim_{h\to 0}(h^{-1}(f(x+h)-f(x))) \tag{5.2.4}$$

Отличие $D\star$-производной Фреше отображения тела $D$ от производной отображения поля $F$ состоит в том, что прежде всего $D\star$-производная Фреше не является линейным отображением. Более точно, $D\star$-производная Фреше является $\star D$-линейным отображением, однако, вообще говоря, не является $D\star$-линейным отображением. Действительно,

$$f(x+h)a-f(x)a=h\frac{df}{d_* x}a+o(h)$$

$$af(x+h)-af(x)=ah\frac{df}{d_* x}+o(h) \tag{5.2.5}$$

Однако, вообще говоря,

$$ah\frac{df}{d_* x}\ne ha\frac{df}{d_* x}$$

Подобная проблема возникает при дифференцировании произведения функций. Согласно определению 5.2.1

$$f(x+h)g(x+h)-f(x)g(x)=h\frac{df(x)g(x)}{d_* x}+o(h) \tag{5.2.6}$$

Выражение в левой части можно представить в виде

$$\begin{aligned}&f(x+h)g(x+h)-f(x)g(x)\\ =&f(x+h)g(x+h)-f(x)g(x+h)+f(x)g(x+h)-f(x)g(x)\\ =&(f(x+h)-f(x))g(x+h)+f(x)(g(x+h)-g(x))\end{aligned} \tag{5.2.7}$$

Согласно определению 5.2.1, мы можем представить (5.2.7) в виде

$$\begin{aligned}&f(x+h)g(x+h)-f(x)g(x)\\ =&\left(h\frac{df(x)}{d_* x}+o(h)\right)\left(g(x)+h\frac{dg(x)}{d_* x}+o(h)\right)\\ &+f(x)\left(h\frac{dg(x)}{d_* x}+o(h)\right)\end{aligned} \tag{5.2.8}$$

Так как, вообще говоря,

$$f(x)h\frac{dg}{d_* x}\ne hf(x)\frac{dg}{d_* x}$$

то естественно ожидать, что функция $f(x)g(x)$ не является $D\star$-дифференцируемой по Фреше.

Таким образом, определение $D\star$-производной Фреше весьма ограничено и не удовлетворяет стандартному определению операции дифференцирования. Чтобы найти решение проблемы дифференцирования, рассмотрим эту проблему с другой стороны.

Пример 5.2.2. Рассмотрим приращение функции $f(x)=x^2$.

$$\begin{aligned} f(x+h)-f(x)&=(x+h)^2-x^2\\ &=xh+hx+h^2\\ &=xh+hx+o(h) \end{aligned}$$

Как мы видим компонента приращения функции $f(x)=x^2$, линейно зависящая от приращения аргумента, имеет вид

$$xh+hx$$

Так как произведение некомутативно, то мы не можем представить приращение функции $f(x+h)-f(x)$ в виде $Ah$ или $hA$, где $A$ не зависит от $h$. Следствием этого является непредсказуемость поведения приращения функции $f(x)=x^2$, когда приращение аргумента стремится к 0. Однако, если бесконечно малая величина $h$ будет бесконечно малой величиной вида $h=ta$, $a\in D$, $t\in R$, $t\to 0$, то ответ становится более определённым

$$(xa+ax)t$$

□

Определение 5.2.3. Пусть $D$ - полное тело характеристики 0.[5.5] Отображение

$$f:D\to D$$

**дифференцируемо по Гато** на множестве $U\subset D$, если в каждой точке $x\in U$ изменение функции $f$ может быть представлено в виде

$$f(x+a)-f(x)=\partial f(x)(a)+o(a)=\frac{\partial f(x)}{\partial x}(a)+o(a) \tag{5.2.9}$$

где **производная Гато $\partial f(x)$ отображения** $f$ - линейное отображение приращения $a$ и $o:D\to D$ - такое непрерывное отображение, что

$$\lim_{a\to 0}\frac{|o(a)|}{|a|}=0$$

□

Замечание 5.2.4. Согласно определению 5.2.3 при заданном $x$ производная Гато $\partial f(x)\in\mathcal{L}(D;D)$. Следовательно, производная Гато отображения $f$ является отображением

$$\partial f:D\to\mathcal{L}(D;D)$$

Выражения $\partial f(x)$ и $\dfrac{\partial f(x)}{\partial x}$ являются разными обозначениями одной и той же функции. Мы будем пользоваться обозначением $\dfrac{\partial f(x)}{\partial x}$, если хотим подчеркнуть, что мы берём производную Гато по переменной $x$. □

Теорема 5.2.5. *Мы можем представить производную Гато отображения $f$ в виде*

$$\partial f(x)(a)=\frac{\partial_{s\cdot 0}f(x)}{\partial x}\ a\ \frac{\partial_{s\cdot 1}f(x)}{\partial x} \tag{5.2.10}$$

Доказательство. Следствие определения 5.2.3 и теоремы 3.1.11. □

[5.5] Определение дано согласно определению [1]-3.1.2, стр. 256.

Определение 5.2.6. Выражение $\dfrac{\partial_{s\cdot p} f(x)}{\partial x}$, $p=0,\ 1$, называется **компонентой производной Гато отображения** $f(x)$. □

Теорема 5.2.7. *Пусть $D$ - тело характеристики $0$. Пусть функция*

$$f:D\to D$$

*дифференцируема по Гато. Тогда*

$$\partial f(x)(ah)=a\partial f(x)(h) \tag{5.2.11}$$

*для любого действительного числа $a$.*

Доказательство. Следствие теоремы 5.1.15 и определения 5.2.3. □

Комбинируя равенство (5.2.11) и определение 5.2.3, мы получим знакомое определение производной Гато

$$\partial f(x)(a)=\lim_{t\to 0,\ t\in R}(t^{-1}(f(x+ta)-f(x))) \tag{5.2.12}$$

Определения производной Гато (5.2.9) и (5.2.12) эквивалентны. На основе этой эквивалентности мы будем говорить, что отображение $f$ дифференцируемо по Гато на множестве $U\subset D$, если в каждой точке $x\in U$ изменение функции $f$ может быть представлено в виде

$$f(x+ta)-f(x)=t\partial f(x)(a)+o(t) \tag{5.2.13}$$

где $o:R\to D$ - такое непрерывное отображение, что

$$\lim_{t\to 0}\frac{|o(t)|}{|t|}=0$$

Если бесконечно малая $a$ в равенстве (5.2.10) является дифференциалом $dx$, то равенство (5.2.10) становится определением **дифференциала Гато**

$$\partial f(x)(dx)=\frac{\partial_{s\cdot 0} f(x)}{\partial x}dx\frac{\partial_{s\cdot 1} f(x)}{\partial x} \tag{5.2.14}$$

Теорема 5.2.8. *Пусть $D$ - тело характеристики $0$. Пусть $\overline{\overline{e}}$ - базис тела $D$ над центром $Z(D)$ тела $D$.* **Стандартное представление производной Гато** (5.2.10) *отображения*

$$f:D\to D$$

*имеет вид*

$$\partial f(x)(a)=\frac{\partial^{ij} f(x)}{\partial x}\ e_i\ a\ e_j \tag{5.2.15}$$

*Выражение* $\dfrac{\partial^{ij} f(x)}{\partial x}$ *в равенстве* (5.2.15) *называется* **стандартной компонентой производной Гато** *отображения $f$.*

Доказательство. Утверждение теоремы является следствием теоремы 3.1.12. □

Теорема 5.2.9. *Пусть $D$ - тело характеристики $0$. Пусть $\overline{\overline{e}}$ - базис тела $D$ над полем $F\subset Z(D)$. Тогда дифференциал Гато отображения*

$$f:D\to D$$

*можно записать в виде*

$$\partial f(x)(a) = a^{i}\ \frac{\partial f^{j}}{\partial x^{i}} e_{j} \tag{5.2.16}$$

*где* $a \in D$ *имеет разложение*

$$a = a^{i}\ e_{i} \qquad a^{i} \in F$$

*относительно базиса* $\overline{\overline{e}}$ *и матрица Якоби отображения* $f$ *имеет вид*

$$\frac{\partial f^{j}}{\partial x^{i}} = \frac{\partial^{kr} f(x)}{\partial x}\ C^{p}_{ki}\ C^{j}_{pr} \tag{5.2.17}$$

Доказательство. Утверждение теоремы является следствием теоремы 3.1.13. □

Чтобы найти конструкцию, подобную производной в коммутативном случае, мы должны выделить производную из дифференциала отображения. Для этого мы должны вынести приращение аргумента за скобки.[5.6] Мы можем вынести $a$ за скобки, опираясь на равенства

$$\frac{\partial_{s\cdot 0} f(x)}{\partial x} a \frac{\partial_{s\cdot 1} f(x)}{\partial x} = \frac{\partial_{s\cdot 0} f(x)}{\partial x} a \frac{\partial_{s\cdot 1} f(x)}{\partial x} a^{-1} a$$
$$\frac{\partial_{s\cdot 0} f(x)}{\partial x} a \frac{\partial_{s\cdot 1} f(x)}{\partial x} = a a^{-1} \frac{\partial_{s\cdot 0} f(x)}{\partial x} a \frac{\partial_{s\cdot 1} f(x)}{\partial x}$$

Определение 5.2.10. Пусть $D$ - полное тело характеристики 0 и $a \in D$. **$D\star$-производная Гато** $\dfrac{\partial f(x)(a)}{\partial_{*} x}$ **отображения** $f : D \to D$ определена равенством

$$\partial f(x)(a) = a \frac{\partial f(x)(a)}{\partial_{*} x} \tag{5.2.18}$$

**$\star D$-производная Гато** $\dfrac{\partial f(x)(a)}{{}_{*}\partial x}$ **отображения** $f : D \to D$ определена равенством

$$\partial f(x)(a) = \frac{\partial f(x)(a)}{{}_{*}\partial x} a \tag{5.2.19}$$

□

Для представления $D\star$-производной Гато мы также будем пользоваться записью

$$\frac{\partial f(x)(a)}{\partial_{*} x} = \frac{\partial}{\partial_{*} x} f(x)(a) = \partial_{*} f(x)(a)$$

Для представления $\star D$-производной Гато мы также будем пользоваться записью

$$\frac{\partial f(x)(a)}{{}_{*}\partial x} = \frac{\partial}{{}_{*}\partial x} f(x)(a) = {}_{*}\partial f(x)(a)$$

Опираясь на принцип двойственности [3]-4.3.8, мы будем изучать $D\star$-производную, имея в виду, что двойственное утверждение справедливо для $\star D$-производной.

[5.6] Это возможно в том случае, если все $\dfrac{{}_{s\cdot 0}\partial f(x)}{\partial x} = e$ либо все $\dfrac{{}_{s\cdot 1}\partial f(x)}{\partial x} = e$. Поэтому определение [1]-3.1.2, стр. 256, приводит нас либо к определению $D\star$-производной Фреше, либо к определению $\star D$-производной Фреше.

Поскольку произведение не коммутативно, понятие дроби не определено в теле. Мы должны явно указать, с какой стороны знаменатель действует на числитель. Эту функцию выполняет символ $*$ в знаменателе дроби. Согласно определению

$$\frac{\partial f}{\partial_* x} = (\partial_* x)^{-1}\partial f$$
$$\frac{\partial f}{{}_*\partial x} = \partial f({}_*\partial x)^{-1}$$

Таким образом, мы можем представить $D\star$-производную Гато функции $f$ как отношение изменения функции к изменению аргумента. Это замечание не распространяется на компоненты дифференциала Гато. В этом случае мы рассматриваем запись $\frac{\partial}{\partial_* x}$ и $\frac{\partial}{{}_*\partial x}$ не как дроби, а как символы операторов.

При этом равенство (5.2.18) можно записать в виде

$$\partial f(x)(dx) = dx(\partial_* x)^{-1}\partial f(x)(dx) \tag{5.2.20}$$

Нетрудно видеть, что в знаменателе дроби в равенстве (5.2.20) мы записали $D\star$-дифференциал переменной $x$, а умножаем дробь на дифференциал переменной $x$.

Теорема 5.2.11. *Пусть $D$ - полное тело характеристики $0$. $D\star$-производная Гато проективна над полем действительных чисел $R$.*

Доказательство. Следствие теоремы 5.2.7 и примера 3.1.24. □

Из теоремы 5.2.11 следует

$$\frac{\partial f(x)(ra)}{\partial_* x} = \frac{\partial f(x)(a)}{\partial_* x} \tag{5.2.21}$$

для любых $r \in R$, $r \neq 0$ и $a \in D$, $a \neq 0$. Следовательно, $D\star$-производная Гато хорошо определена в направлении $a$ над полем $R$, $a \in D$, $a \neq 0$, и не зависит от выбора значения в этом направлении.

Теорема 5.2.12. *Пусть $D$ - полное тело характеристики $0$ и $a \neq 0$. $D\star$-производная Гато и $\star D$-производная Гато отображения $f$ тела $D$ связаны соотношением*

$$\frac{\partial f(x)(a)}{{}_*\partial x} = a\frac{\partial f(x)(a)}{\partial_* x}a^{-1} \tag{5.2.22}$$

Доказательство. Из равенств (5.2.18) и (5.2.19) следует

$$\frac{\partial f(x)(a)}{{}_*\partial x} = \partial f(x)(a)a^{-1} = a\frac{\partial f(x)(a)}{\partial_* x}a^{-1}$$

□

Теорема 5.2.13. *Пусть $D$ - полное тело характеристики $0$. Производная Гато удовлетворяет соотношению*

$$\partial(f(x)g(x))(a) = \partial f(x)(a)\ g(x) + f(x)\ \partial g(x)(a) \tag{5.2.23}$$

Доказательство. Равенство (5.2.23) следует из цепочки равенств

$$\begin{aligned}\partial(f(x)g(x))(a) &= \lim_{t\to 0}(t^{-1}(f(x+ta)g(x+ta)-f(x)g(x)))\\ &= \lim_{t\to 0}(t^{-1}(f(x+ta)g(x+ta)-f(x)g(x+ta)))\\ &+ \lim_{t\to 0}(t^{-1}(f(x)g(x+ta)-f(x)g(x)))\\ &= \lim_{t\to 0}(t^{-1}(f(x+ta)-f(x)))g(x)\\ &+ f(x)\lim_{t\to 0}(t^{-1}(g(x+ta)-g(x)))\end{aligned}$$

основанной на определении (5.2.12). □

Теорема 5.2.14. *Пусть $D$ - полное тело характеристики 0. Допустим производная Гато отображения $f:D\to D$ имеет разложение*

(5.2.24) $$\partial f(x)(a)=\frac{\partial_{s\cdot 0}f(x)}{\partial x}a\frac{\partial_{s\cdot 1}f(x)}{\partial x}$$

*Допустим проиводная Гато отображения $g:D\to D$ имеет разложение*

(5.2.25) $$\partial g(x)(a)=\frac{\partial_{t\cdot 0}g(x)}{\partial x}a\frac{\partial_{t\cdot 1}g(x)}{\partial x}$$

*Компоненты производной Гато отображения $f(x)g(x)$ имеют вид*

(5.2.26) $$\frac{\partial_{s\cdot 0}f(x)g(x)}{\partial x}=\frac{\partial_{s\cdot 0}f(x)}{\partial x}\qquad \frac{\partial_{t\cdot 0}f(x)g(x)}{\partial x}=f(x)\frac{\partial_{t\cdot 0}g(x)}{\partial x}$$

(5.2.27) $$\frac{\partial_{s\cdot 1}f(x)g(x)}{\partial x}=\frac{\partial_{s\cdot 1}f(x)}{\partial x}g(x)\qquad \frac{\partial_{t\cdot 1}f(x)g(x)}{\partial x}=\frac{\partial_{t\cdot 1}g(x)}{\partial x}$$

Доказательство. Подставим (5.2.24) и (5.2.25) в равенство (5.2.23)

(5.2.28) $$\begin{aligned}\partial(f(x)g(x))(a) &= \partial f(x)(a)\ g(x)+f(x)\ \partial g(x)(a)\\ &= \frac{\partial_{s\cdot 0}f(x)}{\partial x}a\frac{\partial_{s\cdot 1}f(x)}{\partial x}g(x)+f(x)\frac{\partial_{t\cdot 0}g(x)}{\partial x}a\frac{\partial_{t\cdot 1}g(x)}{\partial x}\end{aligned}$$

Опираясь на (5.2.28), мы определяем равенства (5.2.26), (5.2.27). □

Теорема 5.2.15. *Пусть $D$ - полное тело характеристики 0. $D\star$-производная Гато удовлетворяет соотношению*

(5.2.29) $$\frac{\partial f(x)g(x)}{\partial_* x}(a)=\frac{\partial f(x)(a)}{\partial_* x}g(x)+a^{-1}f(x)a\frac{\partial g(x)(a)}{\partial_* x}$$

Доказательство. Равенство (5.2.29) следует из цепочки равенств

$$\begin{aligned}\frac{\partial f(x)g(x)}{\partial_* x}(a) &= a^{-1}\partial f(x)g(x)(a)\\ &= a^{-1}(\partial f(x)(a)g(x)+f(x)\partial g(x)(a))\\ &= a^{-1}\partial f(x)(a)g(x)+a^{-1}f(x)aa^{-1}\partial g(x)(a)\\ &= \frac{\partial f(x)(a)}{\partial_* x}g(x)+a^{-1}f(x)a\frac{\partial g(x)(a)}{\partial_* x}\end{aligned}$$

□

Теорема 5.2.16. *Пусть $D$ - полное тело характеристики 0. Либо $D\star$-производная Гато не зависит от направления, либо $D\star$-производная Гато в направлении 0 не определена.*

Доказательство. Утверждение теоремы является следствием теоремы 5.2.11 и теоремы 5.1.16. □

Теорема 5.2.17. *Пусть $D$ - полное тело характеристики $0$. Пусть единичная сфера тела $D$ - компактна. Если $D\star$-производная Гато $\frac{\partial f(x)(a)}{\partial_* x}$ существует в точке $x$ и непрерывна по направлению над полем $R$, то существует норма $\|\partial f(x)\|$ дифференциала Гато.*

Доказательство. Из определения 5.2.10 следует

$$|\partial f(x)(a)| = |a| \left|\frac{\partial f(x)(a)}{\partial_* x}\right| \tag{5.2.30}$$

Из теорем 5.1.18, 5.2.11 следует, что $D\star$-производная Гато непрерывна на единичной сфере. Так как единичная сфера компактна, то множество значений $D\star$-производной Гато функции $f$ в точке $x$ ограниченно

$$\left|\frac{\partial f(x)(a)}{\partial_* x}\right| < F = \sup\left|\frac{\partial f(x)(a)}{\partial_* x}\right|$$

Согласно определению 5.1.13

$$\|\partial f(x)\| = F$$

□

Теорема 5.2.18. *Пусть $D$ - полное тело характеристики $0$. Пусть единичная сфера тела $D$ - компактна. Если $D\star$-производная Гато $\frac{\partial f(x)(a)}{\partial_* x}$ существует в точке $x$ и непрерывна по направлению над полем $R$, то отображение $f$ непрерывно в точке $x$.*

Доказательство. Из теоремы 5.2.17 следует

$$|\partial f(x)(a)| \le \|\partial f(x)\||a| \tag{5.2.31}$$

Из (5.2.9), (5.2.31) следует

$$|f(x+a) - f(x)| < |a|\ \|\partial f(x)\| \tag{5.2.32}$$

Возьмём произвольное $\epsilon > 0$. Положим

$$\delta = \frac{\epsilon}{\|\partial f(x)\|}$$

Тогда из неравенства

$$|a| < \delta$$

следует

$$|f(x+a) - f(x)| \le \|\partial f(x)\|\ \delta = \epsilon$$

Согласно определению 5.1.11 отображение $f$ непрерывно в точке $x$. □

Теорема 5.2.18 имеет интересное обобщение. Если единичная сфера тела $D$ не является компактной, то мы можем рассмотреть компактное множество направлений вместо единичной сферы. В этом случае, мы можем говорить о непрерывности функции $f$ вдоль заданного множества направлений.

## 5.3. Таблица производных Гато отображения тела

Теорема 5.3.1. *Пусть $D$ - полное тело характеристики 0. Тогда для любого $b\in D$*

$$\partial(b)(a)=0 \tag{5.3.1}$$

Доказательство. Непосредственное следствие определения 5.2.3. □

Теорема 5.3.2. *Пусть $D$ - полное тело характеристики 0. Тогда для любых $b$, $c\in D$*

$$\partial(bf(x)c)(a)=b\partial f(x)(a)c \tag{5.3.2}$$

$$\frac{\partial_{s\cdot 0}bf(x)c}{\partial x}=b\frac{\partial_{s\cdot 0}f(x)}{\partial x} \tag{5.3.3}$$

$$\frac{\partial_{s\cdot 1}bf(x)c}{\partial x}=\frac{\partial_{s\cdot 1}f(x)}{\partial x}c \tag{5.3.4}$$

$$\frac{\partial bf(x)c}{\partial_* x}(a)=a^{-1}ba\frac{\partial f(x)(a)}{\partial_* x}c \tag{5.3.5}$$

Доказательство. Непосредственное следствие равенств (5.2.23), (5.2.26), (5.2.27), (5.2.29), так как $\partial b=\partial c=0$. □

Теорема 5.3.3. *Пусть $D$ - полное тело характеристики 0. Тогда для любых $b$, $c\in D$*

$$\partial(bxc)(h)=bhc \tag{5.3.6}$$

$$\frac{\partial_{1\cdot 0}bxc}{\partial x}=b \tag{5.3.7}$$

$$\frac{\partial_{1\cdot 1}bxc}{\partial x}=c \tag{5.3.8}$$

$$\frac{\partial bxc}{\partial_* x}(h)=h^{-1}bhc \tag{5.3.9}$$

Доказательство. Следствие теоремы 5.3.2, когда $f(x)=x$. □

Теорема 5.3.4. *Пусть $D$ - полное тело характеристики 0. Тогда для любого $b\in D$*

$$\partial(xb-bx)(h)=hb-bh \tag{5.3.10}$$

$$\frac{\partial_{1\cdot 0}(xb-bx)}{\partial x}=1 \qquad \frac{\partial_{1\cdot 1}(xb-bx)}{\partial x}=b$$

$$\frac{\partial_{2\cdot 0}(xb-bx)}{\partial x}=-b \qquad \frac{\partial_{2\cdot 1}(xb-bx)}{\partial x}=1$$

$$\frac{\partial(xb-bx)}{\partial_* x}(h)=h^{-1}bhc$$

Доказательство. Следствие теоремы 5.3.2, когда $f(x)=x$. □

Теорема 5.3.5. *Пусть $D$ - полное тело характеристики 0. Тогда*[5.7]

$$\begin{aligned}\partial(x^2)(a) &= xa+ax\\ \frac{\partial x^2}{\partial_* x}(a) &= a^{-1}xa+x\end{aligned}\tag{5.3.11}$$

$$\begin{aligned}&\frac{\partial_{1\cdot 0}x^2}{\partial x}=x \quad \frac{\partial_{1\cdot 1}x^2}{\partial x}=e\\ &\frac{\partial_{2\cdot 0}x^2}{\partial x}=e \quad \frac{\partial_{2\cdot 1}x^2}{\partial x}=x\end{aligned}\tag{5.3.12}$$

Доказательство. (5.3.11) следует из примера 5.2.2 и определения 5.2.10. (5.3.12) следует из примера 5.2.2 и равенства (5.2.14). □

Теорема 5.3.6. *Пусть $D$ - полное тело характеристики 0. Тогда*[5.8]

$$\partial(x^{-1})(h)=-x^{-1}hx^{-1}\tag{5.3.13}$$

$$\frac{\partial x^{-1}}{\partial_* x}(h)=-h^{-1}x^{-1}hx^{-1}$$

$$\frac{\partial_{1\cdot 0}x^{-1}}{\partial x}=-x^{-1} \quad \frac{\partial_{1\cdot 1}x^{-1}}{\partial x}=x^{-1}$$

Доказательство. Подставим $f(x)=x^{-1}$ в определение (5.2.12).

$$\begin{aligned}\partial f(x)(h) &= \lim_{t\to 0,\ t\in R}(t^{-1}((x+th)^{-1}-x^{-1}))\\ &= \lim_{t\to 0,\ t\in R}(t^{-1}((x+th)^{-1}-x^{-1}(x+th)(x+th)^{-1}))\\ &= \lim_{t\to 0,\ t\in R}(t^{-1}(1-x^{-1}(x+th))(x+th)^{-1})\\ &= \lim_{t\to 0,\ t\in R}(t^{-1}(1-1-x^{-1}th)(x+th)^{-1})\\ &= \lim_{t\to 0,\ t\in R}(-x^{-1}h(x+th)^{-1})\end{aligned}\tag{5.3.14}$$

Равенство (5.3.13) следует из цепочки равенств (5.3.14). □

---

[5.7] Утверждение теоремы аналогично примеру VIII, [11], с. 451. Если произведение коммутативно, то равенство (5.3.11) принимает вид

$$\partial(x^2)(h)=2hx$$
$$\frac{dx^2}{dx}=2x$$

[5.8] Утверждение теоремы аналогично примеру IX, [11], с. 451. Если произведение коммутативно, то равенство (5.3.13) принимает вид

$$\partial(x^{-1})(h)=-hx^{-2}$$
$$\frac{dx^{-1}}{dx}=-x^{-2}$$

Теорема 5.3.7. *Пусть $D$ - полное тело характеристики 0. Тогда*[5.9]

$$\partial(xax^{-1})(h) = hax^{-1} - xax^{-1}hx^{-1} \tag{5.3.15}$$

$$\frac{\partial xax^{-1}}{\partial_* x}(h) = ax^{-1} - h^{-1}xax^{-1}hx^{-1}$$

$$\frac{\partial_{1\cdot 0}x^{-1}}{\partial x} = 1 \qquad \frac{\partial_{1\cdot 1}x^{-1}}{\partial x} = ax^{-1}$$

$$\frac{\partial_{2\cdot 0}x^{-1}}{\partial x} = -xax^{-1} \qquad \frac{\partial_{2\cdot 1}x^{-1}}{\partial x} = x^{-1}$$

Доказательство. Равенство (5.3.15) является следствием равенств (5.2.23), (5.3.6), (5.3.15). □

---

[5.9]Если произведение коммутативно, то

$$y = xax^{-1} = a$$

Соответственно производная обращается в 0.

Глава 6

# Дифференцируемые отображения $D$-векторного пространства

## 6.1. Топологическое $D$-векторное пространство

Определение 6.1.1. Пусть $D$ - топологическое тело. $D$-векторное пространство $V$ называется **топологическим $D$-векторным пространством**[6.1], если $V$ наделено топологией, согласующейся со структурой аддитивной группы в $V$, и отображения

$$(a, v) \in D \times V \to av \in V$$
$$(v, a) \in V \times D \to va \in V$$

непрерывны. □

Определение 6.1.2. Отображение

$$f : F \to A$$

множества $F$ в произвольную алгебру $A$ называется **$A$-значной функцией**. Отображение

$$f : F \to V$$

множества $F$ в $D$-векторное пространство $V$ называется **$D$-вектор-функцией**. □

Мы рассматриваем последующие определения в этом разделе для топологического $D$-векторного пространства. Однако определения не изменятся, если мы будем рассматривать топологическое $D^*{}_*$-векторное пространство.

Определение 6.1.3. Пусть $V$ - $D$-векторное пространство размерности $n$. Пусть $\overline{\overline{e}}$ - базис $D^*{}_*$-векторного пространства $V$. Произвольное отображение

$$f : V \to A$$

$D$-векторного пространства $V$ в множество $A$ можно представить в виде **отображения**

$$f' : D^n \to A$$

$n$ **$D$-значных переменных**; функция $f'$ определена равенством

$$f'(a^1, ..., a^n) = f(a^*{}_*e)$$

□

[6.1] Определение дано согласно определению из [5], с. 21

Определение 6.1.4. Отображение

$$f : V \to A$$

топологического $D$-векторного пространства $V$ размерности $n$ в топологическое пространство $A$ называется **непрерывным по совокупности аргументов**, если для произвольной окрестности $U$ образа элемента $\overline{a} = a^*{}_*e \in V$ для каждого $a_i \in D$ существуют такие окрестности $V_i$, что

$$f'(V_1, ..., V_n) \subset U$$

$\square$

Теорема 6.1.5. *Непрерывное отображение*

$$f : V \to A$$

*топологического $D$-векторного пространства $V$ размерности $n$ в топологическое пространство $A$ непрерывна по совокупности аргументов.*

Доказательство. Пусть $b = f(a^*{}_*e) \in A$ и $U$ окрестность точки $b$. Так как отображение $f$ непрерывно, то существует такая окрестность $V$ вектора $\overline{a} = a^*{}_*e$ , что $f(V) \subseteq U$. Так как аддитивная операция непрерывна, то существуют такие окрестности $E_i$ вектора $e_i$ и окрестности $W_i$ элемента $a_i \in D$, что $W_*{}^*E \subset V$. $\square$

Определение 6.1.6. **Норма на $D$-векторном пространстве** $\overline{V}$ над недискретным нормированным телом $D$[6.2] - это отображение

$$v \in V \to \|v\| \in R$$

такое, что

- $\|v\| \geq 0$
- $\|v\| = 0$ равносильно $v = 0$
- $\|v + w\| \leq \|v\| + \|w\|$
- $\|av\| = |a|\|v\|$ для всех $a \in D$ и $v \in V$

$D$-векторное пространство $V$ над недискретным нормированным телом $D$, наделённое структурой, определяемой заданием на $V$ нормы, называется **нормированным $D$-векторным пространством**. $\square$

Инвариантное расстояние на аддитивной группе $D$-векторного пространства $V$

$$d(a, b) = \|a - b\|$$

определяет в $V$ топологию метрического пространства, согласующуюся со структурой $D$-векторного пространства в $V$.

Определение 6.1.7. Пусть $D$ - полное тело характеристики 0. Пусть

$$f : V \to W$$

отображение нормированного $D$-векторного пространства $V$ с нормой $\|x\|_1$ в нормированное $D$-векторное пространство $W$ с нормой $\|y\|_2$. Величина

$$\|f\| = \sup \frac{\|f(x)\|_2}{\|x\|_1}$$

называется **нормой отображения** $f$. $\square$

[6.2] Определение дано согласно определению из [4], гл. IX, §3, п°3

## 6.2. Дифференцируемые отображения $D$-векторного пространства

Определение 6.2.1. Отображение

$$f: V \to W$$

нормированного $D$-векторного пространства $V$ с нормой $\|x\|_1$ в нормированное $D$-векторное пространство $W$ с нормой $\|x\|_2$[6.3] называется **дифференцируемой по Гато** на множестве $U \subset V$, если в каждой точке $x \in U$ изменение отображения $f$ может быть представлено в виде

$$f(x+a) - f(x) = \partial f(x)(a) + o(a) \tag{6.2.1}$$

где **производная Гато** $\partial f(x)$ **отображения** $f$ - линейное отображение приращения $a$ и $o: V \to W$ - такое непрерывное отображение, что

$$\lim_{a \to 0} \frac{\|o(a)\|_2}{\|a\|_1} = 0$$

□

Замечание 6.2.2. Согласно определению 6.2.1 при заданном $x$ производная Гато $\partial f(x) \in \mathcal{L}(D; V; W)$. Следовательно, производная Гато отображения $f$ является отображением

$$\partial f: V \to \mathcal{L}(D; V; W)$$

□

Теорема 6.2.3. *Пусть $D$ - тело характеристики 0. Пусть $V$, $W$ - $D$-векторные пространства. Допустим $\overline{\overline{e}}_V$ - базис в $D^*{}_*$-векторном пространстве $V$ и $h \in V$*

$$\overline{h} = h^*{}_* e_V$$

*Допустим $\overline{\overline{e}}_W$ - базис в $D^*{}_*$-векторном пространстве $W$. Мы можем представить производную Гато отображения $f$ в виде*

$$\partial f(x)(h) = \partial_{s\cdot 0\cdot i} f^j(x)\ h^i\ \partial_{s\cdot 1\cdot i} f^j(x)\ e_{W\cdot j} = \frac{\partial_{s\cdot 0} f^j(x)}{\partial x^i}\ h^i\ \frac{\partial_{s\cdot 1} f^j(x)}{\partial x^i}\ e_{W\cdot j} \tag{6.2.2}$$

Доказательство. Следствие определения 6.2.1 и теоремы 4.1.4. □

Определение 6.2.4. Выражение $\partial_{s\cdot p\cdot i} f^j(x) = \dfrac{\partial_{s\cdot p} f^j(x)}{\partial x^i}$, $p = 0,\ 1$, называется **компонентой производной Гато** отображения $f(x)$. □

Теорема 6.2.5. *Пусть полное тело $D$ характеристики 0 является алгеброй над полем $F \subset Z(D)$. Производная Гато отображения*

$$f: V \to W$$

*нормированного $D$-векторного пространства $V$ в нормированное $D$-векторное пространство $W$ удовлетворяет равенству*

$$\partial f(x)(ph) = p\partial f(x)(h) \quad p \in F$$

[6.3] Если размерность $D$-векторного пространства $V$ равна 1, то мы можем отождествить его с телом $D$. Тогда мы будем говорить об отображении тела $D$ в $D$-векторное пространство $W$. Если размерность $D$-векторного пространства $W$ равна 1, то мы можем отождествить его с телом $D$. Тогда мы будем говорить об отображении $D$-векторного пространства $V$ в тело $D$. В обоих случаях мы будем опускать соответствующие индексы, так как эти индексы имеют единственное значение.

Доказательство. Следствие определений 4.1.1, 6.2.1. □

Теорема 6.2.6. *Пусть $D$ - полное тело характеристики 0. Производная Гато функции*

$$f:V\to W$$

*нормированного $D$-векторного пространства $V$ в нормированное $D$-векторное пространство $W$ удовлетворяет равенству*

$$\partial f(x)(ra)=r\partial f(x)(a) \tag{6.2.3}$$

*для любых $r\in R$, $r\neq 0$ и $a\in V$, $a\neq 0$.*

Доказательство. Следствие теорем 5.1.9, 6.2.5. □

Комбинируя равенство (6.2.3) и определение 6.2.1, мы получим знакомое определение производной Гато

$$\partial f(x)(a)=\lim_{t\to 0}(t^{-1}(f(x+ta)-f(x))) \tag{6.2.4}$$

Определения производной Гато (6.2.1) и (6.2.4) эквивалентны. На основе этой эквивалентности мы будем говорить, что отображение $f$ дифференцируемо по Гато на множестве $U\subset D$, если в каждой точке $x\in U$ изменение функции $f$ может быть представлено в виде

$$f(x+ta)-f(x)=t\partial f(x)(a)+o(t) \tag{6.2.5}$$

где $o:R\to W$ - такое непрерывное отображение, что

$$\lim_{t\to 0}\frac{\|o(t)\|_2}{|t|}=0$$

Если бесконечно малая $a$ является дифференциалом $dx$, то равенство (6.2.2) становится определением **дифференциала Гато**

$$\partial f(x)(dx)=\frac{\partial_{s\cdot 0}f^{j}(x)}{\partial x^{i}}\,dx^{i}\,\frac{\partial_{s\cdot 1}f^{j}(x)}{\partial x^{i}}\,e_{W\cdot j} \tag{6.2.6}$$

Определение 6.2.7. Пусть $D$ - тело характеристики 0. Пусть $V$, $W$ - $D$-векторные пространства. Допустим $\overline{\overline{e}}_V$ - базис в $D^*{}_*$-векторном пространстве $V$ и $h\in V$

$$\overline{h}=h^*{}_*e_V$$

Допустим $\overline{\overline{e}}_W$ - базис в $D^*{}_*$-векторном пространстве $W$. Если мы рассмотрим производную Гато отображения $f$ по переменной $v^{i}$ при условии, что остальные координаты вектора $v$ постояны, то соответствующий аддитивный оператор

$$\frac{\partial f^{j}(v)}{\partial v^{i}}(h^{i})=\lim_{t\to 0}(t^{-1}(f^{j}(v^{1},...,v^{i}+th^{i},...,v^{n})-f^{j}(v))) \tag{6.2.7}$$

называется **частной производной Гато** отображения $f^{j}$ по переменной $v^{i}$. □

Теорема 6.2.8. *Пусть $D$ - тело характеристики 0. Допустим $\overline{\overline{e}}_V$ - базис в $D$-векторном пространстве $V$ и $h\in V$*

$$\overline{h}=h^*{}_*e_V$$

*Допустим $\overline{\overline{e}}_W$ - базис в $D$-векторном пространстве $W$. Пусть частные производные Гато отображения*

$$f:V\to W$$

*непрерывны в области $U \in V$. Тогда на множестве $U$ производная Гато отображения $f$ и частные производные Гато связаны соотношением*

$$\partial f(x)(h) = \left(\frac{\partial f^j(x)}{\partial x^i}(h^i)\right) e_{W\cdot j} \tag{6.2.8}$$

Доказательство. Из теоремы 4.1.2 следует, что мы можем разложить производную Гато отображения $f$ в сумму (6.2.8) частных линейных отображений. Наша задача: выяснить, при каком условии эти частные линейные отображения становятся частными производными Гато.

Согласно определению 6.2.4,

$$\begin{aligned}\partial f(x)(h) &= \lim_{t\to 0}(t^{-1}(f(x+th)-f(x)))\\ &= \lim_{t\to 0}(t^{-1}(f^j(x^1+th^1, x^2+th^2, ..., x^n+th^n)\\ &-f^j(x^1, x^2+th^2, ..., x^n+th^n)\\ &+f^j(x^1, x^2+th^2, ..., x^n+th^n) - ...\\ &-f^j(x^1, ..., x^n)))e_{W\cdot j}\end{aligned} \tag{6.2.9}$$

Из равенств (6.2.7) и (6.2.9) с точностью до бесконечно малой по сравнению с $t$ следует

$$\begin{aligned}\partial f(x)(h) = \lim_{t\to 0}\Bigg(&\frac{\partial f^j(x^1, x^2+th^2, ..., x^n+th^n)}{\partial x^1}(h^1)\\ &+... + \frac{\partial f^j(x^1, ..., x^n)}{\partial x^n}(h^n)\Bigg) e_{W\cdot j}\end{aligned} \tag{6.2.10}$$

Следовательно, для того, чтобы из равенства (6.2.10) следовало равенство (6.2.8), необходимо, чтобы частные производные Гато были непрерывны. □

Определение 6.2.9. Пусть отображение

$$f : V \to W$$

нормированного $D$-векторного пространства $V$ в нормированное $D$-векторное пространство $W$ дифференцируема по Гато на множестве $U \times W$. Допустим мы выбрали базис $\overline{\overline{e}}_V$ в $D^*{}_*$-векторном пространстве $V$. Допустим мы выбрали базис $\overline{\overline{e}}_W$ в $D^*{}_*$-векторном пространстве $W$. **$D^*{}_*$-производная Гато** $\frac{\partial f(x)(a)}{(\partial^*{}_*)x}$ **отображения** $f$ определена равенством

$$\partial f(x)(a) = a^*{}_*\frac{\partial f(x)(a)}{(\partial^*{}_*)x} = a^i\,\frac{\partial f^j(x)(a)}{(\partial^*{}_*)x^i}\, e_{W\cdot j} \tag{6.2.11}$$

Выражение $\frac{\partial f^j(x)(a)}{(\partial^*{}_*)x^i}$ в равенстве (6.2.11) называется **частными $D^*{}_*$-производными Гато отображения** $f^j$ **по переменной** $x^i$. □

Для представления $D^*{}_*$-производной Гато мы также будем пользоваться записью

$$\begin{aligned}\frac{\partial f(x)(a)}{(\partial^*{}_*)x} &= \frac{\partial}{(\partial^*{}_*)x}f(x)(a) &= (\partial^*{}_*)f(x)(a)\\ \frac{\partial f^j(x)(a^i)}{(\partial^*{}_*)x^i} &= \frac{\partial}{(\partial^*{}_*)x^i}f^j(x)(a) &= (\partial_i{}^*{}_*)f^j(x)(a)\end{aligned}$$

Определение 6.2.10. Пусть отображение

$$f:V\to W$$

нормированного $D$-векторного пространства $V$ в нормированное $D$-векторное пространство $W$ дифференцируема по Гато на множестве $U\times W$. Допустим мы выбрали базис $\overline{\overline{e}}_V$ в ${}_*{}^*D$-векторном пространстве $V$. Допустим мы выбрали базис $\overline{\overline{e}}_W$ в ${}_*{}^*D$-векторном пространстве $W$. **${}_*{}^*D$-производная Гато** $\frac{\partial f(x)(a)}{({}_*{}^*\partial)x}$ **отображения** $f$ определена равенством

$$\partial f(x)(a)=\frac{\partial f(x)(a)}{({}_*{}^*\partial)x}{}_*{}^*\overline{a}=e_{W\cdot j}\ \frac{\partial f^j(x)(a)}{({}_*{}^*\partial)x^i}\ a^i \tag{6.2.12}$$

Выражение $\frac{\partial f^j(x)(a)}{({}_*{}^*\partial)x^i}$ в равенстве (6.2.12) называется **частными ${}_*{}^*D$-производными Гато отображения $f^j$ по переменной $x^i$.** □

Для представления ${}_*{}^*D$-производной Гато мы также будем пользоваться записью

$$\begin{aligned}\frac{\partial f(x)(a)}{({}_*{}^*\partial)x}&=\frac{\partial}{({}_*{}^*\partial)x}f(x)(a)&&=({}_*{}^*\partial)f(x)(a)\\ \frac{\partial f^j(x)(a)}{({}_*{}^*\partial)x^i}&=\frac{\partial}{({}_*{}^*\partial)x^i}f^j(x)(a)&&=({}_*{}^*\partial_i)\ f^j(x)(a)\end{aligned}$$

Теорема 6.2.11. *Пусть $D$ - тело характеристики 0. Пусть функция*

$$f:U\to V$$

*нормированного $D$-векторного пространства $U$ в нормированное $D$-векторное пространство $V$ дифференцируема по Гато в точке $x$. Тогда*

$$\partial f(x)(0)=0$$

Доказательство. Следствие определения 6.2.1 и теоремы 4.1.8. □

Теорема 6.2.12. *Пусть $D$ - тело характеристики 0. Пусть $U$ - нормированное $D$-векторное пространство с нормой $\|x\|_1$. Пусть $V$ - нормированное $D$-векторное пространство с нормой $\|y\|_2$. Пусть $W$ - нормированное $D$-векторное пространство с нормой $\|z\|_3$. Пусть отображение*

$$f:U\to V$$

*дифференцируемо по Гато в точке $x$ и норма производной Гато отображения $f$ конечна*

$$\|\partial f(x)\|=F\le\infty \tag{6.2.13}$$

*Пусть отображение*

$$g:V\to W$$

*дифференцируемо по Гато в точке*

$$y=f(x) \tag{6.2.14}$$

*и норма производной Гато отображения $g$ конечна*

$$\|\partial g(y)\|=G\le\infty \tag{6.2.15}$$

*Функция*

$$gf(x)=g(f(x)) \tag{6.2.16}$$

*дифференцируема по Гато в точке $x$*

$$\partial gf(x)(a)=\partial g(f(x))(\partial f(x)(a)) \tag{6.2.17}$$

$$\frac{\partial gf(x)(a)}{(\partial^*{}_*)x}=\frac{\partial f(x)(a)}{(\partial^*{}_*)x}{}^*{}_*\frac{\partial g(f(x))(\partial f(x)(a))}{(\partial^*{}_*)f(x)} \tag{6.2.18}$$

$$\frac{\partial (gf)^{k}(x)(a)}{(\partial^*{}_*)x^{i}}=\frac{\partial f^{j}(x)(a)}{(\partial^*{}_*)x^{i}}\frac{\partial g^{k}(f(x))(\partial f(x)(a))}{(\partial^*{}_*)f^{j}(x)} \tag{6.2.19}$$

$$\frac{\partial_{st\cdot 0}(gf)^{k}(x)}{\partial x^{i}}=\frac{\partial_{s\cdot 0}g^{k}(f(x))}{\partial f^{j}(x)}\frac{\partial_{t\cdot 0}f^{j}(x)}{\partial x^{i}} \tag{6.2.20}$$

$$\frac{\partial_{st\cdot 1}(gf)^{k}(x)}{\partial x^{i}}=\frac{\partial_{t\cdot 1}f^{j}(x)}{\partial x^{i}}\frac{\partial_{s\cdot 1}g^{k}(f(x))}{\partial f^{j}(x)} \tag{6.2.21}$$

Доказательство. Согласно определению 6.2.1

$$g(y+b)-g(y)=\partial g(y)(b)+o_1(b) \tag{6.2.22}$$

где $o_1:V\rightarrow W$ - такое непрерывное отображение, что

$$\lim_{b\rightarrow 0}\frac{\|o_1(b)\|_3}{\|b\|_2}=0$$

Согласно определению 6.2.1

$$f(x+a)-f(x)=\partial f(x)(a)+o_2(a) \tag{6.2.23}$$

где $o_2:U\rightarrow V$ - такое непрерывное отображение, что

$$\lim_{a\rightarrow 0}\frac{\|o_2(a)\|_2}{\|a\|_1}=0$$

Согласно (6.2.23) смещение вектора $x\in U$ в направлении $a$ приводит к смещению вектора $y$ в направлении

$$b=\partial f(x)(a)+o_2(a) \tag{6.2.24}$$

Используя (6.2.14), (6.2.24) в равенстве (6.2.22), мы получим

$$\begin{aligned}&g(f(x+a))-g(f(x))\\ =&g(f(x)+\partial f(x)(a)+o_2(a))-g(f(x))\\ =&\partial g(f(x))(\partial f(x)(a)+o_2(a))-o_1(\partial f(x)(a)+o_2(a))\end{aligned} \tag{6.2.25}$$

Согласно определениям 6.2.1, 4.1.1 из равенства (6.2.25) следует

$$\begin{aligned}&g(f(x+a))-g(f(x))\\ =&\partial g(f(x))(\partial f(x)(a))+\partial g(f(x))(o_2(a))-o_1(\partial f(x)(a)+o_2(a))\end{aligned} \tag{6.2.26}$$

Найдём предел

$$\lim_{a\rightarrow 0}\frac{\|\partial g(f(x))(o_2(a))-o_1(\partial f(x)(a)+o_2(a))\|_3}{\|a\|_1}$$

Согласно определению 6.1.6

$$\begin{aligned}(6.2.27)\quad &\lim_{a\to 0}\frac{\|\partial g(f(x))(o_2(a))-o_1(\partial f(x)(a)+o_2(a))\|_3}{\|a\|_1}\\ &\le\lim_{a\to 0}\frac{\|\partial g(f(x))(o_2(a))\|_3}{\|a\|_1}+\lim_{a\to 0}\frac{\|o_1(\partial f(x)(a)+o_2(a))\|_3}{\|a\|_1}\end{aligned}$$

Из (6.2.15) следует

$$(6.2.28)\quad \lim_{a\to 0}\frac{\|\partial g(f(x))(o_2(a))\|_3}{\|a\|_1}\le G\lim_{a\to 0}\frac{\|o_2(a)\|_2}{\|a\|_1}=0$$

Из (6.2.13) следует

$$\begin{aligned}&\lim_{a\to 0}\frac{\|o_1(\partial f(x)(a)+o_2(a))\|_3}{\|a\|_1}\\ &=\lim_{a\to 0}\frac{\|o_1(\partial f(x)(a)+o_2(a))\|_3}{\|\partial f(x)(a)+o_2(a)\|_2}\lim_{a\to 0}\frac{\|\partial f(x)(a)+o_2(a)\|_2}{\|a\|_1}\\ &\le\lim_{a\to 0}\frac{\|o_1(\partial f(x)(a)+o_2(a))\|_3}{\|\partial f(x)(a)+o_2(a)\|_2}\lim_{a\to 0}\frac{\|\partial f(x)\|\|a\|_1+\|o_2(a)\|_2}{\|a\|_1}\\ &=\lim_{a\to 0}\frac{\|o_1(\partial f(x)(a)+o_2(a))\|_3}{\|\partial f(x)(a)+o_2(a)\|_2}\|\partial f(x)\|\end{aligned}$$

Согласно теореме 6.2.11

$$\lim_{a\to 0}(\partial f(x)(a)+o_2(a))=0$$

Следовательно,

$$(6.2.29)\quad \lim_{a\to 0}\frac{\|o_1(\partial f(x)(a)+o_2(a))\|_3}{\|a\|_1}=0$$

Из равенств (6.2.27), (6.2.28), (6.2.29) следует

$$(6.2.30)\quad \lim_{a\to 0}\frac{\|\partial g(f(x))(o_2(a))-o_1(\partial f(x)(a)+o_2(a))\|_3}{\|a\|_1}=0$$

Согласно определению 6.2.1

$$(6.2.31)\quad gf(x+a)-gf(x)=\partial gf(x)(a)+o(a)$$

где $o:U\to W$ - такое непрерывное отображение, что

$$\lim_{a\to 0}\frac{\|o(a)\|_3}{\|a\|_1}=0$$

Равенство (6.2.17) следует из (6.2.26), (6.2.30), (6.2.31).

Из равенства (6.2.17) и определения 6.2.9 следует

$$\begin{aligned}(6.2.32)\quad a^*{}_*\frac{\partial gf(x)(a)}{(\partial^*{}_*)x}&=\partial f(x)(a)^*{}_*\frac{\partial g(f(x))(\partial f(x)(a)}{(\partial^*{}_*)f(x)}\\ &=a^*{}_*\frac{\partial f(x)(a)}{(\partial^*{}_*)x}{}^*{}_*\frac{\partial g(f(x))(\partial f(x)(a)}{(\partial^*{}_*)f(x)}\end{aligned}$$

Так как приращение $a$ произвольно, то равенство (6.2.18) следует из равенства (6.2.32).

Из равенства (6.2.17) и теоремы 6.2.3 следует

$$
\begin{aligned}
&\frac{\partial_{st\cdot 0} gf^{\boldsymbol{k}}(x)}{\partial x^{\boldsymbol{i}}}\ a^{\boldsymbol{i}}\ \frac{\partial_{st\cdot 1} gf^{\boldsymbol{k}}(x)}{\partial x^{\boldsymbol{i}}}\ e_{W\cdot\boldsymbol{k}}\\
&=\frac{\partial_{s\cdot 0} g^{\boldsymbol{k}}(f(x))}{\partial f^{\boldsymbol{j}}(x)}\ \partial f^{\boldsymbol{j}}(x)(a)\ \frac{\partial_{s\cdot 1} g^{\boldsymbol{k}}(f(x))}{\partial f^{\boldsymbol{j}}(x)}\ e_{W\cdot\boldsymbol{k}}\\
&=\frac{\partial_{s\cdot 0} g^{\boldsymbol{k}}(f(x))}{\partial f^{\boldsymbol{j}}(x)}\frac{\partial_{t\cdot 0} f^{\boldsymbol{j}}(x)}{\partial x^{\boldsymbol{i}}}\ a^{\boldsymbol{i}}\ \frac{\partial_{t\cdot 1} f^{\boldsymbol{j}}(x)}{\partial x^{\boldsymbol{i}}}\ \frac{\partial_{s\cdot 1} g^{\boldsymbol{k}}(f(x))}{\partial f^{\boldsymbol{j}}(x)}\ e_{W\cdot\boldsymbol{k}}
\end{aligned}
\tag{6.2.33}
$$

(6.2.20), (6.2.21) следуют из равенства (6.2.33). □

Глава 7

# Алгебра кватернионов

## 7.1. Линейная функция комплексного поля

Теорема 7.1.1 (Уравнения Коши-Римана). *Рассмотрим поле комплексных чисел $C$ как двумерную алгебру над полем действительных чисел. Положим*

$$e_{C\cdot 0} = 1 \quad e_{C\cdot 1} = i \tag{7.1.1}$$

*базис алгебры $C$. Тогда в этом базисе произведение имеет вид*

$$e_{C\cdot 1}^2 = -e_{C\cdot 0} \tag{7.1.2}$$

*и структурные константы имеют вид*

$$\begin{aligned} C_{C\cdot}{}^{0}_{00} &= 1 & C_{C\cdot}{}^{1}_{01} &= \ \ 1 \\ C_{C\cdot}{}^{1}_{10} &= 1 & C_{C\cdot}{}^{0}_{11} &= -1 \end{aligned} \tag{7.1.3}$$

*Матрица линейной функции*

$$y^i = x^j f^i_j$$

*поля комплексных чисел над полем действительных чисел удовлетворяет соотношению*

$$f^0_0 = f^1_1 \tag{7.1.4}$$

$$f^1_0 = -f^0_1 \tag{7.1.5}$$

Доказательство. Равенства (7.1.2) и (7.1.3) следуют из равенства $i^2 = -1$. Пользуясь равенством (3.1.17) получаем соотношения

$$f^0_0 = f^{kr} C_{C\cdot}{}^{p}_{k0} C_{C\cdot}{}^{0}_{pr} = f^{0r} C_{C\cdot}{}^{0}_{00} C_{C\cdot}{}^{0}_{0r} + f^{1r} C_{C\cdot}{}^{1}_{10} C_{C\cdot}{}^{0}_{1r} = f^{00} - f^{11} \tag{7.1.6}$$

$$f^1_0 = f^{kr} C_{C\cdot}{}^{p}_{k0} C_{C\cdot}{}^{1}_{pr} = f^{0r} C_{C\cdot}{}^{0}_{00} C_{C\cdot}{}^{1}_{0r} + f^{1r} C_{C\cdot}{}^{1}_{10} C_{C\cdot}{}^{1}_{1r} = f^{01} + f^{10} \tag{7.1.7}$$

$$f^0_1 = f^{kr} C_{C\cdot}{}^{p}_{k1} C_{C\cdot}{}^{0}_{pr} = f^{0r} C_{C\cdot}{}^{1}_{01} C_{C\cdot}{}^{0}_{1r} + f^{1r} C_{C\cdot}{}^{0}_{11} C_{C\cdot}{}^{0}_{0r} = -f^{01} - f^{10} \tag{7.1.8}$$

$$f^1_1 = f^{kr} C_{C\cdot}{}^{p}_{k1} C_{C\cdot}{}^{1}_{pr} = f^{0r} C_{C\cdot}{}^{1}_{01} C_{C\cdot}{}^{1}_{1r} + f^{1r} C_{C\cdot}{}^{0}_{11} C_{C\cdot}{}^{1}_{0r} = f^{00} - f^{11} \tag{7.1.9}$$

Из равенств (7.1.6) и (7.1.9) следует (7.1.4). Из равенств (7.1.7) и (7.1.8) следует (7.1.5). □

Теорема 7.1.2 (Уравнения Коши-Римана). *Если матрица*

$$\begin{pmatrix} \dfrac{\partial y^0}{\partial x^0} & \dfrac{\partial y^0}{\partial x^1} \\ \dfrac{\partial y^1}{\partial x^0} & \dfrac{\partial y^1}{\partial x^1} \end{pmatrix}$$

*является матрицей Якоби функции комплексного переменного*

$$x = x^0 + x^1 i \to y = y^0(x^0, x^1) + y^1(x^0, x^1) i$$

*над полем действительных чисел, то*

$$\begin{aligned}\frac{\partial y^1}{\partial x^0} &= -\frac{\partial y^0}{\partial x^1}\\ \frac{\partial y^0}{\partial x^0} &= \frac{\partial y^1}{\partial x^1}\end{aligned} \tag{7.1.10}$$

Доказательство. Следствие теоремы 7.1.1. □

Теорема 7.1.3. *Производная функции комплексного переменного удовлетворяет равенству*

$$\frac{\partial y}{\partial x^0} + i\frac{\partial y}{\partial x^1} = 0 \tag{7.1.11}$$

Доказательство. Равенство

$$\frac{\partial y^0}{\partial x^0} + i\frac{\partial y^1}{\partial x^0} + i\frac{\partial y^0}{\partial x^1} - \frac{\partial y^1}{\partial x^1} = 0$$

следует из равенств (7.1.10). □

Равенство (7.1.11) эквивалентно равенству

$$\begin{pmatrix}1 & i\end{pmatrix}\begin{pmatrix}\frac{\partial y^0}{\partial x^0} & \frac{\partial y^0}{\partial x^1}\\ \frac{\partial y^1}{\partial x^0} & \frac{\partial y^1}{\partial x^1}\end{pmatrix}\begin{pmatrix}1\\ i\end{pmatrix} = 0 \tag{7.1.12}$$

## 7.2. Алгебра кватернионов

В этой статье я рассматриваю множество кватернионных алгебр, определённых в [8].

Определение 7.2.1. Пусть $F$ - поле. Расширение $F(i,j,k)$ поля $F$ называется **алгеброй** $E(F,a,b)$ **кватернионов над полем** $F$[7.1], если произведение в алгебре $E$ определено согласно правилам

(7.2.1)

| | $i$ | $j$ | $k$ |
|---|---|---|---|
| $i$ | $a$ | $k$ | $aj$ |
| $j$ | $-k$ | $b$ | $-bi$ |
| $k$ | $-aj$ | $bi$ | $-ab$ |

где $a$, $b \in F$, $ab \neq 0$. □

Элементы алгебры $E(F,a,b)$ имеют вид

$$x = x^0 + x^1 i + x^2 j + x^3 k$$

где $x^i \in F$, $i = 0, 1, 2, 3$. Кватернион

$$\overline{x} = x^0 - x^1 i - x^2 j - x^3 k$$

[7.1] Я буду следовать определению из [8].

называется сопряжённым кватерниону $x$. Мы определим **норму кватерниона** $x$ равенством

$$|x|^2 = x\overline{x} = (x^0)^2 - a(x^1)^2 - b(x^2)^2 + ab(x^3)^2 \tag{7.2.2}$$

Из равенства (7.2.2) следует, что $E(F,a,b)$ является алгеброй с делением только когда $a < 0$, $b < 0$. Тогда мы можем пронормировать базис так, что $a = -1$, $b = -1$.

Мы будем обозначать символом $E(F)$ алгебру $E(F,-1,-1)$ кватернионов с делением над полем $F$. Мы будем полагать $H = E(R,-1,-1)$. Произведение в алгебре кватернионов $E(F)$ определено согласно правилам

(7.2.3)

| | $i$ | $j$ | $k$ |
|---|---|---|---|
| $i$ | $-1$ | $k$ | $-j$ |
| $j$ | $-k$ | $-1$ | $i$ |
| $k$ | $j$ | $-i$ | $-1$ |

В алгебре $E(F)$ норма кватерниона имеет вид

$$|x|^2 = x\overline{x} = (x^0)^2 + (x^1)^2 + (x^2)^2 + (x^3)^2 \tag{7.2.4}$$

При этом обратный элемент имеет вид

$$x^{-1} = |x|^{-2}\overline{x} \tag{7.2.5}$$

Внутренний автоморфизм алгебры кватернионов $H$[7.2]

$$\begin{gathered} p \to qpq^{-1} \\ q(ix + jy + kz)q^{-1} = ix' + jy' + kz' \end{gathered} \tag{7.2.6}$$

описывает вращение вектора с координатами $x$, $y$, $z$. Если $q$ записан в виде суммы скаляра и вектора

$$q = \cos\alpha + (ia + jb + kc)\sin\alpha \quad a^2 + b^2 + c^2 = 1$$

то (7.2.6) описывает вращение вектора $(x,y,z)$ вокруг вектора $(a,b,c)$ на угол $2\alpha$.

## 7.3. Линейная функция алгебры кватернионов

Теорема 7.3.1. *Положим*

$$e_0 = 1 \quad e_1 = i \quad e_2 = j \quad e_3 = k \tag{7.3.1}$$

*базис алгебры кватернионов $H$. Тогда в базисе (7.3.1) структурные константы имеют вид*

$$\begin{array}{llll}
C^0_{00} = 1 & C^1_{01} = \phantom{-}1 & C^2_{02} = \phantom{-}1 & C^3_{03} = \phantom{-}1 \\
C^1_{10} = 1 & C^0_{11} = -1 & C^3_{12} = \phantom{-}1 & C^2_{13} = -1 \\
C^2_{20} = 1 & C^3_{21} = -1 & C^0_{22} = -1 & C^1_{23} = \phantom{-}1 \\
C^3_{30} = 1 & C^2_{31} = \phantom{-}1 & C^1_{32} = -1 & C^0_{33} = -1
\end{array}$$

[7.2] См. [10], с. 643.

Доказательство. Значение структурных констант следует из таблицы умножения (7.2.3). □

Так как вычисления в этом разделе занимают много места, я собрал в одном месте ссылки на теоремы в этом разделе.

**Теорема 7.3.2, страница 72:** определение координат линейного отображения алгебры кватернионов $H$ через стандартные компоненты этого отображения.

**Равенство** (7.3.22)**, страница 75:** матричная форма зависимости координат линейного отображения алгебры кватернионов $H$ от стандартных компонент этого отображения.

**Равенство** (7.3.23)**, страница 75:** матричная форма зависимости стандартных компонент линейного отображения алгебры кватернионов $H$ от координат этого отображения.

**Теорема 7.3.4, страница 78:** зависимость стандартных компонент линейного отображения алгебры кватернионов $H$ от координат этого отображения.

Теорема 7.3.2. *Стандартные компоненты линейной функции алгебры кватернионов $H$ относительно базиса* (7.3.1) *и координаты соответствующего линейного преобразования удовлетворяют соотношениям*

$$\left\{\begin{array}{l} f_0^0 = f^{00} - f^{11} - f^{22} - f^{33} \\ f_1^1 = f^{00} - f^{11} + f^{22} + f^{33} \\ f_2^2 = f^{00} + f^{11} - f^{22} + f^{33} \\ f_3^3 = f^{00} + f^{11} + f^{22} - f^{33} \end{array}\right. \tag{7.3.2}$$

$$\left\{\begin{array}{l} f_0^1 = f^{01} + f^{10} + f^{23} - f^{32} \\ f_1^0 = -f^{01} - f^{10} + f^{23} - f^{32} \\ f_2^3 = -f^{01} + f^{10} - f^{23} - f^{32} \\ f_3^2 = f^{01} - f^{10} - f^{23} - f^{32} \end{array}\right. \tag{7.3.3}$$

$$\left\{\begin{array}{l} f_0^2 = f^{02} - f^{13} + f^{20} + f^{31} \\ f_1^3 = f^{02} - f^{13} - f^{20} - f^{31} \\ f_2^0 = -f^{02} - f^{13} - f^{20} + f^{31} \\ f_3^1 = -f^{02} - f^{13} + f^{20} - f^{31} \end{array}\right. \tag{7.3.4}$$

$$\left\{\begin{array}{l} f_0^3 = f^{03} + f^{12} - f^{21} + f^{30} \\ f_1^2 = -f^{03} - f^{12} - f^{21} + f^{30} \\ f_2^1 = f^{03} - f^{12} - f^{21} - f^{30} \\ f_3^0 = -f^{03} + f^{12} - f^{21} - f^{30} \end{array}\right. \tag{7.3.5}$$

Доказательство. Пользуясь равенством (3.1.17) получаем соотношения

$$\begin{aligned} f_0^0 &= f^{kr} C_{k0}^p C_{pr}^0 \\ &= f^{00} C_{00}^0 C_{00}^0 + f^{11} C_{10}^1 C_{11}^0 + f^{22} C_{20}^2 C_{22}^0 + f^{33} C_{30}^3 C_{33}^0 \\ &= f^{00} - f^{11} - f^{22} - f^{33} \end{aligned} \tag{7.3.6}$$

$$\begin{aligned} f_0^1 &= f^{kr} C_{k0}^p C_{pr}^1 \\ &= f^{01} C_{00}^0 C_{01}^1 + f^{10} C_{10}^1 C_{10}^1 + f^{23} C_{20}^2 C_{23}^1 + f^{32} C_{30}^3 C_{32}^1 \\ &= f^{01} + f^{10} + f^{23} - f^{32} \end{aligned} \tag{7.3.7}$$

$$\begin{aligned} f_0^2 &= f^{kr} C_{k0}^p C_{pr}^2 \\ &= f^{02} C_{00}^0 C_{02}^2 + f^{13} C_{10}^1 C_{13}^2 + f^{20} C_{20}^2 C_{20}^2 + f^{31} C_{30}^3 C_{31}^2 \\ &= f^{02} - f^{13} + f^{20} + f^{31} \end{aligned} \tag{7.3.8}$$

$$\begin{aligned} f_0^3 &= f^{kr} C_{k0}^p C_{pr}^3 \\ &= f^{03} C_{00}^0 C_{03}^3 + f^{12} C_{10}^1 C_{12}^3 + f^{21} C_{20}^2 C_{21}^3 + f^{30} C_{30}^3 C_{30}^3 \\ &= f^{03} + f^{12} - f^{21} + f^{30} \end{aligned} \tag{7.3.9}$$

$$\begin{aligned} f_1^0 &= f^{kr} C_{k1}^p C_{pr}^0 \\ &= f^{01} C_{01}^1 C_{11}^0 + f^{10} C_{11}^0 C_{00}^0 + f^{23} C_{21}^3 C_{33}^0 + f^{32} C_{31}^2 C_{22}^0 \\ &= -f^{01} - f^{10} + f^{23} - f^{32} \end{aligned} \tag{7.3.10}$$

$$\begin{aligned} f_1^1 &= f^{kr} C_{k1}^p C_{pr}^1 \\ &= f^{00} C_{01}^1 C_{10}^1 + f^{11} C_{11}^0 C_{01}^1 + f^{22} C_{21}^3 C_{32}^1 + f^{33} C_{31}^2 C_{23}^1 \\ &= f^{00} - f^{11} + f^{22} + f^{33} \end{aligned} \tag{7.3.11}$$

$$\begin{aligned} f_1^2 &= f^{kr} C_{k1}^p C_{pr}^2 \\ &= f^{03} C_{01}^1 C_{13}^2 + f^{12} C_{11}^0 C_{02}^2 + f^{21} C_{21}^3 C_{31}^2 + f^{30} C_{31}^2 C_{20}^2 \\ &= -f^{03} - f^{12} - f^{21} + f^{30} \end{aligned} \tag{7.3.12}$$

$$\begin{aligned} f_1^3 &= f^{kr} C_{k1}^p C_{pr}^3 \\ &= f^{02} C_{01}^1 C_{12}^3 + f^{13} C_{11}^0 C_{03}^3 + f^{20} C_{21}^3 C_{30}^3 + f^{31} C_{31}^2 C_{21}^3 \\ &= f^{02} - f^{13} - f^{20} - f^{31} \end{aligned} \tag{7.3.13}$$

$$\begin{aligned}f_2^0 &= f^{kr}C_{k2}^pC_{pr}^0\\ &= f^{02}C_{02}^2C_{22}^0+f^{13}C_{12}^3C_{33}^0+f^{20}C_{22}^0C_{00}^0+f^{31}C_{32}^1C_{11}^0\\ &= -f^{02}-f^{13}-f^{20}+f^{31}\end{aligned}\tag{7.3.14}$$

$$\begin{aligned}f_2^1 &= f^{kr}C_{k2}^pC_{pr}^1\\ &= f^{03}C_{02}^2C_{23}^1+f^{12}C_{12}^3C_{32}^1+f^{21}C_{22}^0C_{01}^1+f^{30}C_{32}^1C_{10}^1\\ &= f^{03}-f^{12}-f^{21}-f^{30}\end{aligned}\tag{7.3.15}$$

$$\begin{aligned}f_2^2 &= f^{kr}C_{k2}^pC_{pr}^2\\ &= f^{00}C_{02}^2C_{20}^2+f^{11}C_{12}^3C_{31}^2+f^{22}C_{22}^0C_{02}^2+f^{33}C_{32}^1C_{13}^2\\ &= f^{00}+f^{11}-f^{22}+f^{33}\end{aligned}\tag{7.3.16}$$

$$\begin{aligned}f_2^3 &= f^{kr}C_{k2}^pC_{pr}^3\\ &= f^{01}C_{02}^2C_{21}^3+f^{10}C_{12}^3C_{30}^3+f^{23}C_{22}^0C_{03}^3+f^{32}C_{32}^1C_{12}^3\\ &= -f^{01}+f^{10}-f^{23}-f^{32}\end{aligned}\tag{7.3.17}$$

$$\begin{aligned}f_3^0 &= f^{kr}C_{k3}^pC_{pr}^0\\ &= f^{03}C_{03}^3C_{33}^0+f^{12}C_{13}^2C_{22}^0+f^{21}C_{23}^1C_{11}^0+f^{30}C_{33}^0C_{00}^0\\ &= -f^{03}+f^{12}-f^{21}-f^{30}\end{aligned}\tag{7.3.18}$$

$$\begin{aligned}f_3^1 &= f^{kr}C_{k3}^pC_{pr}^1\\ &= f^{02}C_{03}^3C_{32}^1+f^{13}C_{13}^2C_{23}^1+f^{20}C_{23}^1C_{10}^1+f^{31}C_{33}^0C_{01}^1\\ &= -f^{02}-f^{13}+f^{20}-f^{31}\end{aligned}\tag{7.3.19}$$

$$\begin{aligned}f_3^2 &= f^{kr}C_{k3}^pC_{pr}^2\\ &= f^{01}C_{03}^3C_{31}^2+f^{10}C_{13}^2C_{20}^2+f^{23}C_{23}^1C_{13}^2+f^{32}C_{33}^0C_{02}^2\\ &= f^{01}-f^{10}-f^{23}-f^{32}\end{aligned}\tag{7.3.20}$$

$$\begin{aligned}f_3^3 &= f^{kr}C_{k3}^pC_{pr}^3\\ &= f^{00}C_{03}^3C_{30}^3+f^{11}C_{13}^2C_{21}^3+f^{22}C_{23}^1C_{12}^3+f^{33}C_{33}^0C_{03}^3\\ &= f^{00}+f^{11}+f^{22}-f^{33}\end{aligned}\tag{7.3.21}$$

Уравнения (7.3.6), (7.3.11), (7.3.16), (7.3.21) формируют систему линейных уравнений (7.3.2).

Уравнения (7.3.7), (7.3.10), (7.3.17), (7.3.20) формируют систему линейных уравнений (7.3.3).

Уравнения (7.3.8), (7.3.13), (7.3.14), (7.3.19) формируют систему линейных уравнений (7.3.4).

Уравнения (7.3.9), (7.3.12), (7.3.15), (7.3.18) формируют систему линейных уравнений (7.3.5). □

Теорема 7.3.3. *Рассмотрим алгебру кватернионов $H$ с базисом (7.3.1). Стандартные компоненты линейной функции над полем $F$ и координаты этой функции над полем $F$ удовлетворяют соотношениям*

$$\begin{pmatrix} f_{0}^{0} & f_{1}^{0} & f_{2}^{0} & f_{3}^{0} \\ f_{1}^{1} & -f_{1}^{0} & f_{3}^{1} & -f_{2}^{1} \\ f_{2}^{2} & -f_{3}^{2} & -f_{0}^{2} & f_{1}^{2} \\ f_{3}^{3} & f_{2}^{3} & -f_{1}^{3} & -f_{0}^{3} \end{pmatrix} = \begin{pmatrix} 1 & -1 & -1 & -1 \\ 1 & -1 & 1 & 1 \\ 1 & 1 & -1 & 1 \\ 1 & 1 & 1 & -1 \end{pmatrix} \begin{pmatrix} f^{00} & -f^{01} & -f^{02} & -f^{03} \\ f^{11} & f^{10} & f^{13} & -f^{12} \\ f^{22} & -f^{23} & f^{20} & f^{21} \\ f^{33} & f^{32} & -f^{31} & f^{30} \end{pmatrix} \tag{7.3.22}$$

$$\begin{pmatrix} f^{00} & -f^{01} & -f^{02} & -f^{03} \\ f^{11} & f^{10} & f^{13} & -f^{12} \\ f^{22} & -f^{23} & f^{20} & f^{21} \\ f^{33} & f^{32} & -f^{31} & f^{30} \end{pmatrix} = \frac{1}{4} \begin{pmatrix} 1 & 1 & 1 & 1 \\ -1 & -1 & 1 & 1 \\ -1 & 1 & -1 & 1 \\ -1 & 1 & 1 & -1 \end{pmatrix} \begin{pmatrix} f_{0}^{0} & f_{1}^{0} & f_{2}^{0} & f_{3}^{0} \\ f_{1}^{1} & -f_{1}^{0} & f_{3}^{1} & -f_{2}^{1} \\ f_{2}^{2} & -f_{3}^{2} & -f_{0}^{2} & f_{1}^{2} \\ f_{3}^{3} & f_{2}^{3} & -f_{1}^{3} & -f_{0}^{3} \end{pmatrix} \tag{7.3.23}$$

*где*

$$\begin{pmatrix} 1 & -1 & -1 & -1 \\ 1 & -1 & 1 & 1 \\ 1 & 1 & -1 & 1 \\ 1 & 1 & 1 & -1 \end{pmatrix}^{-1} = \frac{1}{4} \begin{pmatrix} 1 & 1 & 1 & 1 \\ -1 & -1 & 1 & 1 \\ -1 & 1 & -1 & 1 \\ -1 & 1 & 1 & -1 \end{pmatrix}$$

Доказательство. Запишем систему линейных уравнений (7.3.2) в виде произведения матриц

$$\begin{pmatrix} f_{0}^{0} \\ f_{1}^{1} \\ f_{2}^{2} \\ f_{3}^{3} \end{pmatrix} = \begin{pmatrix} 1 & -1 & -1 & -1 \\ 1 & -1 & 1 & 1 \\ 1 & 1 & -1 & 1 \\ 1 & 1 & 1 & -1 \end{pmatrix} \begin{pmatrix} f^{00} \\ f^{11} \\ f^{22} \\ f^{33} \end{pmatrix} \tag{7.3.24}$$

Запишем систему линейных уравнений (7.3.3) в виде произведения матриц

$$
\begin{pmatrix} f_0^1 \\ f_1^0 \\ f_2^3 \\ f_3^2 \end{pmatrix} = \begin{pmatrix} 1 & 1 & 1 & -1 \\ -1 & -1 & 1 & -1 \\ -1 & 1 & -1 & -1 \\ 1 & -1 & -1 & -1 \end{pmatrix} \begin{pmatrix} f^{01} \\ f^{10} \\ f^{23} \\ f^{32} \end{pmatrix} \tag{7.3.25}
$$

Из равенства (7.3.25) следует

$$
\begin{pmatrix} -f_1^0 \\ f_1^0 \\ f_2^3 \\ -f_3^2 \end{pmatrix} = \begin{pmatrix} -1 & -1 & -1 & 1 \\ -1 & -1 & 1 & -1 \\ -1 & 1 & -1 & -1 \\ -1 & 1 & 1 & 1 \end{pmatrix} \begin{pmatrix} f^{01} \\ f^{10} \\ f^{23} \\ f^{32} \end{pmatrix}
= \begin{pmatrix} 1 & -1 & 1 & 1 \\ 1 & -1 & -1 & -1 \\ 1 & 1 & 1 & -1 \\ 1 & 1 & -1 & 1 \end{pmatrix} \begin{pmatrix} -f^{01} \\ f^{10} \\ -f^{23} \\ f^{32} \end{pmatrix}
$$

$$
\begin{pmatrix} f_1^0 \\ -f_1^0 \\ -f_3^2 \\ f_2^3 \end{pmatrix} = \begin{pmatrix} 1 & -1 & -1 & -1 \\ 1 & -1 & 1 & 1 \\ 1 & 1 & -1 & 1 \\ 1 & 1 & 1 & -1 \end{pmatrix} \begin{pmatrix} -f^{01} \\ f^{10} \\ -f^{23} \\ f^{32} \end{pmatrix} \tag{7.3.26}
$$

Запишем систему линейных уравнений (7.3.4) в виде произведения матриц

$$
\begin{pmatrix} f_0^2 \\ f_1^3 \\ f_2^0 \\ f_3^1 \end{pmatrix} = \begin{pmatrix} 1 & -1 & 1 & 1 \\ 1 & -1 & -1 & -1 \\ -1 & -1 & -1 & 1 \\ -1 & -1 & 1 & -1 \end{pmatrix} \begin{pmatrix} f^{02} \\ f^{13} \\ f^{20} \\ f^{31} \end{pmatrix} \tag{7.3.27}
$$

Из равенства (7.3.27) следует

$$\begin{pmatrix}-f_{0}^{2}\\-f_{1}^{3}\\f_{2}^{0}\\f_{3}^{1}\end{pmatrix}=\begin{pmatrix}-1&1&-1&-1\\-1&1&1&1\\-1&-1&-1&1\\-1&-1&1&-1\end{pmatrix}\begin{pmatrix}f^{02}\\f^{13}\\f^{20}\\f^{31}\end{pmatrix}$$

$$=\begin{pmatrix}1&1&-1&1\\1&1&1&-1\\1&-1&-1&-1\\1&-1&1&1\end{pmatrix}\begin{pmatrix}-f^{02}\\f^{13}\\f^{20}\\-f^{31}\end{pmatrix}$$

$$\begin{pmatrix}f_{2}^{0}\\f_{3}^{1}\\-f_{0}^{2}\\-f_{1}^{3}\end{pmatrix}=\begin{pmatrix}1&-1&-1&-1\\1&-1&1&1\\1&1&-1&1\\1&1&1&-1\end{pmatrix}\begin{pmatrix}-f^{02}\\f^{13}\\f^{20}\\-f^{31}\end{pmatrix}\tag{7.3.28}$$

Запишем систему линейных уравнений (7.3.5) в виде произведения матриц

$$\begin{pmatrix}f_{0}^{3}\\f_{1}^{2}\\f_{2}^{1}\\f_{3}^{0}\end{pmatrix}=\begin{pmatrix}1&1&-1&1\\-1&-1&-1&1\\1&-1&-1&-1\\-1&1&-1&-1\end{pmatrix}\begin{pmatrix}f^{03}\\f^{12}\\f^{21}\\f^{30}\end{pmatrix}\tag{7.3.29}$$

Из равенства (7.3.29) следует

$$\begin{pmatrix}-f_{0}^{3}\\f_{1}^{2}\\-f_{2}^{1}\\f_{3}^{0}\end{pmatrix}=\begin{pmatrix}-1&-1&1&-1\\-1&-1&-1&1\\-1&1&1&1\\-1&1&-1&-1\end{pmatrix}\begin{pmatrix}f^{03}\\f^{12}\\f^{21}\\f^{30}\end{pmatrix}$$

$$=\begin{pmatrix}1&1&1&-1\\1&1&-1&1\\1&-1&1&1\\1&-1&-1&-1\end{pmatrix}\begin{pmatrix}-f^{03}\\-f^{12}\\f^{21}\\f^{30}\end{pmatrix}$$

$$\begin{pmatrix} f_3^0 \\ -f_2^1 \\ f_1^2 \\ -f_0^3 \end{pmatrix} = \begin{pmatrix} 1 & -1 & -1 & -1 \\ 1 & -1 & 1 & 1 \\ 1 & 1 & -1 & 1 \\ 1 & 1 & 1 & -1 \end{pmatrix} \begin{pmatrix} -f^{03} \\ -f^{12} \\ f^{21} \\ f^{30} \end{pmatrix} \tag{7.3.30}$$

Мы объединяем равенства (7.3.24), (7.3.26), (7.3.28), (7.3.30) в равенстве (7.3.22). □

Теорема 7.3.4. *Стандартные компоненты линейной функции алгебры кватернионов $H$ относительно базиса* (7.3.1) *и координаты соответствующего линейного преобразования удовлетворяют соотношениям*

$$\left\{\begin{array}{l} 4f^{00} = f_0^0 + f_1^1 + f_2^2 + f_3^3 \\ 4f^{11} = -f_0^0 - f_1^1 + f_2^2 + f_3^3 \\ 4f^{22} = -f_0^0 + f_1^1 - f_2^2 + f_3^3 \\ 4f^{33} = -f_0^0 + f_1^1 + f_2^2 - f_3^3 \end{array}\right. \tag{7.3.31}$$

$$\left\{\begin{array}{l} 4f^{10} = -f_1^0 + f_1^0 - f_3^2 + f_2^3 \\ 4f^{01} = -f_1^0 + f_1^0 + f_3^2 - f_2^3 \\ 4f^{32} = -f_1^0 - f_1^0 - f_3^2 - f_2^3 \\ 4f^{23} = f_1^0 + f_1^0 - f_3^2 - f_2^3 \end{array}\right. \tag{7.3.32}$$

$$\left\{\begin{array}{l} 4f^{20} = -f_2^0 + f_3^1 + f_0^2 - f_1^3 \\ 4f^{31} = f_2^0 - f_3^1 + f_0^2 - f_1^3 \\ 4f^{02} = -f_2^0 - f_3^1 + f_0^2 + f_1^3 \\ 4f^{13} = -f_2^0 - f_3^1 - f_0^2 - f_1^3 \end{array}\right. \tag{7.3.33}$$

$$\left\{\begin{array}{l} 4f^{30} = -f_3^0 - f_2^1 + f_1^2 + f_0^3 \\ 4f^{21} = -f_3^0 - f_2^1 - f_1^2 - f_0^3 \\ 4f^{12} = f_3^0 - f_2^1 - f_1^2 + f_0^3 \\ 4f^{03} = -f_3^0 + f_2^1 - f_1^2 + f_0^3 \end{array}\right. \tag{7.3.34}$$

Доказательство. Системы линейных уравнений (7.3.31), (7.3.32), (7.3.33), (7.3.34) получены в результате перемножения матриц в равенстве (7.3.23). □

## 7.4. Дифференцируемое отображение тела кватернионов

Теорема 7.4.1. *Если матрица* $\left(\frac{\partial y^{i}}{\partial x^{j}}\right)$ *является матрица Якоби функции* $x \to y$ *тела кватернионов над полем действительных чисел, то*

$$\begin{cases}\frac{\partial y^{0}}{\partial x^{0}}=\frac{\partial^{00}y}{\partial x}-\frac{\partial^{11}y}{\partial x}-\frac{\partial^{22}y}{\partial x}-\frac{\partial^{33}y}{\partial x}\\ \frac{\partial y^{1}}{\partial x^{1}}=\frac{\partial^{00}y}{\partial x}-\frac{\partial^{11}y}{\partial x}+\frac{\partial^{22}y}{\partial x}+\frac{\partial^{33}y}{\partial x}\\ \frac{\partial y^{2}}{\partial x^{2}}=\frac{\partial^{00}y}{\partial x}+\frac{\partial^{11}y}{\partial x}-\frac{\partial^{22}y}{\partial x}+\frac{\partial^{33}y}{\partial x}\\ \frac{\partial y^{3}}{\partial x^{3}}=\frac{\partial^{00}y}{\partial x}+\frac{\partial^{11}y}{\partial x}+\frac{\partial^{22}y}{\partial x}-\frac{\partial^{33}y}{\partial x}\end{cases}\tag{7.4.1}$$

$$\begin{cases}\frac{\partial y^{1}}{\partial x^{0}}=\frac{\partial^{01}y}{\partial x}+\frac{\partial^{10}y}{\partial x}+\frac{\partial^{23}y}{\partial x}-\frac{\partial^{32}y}{\partial x}\\ \frac{\partial y^{0}}{\partial x^{1}}=-\frac{\partial^{01}y}{\partial x}-\frac{\partial^{10}y}{\partial x}+\frac{\partial^{23}y}{\partial x}-\frac{\partial^{32}y}{\partial x}\\ \frac{\partial y^{3}}{\partial x^{2}}=-\frac{\partial^{01}y}{\partial x}+\frac{\partial^{10}y}{\partial x}-\frac{\partial^{23}y}{\partial x}-\frac{\partial^{32}y}{\partial x}\\ \frac{\partial y^{2}}{\partial x^{3}}=\frac{\partial^{01}y}{\partial x}-\frac{\partial^{10}y}{\partial x}-\frac{\partial^{23}y}{\partial x}-\frac{\partial^{32}y}{\partial x}\end{cases}\tag{7.4.2}$$

$$\begin{cases}\frac{\partial y^{2}}{\partial x^{0}}=\frac{\partial^{02}y}{\partial x}-\frac{\partial^{13}y}{\partial x}+\frac{\partial^{20}y}{\partial x}+\frac{\partial^{31}y}{\partial x}\\ \frac{\partial y^{3}}{\partial x^{1}}=\frac{\partial^{02}y}{\partial x}-\frac{\partial^{13}y}{\partial x}-\frac{\partial^{20}y}{\partial x}-\frac{\partial^{31}y}{\partial x}\\ \frac{\partial y^{0}}{\partial x^{2}}=-\frac{\partial^{02}y}{\partial x}-\frac{\partial^{13}y}{\partial x}-\frac{\partial^{20}y}{\partial x}+\frac{\partial^{31}y}{\partial x}\\ \frac{\partial y^{1}}{\partial x^{3}}=-\frac{\partial^{02}y}{\partial x}-\frac{\partial^{13}y}{\partial x}+\frac{\partial^{20}y}{\partial x}-\frac{\partial^{31}y}{\partial x}\end{cases}\tag{7.4.3}$$

$$\begin{cases}\frac{\partial y^{3}}{\partial x^{0}}=\frac{\partial^{03}y}{\partial x}+\frac{\partial^{12}y}{\partial x}-\frac{\partial^{21}y}{\partial x}+\frac{\partial^{30}y}{\partial x}\\ \frac{\partial y^{2}}{\partial x^{1}}=-\frac{\partial^{03}y}{\partial x}-\frac{\partial^{12}y}{\partial x}-\frac{\partial^{21}y}{\partial x}+\frac{\partial^{30}y}{\partial x}\\ \frac{\partial y^{1}}{\partial x^{2}}=\frac{\partial^{03}y}{\partial x}-\frac{\partial^{12}y}{\partial x}-\frac{\partial^{21}y}{\partial x}-\frac{\partial^{30}y}{\partial x}\\ \frac{\partial y^{0}}{\partial x^{3}}=-\frac{\partial^{03}y}{\partial x}+\frac{\partial^{12}y}{\partial x}-\frac{\partial^{21}y}{\partial x}-\frac{\partial^{30}y}{\partial x}\end{cases}\tag{7.4.4}$$

$$
\begin{cases}
4\dfrac{\partial^{00}y}{\partial x}=\dfrac{\partial y^0}{\partial x^0}+\dfrac{\partial y^1}{\partial x^1}+\dfrac{\partial y^2}{\partial x^2}+\dfrac{\partial y^3}{\partial x^3}\\
4\dfrac{\partial^{11}y}{\partial x}=-\dfrac{\partial y^0}{\partial x^0}-\dfrac{\partial y^1}{\partial x^1}+\dfrac{\partial y^2}{\partial x^2}+\dfrac{\partial y^3}{\partial x^3}\\
4\dfrac{\partial^{22}y}{\partial x}=-\dfrac{\partial y^0}{\partial x^0}+\dfrac{\partial y^1}{\partial x^1}-\dfrac{\partial y^2}{\partial x^2}+\dfrac{\partial y^3}{\partial x^3}\\
4\dfrac{\partial^{33}y}{\partial x}=-\dfrac{\partial y^0}{\partial x^0}+\dfrac{\partial y^1}{\partial x^1}+\dfrac{\partial y^2}{\partial x^2}-\dfrac{\partial y^3}{\partial x^3}
\end{cases}
\tag{7.4.5}
$$

$$
\begin{cases}
4\dfrac{\partial^{10}y}{\partial x}=-\dfrac{\partial y^0}{\partial x^1}+\dfrac{\partial y^1}{\partial x^0}-\dfrac{\partial y^2}{\partial x^3}+\dfrac{\partial y^3}{\partial x^2}\\
4\dfrac{\partial^{01}y}{\partial x}=-\dfrac{\partial y^0}{\partial x^1}+\dfrac{\partial y^1}{\partial x^0}+\dfrac{\partial y^2}{\partial x^3}-\dfrac{\partial y^3}{\partial x^2}\\
4\dfrac{\partial^{32}y}{\partial x}=-\dfrac{\partial y^0}{\partial x^1}-\dfrac{\partial y^1}{\partial x^0}-\dfrac{\partial y^2}{\partial x^3}-\dfrac{\partial y^3}{\partial x^2}\\
4\dfrac{\partial^{23}y}{\partial x}=\dfrac{\partial y^0}{\partial x^1}+\dfrac{\partial y^1}{\partial x^0}-\dfrac{\partial y^2}{\partial x^3}-\dfrac{\partial y^3}{\partial x^2}
\end{cases}
\tag{7.4.6}
$$

$$
\begin{cases}
4\dfrac{\partial^{20}y}{\partial x}=-\dfrac{\partial y^0}{\partial x^2}+\dfrac{\partial y^1}{\partial x^3}+\dfrac{\partial y^2}{\partial x^0}-\dfrac{\partial y^3}{\partial x^1}\\
4\dfrac{\partial^{31}y}{\partial x}=+\dfrac{\partial y^0}{\partial x^2}-\dfrac{\partial y^1}{\partial x^3}+\dfrac{\partial y^2}{\partial x^0}-\dfrac{\partial y^3}{\partial x^1}\\
4\dfrac{\partial^{02}y}{\partial x}=-\dfrac{\partial y^0}{\partial x^2}-\dfrac{\partial y^1}{\partial x^3}+\dfrac{\partial y^2}{\partial x^0}+\dfrac{\partial y^3}{\partial x^1}\\
4\dfrac{\partial^{13}y}{\partial x}=-\dfrac{\partial y^0}{\partial x^2}-\dfrac{\partial y^1}{\partial x^3}-\dfrac{\partial y^2}{\partial x^0}-\dfrac{\partial y^3}{\partial x^1}
\end{cases}
\tag{7.4.7}
$$

$$
\begin{cases}
4\dfrac{\partial^{30}y}{\partial x}=-\dfrac{\partial y^0}{\partial x^3}-\dfrac{\partial y^1}{\partial x^2}+\dfrac{\partial y^2}{\partial x^1}+\dfrac{\partial y^3}{\partial x^0}\\
4\dfrac{\partial^{21}y}{\partial x}=-\dfrac{\partial y^0}{\partial x^3}-\dfrac{\partial y^1}{\partial x^2}-\dfrac{\partial y^2}{\partial x^1}-\dfrac{\partial y^3}{\partial x^0}\\
4\dfrac{\partial^{12}y}{\partial x}=\dfrac{\partial y^0}{\partial x^3}-\dfrac{\partial y^1}{\partial x^2}-\dfrac{\partial y^2}{\partial x^1}+\dfrac{\partial y^3}{\partial x^0}\\
4\dfrac{\partial^{03}y}{\partial x}=-\dfrac{\partial y^0}{\partial x^3}+\dfrac{\partial y^1}{\partial x^2}-\dfrac{\partial y^2}{\partial x^1}+\dfrac{\partial y^3}{\partial x^0}
\end{cases}
\tag{7.4.8}
$$

Доказательство. Следствие теоремы 7.3.2. □

Теорема 7.4.2. *Отображение кватернионов*

$$f(x)=\overline{x}$$

*имеет производную Гато*

$$\partial(\overline{x})(h)=-\frac{1}{2}(h+ihi+jhj+khk)\tag{7.4.9}$$

Доказательство. Матрица Якоби отображения $f$ имеет вид

$$\begin{pmatrix} 1 & 0 & 0 & 0 \\ 0 & -1 & 0 & 0 \\ 0 & 0 & -1 & 0 \\ 0 & 0 & 0 & -1 \end{pmatrix} \tag{7.4.10}$$

Из равенств (7.4.5) следует

$$\frac{\partial^{00}y}{\partial x} = \frac{\partial^{11}y}{\partial x} = \frac{\partial^{22}y}{\partial x} = \frac{\partial^{33}y}{\partial x} = -\frac{1}{2} \tag{7.4.11}$$

Из равенств (7.4.6), (7.4.7), (7.4.8) следует

$$\frac{\partial^{ij}y}{\partial x} = 0 \quad i \neq j \tag{7.4.12}$$

Равенство (7.4.9) следует из равенств (5.2.15), (7.4.11), (7.4.12). □

Теорема 7.4.3. *Норма кватерниона удовлетворяет равенствам*[7.3]

$$\partial(|x|^2) \circ dx = x\overline{dx} + \overline{x}dx \tag{7.4.13}$$

$$\partial(|x|)(h) = \frac{1}{2\sqrt{x\overline{dx} + \overline{x}dx}} \tag{7.4.14}$$

Доказательство. Равенство (7.4.13) следует из теоремы 7.4.2 и равенств (7.2.4), (5.2.23). Так как $x\overline{dx} + \overline{x}dx \in R$, то равенство (7.4.14) следует из равенства (7.4.13). □

Теорема 7.4.4. *В теле кватернионов над полем действительных чисел $H$ с базисом (7.3.1) стандартное $D\star$-представление дифференциала Гато отображения $x^2$ имеет вид*

$$\partial x^2(a) = (x + x_1)a + x_2 ai + x_3 aj + x_4 ak$$

Доказательство. Следствие определения 5.2.10. □

Теорема 7.4.5. *Мы можем отождествить кватернион*

$$a = a^0 + a^1 i + a^2 j + a^3 k \tag{7.4.15}$$

*и матрицу*

$$J_a = \begin{pmatrix} a^0 & -a^1 & -a^2 & -a^3 \\ a^1 & a^0 & -a^3 & a^2 \\ a^2 & a^3 & a^0 & -a^1 \\ a^3 & -a^2 & a^1 & a^0 \end{pmatrix} \tag{7.4.16}$$

Доказательство. Произведение кватернионов (7.4.15) и

$$x = x^0 + x^1 i + x^2 j + x^3 k$$

[7.3] Утверждение теоремы аналогично примеру X, [11], с. 455.

имеет вид

$$\begin{aligned}ax = a^0x^0 - a^1x^1 - a^2x^2 - a^3x^3 + (a^0x^1 + a^1x^0 + a^2x^3 - a^3x^2)i\\ + (a^0x^2 + a^2x^0 + a^3x^1 - a^1x^3)j + (a^0x^3 + a^3x^0 + a^1x^2 - a^2x^1)k\end{aligned}$$

Следовательно, функция $f_a(x) = ax$ имеет матрицу Якоби (7.4.16). Очевидно, что $f_a \circ f_b = f_{ab}$. Аналогичное равенство верно для матриц

$$\begin{pmatrix} a^0 & -a^1 & -a^2 & -a^3 \\ a^1 & a^0 & -a^3 & a^2 \\ a^2 & a^3 & a^0 & -a^1 \\ a^3 & -a^2 & a^1 & a^0 \end{pmatrix} \begin{pmatrix} b^0 & -b^1 & -b^2 & -b^3 \\ b^1 & b^0 & -b^3 & b^2 \\ b^2 & b^3 & b^0 & -b^1 \\ b^3 & -b^2 & b^1 & b^0 \end{pmatrix}$$

$$= \begin{pmatrix} \begin{matrix} a^0b^0 - a^1b^1 \\ -a^2b^2 - a^3b^3 \end{matrix} & \begin{matrix} -a^0b^1 - a^1b^0 \\ -a^2b^3 + a^3b^2 \end{matrix} & \begin{matrix} -a^0b^2 + a^1b^3 \\ -a^2b^0 - a^3b^1 \end{matrix} & \begin{matrix} -a^0b^3 - a^1b^2 \\ +a^2b^1 - a^3b^0 \end{matrix} \\ \\ \begin{matrix} a^0b^1 + a^1b^0 \\ +a^2b^3 - a^3b^2 \end{matrix} & \begin{matrix} a^0b^0 - a^1b^1 \\ -a^2b^2 - a^3b^3 \end{matrix} & \begin{matrix} -a^0b^3 - a^1b^2 \\ +a^2b^1 - a^3b^0 \end{matrix} & \begin{matrix} a^0b^2 - a^1b^3 \\ +a^2b^0 + a^3b^1 \end{matrix} \\ \\ \begin{matrix} a^0b^2 - a^1b^3 \\ +a^2b^0 + a^3b^1 \end{matrix} & \begin{matrix} a^0b^3 + a^1b^2 \\ -a^2b^1 + a^3b^0 \end{matrix} & \begin{matrix} a^0b^0 - a^1b^1 \\ -a^2b^2 - a^3b^3 \end{matrix} & \begin{matrix} -a^0b^1 - a^1b^0 \\ -a^2b^3 - a^3b^2 \end{matrix} \\ \\ \begin{matrix} a^0b^3 + a^1b^2 \\ -a^2b^1 + a^3b^0 \end{matrix} & \begin{matrix} -a^0b^2 + a^1b^3 \\ -a^2b^0 - a^3b^1 \end{matrix} & \begin{matrix} a^0b^1 + a^1b^0 \\ +a^2b^3 - a^3b^2 \end{matrix} & \begin{matrix} a^0b^0 - a^1b^1 \\ -a^2b^2 - a^3b^3 \end{matrix} \end{pmatrix}$$

Следовательно, мы можем отождествить кватернион $a$ и матрицу $J_a$. □

Глава 8

# Производная второго порядка отображения тела

## 8.1. Производная второго порядка отображения тела

Пусть $D$ - нормированное тело характеристики 0. Пусть

$$f:D\to D$$

функция, дифференцируемая по Гато. Согласно замечанию 5.2.4 производная Гато является отображением

$$\partial f:D\to\mathcal{L}(D;D)$$

Согласно теоремам 3.1.5, 3.1.6 и определению 5.1.13 множество $\mathcal{L}(D;D)$ является нормированным $D$-векторным пространством. Следовательно, мы можем рассмотреть вопрос, является ли отображение $\partial f$ дифференцируемым по Гато.

Согласно определению 6.2.1

$$\partial f(x+a_2)(a_1)-\partial f(x)(a_1)=\partial(\partial f(x)(a_1))(a_2)+o_2(a_2) \tag{8.1.1}$$

где $o_2:D\to\mathcal{L}(D;D)$ - такое непрерывное отображение, что

$$\lim_{a_2\to 0}\frac{\|o_2(a_2)\|}{|a_2|}=0$$

Согласно определению 6.2.1 отображение $\partial(\partial f(x)(a_1))(a_2)$ линейно по переменной $a_2$. Из равенства (8.1.1) следует, что отображение $\partial(\partial f(x)(a_1))(a_2)$ линейно по переменной $a_1$.

Определение 8.1.1. Полилинейное отображение

$$\partial^2 f(x)(a_1;a_2)=\frac{\partial^2 f(x)}{\partial x^2}(a_1;a_2)=\partial(\partial f(x)(a_1))(a_2) \tag{8.1.2}$$

называется **производной Гато второго порядка** отображения $f$. □

Замечание 8.1.2. Согласно определению 8.1.1 при заданном $x$ производная Гато второго порядка $\partial^2 f(x)\in\mathcal{L}(D,D;D)$. Следовательно, производная Гато второго порядка отображения $f$ является отображением

$$\partial^2 f:D\to\mathcal{L}(D,D;D)$$

Теорема 8.1.3. *Мы можем представить производную Гато второго порядка отображения $f$ в виде*

$$\partial^2 f(x)(a_1;a_2)=\frac{\partial^2_{s\cdot 0}f(x)}{\partial x^2}\sigma_s(a_1)\frac{\partial^2_{s\cdot 1}f(x)}{\partial x^2}\sigma_s(a_2)\frac{\partial^2_{s\cdot 2}f(x)}{\partial x^2} \tag{8.1.3}$$

Доказательство. Следствие определения 8.1.1 и теоремы 4.2.3. □

Определение 8.1.4. Мы будем называть выражение[8.1] $\frac{\partial^2_{s\cdot p} f(x)}{\partial x^2}$, $p=0,1,2$, **компонентой производной Гато второго порядка** отображения $f(x)$. □

По индукции, предполагая, что определена производная Гато $\partial^{n-1}f(x)$ порядка $n-1$, мы определим

$$\partial^n f(x)(a_1;...;a_n)=\frac{\partial^n f(x)}{\partial x^n}(a_1;...;a_n)=\partial(\partial^{n-1}f(x)(a_1;...;a_{n-1}))(a_n) \tag{8.1.4}$$

**производную Гато порядка** $n$ отображения $f$. Мы будем также полагать $\partial^0 f(x)=f(x)$.

## 8.2. Ряд Тейлора

Рассмотрим многочлен одной переменной над телом $D$ степени $n$, $n>0$. Нас интересует структура одночлена $p_k(x)$ многочлена степени $k$.

Очевидно, что одночлен степени 0 имеет вид $a_0$, $a_0\in D$. Пусть $k>0$. Докажем, что

$$p_k(x)=p_{k-1}(x)xa_k$$

где $a_k\in D$. Действительно, последний множитель одночлена $p_k(x)$ является либо $a_k\in D$, либо имеет вид $x^l$, $l\ge 1$. В последнем случае мы положим $a_k=1$. Множитель, предшествующий $a_k$, имеет вид $x^l$, $l\ge 1$. Мы можем представить этот множитель в виде $x^{l-1}x$. Следовательно, утверждение доказано.

В частности, одночлен степени 1 имеет вид $p_1(x)=a_0xa_1$.

Не нарушая общности, мы можем положить $k=n$.

Теорема 8.2.1. *Для произвольного $m>0$ справедливо равенство*

$$\begin{aligned}
&\partial^m(f(x)x)(h_1;...;h_m)\\
&=\partial^m f(x)(h_1;...;h_m)x+\partial^{m-1}f(x)(h_1;...;h_{m-1})h_m\\
&+\partial^{m-1}f(x)(\widehat{h_1};...;h_{m-1};h_m)h_1+...\\
&+\partial^{m-1}f(x)(h_1;...;\widehat{h_{m-1}};h_m)h_{m-1}
\end{aligned} \tag{8.2.1}$$

*где символ $\widehat{h^i}$ означает отсутствие переменной $h^i$ в списке.*

Доказательство. Для $m=1$ - это следствие равенства (5.2.23)

$$\partial(f(x)x)(h_1)=\partial f(x)(h_1)x+f(x)h_1$$

Допустим, (8.2.1) справедливо для $m-1$. Тогда

$$\begin{aligned}
&\partial^{m-1}(f(x)x)(h_1;...;h_{m-1})\\
&=\partial^{m-1}f(x)(h_1;...;h_{m-1})x+\partial^{m-2}f(x)(h_1;...;h_{m-2})h_{m-1}\\
&+\partial^{m-2}f(x)(\widehat{h_1},...,h_{m-2},h_{m-1})h_1+...\\
&+\partial^{m-2}f(x)(h_1,...,\widehat{h_{m-2}},h_{m-1})h_{m-2}
\end{aligned} \tag{8.2.2}$$

[8.1] Мы полагаем

$$\frac{\partial^2_{s\cdot p}f(x)}{\partial x^2}=\frac{\partial^2_{s\cdot p}f(x)}{\partial x\partial x}$$

Пользуясь равенствами (5.2.23) и (5.3.2) получим

$$
\begin{aligned}
&\partial^m(f(x)x)(h_1;...;h_{m-1};h_m)\\
=&\partial^m f(x)(h_1;...;h_{m-1};h_m)x+\partial^{m-1}f(x)(h_1;...;h_{m-2};h_{m-1})h_m\\
&+\partial^{m-1}f(x)(h_1;...;h_{m-2};\widehat{h_{m-1}};h_m)h_{m-1}\\
&+\partial^{m-2}f(x)(\widehat{h_1},...,h_{m-2},h_{m-1};h_m)h_1+...\\
&+\partial^{m-2}f(x)(h_1,...,\widehat{h_{m-2}},h_{m-1};h_m)h_{m-2}
\end{aligned}
\tag{8.2.3}
$$

Равенства (8.2.1) и (8.2.3) отличаются только формой записи. Теорема доказана. □

Теорема 8.2.2. *Производная Гато $\partial^m p_n(x)(h_1,...,h_m)$ является симметричным многочленом по переменным $h_1$, ..., $h_m$.*

Доказательство. Для доказательства теоремы мы рассмотрим алгебраические свойства производной Гато и дадим эквивалентное определение. Мы начнём с построения одночлена. Для произвольного одночлена $p_n(x)$ мы построим симметричный многочлен $r_n(x)$ согласно следующим правилам

- Если $p_1(x)=a_0xa_1$, то $r_1(x_1)=a_0x_1a_1$
- Если $p_n(x)=p_{n-1}(x)a_n$, то

$$r_n(x_1,...,x_n)=r_{n-1}(x_{[1},...,x_{n-1})x_{n]}a_n$$

  где квадратные скобки выражают симметризацию выражения по переменным $x_1$, ..., $x_n$.

Очевидно, что

$$p_n(x)=r_n(x_1,...,x_n)\quad x_1=...=x_n=x$$

Мы определим производную Гато порядка $k$ согласно правилу

$$\partial^k p_n(x)(h_1,...,h_k)=r_n(h_1,...,h_k,x_{k+1},...,x_n)\quad x_{k+1}=x_n=x \tag{8.2.4}$$

Согласно построению многочлен $r_n(h_1,...,h_k,x_{k+1},...,x_n)$ симметричен по переменным $h_1$, ..., $h_k$, $x_{k+1}$, ..., $x_n$. Следовательно, многочлен (8.2.4) симметричен по переменным $h_1$, ..., $h_k$.

При $k=1$ мы докажем, что определение (8.2.4) производной Гато совпадает с определением (5.2.18).

Для $n=1$, $r_1(h_1)=a_0h_1a_1$. Это выражение совпадает с выражением производной Гато в теореме 5.3.3.

Пусть утверждение справедливо для $n-1$. Справедливо равенство

$$r_n(h_1,x_2,...,x_n)=r_{n-1}(h_1,x_{[2},...,x_{n-1})x_{n]}a_n+r_{n-1}(x_2,...,x_n)h_1a_n \tag{8.2.5}$$

Положим $x_2=...=x_n=x$. Согласно предположению индукции, из равенств (8.2.4), (8.2.5) следует

$$r_n(h_1,x_2,...,x_n)=\partial p_{n-1}(x)(h_1)xa_n+p_{n-1}(x)h_1a_n$$

Согласно теореме 8.2.1

$$r_n(h_1,x_2,...,x_n)=\partial p_n(x)(h_1)$$

что доказывает равенство (8.2.4) для $k=1$.

Докажем теперь, что определение (8.2.4) производной Гато совпадает с определением (8.1.4) для $k>1$.

Пусть равенство (8.2.4) верно для $k-1$. Рассмотрим произвольное слагаемое многочлена $r_n(h_1,...,h_{k-1},x_k,...,x_n)$. Отождествляя переменные $h_1$, ..., $h_{k-1}$ с элементами тела $D$, мы рассмотрим многочлен

$$R_{n-k}(x_k,...,x_n)=r_n(h_1,...,h_{k-1},x_k,...,x_n) \tag{8.2.6}$$

Положим $P_{n-k}(x)=R_{n-k}(x_k,...,x_n)$, $x_k=...=x_n=x$. Следовательно

$$P_{n-k}(x)=\partial^{k-1}p_n(x)(h_1;...;h_{k-1})$$

Согласно определению (8.1.4) производной Гато

$$\begin{aligned}\partial P_{n-k}(x)(h_k)=&\partial(\partial^{k-1}p_n(x)(h_1;...;h_{k-1}))(h_k)\\=&\partial^k p_n(x)(h_1;...;h_{k-1};h_k)\end{aligned} \tag{8.2.7}$$

Согласно определению производной Гато (8.2.4)

$$\partial P_{n-k}(x)(h_k)=R_{n-k}(h_k,x_{k+1},...,x_n)\quad x_{k+1}=x_n=x \tag{8.2.8}$$

Согласно определению (8.2.6), из равенства (8.2.8) следует

$$\partial P_{n-k}(x)(h_k)=r_n(h_1,...,h_k,x_{k+1},...,x_n)\quad x_{k+1}=x_n=x \tag{8.2.9}$$

Из сравнения равенств (8.2.7) и (8.2.9) следует

$$\partial^k p_n(x)(h_1;...;h_k)=r_n(h_1,...,h_k,x_{k+1},...,x_n)\quad x_{k+1}=x_n=x$$

Следовательно равенсто (8.2.4) верно для любых $k$ и $n$.

Утверждение теоремы доказано. □

Теорема 8.2.3. *Для произвольного $n\ge 0$ справедливо равенство*

$$\partial^{n+1}p_n(x)(h_1;...;h_{n+1})=0 \tag{8.2.10}$$

Доказательство. Так как $p_0(x)=a_0$, $a_0\in D$, то при $n=0$ теорема является следствием теоремы 5.3.1. Пусть утверждение теоремы верно для $n-1$. Согласно теореме 8.2.1 при условии $f(x)=p_{n-1}(x)$ мы имеем

$$\begin{aligned}\partial^{n+1}p_n(x)(h_1;...;h_{n+1})=&\partial^{n+1}(p_{n-1}(x)xa_n)(h_1;...;h_{n+1})\\=&\partial^{n+1}p_{n-1}(x)(h_1;...;h_m)xa_n\\&+\partial^n p_{n-1}(x)(h_1;...;h_{m-1})h_m a_n\\&+\partial^n p_{n-1}(x)(\widehat{h_1};...;h_{m-1};h_m)h_1a_n+...\\&+\partial^{m-1}p_{n-1}(x)(h_1;...;\widehat{h_{m-1}};h_m)h_{m-1}a_n\end{aligned}$$

Согласно предположению индукции все одночлены равны 0. □

Теорема 8.2.4. *Если $m<n$, то справедливо равенство*

$$\partial^m p_n(0)(h)=0 \tag{8.2.11}$$

Доказательство. Для $n=1$ справедливо равенство

$$\partial^0 p_1(0)=a_0xa_1=0$$

Допустим, утверждение справедливо для $n-1$. Тогда согласно теореме 8.2.1

$$\begin{aligned}&\partial^m(p_{n-1}(x)xa_n)(h_1;...;h_m)\\=&\partial^m p_{n-1}(x)(h_1;...;h_m)xa_n+\partial^{m-1}p_{n-1}(x)(h_1;...;h_{m-1})h_ma_n\\&+\partial^{m-1}p_{n-1}(x)(\widehat{h_1};...;h_{m-1};h_m)h_1a_n+...\\&+\partial^{m-1}p_{n-1}(x)(h_1;...;\widehat{h_{m-1}};h_m)h_{m-1}a_n\end{aligned}$$

Первое слагаемое равно 0 так как $x=0$. Так как $m-1<n-1$, то остальные слагаемые равны 0 согласно предположению индукции. Утверждение теоремы доказано. □

Если $h_1=...=h_n=h$, то мы положим

$$\partial^n f(x)(h)=\partial^n f(x)(h_1;...;h_n)$$

Эта запись не будет приводить к неоднозначности, так как по числу аргументов ясно, о какой функции идёт речь.

Теорема 8.2.5. *Для произвольного $n>0$ справедливо равенство*

$$\partial^n p_n(x)(h)=n!p_n(h) \tag{8.2.12}$$

Доказательство. Для $n=1$ справедливо равенство

$$\partial p_1(x)(h)=\partial(a_0xa_1)(h)=a_0ha_1=1!p_1(h)$$

Допустим, утверждение справедливо для $n-1$. Тогда согласно теореме 8.2.1

$$\begin{aligned}\partial^n p_n(x)(h)=&\partial^n p_{n-1}(x)(h)xa_n+\partial^{n-1}p_{n-1}(x)(h)ha_n\\&+...+\partial^{n-1}p_{n-1}(x)(h)ha_n\end{aligned} \tag{8.2.13}$$

Первое слагаемое равно 0 согласно теореме 8.2.3. Остальные $n$ слагаемых равны, и согласно предположению индукции из равенства (8.2.13) следует

$$\partial^n p_n(x)(h)=n\partial^{n-1}p_{n-1}(x)(h)ha_n=n(n-1)!p_{n-1}(h)ha_n=n!p_n(h)$$

Следовательно, утверждение теоремы верно для любого $n$. □

Пусть $p(x)$ - многочлен степени $n$.[8.2]

$$p(x)=p_0+p_{1i_1}(x)+...+p_{ni_n}(x)$$

Мы предполагаем сумму по индексу $i_k$, который нумерует слагаемые степени $k$. Согласно теоремам 8.2.3, 8.2.4, 8.2.5

$$\partial^k p(0)(x)=k!p_{ki_k}(x)$$

Следовательно, мы можем записать

$$p(x)=p_0+(1!)^{-1}\partial p(0)(x)+(2!)^{-1}\partial^2 p(0)(x)+...+(n!)^{-1}\partial^n p(0)(x)$$

Это представление многочлена называется **формулой Тейлора для многочлена**. Если рассмотреть замену переменных $x=y-y_0$, то рассмотренное построение остаётся верным для многочлена

$$p(y)=p_0+p_{1i_1}(y-y_0)+...+p_{ni_n}(y-y_0)$$

---

[8.2] Я рассматриваю формулу Тейлора для многочлена по аналогии с построением формулы Тейлора в [7], с. 246.

откуда следует

$$p(y)=p_0+(1!)^{-1}\partial p(y_0)(y-y_0)+(2!)^{-1}\partial^2 p(y_0)(y-y_0)+...+(n!)^{-1}\partial^n p(y_0)(y-y_0)$$

Предположим, что функция $f(x))$ в точке $x_0$ дифференцируема в смысле Гато до любого порядка.[8.3]

Теорема 8.2.6. *Если для функции $f(x)$ выполняется условие*

$$f(x_0)=\partial f(x_0)(h)=...=\partial^n f(x_0)(h)=0 \tag{8.2.14}$$

*то при $t\to 0$ выражение $f(x+th)$ является бесконечно малой порядка выше $n$ по сравнению с $t$*

$$f(x_0+th)=o(t^n)$$

Доказательство. При $n=1$ это утверждение следует из равенства (5.2.13).

Пусть утверждение справедливо для $n-1$. Для отображения

$$f_1(x)=\partial f(x)(h)$$

выполняется условие

$$f_1(x_0)=\partial f_1(x_0)(h)=...=\partial^{n-1} f_1(x_0)(h)=0$$

Согласно предположению индукции

$$f_1(x_0+th)=o(t^{n-1})$$

Тогда равенство (5.2.12) примет вид

$$o(t^{n-1})=\lim_{t\to 0,\ t\in R}(t^{-1}f(x+th))$$

Следовательно,

$$f(x+th)=o(t^n)$$

□

Составим многочлен

$$p(x)=f(x_0)+(1!)^{-1}\partial f(x_0)(x-x_0)+...+(n!)^{-1}\partial^n f(x_0)(x-x_0) \tag{8.2.15}$$

Согласно теореме 8.2.6

$$f(x_0+t(x-x_0))-p(x_0+t(x-x_0))=o(t^n)$$

Следовательно, полином $p(x)$ является хорошей апроксимацией отображения $f(x)$.

Если отображение $f(x)$ имеет производную Гато любого порядка, то переходя к пределу $n\to\infty$, мы получим разложение в ряд

$$f(x)=\sum_{n=0}^{\infty}(n!)^{-1}\partial^n f(x_0)(x-x_0)$$

который называется **рядом Тейлора**.

[8.3] Я рассматриваю построение ряда Тейлора по аналогии с построением ряда Тейлора в [7], с. 248, 249.

## 8.3. Интеграл

Понятие интеграла имеет разные аспекты. В этом разделе мы рассмотрим интегрирование, как операцию, обратную дифференцированию. По сути дела, мы рассмотрим процедуру решения обыкновенного дифференциального уравнения

$$\partial f(x)(h) = F(x; h)$$

Пример 8.3.1. Я начну с примера дифференциального уравнения над полем действительных чисел.

$$y' = 3x^2 \tag{8.3.1}$$

$$x_0 = 0 \quad y_0 = 0 \tag{8.3.2}$$

Последовательно дифференцируя равенство (8.3.1), мы получаем цепочку уравнений

$$y'' = 6x \tag{8.3.3}$$

$$y''' = 6 \tag{8.3.4}$$

$$y^{(n)} = 0 \qquad n > 3 \tag{8.3.5}$$

Из уравнений (8.3.1), (8.3.2), (8.3.3), (8.3.4), (8.3.5) следует разложение в ряд Тейлора

$$y = x^3$$

□

Пример 8.3.2. Рассмотрим аналогичное дифференциальное уравнение над телом

$$\partial(y)(h) = hx^2 + xhx + x^2h \tag{8.3.6}$$

$$x_0 = 0 \quad y_0 = 0 \tag{8.3.7}$$

Последовательно дифференцируя равенство (8.3.6), мы получаем цепочку уравнений

$$\begin{aligned}&\partial^2(y)(h_1; h_2)\\ =&h_1h_2x + h_1xh_2 + h_2h_1x + xh_1h_2 + h_2xh_1 + xh_2h_1\end{aligned} \tag{8.3.8}$$

$$\begin{aligned}&\partial^3(y)(h_1; h_2; h_3)\\ =&h_1h_2h_3 + h_1h_3h_2 + h_2h_1h_3 + h_3h_1h_2 + h_2h_3h_1 + h_3h_2h_1\end{aligned} \tag{8.3.9}$$

$$\partial^n(y)(h_1; ...; h_n) = 0 \quad n > 3 \tag{8.3.10}$$

Из уравнений (8.3.6), (8.3.7), (8.3.8), (8.3.9), (8.3.10) следует разложение в ряд Тейлора

$$y = x^3$$

□

Замечание 8.3.3. Дифференциальное уравнение

$$\partial(y)(h) = 3hx^2 \tag{8.3.11}$$

$$x_0 = 0 \quad y_0 = 0 \tag{8.3.12}$$

так же приводит к решению $y = x^3$. Очевидно, что это отображение не удовлетворяет дифференциальному уравнению. Однако, вопреки теореме 8.2.2 вторая производная не является симметричным многочленом. Это говорит о том, что уравнение (8.3.11) не имеет решений. □

Пример 8.3.4. Очевидно, если функция удовлетворяет дифференциальному уравнению

$$\partial(y)(h) = f_{s\cdot 0}\ h\ f_{s\cdot 1} \tag{8.3.13}$$

то вторая производная Гато

$$\partial^2 f(x)(h_1; h_2) = 0$$

Следовательно, если задано начальное условие $y(0) = 0$, то дифференциальное уравнение (8.3.13) имеет решение

$$y = f_{s\cdot 0} x f_{s\cdot 1}$$

□

## 8.4. Экспонента

В этом разделе мы рассмотрим одну из возможных моделей построения экспоненты.

В поле мы можем определить экспоненту как решение дифференциального уравнения

$$y' = y \tag{8.4.1}$$

Очевидно, что мы не можем записать подобное уравнения для тела. Однако мы можем воспользоваться равенством

$$\partial(y)(h) = y'h \tag{8.4.2}$$

Из уравнений (8.4.1), (8.4.2) следует

$$\partial(y)(h) = yh \tag{8.4.3}$$

Это уравнение уже ближе к нашей цели, однако остаётся открытым вопрос в каком порядке мы должны перемножать $y$ и $h$. Что бы ответить на этот вопрос, мы изменим запись уравнения

$$\partial(y)(h) = \frac{1}{2}(yh + hy) \tag{8.4.4}$$

Следовательно, наша задача - решить дифференциальное уравнение (8.4.4) при начальном условии $y(0) = 1$.

Для формулировки и доказательства теоремы 8.4.1 я введу следующее обозначение. Пусть

$$\sigma = \begin{pmatrix} y & h_1 & ... & h_n \\ \sigma(y) & \sigma(h_1) & ... & \sigma(h_n) \end{pmatrix}$$

перестановка кортежа переменных

$$\begin{pmatrix} y & h_1 & ... & h_n \end{pmatrix}$$

Обозначим $p_\sigma(h_i)$ позицию, которую занимает переменная $h_i$ в кортеже

$$\begin{pmatrix} \sigma(y) & \sigma(h_1) & ... & \sigma(h_n) \end{pmatrix}$$

Например, если перестановка $\sigma$ имеет вид

$$\begin{pmatrix} y & h_1 & h_2 & h_3 \\ h_2 & y & h_3 & h_1 \end{pmatrix}$$

то следующие кортежи равны

$$\begin{aligned}\begin{pmatrix}\sigma(y) & \sigma(h_1) & \sigma(h_2) & \sigma(h_3)\end{pmatrix} &= \begin{pmatrix} h_2 & y & h_3 & h_1\end{pmatrix}\\ &= \begin{pmatrix} p_\sigma(h_2) & p_\sigma(y) & p_\sigma(h_3) & p_\sigma(h_1)\end{pmatrix}\end{aligned}$$

Теорема 8.4.1. *Производная Гато порядка $n$ функции $y$, удовлетворяющей дифференциальному уравнению (8.4.4) имеет вид*

$$\partial^n(y)(h_1, ..., h_n) = \frac{1}{2^n}\sum_\sigma \sigma(y)\sigma(h_1)...\sigma(h_n) \tag{8.4.5}$$

*где сумма выполнена по перестановкам*

$$\sigma = \begin{pmatrix} y & h_1 & ... & h_n \\ \sigma(y) & \sigma(h_1) & ... & \sigma(h_n)\end{pmatrix}$$

*множества переменных $y$, $h_1$, ..., $h_n$. Перестановка $\sigma$ обладает следующими свойствами*

(1) *Если существуют $i$, $j$, $i \neq j$, такие, что $p_\sigma(h_i)$ располагается в произведении (8.4.5) левее $p_\sigma(h_j)$ и $p_\sigma(h_j)$ располагается левее $p_\sigma(y)$, то $i < j$.*

(2) *Если существуют $i$, $j$, $i \neq j$, такие, что $p_\sigma(h_i)$ располагается в произведении (8.4.5) правее $p_\sigma(h_j)$ и $p_\sigma(h_j)$ располагается правее $p_\sigma(y)$, то $i > j$.*

Доказательство. Мы докажем это утверждение индукцией. Для $n = 1$ утверждение верно, так как это дифференциальное уравнение (8.4.4). Пусть утверждение верно для $n = k - 1$. Следовательно

$$\partial^{k-1}(y)(h_1, ..., h_{k-1}) = \frac{1}{2^{k-1}}\sum_\sigma \sigma(y)\sigma(h_1)...\sigma(h_{k-1}) \tag{8.4.6}$$

где сумма выполнена по перестановкам

$$\sigma = \begin{pmatrix} y & h_1 & ... & h_{k-1} \\ \sigma(y) & \sigma(h_1) & ... & \sigma(h_{k-1})\end{pmatrix}$$

множества переменных $y$, $h_1$, ..., $h_{k-1}$. Перестановка $\sigma$ удовлетворяет условиям (1), (2), сформулированным в теореме. Согласно определению (8.1.4) производная Гато‘ порядка $k$ имеет вид

$$\begin{aligned}\partial^k(y)(h_1, ..., h_k) =&\partial(\partial^{k-1}(y)(h_1, ..., h_{k-1}))(h_k)\\ =&\frac{1}{2^{k-1}}\partial\left(\sum_\sigma \sigma(y)\sigma(h_1)...\sigma(h_{k-1})\right)(h_k) \end{aligned}\tag{8.4.7}$$

Из равенств (8.4.4), (8.4.7) следует

$$
\begin{aligned}
&\partial^k(y)(h_1,...,h_k)\\
&=\frac{1}{2^{k-1}}\frac{1}{2}\left(\sum_\sigma \sigma(yh_k)\sigma(h_1)...\sigma(h_{k-1})+\sum_\sigma \sigma(h_ky)\sigma(h_1)...\sigma(h_{k-1})\right)
\end{aligned}
\tag{8.4.8}
$$

Нетрудно видеть, что произвольная перестановка $\sigma$ из суммы (8.4.8) порождает две перестановки

$$
\begin{aligned}
\tau_1 &= \begin{pmatrix} y & h_1 & ... & h_{k-1} & h_k\\ \tau_1(y) & \tau_1(h_1) & ... & \tau_1(h_{k-1}) & \tau_1(h_k)\end{pmatrix}\\
&= \begin{pmatrix} h_ky & h_1 & ... & h_{k-1}\\ \sigma(h_ky) & \sigma(h_1) & ... & \sigma(h_{k-1})\end{pmatrix}\\
\tau_2 &= \begin{pmatrix} y & h_1 & ... & h_{k-1} & h_k\\ \tau_2(y) & \tau_2(h_1) & ... & \tau_2(h_{k-1}) & \tau_2(h_k)\end{pmatrix}\\
&= \begin{pmatrix} yh_k & h_1 & ... & h_{k-1}\\ \sigma(yh_k) & \sigma(h_1) & ... & \sigma(h_{k-1})\end{pmatrix}
\end{aligned}
\tag{8.4.9}
$$

Из (8.4.8) и (8.4.9) следует

$$
\begin{aligned}
&\partial^k(y)(h_1,...,h_k)\\
&=\frac{1}{2^k}\left(\sum_{\tau_1}\tau_1(y)\tau_1(h_1)...\tau_1(h_{k-1})\tau_1(h_k)\right.\\
&\left.+\sum_{\tau_2}\tau_2(y)\tau_2(h_1)...\tau_2(h_{k-1})\tau_2(h_k)\right)
\end{aligned}
\tag{8.4.10}
$$

В выражении (8.4.10) $p_{\tau_1}(h_k)$ записано непосредственно перед $p_{\tau_1}(y)$. Так как $k$ - самое большое значение индекса, то перестановка $\tau_1$ удовлетворяет условиям (1), (2), сформулированным в теореме. В выражении (8.4.10) $p_{\tau_2}(h_k)$ записано непосредственно после $p_{\tau_2}(y)$. Так как $k$ - самое большое значение индекса, то перестановка $\tau_2$ удовлетворяет условиям (1), (2), сформулированным в теореме.

Нам осталось показать, что в выражении (8.4.10) перечислены все перестановки $\tau$, удовлетворяющие условиям (1), (2), сформулированным в теореме. Так как $k$ - самый большой индекс, то согласно условиям (1), (2), сформулированным в теореме, $\tau(h_k)$ записано непосредственно перед или непосредственно после $\tau(y)$. Следовательно, любая перестановка $\tau$ имеет либо вид $\tau_1$, либо вид $\tau_2$. Пользуясь равенством (8.4.9), мы можем для заданной перестановки $\tau$ найти соответствующую перестановку $\sigma$. Следовательно, утверждение теоремы верно для $n=k$. Теорема доказана. □

Теорема 8.4.2. *Решением дифференциального уравнения* (8.4.4) *при начальном условии* $y(0)=1$ *является экспонента* $y=e^x$ *которая имеет следующее разложение в ряд Тейлора*

$$
e^x=\sum_{n=0}^\infty\frac{1}{n!}x^n
\tag{8.4.11}
$$

Доказательство. Производная Гато порядка $n$ содержит $2^n$ слагаемых. Действительно, производная Гато порядка 1 содержит 2 слагаемых, и каждое дифференцирование увеличивает число слагаемых вдвое. Из начального условия $y(0)=1$ и теоремы 8.4.1 следует, что производная Гато порядка $n$ искомого решения имеет вид

$$\partial^n(y(0))(h,...,h)=1 \tag{8.4.12}$$

Из равенства (8.4.12) следует разложение (8.4.11). в ряд Тейлора. □

Теорема 8.4.3. *Равенство*

$$e^{a+b}=e^ae^b \tag{8.4.13}$$

*справедливо тогда и только тогда, когда*

$$ab=ba \tag{8.4.14}$$

Доказательство. Для доказательства теоремы достаточно рассмотреть ряды Тейлора

$$e^a=\sum_{n=0}^\infty\frac{1}{n!}a^n \tag{8.4.15}$$

$$e^b=\sum_{n=0}^\infty\frac{1}{n!}b^n \tag{8.4.16}$$

$$e^{a+b}=\sum_{n=0}^\infty\frac{1}{n!}(a+b)^n \tag{8.4.17}$$

Перемножим выражения (8.4.15) и (8.4.16). Сумма одночленов порядка 3 имеет вид

$$\frac{1}{6}a^3+\frac{1}{2}a^2b+\frac{1}{2}ab^2+\frac{1}{6}b^3 \tag{8.4.18}$$

и не совпадает, вообще говоря, с выражением

$$\frac{1}{6}(a+b)^3=\frac{1}{6}a^3+\frac{1}{6}a^2b+\frac{1}{6}aba+\frac{1}{6}ba^2+\frac{1}{6}ab^2+\frac{1}{6}bab+\frac{1}{6}b^2a+\frac{1}{6}b^3 \tag{8.4.19}$$

Доказательство утверждения, что (8.4.13) следует из (8.4.14) тривиально. □

Смысл теоремы 8.4.3 становится яснее, если мы вспомним, что существует две модели построения экспоненты. Первая модель - это решение дифференциального уравнения (8.4.4). Вторая - это изучение однопараметрической группы преобразований. В случае поля обе модели приводят к одной и той же функции. Я не могу этого утверждать сейчас в общем случае. Это вопрос отдельного исследования. Но если вспомнить, что кватернион является аналогом преобразования трёхмерного пространства, то утверждение теоремы становится очевидным.

Глава 9

# Производная второго порядка отображения $D$-векторного пространства

## 9.1. Производная второго порядка отображения $D$-векторного пространства

Пусть $D$ - нормированное тело характеристики 0. Пусть $V$ - нормированное $D$-векторное пространство с нормой $\|y\|_V$. Пусть $W$ - нормированное $D$-векторное пространство с нормой $\|z\|_W$. Пусть

$$f : V \to W$$

функция, дифференцируемая по Гато. Согласно замечанию 6.2.2 производная Гато является отображением

$$\partial f : V \to \mathcal{L}(D; V; W)$$

Согласно теоремам 3.1.5, 3.1.6 и определению 6.1.7 множество $\mathcal{L}(D; V; W)$ является нормированным $D$-векторным пространством. Следовательно, мы можем рассмотреть вопрос, является ли отображение $\partial f$ дифференцируемым по Гато.

Согласно определению 6.2.1

$$\partial f(x + a_2)(a_1) - \partial f(x)(a_1) = \partial(\partial f(x)(a_1))(a_2) + o_2(a_2) \tag{9.1.1}$$

где $o_2 : V \to \mathcal{L}(D; V; W)$ - такое непрерывное отображение, что

$$\lim_{a_2 \to 0} \frac{\|o_2(a_2)\|_W}{\|a_2\|_V} = 0$$

Согласно определению 6.2.1 отображение $\partial(\partial f(x)(a_1))(a_2)$ линейно по переменной $a_2$. Из равенства (9.1.1) следует, что отображение $\partial(\partial f(x)(a_1))(a_2)$ линейно по переменной $a_1$.

Определение 9.1.1. Полилинейное отображение

$$\partial^2 f(x)(a_1; a_2) = \partial(\partial f(x)(a_1))(a_2) \tag{9.1.2}$$

называется **производной Гато второго порядка** отображения $f$. □

Замечание 9.1.2. Согласно определению 9.1.1 при заданном $x$ производная Гато второго порядка $\partial^2 f(x) \in \mathcal{L}(D; V, V; W)$. Следовательно, производная Гато второго порядка отображения $f$ является отображением

$$\partial^2 f : V \to \mathcal{L}(D; V, V; W)$$

Теорема 9.1.3. *Допустим $\overline{\overline{e}}_V$ - базис в $D$-векторном пространстве $V$. Допустим $\overline{\overline{e}}_W$ - базис в $D$-векторном пространстве $W$. Производная Гато* *второго порядка отображения $f$ относительно базиса $\overline{\overline{e}}_V$ и базиса $\overline{\overline{e}}_W$ имеет вид*

$$\partial^2 f(x)(a_1;a_2) = \frac{\partial^2_{s\cdot 0} f^{j}(x)}{\partial x^{i_1}\partial x^{i_2}}\ \sigma_s(a_1^{i_1})\ \frac{\partial^2_{s\cdot 1} f^{j}(x)}{\partial x^{i_1}\partial x^{i_2}}\ \sigma_s(a_2^{i_2})\ \frac{\partial^2_{s\cdot 2} f^{j}(x)}{\partial x^{i_1}\partial x^{i_2}} e_{W\cdot j} \tag{9.1.3}$$

Доказательство. Следствие определения 9.1.1 и теоремы 4.2.3. □

Определение 9.1.4. Мы будем называть выражение $\dfrac{\partial^2_{s\cdot p} f^{j}(x)}{\partial x^{i_1}\partial x^{i_2}}$, $p=0,1,2$, **компонентой производной Гато второго порядка** отображения $f(x)$. □

Определение 9.1.5. Пусть $D$ - тело характеристики 0. Допустим $\overline{\overline{e}}_V$ - базис в $D$-векторном пространстве $V$. Допустим $\overline{\overline{e}}_W$ - базис в $D$-векторном пространстве $W$. Частная производная Гато

$$\frac{\partial^2 f^{k}(v)}{\partial v^{j}\partial v^{i}}(h^{i},h^{j}) = \frac{\partial}{\partial v^{j}}\left(\frac{\partial f^{k}(v)}{\partial v^{i}}(h^{i})\right)(h^{j}) \tag{9.1.4}$$

называется **смешанной частной производной Гато** отображения $f^{k}$ по переменным $v^{i}$, $v^{j}$. □

Теорема 9.1.6. *Пусть $D$ - тело характеристики 0.*[9.1] *Допустим $\overline{\overline{e}}_V$ - базис в $D$-векторном пространстве $V$. Допустим $\overline{\overline{e}}_W$ - базис в $D$-векторном пространстве $W$. Пусть частные производные Гато отображения*

$$f:V\to W$$

*непрерывны и дифференцируемы в смысле Гато в области $U\subset V$. Пусть смешанные частные производные Гато второго порядка непрерывны в области $U\subset V$. Тогда на множестве $U$ смешанные частные производные Гато удовлетворяют равенству*

$$\frac{\partial^2 f^{k}(v)}{\partial v^{j}\partial v^{i}}(h^{i},h^{j}) = \frac{\partial^2 f^{k}(v)}{\partial v^{i}\partial v^{j}}(h^{j},h^{i}) \tag{9.1.5}$$

Доказательство. Для доказательства теоремы мы будем пользоваться обозначением $f^{k}(v)=f^{k}(v^{i},v^{j})$, подразумевая, что остальные переменные имеют заданные значения. Утверждение теоремы следует из следующей цепочки

[9.1] Доказательство теоремы сделано по аналогии с доказательством теоремы [7], с. 405, 406.

равенств

$$\begin{aligned}
&\frac{\partial^2 f(v^i,v^j)}{\partial v^i\partial v^j}(h^j,h^i)\\
=&\frac{\partial}{\partial v^i}\left(\frac{\partial f(v^i,v^j)}{\partial v^j}(h^j)\right)(h^i)\\
=&\lim_{t\to 0} t^{-1}\left(\frac{\partial f(v^i+h^it,v^j)}{\partial v^j}(h^j)-\frac{\partial f(v^i,v^j)}{\partial v^j}(h^j)\right)\\
=&\lim_{t\to 0} t^{-2}\left(f(v^i+h^it,v^j+h^jt)-f(v^i+h^it,v^j)-f(v^i,v^j+h^jt)+f(v^i,v^j)\right)\\
=&\lim_{t\to 0} t^{-2}\left(f(v^i+h^it,v^j+h^jt)-f(v^i,v^j+h^jt)-f(v^i+h^it,v^j)+f(v^i,v^j)\right)\\
=&\lim_{t\to 0} t^{-1}\left(\frac{\partial f(v^i,v^j+h^jt)}{\partial v^i}(h^i)-\frac{\partial f(v^i,v^j)}{\partial v^i}(h^i)\right)\\
=&\frac{\partial}{\partial v^j}\left(\frac{\partial f(v^i,v^j)}{\partial v^i}(h^i)\right)(h^j)\\
=&\frac{\partial^2 f(v^i,v^j)}{\partial v^j\partial v^i}(h^i,h^j)
\end{aligned}$$

□

По индукции, предполагая, что определена производная Гато $\partial^{n-1}f(x)$ порядка $n-1$, мы определим

$$\partial^n f(x)(a_1;...;a_n)=\partial(\partial^{n-1}f(x)(a_1;...;a_{n-1}))a_n$$

**производную Гато порядка** $n$ отображения $f$. Мы будем также полагать $\partial^0 f(x)=f(x)$.

## 9.2. Однородное отображение

Определение 9.2.1. Отображение

$$f:V\to W$$

$D$-векторного пространства $V$ в $D$-векторное пространство $W$ называется **однородным отображением степени** $k$ **над полем** $F$, если для любого $a\in F$

(9.2.1) $$f(av)=a^kf(v)$$

□

Теорема 9.2.2 (теорема Эйлера). *Пусть отображение*

$$f:V\to W$$

*$D$-векторного пространства $V$ в $D$-векторное пространство $W$ однородно степени $k$ над полем $F\subset Z(D)$.*

(9.2.2) $$\partial f(v)(v)=kf(v)$$

Доказательство. Продифференцируем равенство (9.2.1) по $a$[9.2]

$$\frac{d}{da}(f(av))\Delta a=\frac{d}{da}(a^kf(v))\Delta a \tag{9.2.3}$$

Из теоремы 6.2.12 следует

$$\begin{aligned}\frac{d}{da}(f(av))\Delta a&=\partial f(av)(\partial(av)(\Delta a))\\&=\partial f(av)\left(\frac{d(av)}{da}\Delta a\right)\\&=\partial f(av)(v\Delta a)\end{aligned} \tag{9.2.4}$$

Из теоремы 6.2.5 и равенства (9.2.4) следует

$$\frac{d}{da}(f(av))\Delta a=\partial f(av)(v)\Delta a \tag{9.2.5}$$

Очевидно

$$\frac{d}{da}(a^kf(v))=ka^{k-1}f(v) \tag{9.2.6}$$

Подставим (9.2.5) и (9.2.6) в (9.2.3)

$$\partial f(av)(v)\Delta a=ka^{k-1}f(v)\Delta a \tag{9.2.7}$$

Так как $\Delta a\ne 0$, то (9.2.2) следует из (9.2.7) при условии $a=1$. □

Пусть отображение

$$f:V\to W$$

$D$-векторного пространства $V$ в $D$-векторное пространство $W$ однородно степени $k$ над полем $F\subset Z(D)$. Из теоремы 6.2.5 следует

$$\partial f(av)(av)=a\partial f(av)(v) \tag{9.2.8}$$

Из равенств (9.2.2), (9.2.1) следует

$$\partial f(av)(av)=kf(av)=ka^kf(v)=a^k\partial f(v)(v) \tag{9.2.9}$$

Из равенств (9.2.8), (9.2.9) следует

$$\partial f(av)(v)=a^{k-1}\partial f(v)(v) \tag{9.2.10}$$

Однако из равенства (9.2.10) не следует, что отображение

$$\partial f:V\to\mathcal{A}(D;V;W)$$

однородно степени $k-1$ над полем $F$.

[9.2] Для отображения поля справедливо

$$\partial f(x)(h)=\frac{df(x)}{dx}h$$

Глава 10

# Список литературы

Глава 11

# Предметный указатель

Глава 12

# Специальные символы и обозначения